\documentclass[11pt, a4paper,twoside]{article}
\usepackage{amsfonts,amsmath,amssymb,amsthm}
\usepackage{latexsym,amscd,mathrsfs}
\usepackage{amsbsy}

\newcommand{\mput}{\multiput}
\newcommand{\bcen}{\begin{center}}      \newcommand{\ecen}{\end{center}}
\def\az{\alpha}

\def\dz{\delta}
\def\lz{\lambda}

\def\A{{\cal A}}

\def\K{{\cal K}}
\def\rK{ k }
\def\M{{\cal M}}
\def\I{{\cal I}}
\def\J{{\cal J}}
\def\L{{\cal L}}

\def\P{{\cal P}}
\def\Q{{\cal Q}}
\def\l{\underline{l}}
\def\n{\underline{n}}
\def\m{\underline{m}}
\def\e{\underline{e}}

\def\rad{\mbox{rad}}

\def\wt{\widetilde}
\def\lra{\longrightarrow}

\def\Hom{\mbox{Hom}}

\def\dim{\mbox{dim}}
\def\mod{\mbox{mod}}
\def\ind{\mbox{ind}}
\def\Im{\mbox{Im}}
\def\Ker{\mbox{Ker}}
\def\Coker{\mbox{Coker}}

\def\soc{\mbox{soc}}
\def\rank{\mbox{rank}}
\def\s{\scriptstyle}
\def\c{c_{\big(\stackrel{u}{{\s\I}{\s\L}}\big)
\big(\stackrel{v}{{\s\L}{\s\J}}\big)}
^{\big(\stackrel{w}{{\s\I}{\s\J}}\big)}}
\def\kg{\hspace*{0.65cm}}
\begin{document}
\thispagestyle{empty} \vspace*{6.5cm}

\begin{center}
{\bf \Huge Tameness and Homogeneity
} \vspace{2cm}
\begin{center}
\Large Dedicated to Professor C.M. Ringel on the occasion of\\ his
60'th birthday
\end{center}
\end{center}
 \vspace{3cm}
\begin{center} {\bf \huge Zhang Yingbo\\[4mm]
Xu Yunge }
\end{center}

\newpage
\thispagestyle{empty} \tableofcontents \thispagestyle{empty}
\newpage

\bigskip
\thispagestyle{empty}
\begin{abstract}

Let $\Lambda$ be a finite-dimensional algebra over an algebraically closed field, then $\Lambda$ is either tame
or wild. Is there any homological  description in terms of AR-translations on tameness? Or equivalently, is
there any  combinatorial description in terms of AR-quivers? The answer is yes. In the present paper we prove
the following main theorem: ``$\Lambda$ is tame if and only if almost all modules are isomorphic to their
Auslander-Reiten translations, if and only if they lie in homogeneous tubes".

The method used in the paper is bimodule problem, which has been studied by several authors. We treat bimodule
problems and their reductions in terms of matrices and generalized Jordan forms respectively. In particular, we
introduce a concept of freely parameterized bimodule problems corresponding to layered bocses in order to define
the reductions of bimodule problems. Moreover, we list all the possibilities of the differential of the first
arrow of a layered bocs such that the induced bocs is still layered and preserves all the free parameters,
especially in wild case.

\bigskip



\end{abstract}

\newpage
\setcounter{page}{1}

\baselineskip=16pt
\parskip=4pt


\begin{center} \section{Introduction} \end{center}

\subsection{AR-sequences and AR-quivers}

\thispagestyle{empty}
 \kg Let $\Lambda$ be a finite-dimensional algebra
(associative, with 1) over an algebraically closed field $\rK$,
and $\Lambda$-mod be the category of finite-dimensional left
$\Lambda$-modules. Auslander and Reiten have defined a concept of
almost split sequence (i.e. AR-sequence) in a module category
\cite{AR}.  If we take any non-projective indecomposable
$\Lambda$-module $N$ (resp. non-injective module $M$), there
exists an almost split sequence ending at $N$ (resp. starting at
$M$):
$$
(e) \quad 0\lra M\stackrel{\iota}{\lra}E\stackrel{\pi}{\lra}N\lra
0
$$
where $M\cong DTr(N)$, the Auslander-Reiten translation of $N$, or
AR-translation of $N$ for short. $DTr:
\Lambda\mbox{-}\underline{\mbox{mod}}\rightarrow\Lambda\mbox{-}\overline{\mbox{mod}}$
is a functor, which maps indecomposable non-projective
$\Lambda$-modules to non-injective $\Lambda$-modules \cite{AR}.
AR-sequences yield a combinatorial description of $\Lambda$-mod.
Namely,  we draw a vertex $[M]$  corresponding to an iso-class of
indecomposable module $M$, and an arrow
$[L]\stackrel{[\zeta]}{\lra}[M]$ representing that there is an
irreducible map $\zeta: L\rightarrow M$, then the obtained quiver
is called an Auslander-Reiten quiver of $\Lambda$, or AR-quiver of
$\Lambda$ for short. If  we have
$$E=\oplus_{i=1}^lE_i, \ \iota=(\iota_1,
\iota_2,\cdots, \iota_l)\ {\rm and}\
\pi=(\pi_1,\pi_2,\cdots,\pi_l)^T$$  in the sequence $(e)$, where
$E_i$ are indecomposable, $T$ stands for the transpose of a
matrix, then $\iota_i: M\rightarrow E_i$ and $\pi_i:
E_i\rightarrow N$ are irreducible maps between indecomposables.
Conversely, if $M\stackrel{\iota'}{\longrightarrow}E'$ (resp.
$E'\stackrel{\pi'}{\longrightarrow} N$) with $E'$ being
indecomposable is an irreducible map, there exists some module
$E''$ and morphism $M\stackrel{\iota''}{\longrightarrow}E''$
(resp. $E''\stackrel{\pi''}{\longrightarrow} N$) such that the
sequence
$$(e') \quad \begin{CD} 0@>>> M @>(\iota',\; \iota'')>>E'\oplus
E'' @>{\pi'\choose \pi''}>> N @>>> 0 \end{CD}$$ is equivalent to
the almost split sequence $(e)$.
\medskip

 Furthermore, the shape of stable components of
AR-quivers has been completely described in \cite{Rie}, \cite{HPR}
and \cite{Z}. The  simplest stable components which mostly deserve
to pay attention to are so-called homogeneous tubes consisting of
 vertices $[M_i]$ and arrows \hspace{5mm}
  {\unitlength=1mm
\special{em:linewidth 0.4pt} \linethickness{0.4pt}
\begin{picture}(28.00, 5.00)
\put(6.00,2.00){\vector(1,0){8}}
\put(14.00,0.50){\vector(-1,0){8}}

\put(1.00,1.00){\makebox(0,0)[cc]{$[M_i]$}}
\put(21.00,1.00){\makebox(0,0)[cc]{$[M_{i+1}]$}}
\put(10.00,3.00){\makebox(0,0)[cc]{$\scriptscriptstyle
[\iota_i]$}}
\put(10.00,-0.70){\makebox(0,0)[cc]{$\scriptscriptstyle [\pi_i]$}}
\end{picture},} where $i$ is any positive integer.
 Thus almost
split sequences lying in a homogeneous tube have the shape of
$$
0\lra M_1\stackrel{\iota_1}{\lra}M_2\stackrel{\pi_1}{\lra}M_1\lra
0 $$
$$
\begin{CD}
 0@>>> M_i @>(\pi_{i-1},\;
\iota_i)>> M_{i-1}\oplus M_{i+1}
@>{\iota_{i-1}\choose\pi_i}>>M_i@>>> 0
\end{CD}$$ for $i \ge 2$, and $DTr(M_i)\cong M_i$ for $i=1,2,\cdots$.

\subsection{Tameness and wildness}
\kg The well known Drozd's  theorem tells us that a
finite-dimensional algebra $\Lambda$ over an algebraically closed
field $\rK$ is either of tame representation type or of wild
representation type \cite{D}.

{\bf Definition 1.2.1} \cite{D, CB1, DS} A finite-dimensional
$k$-algebra $\Lambda$ is of tame representation type, if for any
positive integer $d$, there are a finite number of localizations
$R_i=\rK[x,f_i(x)^{-1}]$ of $\rK[x]$ and $\Lambda$-$R_i$-bimodules
$T_i$ which are free  as right $R_i$-modules, such that almost all
(except finitely many) indecomposable $\Lambda$-modules of
dimension at most $d$ are isomorphic to
$$T_i\otimes_{R_i}R_i/(x-\lz)^{m},$$ for some  $\lz\in \rK,$
$f_i(\lz)\ne 0$, and some positive integer $m$.

{\bf Definition 1.2.2} \cite{D, CB1} A finite-dimensional
$\rK$-algebra $\Lambda$ is  of wild representation type if there
is a finitely generated $\Lambda$-$\rK\langle x,y\rangle$-bimodule
$T$, which is free as a right $\rK\langle x,y\rangle$-module, such
that the functor $$T\otimes_{\rK\langle x,y\rangle}-: k\langle
x,y\rangle\mbox{-mod}\rightarrow \Lambda\mbox{-mod}$$ preserves
indecomposability and isomorphism classes.

\medskip
Crawley-Boevey once proposed to consider  generic modules instead
of dealing with families of modules \cite{CB4}. A $\Lambda$-module
$G$ is called generic,  if $G$ is of infinite length over
$\Lambda$, but of finite endolength (the length over its own
endomorphism ring) and $G$ is indecomposable. Then $\Lambda$ is
tame, if and only if $\Lambda$ is generically tame. Generically
tame type means for any positive integer $d$, there are only
finitely many generic modules of endolength at most $d$. Based on
Crawley-Boevey's work on tame algebras, Krause presents a new
description of tameness in terms of Ziegler spectrum, functor
category and  model theory. More precisely, he proves that a
finite-dimensional algebra over an algebraically closed field is
of tame representation type, if and only if every generic
$\Lambda$-module appears as the only generic module on the Ziegler
closure of a homogeneous tube \cite{Kr}.

In summary, the definition of tameness given by Drozd, and the
 description in terms of generic modules by Crawley-Boevey,
 as well as the description in terms of Ziegler spectrum by Krause,  all
involve infinite-dimensional modules outside $\Lambda$-mod. A
natural question arises:
 is there any internal description of tameness, which only involves
finite-dimensional modules? Moreover, is there any combinatorial
description of tameness in terms of AR-quivers?

\subsection{The main theorem: tameness and
homogeneity}

\kg A result due to Crawley-Boevey may lead to an internal
description of tameness.

{\bf Theorem 1.3.1} \cite{CB1}: Let $\Lambda$ be a
finite-dimensional algebra over an algebraically closed field. If
$\Lambda$ is of tame representation type, then almost all
$\Lambda$-modules are isomorphic to their Auslander-Rieten
translations. And then lie in homogeneous tubes.

Here ``almost all modules'' means ``for any fixed positive integer
$d$, all but a finite number of isomorphism classes of
indecomposable $\Lambda$-modules  of dimension at most $d$.''

The proof was based on the method of bocses. In fact, a minimal
bocs possesses such a property. It is also conjectured in  the
same paper that the property may possibly describe tameness. Since
then, several experts have tried to prove the converse of the
theorem, expecting
 that the same method of bocses would work, that is, none of the
minimally wild bocses would have such a property. Unfortunately,
this was proved to be wrong. A counterexample  was constructed  in
\cite{ZLB}, and later a great number of wild  categories having
such a property were given in \cite{BCLZ} and \cite{V}.

Another approach to testify the converse of Crawley-Boevey's
theorem is to prove Ringel's conjecture. A wild algebra $\Lambda$
is called $\tau$-wild if there exist infinitely many
non-isomorphic $\tau$-variant (i.e. non-homogeneous)
indecomposable modules of dimension $d$, for a certain  positive
integer $d$, where $\tau=DTr$. Han Yang introduces the concept of
controlled wild algebras (see \cite{H} and \cite{R3}). And Ringel
conjectures that all wild algebras in the meaning  of Drozd are
  controlled wild with finite controlled index. Nagase shows
  that all controlled wild
algebras with finite controlling index are $\tau$-wild \cite{N2}.
Therefore, if Ringel's conjecture holds, so does the converse
 of Crawley-Boevey's theorem. The covering criterion for an algebra
to be controlled wild,  given by Han Yang in \cite{H} is very
effective. It seems to be impossible to find a concrete wild
algebra which is not controlled wild. However, Ringel's conjecture
is still open up to now.

In this paper we will prove the converse of Crawley-Boevey's
theorem: if almost all $\Lambda$-modules are isomorphic to their
$AR$-translations, then $\Lambda$ is tame.

{\bf Main theorem 1.3.2 } {\it Let $\Lambda$ be a
finite-dimensional algebra over an algebraically closed field.
Then $\Lambda$ is tame, if and only if almost all
$\Lambda$-modules are isomorphic to their AR-translations, and
then if and only if they lie in homogeneous tubes.}

\subsection{The outline of the proof}

\kg Our argument relies on the methods of bimodule problems \cite{S} and bocses \cite{Ro, D,CB1}. The notion of
bocses was introduced by Rojter in \cite{Ro} in order to apply formulation of the matrix problems to
representation theory. There are several formulations of matrix problems \cite{ZZ} such as bocses, differential
graded categories \cite{Ro}, differential biquivers,
 relatively projective categories
\cite{BK}, bimodule problems \cite{CB2}, lift categories
\cite{CB3}, etc. The following picture is suggested by
Crawley-Boevey.

\begin{minipage}{14cm}
$$
\begin{aligned}
& \{\mbox{matrix problems given by ``layered" or ``almost
free" bocses \cite{D,CB1}}\} \\
\supset & \{\mbox{matrix problems of the form }
(\K,\M,d), \mbox{ where } \K \mbox{ is an algebra, } \M \mbox{ an} \\
& \   \K\mbox{-}\K  \mbox{-bimodule}, \  d \mbox{ a derivation,
which is called a bimodule problem \cite{CB3}}\}\\
= & \{\mbox{matrix problems of the form }
(\K,\M,H), \mbox{ where } \K \mbox{ is an upper triangular} \\
& \mbox{ matrix algebra, } \M \mbox{  a matrix }
\K\mbox{-}\K\mbox{-bimodule},\  H \mbox{ is a fixed matrix, which} \\
& \  \mbox{determins a
derivation (see 2.2)}\} \\
 \supset &\{\mbox{matrix problem given by bipartite bimodule problems (see
 10.1)}\}\\
 \supset &\{\mbox{matrix problems of the form } P_1(\Lambda)\mbox{ (see
 2.1)}\}.
 \end{aligned}
 $$
 \end{minipage}

\vspace{6mm}
 \noindent Throughout the paper, we concentrate on the third set of matrix problems, which are
convenient for the calculations in the proof of the main theorem. On the other hand, such matrix problems are
helpful to understand the notion of bocses, since their structure is very concrete.

The reductions of matrix problems given by layered bocses are
well-established by the same authors, which seem to be elegant for
theoretical purposes. Among all the reductions, Belitskii's
reduction algorithm is very effective for reducing an individual
matrix to a canonical form, which may be considered as a
generalized Jordan form,  under some admissible transformations.
 In this paper, we will describe the
reductions in terms of generalized Jordan forms for a bimodule problem given in the third set, which is
 based only on linear algebra.

 The idea of our proof of the main theorem is as follows. Suppose an algebra $\Lambda$ is of wild
representation type and has homogeneous property, which means that
the bimodule problem $(\K,\M)$, and equivalently the Drozd bocs
$\mathfrak{A}$, corresponding to $\Lambda$ is also wild and
homogeneous. The homogeneous property is inherited to any induced
structure given by reductions. Then we are restricted to dealing
with one of those configurations listed in the wild theorem 5.1
obtained at some stage of reductions. For some easy cases, using
the method proposed by Bautista in [B], we can  directly construct
infinitely many iso-classes of non-homogeneous modules of a fixed
dimension. For the most difficult case, we will construct a
non-homogeneous parameterized matrix based on the bipartite
property of $(\K,\M)$. And thereby we show that $(\K,\M)$, as well
as $\Lambda$, has no homogeneous property, which is a
contradiction to the assumption.

\subsection{Notations.}

\kg In this paper $\rK$ represents an algebraically closed field.
By an algebra we mean an associative finite-dimensional
$\rK$-algebra with an identity, which we may assume (without loss
of generality) to be basic and connected. By a module over an
algebra $\Lambda$ we mean a finite-dimensional left
$\Lambda$-module. And by $\Lambda$-mod we denote the category of
finite-dimensional left $\Lambda$-modules.

 Let
$Q=(Q_0, Q_1)$ be a quiver, where $Q_0$ is the set of vertices and
$Q_1$ is the set of arrows. For any arrow $\az: i\rightarrow j$ in
$Q_1$,  $s(\az)$ stands for the starting vertex of $\az$ and
$e(\az)$ for the ending one. Let $p$ be any path in $Q$, $s(p)$
and $e(p)$ stand also for the starting and ending vertices of $p$
respectively. We write the composition of paths from left to
right, i.e.
$$ p\cdot q=\left\{\begin{array}{cl} pq & \mbox{if } e(p)=s(q),\\
0 & \mbox{if } e(p)\neq s(q).
\end{array}\right.$$
 According to a theorem due to Gabriel, for any basic finite-dimensional
 algebra $\Lambda$ there exists some
quiver $Q$ and some admissible ideal $I$ of $\rK Q$, such that
$\Lambda\cong \rK Q/I$.

Let us denote the set of natural numbers by $\mathbb{N}$, the set
of $m\times n$ matrices over $\rK$ by $M_{m\times n}(\rK)$, or
 $M_n(\rK)$ whenever $m=n$. We denote by $E_{ij}$ the matrix with
the $(i,j)$-entry $1$ and others zero,
  and by $I_n$ the $n\times n$ identity matrix.

 We assume that all the categories and functors under the
 consideration are $\rK$-linear.

\newpage

\begin{center} \section{Bimodule problems} \end{center}

\kg This section is devoted to describing a special class of
bimodule problems of the form $(\K,\M,H)$ listed in the third set
of subsection 1.4. We will give the definition of $(\K,\M,H)$,
calculate the triangular bases of $\M$ and $\rad \K$, define Weyr
matrices, and give the reductions for a triple of $(\K,\M,H)$.

\bigskip
\subsection{The category $P_1(\Lambda)$}

\kg Let $\Lambda$ be a finite-dimensional  $\rK$-algebra, which is
basic,  and $\Lambda$-proj be the full subcategory of
$\Lambda$-mod consisting of projective $\Lambda$-modules. Let us
recall the category
$$
P(\Lambda)=\{P_1\stackrel{\alpha}{\lra}P_0\mid P_1,P_0\in
\Lambda\mbox{-proj}, \alpha\in\Hom_{\Lambda}(P_1, P_0)\},
$$
whose objects are morphisms between $\Lambda$-projectives. If
$(P_1\stackrel{\alpha}{\lra}P_0),
(Q_1\stackrel{\beta}{\lra}Q_0)\in P(\Lambda)$, the morphisms from
$(P_1\stackrel{\alpha}{\lra}P_0)$ to $
(Q_1\stackrel{\beta}{\lra}Q_0)$ are pairs $(f_1, f_0)$ with
$f_i\in
 {\rm Hom}_{\Lambda}(P_i, Q_i)$, $i=1,0$, such that the following diagram
 commutes:
$$
\begin{CD} P_1 @>\alpha>> P_0\\
@Vf_1VV   @VVf_0V\\
 Q_1 @>\beta>> Q_0
\end{CD}
$$
The compositions and additions of the morphisms are given
componentwise.

 There are a full subcategory
$$
P_1(\Lambda)=\{(P_1\stackrel{\alpha}{\lra}P_0)| \ \ \Im(\alpha)
\subseteq \rad(P_0) \}\subset P(\Lambda),
$$
and a full subcategory
$$
P_2(\Lambda)=\{(P_1\stackrel{\alpha}{\lra}P_0)\in P_1(\Lambda)| \
\ \Ker(\alpha) \subseteq \rad(P_1) \}\subset P_1(\Lambda).
$$
It is well known that the functor $Cok: P_2(\Lambda)\rightarrow
\Lambda$-mod is a representation equivalence, i.e. it is dense,
full and reflects isomorphisms.

The category $P_1(\Lambda)$ can be described explicitly in terms
of matrices. Let $J$ be the Jacobson radical of $\Lambda$ with
$J^m\neq 0$, $J^{m+1}=0$, and $\Lambda/J\cong \rK
e_1\times\cdots\times\rK e_s$ with $\{e_1,e_2,\cdots,e_s\}$ being
a complete set of orthogonal primitive idempotents of $\Lambda$.
Let us first  take  a $\rK$-basis of $e_iJ^me_j$, $1\le i, j\le
s$. Suppose we have already had a basis of $J^l$, and then extend
it to a $\rK$-basis
 of $J^{l-1}$, such that the images of the added new elements
 form a $\rK$-basis of $e_i(J^{l-1}/J^l)e_j$,
for $1\le i, j\le s$. The obtained $\rK$-basis of $J$ is denoted
by
$$B = \{\zeta_1,\zeta_2,\cdots,\zeta_{t-s} \}$$
 with a fixed order given above, where $t=\dim\Lambda$. Moreover, by
adding the primitive idempotents $e_1,e_2,\ldots,e_s$ after
$\zeta_{t-s}$, we obtain an ordered $\rK$-basis of $\Lambda$.

Given any $a\in \Lambda$, we   define a map $\bar{a}:
\Lambda\rightarrow \Lambda$ given by $\bar{a}(b)=ab$, which is a
$\rK$-linear transformation of the $\rK$-vector space $\Lambda$.
Then $\zeta_i$ and $e_j$ will correspond to some $t\times t$
matrices $\bar{\zeta}_i$ and $\bar{e}_j$ under the ordered basis
respectively, which are upper triangular because of the order of
the basis. Thus  a matrix algebra ${\widetilde \Lambda}$ follows,
which is  upper triangular and  isomorphic to $\Lambda$, with a
$\rK$-basis $\{\bar{\zeta}_1, \cdots, \bar{\zeta}_{t-s},
\bar{e}_1, \cdots, \bar{e}_s\}$. And $\widetilde{\Lambda}$ is
usually called a {\it left regular representation} of $\Lambda$.

{\bf Example 1.} Let $Q= $ {\unitlength=1mm
\begin{picture}(20,8) \put(10,1){\circle*{0.6}}
\put(7,1){\circle{4}} \put(13,1){\circle{4}}
\put(10,4){\circle{4}} \put(9,2){\vector(0,1){0}}
\put(11,2){\vector(1,0){0}} \put(11,2){\vector(0,1){0}}
\put(2,0){$a$} \put(9,6){$b$} \put(16,0){$c$}
\end{picture}}
be a quiver, $\rK Q$ be the path algebra of $Q$,
 $$I=\langle a^2,
b^2, c^2, ab,\\  ba, bc,cb,ca,ac\rangle$$ be an ideal of $\rK Q$
generated by the elements in the bracket, and $\Lambda=\rK Q/I$.
We denote
 the residue classes of $\{e,a,b,c\}$ in
$\Lambda$ still by $\{e,a,b,c\}$ respectively. Then an ordered
$\rK$-basis $\{ c,b,a,e\}$ of $\Lambda$ yields
\begin{eqnarray*}
\widetilde{\Lambda}&=&\left\{\left(
                          \begin{array}{cccc}
                          s_1 & 0 & 0  &s_4  \\
                        & s_1 & 0 &  s_3  \\
                        & & s_1  & s_2  \\
                        & &  & s_1
                          \end{array}
                      \right)  \, \mid
\forall s_i\in\rK\right\}, \\
\rad\widetilde{\Lambda}&=&\left\{ \left(
\begin{array}{cccc}
0 & 0 & 0  & m_c  \\
  & 0 & 0  & m_b  \\
  &   & 0  & m_a  \\
  &   &   & 0
\end{array}
\right)\mid \forall m_a,m_b,m_c\in\rK\right\}.
\end{eqnarray*}

{\bf Example 2.} Let $Q=$ {\unitlength=1.2mm \special{em:linewidth
0.4pt} \linethickness{0.4pt}
\begin{picture}(32.00, 8.00)
\put(5.00,2.00){\circle*{1.00}} \put(16.50,2.00){\circle*{1.00}}
\put(28.50,2.00){\circle*{1.00}}
\put(16.00,2.50){\vector(-1,0){10}}
\put(16.00,1.50){\vector(-1,0){10}}
\put(17.00,2.00){\vector(1,0){10}}

\put(28.50,0.00){\makebox(0,0)[cc]{$3$}}
\put(16.50,0.00){\makebox(0,0)[cc]{$2$}}
\put(5.00,0.00){\makebox(0,0)[cc]{$1$}}
\put(10.00,4.00){\makebox(0,0)[cc]{$a$}}
\put(10.00,0.00){\makebox(0,0)[cc]{$b$}}
\put(22.00,3.00){\makebox(0,0)[cc]{$c$}}
\end{picture}}
be a quiver and $\Lambda=\rK Q$ be the path algebra  of $Q$. If we
choose an ordered $\rK$-basis
   $\{ c,b,a,e_1,e_2,e_3\}$ of $\Lambda$, then
 $$
 \widetilde{\Lambda}=\left\{\left(
                          \begin{array}{cccccc}
                          s_2 & 0 & 0 & 0  & 0& s_6  \\
                        & s_2 & 0 & s_5& 0 & 0 \\
                        & & s_2 & s_4 & 0 &0 \\
                        & &  & s_1 & 0 & 0 \\
                        & & &  & s_2 & 0 \\
                        & & & & & s_3
                          \end{array}
                      \right)  \mid
\forall s_i\in\rK\right\}; $$
$$
\rad\widetilde{\Lambda}= \left\{\left(
\begin{array}{cccccc}
0 & 0 & 0 & 0 & 0 & m_c \\
  & 0 & 0 & m_b &0 & 0 \\
  &   & 0 & m_a & 0 & 0 \\
  &   &   & 0 & 0 & 0 \\
  &   &   &   & 0 & 0 \\
  &   &   &   &   & 0
\end{array}
\right)\mid \forall m_a,m_b,m_c\in\rK\right\} .
$$

{\bf Example 3.} \cite{R1, D} Let $Q= $ {\unitlength=1mm
\begin{picture}(20,4) \put(10,1){\circle*{0.6}}
\put(7,1){\circle{4}} \put(13,1){\circle{4}}
\put(9,2){\vector(0,1){0}} \put(11,2){\vector(0,1){0}}
\put(2,0){$a$} \put(16,0){$b$}
\end{picture}}
be a quiver, $I=\langle a^2,ba-\alpha ab, ab^2 , b^3\rangle$ be an
ideal of $\rK Q$ with $\alpha\in \rK$ being a fixed non-zero
constant, and let  algebra $\Lambda=\rK Q/I$. Denote  the residue
classes of $e, a, b$ in $\Lambda$ still by $e, a,b $ respectively.
Moreover let us set $ c=b^2, d=ab$. Then an ordered $\rK$-basis
$\{d, c, b, a, e\}$ of $\Lambda$ yields
$$
\widetilde{\Lambda}=\left\{\left(
\begin{array}{ccccc}
s_1 &  0  & s_2 &  \alpha s_3  & s_5 \\
       & s_1    & s_3    & 0   & s_4 \\
       &        & s_1    & 0   & s_3 \\
       &        &        & s_1 & s_2 \\
       &        &        &     & s_1
\end{array}
\right)\mid \forall s_i\in\rK\right\},
$$
$$
\rad\widetilde{\Lambda}=\left\{\left(
\begin{array}{ccccc}
0 & 0 &  m_a    & \az m_b & m_d \\
  & 0 &  m_b    & 0 & m_c \\
  &   &   0   & 0 & m_b \\
  &   &       & 0 & m_a \\
  &   &       &   & 0
\end{array}
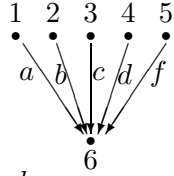
\right)\mid \forall m_a,m_b,m_c,m_d\in\rK\right\}.
$$

\medskip
{\bf Example 4.} Let $Q=$ {\unitlength=1mm
\begin{picture}(23.00, 20.00)
\put(1.00,15.00){\circle*{1.00}} \put(6.00,15.00){\circle*{1.00}}
\put(11.00,15.00){\circle*{1.00}}
\put(16.00,15.00){\circle*{1.00}}
\put(21.00,15.00){\circle*{1.00}} \put(11.00,1.00){\circle*{1.00}}

\put(2.00,14.00){\vector(2,-3){8}}
\put(6.50,14.00){\vector(1,-3){4}}
\put(11.00,14.00){\vector(0,-1){12}}
\put(16.00,14.00){\vector(-1,-3){4}}
\put(21.00,14.00){\vector(-2,-3){8}}

\put(1.00,18.00){\makebox(0,0)[cc]{$1$}}
\put(6.00,18.00){\makebox(0,0)[cc]{$2$}}
\put(11.00,18.00){\makebox(0,0)[cc]{$3$}}
\put(16.00,18.00){\makebox(0,0)[cc]{$4$}}
\put(21.00,18.00){\makebox(0,0)[cc]{$5$}}
\put(11.00,-1.50){\makebox(0,0)[cc]{$6$}}

\put(2.50,10.00){\makebox(0,0)[cc]{$a$}}
\put(7.0,10.00){\makebox(0,0)[cc]{$b$}}
\put(12.0,10.00){\makebox(0,0)[cc]{$c$}}
\put(15.50,10.00){\makebox(0,0)[cc]{$d$}}
\put(20.00,10.00){\makebox(0,0)[cc]{$f$}}

\end{picture}}
be a quiver and $\Lambda=\rK Q$ be the  path algebra of $Q$. Then
an $\rK$-basis
   $\{f,d, c,b,a,e_1,e_2,e_3, e_4, e_5, e_6\}$ of $\Lambda$ yields
 $$
 \widetilde{\Lambda}=\left\{\left(
                          \begin{array}{ccccccccccc}
                    s_5 & 0   & 0   & 0   & 0 & 0 & 0 & 0 & 0 & 0 & s_{11}  \\
                        & s_4 & 0   & 0   & 0 & 0 & 0 & 0 & 0 & 0 & s_{10} \\
                        &   & s_3 & 0   & 0 & 0 & 0 & 0 & 0 & 0 & s_9 \\
                        &   &     & s_2 & 0 & 0 & 0 & 0 & 0 & 0 & s_8\\
                        &   &     &     & s_1 & 0 & 0 & 0 & 0 & 0 & s_7 \\
                        &   &      &    &     & s_1& 0 & 0 & 0 & 0 &0\\
                        &   &      &    &     &  & s_2 & 0 & 0 & 0 &0\\
                        &   &      &    &     &  &     & s_3 & 0 & 0 &0\\
                        &   &      &    &     &  &     &     & s_4 & 0 &0\\
                        &   &      &    &     &  &  &  &  & s_5 &0\\
                        &   &      &    &     & &  &  &  &  &s_6
                                                   \end{array}
                      \right)  \mid
\forall s_i\in\rK\right\}; $$
$$
\rad\widetilde{\Lambda}= \left\{\left(
\begin{array}{ccccccccccc}
  0 & 0   & 0   & 0   & 0 & 0 & 0 & 0 & 0 & 0 & m_f  \\
    & 0   & 0   & 0   & 0 & 0 & 0 & 0 & 0 & 0 & m_d \\
    &     & 0   & 0   & 0 & 0 & 0 & 0 & 0 & 0 & m_c \\
    &     &     & 0   & 0 & 0 & 0 & 0 & 0 & 0 & m_b\\
    &     &      &    & 0 & 0 & 0 & 0 & 0 & 0 & m_a \\
    &     &      &    &   & 0 & 0 & 0 & 0 & 0 & 0  \\
    &     &      &    &   &   & 0 & 0 & 0 & 0 & 0\\
    &     &      &    &   &   &   & 0 & 0 & 0 & 0 \\
    &     &      &    &   &   &  &    & 0 & 0 & 0 \\
    &     &      &    &   &   &  &    &   & 0 & 0 \\
    &     &      &    &   &   &  &    &   &    & 0
\end{array}
\right)\mid \forall m_a,m_b,m_c,m_d,m_f\in\rK\right\} .
$$

\subsection{The bimodule problem  $(\K,\M, H)$}

\kg {\bf Definition 2.2.0} [CB2] By a bimodule problem we mean a
triple $(K,M,d)$, where $K$ is a $k$-algebra, $M$ is a
$K$-$K$-bimodule, and $d:K\rightarrow M$ is a derivation.
\hfill$\square$

{\bf Proposition 2.2.0} (given by Crawley-Boevey) Let $(K,M,d)$ be
a bimodule problem, then there exist

(1) a positive integer $t$;

(2) an upper triangular matrix algebra $\K \subset M_t(k)$ and  an
algebra isomorphism $\varphi: K \rightarrow \K$;

(3) a $\K$-$\K$-bimodule $\M \subset M_t(\rK)$, and  an
isomorphism $\psi:  M \rightarrow \M$, such that $\psi
(xm)=\varphi(x) \psi (m)$, $\psi (mx)=\psi (m) \varphi (x)$ for
any $x \in K, m \in M$;

(4) a fixed matrix $H \in M_t(\rK)$,  such that $\psi (d(x))=
\varphi (x) H-H \varphi (x)$ for any $x \in M$.

{\bf Proof.} Let $N=M \oplus K$ as vector spaces. We define a
$K$-$K$-bimodule structure on $N$ given by $y(m,x)=(ym,yx)$, and
$(m,x)y=(my-xd(y),xy)$ for any $m \in M, x,y \in K$. If we choose
$h=(0,1) \in N$, then
$$xh-hx=x(0,1)-(0,1)x=(0,x)-(-d(x),x)=(d(x),0) \in M.$$
Let $A=K \oplus N$ be an algebra with multiplication
$(x,n)(x',n')=(xx',xn'+nx')$. It is clear that $M \subset N
\subset A$ and $K \subset A$. Since $N^2=0$, the idea $N \subset
rad(A)$. Thus any complete set of orthogonal primitive idempotents
of $K$ is also that of $A$. Suppose that $dim(A)=t$, choose a
suitable $k$-basis of $rad(A)=rand(K)\oplus N$ according to the
method given in 2.1, we can embed $A$ as a upper triangular
subalgebra ${\widetilde{A}} \subset M_t(k)$, i.e there is an
algebra isomorphism $\theta: A \rightarrow \widetilde{A}$. Denote
by $\varphi$ the restriction of $\theta$ on $K$, by $\psi$ that on
$M$. And let $\K=\varphi (K)$, $\M=\psi (M)$, $H=\theta (h)$, the
proof is completed. \hfill$\square$

Proposition 2.2.0 suggests the following definition of bimodule
problems in terms of matrices  in order to do some calculations.

\kg {\bf Definition 2.2.1} A bimodule problem $(\K,\M, H)$ given
by matrices consists of the following datum:
\begin{itemize}
\item[\textrm{I}.] A linearly ordered set of integers
$T=\{1,2,\cdots,t\}$ and an equivalent relation $\sim$ on $T$,
where
$$T/\!\sim \,=\{\I_1,\I_2,\cdots,\I_s\}.$$

\item[\textrm{II}.] An upper triangular matrix algebra
$\K=\{(s_{ij})_{t\times t}
 \in M_{t}(k)\}$, where $s_{ii}=s_{jj}$ when  $i\sim j$; $
 s_{ij}$, when
  $i<j$, satisfy the following system of homogeneous linear
  equations  in indeterminates $x_{ij}$:
 $$
\sum_{\I\ni i<j\in\J}c_{ij}^{l} x_{ij}=0,
$$
where  the coefficients $c_{ij}^{l}\in\rK$, the equations are
indexed by $1\le l\le q_{\I\J}$ for some $q_{\I\J}\in \mathbb{N}$,
and for each pair $$(\I,\J)\in (T/\!\sim)\times (T/\!\sim).$$

\item[\textrm{III}.] A $\K$-$\K$-bimodule $\M=\{(m_{ij})_{t\times
t}\in M_{t}(k)\}$, where $m_{ij}$ satisfy the following
 system of homogeneous linear equations in indeterminates $z_{ij}$:
$$
 \sum_{(i,j)\in\I\times\J}d_{ij}^{l} z_{ij}=0,
$$
where  the coefficients $d_{ij}^{l}\in\rK$, the equations are
indexed by $1\le l\le q'_{\I\J}$ for some $q'_{\I\J}\in
\mathbb{N}$, and for each pair
$$(\I,\J)\in (T/\!\sim)\times
(T/\!\sim).$$

\item[\textrm{IV}.] A fixed matrix $H=(h_{ij})_{t\times t}\in
M_{t}(k)$, where $h_{ij}=0$ when $i\nsim j$. And a derivation $d:
\K\rightarrow \M$ is given by
$$
d(S)=SH-HS
$$
for any $S\in \K$.  \hfill$\square$
\end{itemize}

It is obvious that $(\K,\M,H)$ is a bimodule problem in the sense
of Definition 2.2.0. We stress that  definition 2.2.1 is a
modification of definition 1.1 of \cite{S}.  Sometimes we  write
$(\K,\M)$ instead of $(\K,\M, H)$ for short.

{\bf Corollary 2.2.1} $\{E_{\I}=\sum_{i\in \I}E_{ii}\mid \forall \
\I\in T/\!\sim\}$
 is a complete set of primitive idempotents of $\K$. And
 $d(E_{\I})=0$, $\forall\; \I\in T/\sim$.\hfill$\square$

Sometimes we write $s_{\I}=s_{ii}$ for any $i\in \I$.

{\bf Proposition 2.2.1} Let $\Lambda$ be a $\rK$-algebra of
dimension t having s pairwise orthogonal primitive idempotents,
let $\widetilde{\Lambda}$ be the left regular representation of
$\Lambda$ defined in 2.1. Set
$$
{\cal K} =
  \left(
        \begin{array}{cc}
         {\widetilde \Lambda}&0\\ 0&{\widetilde \Lambda}
        \end{array}
  \right), \ \
{\cal M} =
  \left(
        \begin{array}{cc}
        0&\rad {\widetilde \Lambda}\\ 0&0
        \end{array}
  \right), \ H=(0).
$$ and
$$T=\{1,2,\cdots,t; t+1, t+2, \cdots, 2t\},$$ is the set of row (or
column) indices of $\K$,
$$T/\!\sim=\{\I_1,\cdots,\I_s; \I_{s+1},\cdots,\I_{2s}\},$$  where
$ \I_l=\{i | s_{ii}=s_{i_l i_l}\ \mbox {for a fixed}\ i_l \in
T\}.$ Then $(\K, \M, H)$ is a bimodule problem. \hfill$\square$

The rest of this subsection is devoted to illustrating the
representation category $Mat(\K,$ $\M)$ of a bimodule problem
$(\K, \M, H)$, whose objects and morphisms are both given by
matrices. The category is coincide with that given in \cite{CB2}.

 A vector $\n=(n_1,n_2,\cdots,n_t)\in \mathbb{N}^t$ is called a {\it size vector}
  of $(T,\sim)$ provided that $n_i=n_j$ when $i\sim j$ in $T$.  And
$n=\sum\limits_{i=1}^tn_i$ is called the {\it size of $\n$}. On
the other hand,  if we denote $n_i$ by $n_{\I}$ $\forall\ i\in
\I$, then $\underline{d}=(n_{\I_1}, n_{\I_2},\cdots,n_{\I_s})$ is
called a {\it dimension vector} of $(T, \sim)$, and
$d=\sum\limits_{\I\in T/\!\sim}n_{\I}$ is called the {\it
dimension of $\underline{d}$}.

 Let $\n$ be any size vector,
 $$
 \M_{\n}=
 \{ (N_{ij})_{t\times t}\mid N_{ij} \mbox{ are } n_i\times n_j\
\mbox{ blocks satisfying the  equation system \textrm{III} } \}.$$
Then any partitioned matrix
  $N\in \M_{\n}$ is called a
 {\it representation}, or a {\it matrix} over  $(\K,\M)$
of  size vector $\n$.

 If $\m$
is also  a size vector of $(T,\sim)$, $M\in \M_{\m}$. Let
$$
\K_{\m\times\n} =\{(S_{ij})_{t\times t}\mid S_{ij} \mbox{ are }
m_i\times n_j \mbox{ blocks }\},
$$
where $S_{ij}=0$ if $i>j$; $ S_{ii}=S_{jj}$ if $i\sim j$; $S_{ij}$
satisfy the equation system {II} if $i<j$. Then a partitioned
matrix $S\in \K_{\m\times\n}$ is called a {\it morphism} from $M$
to $N$, provided that
$$
MS-SN=d(S),
$$
or equivalently
$$
(M+H_{\m})S=S(N+H_{\n}),
$$
where if $H=(h_{ij})_{t\times t}$, then
$H_{\n}=(H_{ij})_{\n\times\n}$ is a partitioned matrix with
 $H_{ij}=0$ if $i\not\sim j$, and
$H_{ij}=h_{ij}I_{n_i}$ if  $i\sim j$. $H_{\m}$ is defined
similarly for the size vector $\m$. We denote the partitioned
diagonal part of $S$ by $S_0$.

A morphism $S: M\rightarrow N$ is called an {\it  isomorphism}
provided that $S$ is invertible. In this case, $M$ and $N$ are
said to be {\it isomorphic}, and  denoted by $M\simeq N$.

{\bf Corollary 2.2.2} (1) If $S\in \K_{\m\times\n}$, and $S'\in
\K_{\n\times\m}$ such that $SS'=I$, $S'S=I$, then $S_0S'_0=I$,
$S'_0S_0=I$.

(2) If $M\simeq N$, then $\m=\n$.

{\bf Proof.} (1) The reason is that $S$ and $S'$ are both
partitioned upper triangular.

(2) If $S: M\rightarrow N$ and $S': N\rightarrow M$ are morphisms,
such that $SS'=id_{M}$, and $S'S=id_{N}$, then
$S_{\I}S'_{\I}=I_{m_{\I}}$ and $S'_{\I}S_{\I}=I_{n_{\I}}$ for any
$\I\in T/\!\!\sim$. Thus $m_{\I}=n_{\I}$, consequently $\m=\n$ as
desired. \hfill$\square$

Finally we fix some notations of various indices used in this
paper by an example. The equivalent classes in $T/\!\sim$ are
denoted by italic English letters $\I,\J,\P,\Q,\cdots$; the
positive integers in $T$ by lower-case English $i,j,p,q,\cdots$;
the usual row and column indices of a matrix by lower-case Greek
$\alpha, \beta, \gamma,\cdots$. And we write $\alpha\in i,
\beta\in j$ whenever the $(\alpha,\beta)$-entry is sitting at the
$(i,j)$-block.

{\bf Example.} See Example 2 of 2.1 and Proposition 2.2.1. In the
corresponding bimodule problem of $\Lambda$, $T=\{1,2,3, \cdots,
12\}$,
$$\I_1=\{1,2,3,5\},\ \I_2=\{4\},\ \I_3=\{6\},\ \I_4=\{7,8,9,11\},\
\I_5=\{10\},\ \I_6=\{12\}.$$ Let
$$\m=(2,2,2,0,2,0;  0,0,0,2,0,1)$$
 be a size vector of
$(T,\sim)$. Let us take $$i_l=l,\ \  l=1,2,\cdots, 12;\ \ {\rm and}\ \ \alpha_j=j, \ \ j=1, \cdots, 11.$$ Then
$\alpha_1, \alpha_2\in i_1$; $\alpha_3, \alpha_4\in i_2$; $\alpha_5, \alpha_6\in i_3$; $\alpha_7, \alpha_8\in
i_5$; $\alpha_9, \alpha_{10}\in i_{10}$; $\alpha_{11}\in i_{12}$. The following is an $11\times 11$ matrix $M\in
{\rm Mat}(\K,\M)$ of size vector $\m$, and a morphism $S\in {\rm End}_{(\K,\M)}(M)$.
\begin{center} \unitlength=0.8mm
\begin{picture}(120,60)
\put(0,0){\framebox(55,55)} \mput(10,0)(10,0){5}{\line(0,1){55}}
\mput(0,45)(0,-10){5}{\line(1,0){55}} \put(41,26){$0$}
\put(46,26){$0$} \put(41,31){$0$} \put(46,31){$1$}
\put(41,36){$1$} \put(46,36){$0$} \put(41,41){$0$}
\put(46,41){$\lambda$} \put(52,46){$0$} \put(52,51){$1$}
\put(-10,27){$M\!\!=$} \put(56,0){,}

\thicklines\put(40,0){\line(0,1){55}} \put(0,15){\line(1,0){55}}
\put(105,0){\line(0,1){55}} \put(65,15){\line(1,0){55}}

\thinlines \put(58,27){$S\!\!=$} \put(65,0){\framebox(55,55)}
\mput(75,0)(10,0){5}{\line(0,1){55}}
\mput(65,45)(0,-10){5}{\line(1,0){55}}

\put(66,51){$ s_{1}$} \put(72,51){$0$} \put(71,46){$ s_{2}$}
\put(67,46){$0$} \put(76,41){$ s_{1}$} \put(82,41){$0$}
\put(81,36){$ s_{2}$} \put(77,36){$0$} \put(86,31){$ s_{1}$}
\put(92,31){$0$} \put(91,26){$ s_{2}$} \put(87,26){$0$}
\put(96,21){$ s_{1}$} \put(102,21){$0$} \put(101,16){$ s_{2}$}
\put(97,16){$0$} \put(106,11){$ s_{2}$} \put(112,11){$0$}
\put(111,6){$ s_{1}$} \put(107,6){$0$} \put(116,1){$ s_{1}$}
\put(121,49){$i_1$} \put(121,39){$i_2$} \put(121,29){$i_3$}
\put(121,19){$i_5$} \put(121,9){$i_{10}$} \put(121,1){$i_{12}$}
\end{picture}
\end{center}

\subsection{The triangular basis of $(\K,\M)$}

\kg  Definition 2.2.1 shows that $\M$ and rad$\K$  are solution
spaces of the equation systems  \textrm{III} and \textrm{II}
respectively.  Some nice triangular bases of the spaces will be
chosen in this subsection.

{\bf Definition 2.3.1} Let $M=(m_{ij})_{t\times t}$ be a matrix.
An order on the  indices is defined as follows:  $(i,j)\prec
(i',j')$ provided that $i>i'$ or $i=i'$, $j<j'$. The order is also
valid on the indices of blocks of a $t\times t$ partitioned
matrix.

In fact there are many possibilities to define an order on the indices of a $t\times t$ matrix. For example,
$(i,j)\prec(i',j')$ if $j<j'$, or $j=j'$, $i>i'$. Another example: $(i,j)\prec(i',j')$ if $$i-j>i'-j',\ \ {\rm
or}\ \ i-j=i'-j',\ \ i>i'.$$ The principle is to ensure that (the index of  $m_{ij}s_{jj'})\succ (\mbox{ that of
}m_{ij}$) when $j<j'$, and ( the index of $s_{i'i}m_{ij})\succ (\mbox{ that of } m_{ij}$) when $i'<i$ for any
entry $m_{ij}$ of $M$ in $\M$, and $s_{jj'}$, $s_{i'i}$ of $S$ in $\K$. In this
 paper  we  use mainly the order given by Definition 2.3.1,
unless otherwise stated.

{\bf Lemma 2.3.1} Let $\sum_{j=1}^n d_{ij}z_j=0, i=1,2,\cdots, m$,
be a system of homogeneous linear equations. Then there exists a
unique choice of free indeterminates $z_{p_1}, \cdots, z_{p_r}$
such that
$$
\left\{\begin {array}
{ccccl} z_{j_1}&&&&=\sum_{l=1}^ra_{1l}z_{p_l}\\ &z_{j_2}&&&=\sum_{l=1}^ra_{2l}z_{p_l}\\
&&\ddots&&\\ &&&z_{j_{n-r}}&=\sum_{l=1}^ra_{n-r,l}z_{p_l}
\end{array} \right.
$$
where $p_1<p_2<\cdots<p_r$, $j_1<j_2<\cdots<j_{n-r}$,  and $p_l <
j_k$, whenever $a_{kl}\neq 0$ for each $1\leqslant k\leqslant
n-r$.

{\bf Proof.} $n=1$ is trivial. Suppose that the assertion is true
for $(n-1)$ indeterminates, and we are in the case of $n$
indeterminates. If $z_1=0$, then replacing $z_1$ by $0$, we obtain
a new system with $(n-1)$ indeterminates. If $z_1\neq 0$, then
$z_1$ is taken as a free indeterminate. Suppose that
$z_1,z_2,\cdots z_{j_1-1}$ are free, but $z_{j_1}=\sum
_{l=1}^{j_1-1}a_{1l}z_l$, then substituting $z_{j_1}$ for $\sum
_{l=1}^{j_1-1}a_{1l}z_l$ we obtain a new system with $(n-1) $
indeterminates. Thus the assertion follows by induction.
\hfill$\square$

The lemma  allows us to choose some nice basis for $\M$ (resp.
$rad \K$) according to the order given in Definition 2.3.1. Let
$$z_{p_{\I\J}^1,q_{\I\J}^1}, \cdots , z_{p_{\scriptscriptstyle
\I\J}^{r_{\scriptscriptstyle  \I\J}},q_{\scriptscriptstyle \I\J}^{r_{\scriptscriptstyle \I\J}}}$$ be all the
free indeterminates of the equation system \textrm{III} of Definition 2.2.1 for any fixed $(\I,\J)$, such that
$(p^1_{\I\J}, q^1_{\I\J})\prec\cdots\prec (p^{r_{\scriptscriptstyle \I\J}}_{\scriptscriptstyle
\I\J},q^{r_{\scriptscriptstyle \I\J}}_{\scriptscriptstyle \I\J} )$ and
$$
z_{ij}= \sum _{w=1}^{r_{\I\J}}a_{\I\J}^w(i,j)
z_{p_{\I\J}^wq_{\I\J}^w},
$$
where $(p^w_{\I\J},q^w_{\I\J})\prec(i,j)$ whenever
$a_{\I\J}^w(i,j)\neq 0 $, for $i\in \I$, $j\in \J$.  Define a set
\begin{equation}\begin{array}{l}\displaystyle
A_{\I\J}=\{\rho_{\I\J}^w=E_{p_{\I\J}^w,q_{\I\J}^w}+\sum_{(i,j)\in
\I\times \J}a_{\I\J}^w(i,j) E_{ij}\mid 1\leqslant w\leqslant
r_{\I\J} \}, \mbox{ and } \\[+2ex]\displaystyle A=\bigcup_{\I,\J\in
T/\!\sim}A_{\I\J}.\end{array}
\end{equation}
Then $A_{\I\J}$  is  a  $\rK$-basis of the solution space
$E_{\I}\M E_{\J}$  of   equation system \textrm{III} of Definition
2.2.1 for any fixed pair $(\I,\J)$. Moreover $A$ is a $\rK$-basis
of $\M$.

\medskip
Similarly we can choose a $\rK$-basis $B_{\I\J}$ of the solution
space $E_{\I}(rad \K)E_{\J}$ of the equation system  \textrm{II}
of Definition 2.2.1 for any fixed pair $(\I,\J)$, and obtain a
$\rK$-basis of $rad\K$:
\begin{equation}\begin{array}{l}\displaystyle
B_{\I\J}=\{\zeta_{\I\J}^w=E_{\bar{p}_{\I\J}^w\bar{q}_{\I\J}^w}+\sum
_{\I\ni i<j\in \J}b_{\I\J}^w(i,j) E_{ij} \mid 1\leqslant w
\leqslant \bar r_{\I\J}\}, \mbox{ and }\\[+2ex]\displaystyle
B=\bigcup_{\I,\J\in T/\!\sim}B_{\I\J}.\end{array}
\end{equation}

The pair $(A,B)$ is said to be a {\it triangular basis}  of
bimodule problem $(\K,\M)$.

\subsection{An order on triangular basis}

\kg {\bf Definition 2.4.1} There is a natural linear order on $A:$
$\rho_{\I\J}^w \prec \rho_{\I'\J'}^{w'}$ provided $(p_{\I\J}^w,
q_{\I\J}^w)\prec(p_{\I'\J'}^{w'}, q_{\I'\J'}^{w'})$. Let us take
$$(p,q)=(p_{\P\Q}^{w_0}, q_{\P\Q}^{w_0})=\min \{(p_{\I\J}^w,
q_{\I\J}^w) \mid w=1, \cdots , r_{\I\J},  \forall  (\I,\J)\in
(T/\!\!\sim)\times (T/\!\!\sim)\},$$
 and $\rho =\rho_{\P\Q}^{w_0}$,
then $p\in \P, q\in \Q$. Similarly there is also a linear order on
$B$.

{\bf Proposition 2.4.1} Let $(A,B)$ be the triangular basis of a
bimodule problem $(\K,\M,H)$. If the left and right module actions
are
$$
\left\{\begin{array}{l} \zeta_{\I\L}^u\rho_{\L\J}^v=\sum\limits_w
~_l\c\rho_{\I\J}^w, \\[8mm]
\rho_{\I\L}^u\zeta_{\L\J}^v=\sum\limits_w ~_r\c\rho_{\I\J}^w
\end{array}\right.
$$
respectively, with the structure constants $~_lc, ~_rc\in \rK$,
then $\rho_{\I\J}^w\succ \rho_{\L\J}^v$ for the first case, and
$\rho_{\I\J}^w\succ \rho_{\I\L}^u$ for the second case.

{\bf Proof.} We will only prove the first one, and the second one
can be obtained similarly.
$$
\zeta_{\I\L}^u \rho_{\L\J}^v=
(E_{\bar{p}_{\I\L}^u\bar{q}_{\I\L}^u} + \sum_{\I\ni i<l\in \L}
b_{\I\L}^u(i,l) E_{il}) (E_{p_{\L\J}^vq_{\L\J}^v} +
\sum_{(l',j)\in \L\times\J} a_{\L\J}^v(l',j) E_{l'j})
$$
 by Formulae (1) and (2) in 2.3. We may write
$$\zeta_{\I\L}^u=\sum\limits_{\I\ni i<l\in \L} b_{\I\L}^u(i,l)
E_{il}\ \  {\rm with}\ \ b^u_{\I\L}(\overline{p}_{\I\L}^u,\overline{q}_{\I\L}^u)=1$$ and
$$\rho_{\L\J}^v=\sum\limits_{(l',j)\in \L\times\J} a_{\L\J}^v(l',j) E_{l'j}\ \ {\rm with}\ \
a^v_{\L\J}(p_{\L\J}^v,q_{\L\J}^v)=1$$ for simplicity. Then
$$
E_{il}E_{l'j}=\left\{\begin{array}{cl} E_{ij} & \mbox{when }
l=l',\\ 0& \mbox{otherwise,} \end{array}\right.
$$
 since $\zeta_{\I\L}^u$ is an upper triangular nilpotent
 matrix, $i<l$.
  Thus $(i,j)\succ (l',j)\succeq (p_{\L\J}^v, q_{\L\J}^v)$, i.e.
the indices of $E_{ij}$ in the expression of  $\rho_{\I\J}^w$ in
the first formula are all greater than $(p_{\L\J}^v, q_{\L\J}^v)$.
Therefore $\rho_{\I\J}^w\succ\rho_{\L\J}^v$ as required.
\hfill$\square$

Let $C,D$ be two partitioned matrices. We write
$$C\equiv_{\prec(p,q)}D\quad (\mbox{resp. } C\equiv_{\preceq(p,q)}D)$$ if
 $(i,j)$-blocks of $C,D$ are
the same for all $(i,j)\prec(p,q)$ (resp. $(i,j)\preceq(p,q)$).

{\bf Corollary 2.4.1 } Let $(\K,\M,H)$ be a bimodule problem, with
a pair $(p,q)$ given by Definition 2.4.1. Then
$$
SH\equiv_{\prec(p,q)}HS
$$
for any $S\in \K$.

{\bf Proof.} Since $d(S)=SH-HS\in\M$ by  IV of Definition 2.2.1,
and all the entries, which are smaller than $(p,q)$ in a matrix of
$\M$, equal zero, so that $SH-HS\equiv_{\prec(p,q)}0,$ i.e. $
SH\equiv_{\prec(p,q)}HS $. \hfill$\square$

\medskip
 Let $D=Hom _{\rK}(-,\rK)$ be the usual dual functor. Let
 $\M^*=D\M$, and
 $(\rad\K)^*=D(\rad\K)$ be the dual space of $\M$ and $\rad\K$
 respectively. Then we write their dual bases respectively by
\begin{equation}
A^*=\{(\rho_{\I\J}^w)^* \mid 1\leqslant w\leqslant r_{\I\J},
\forall\ (\I,\J) \in (T/\!\sim)\times (T/\!\sim)\} \label{(3)}
\end{equation}
\begin{equation}
B^*=\{(\zeta_{\I\J}^w)^* \mid 1\leqslant w\leqslant
\overline{r}_{\I\J},\forall\ (\I,\J) \in
(T/\!\sim)\times(T/\!\sim)\} \label{(4)}
\end{equation}
thus $(\rho_{\I\J}^w)^*$  can be regarded as the coefficient
functions of $\rho_{\I\J}^w$ in $\M$ and $(\zeta_{\I\J}^w)^*$ as
those of $\zeta_{\I\J}^w$ in $\rad\K$, which  yield a dual
structure of the bimodule problem  (see 3.6 below). The linear
order on $A$ is transferred to that on $A^*$, namely,
$(\rho_{\I\J}^w)^*\prec (\rho_{\I'\J'}^{w'})^*$ provided
$\rho_{\I\J}^w\prec\rho_{\I'\J'}^{w'}$. And $B^*$ is also ordered
such that $(\zeta_{\I\J}^w)^*\prec(\zeta_{\I'\J'}^{w'})^*$
provided $\zeta_{\I\J}^w\prec\zeta_{\I'\J'}^{w'}$.

{\bf Example.} See Example 3 of  2.1, and Proposition 2.2.1. In
the corresponding bimodule problem of $\Lambda$,
$$\rho^1=E_{4, 10}+E_{18}=a,\ \ \rho^2=E_{3, 10}+E_{28}+\alpha
E_{19}-b,\ \ \rho^3=E_{2, 10}=c,\ \ \rho^4=E_{1, 10}=d;$$ and
$$(\rho^1)^*=a^*,\ \ (\rho^2)^*=b^*,\ \ (\rho^3)^*=c^*,\ \
(\rho^4)^*=d^*,$$ if the $k$-basis of $rad(\Lambda)$, $a,b,c,d$
still stand for the $k$-basis of $rad(\widetilde {\Lambda})$, as
well as$ \left(
        \begin{array}{cc}
        0&\rad {\widetilde \Lambda}\\ 0&0
        \end{array}
  \right)$  with $a\prec b\prec
c\prec d$, and $a^*, b^*, c^*, d^*$ are dual basis.

\subsection{Weyr matrices}

\kg Let $$ J(\lambda)=J_d(\lambda)^{e_d}\oplus J_{d-1}(\lambda)^{e_{d-1}}\oplus \cdots\oplus J_1(\lambda)^{e_1}
$$ be a direct sum of Jordan blocks with  a common  eigenvalue $\lambda$, where $e_j$ stands for the number of
the summands of $J_j(\lambda)$, write
$$m_j=e_d+e_{d-1}+\cdots+e_j$$
for $j=1,2,\cdots, d$. Then $m_1\geq m_2\geq\cdots\geq m_d$, and
$e_j=m_j-m_{j+1}$, $m_{d+1}=0$. The following partitioned matrix
$W_{\lambda}$ is called a {\it Weyr matrix of eigenvalue
$\lambda$}:
$$
W_{\lz}=\left(
\begin{array}{cccccc}
\lz I_{m_1} & W_{12} &0&  \cdots & 0&0 \\
  & \lz I_{m_2} & W_{23}&\cdots &0&0 \\
  && \lz I_{m_3} &\ddots & 0 &0\\
&&&\ddots&\ddots&\vdots\\
  & & & & \lz I_{m_{d-1}}& W_{d-1,d} \\
 & & & && \lz I_{m_d}
\end{array}
\right)_{d\times d}, \ \ W_{j,j+1}=\left(
\begin{array}{c}
I_{m_{j+1}}\\ 0
\end{array}
\right)_{m_j\times m_{j+1}}.
$$

{\bf Lemma 2.5.1} $W_{\lz}$ is obtained from $J(\lambda)$ by a
series of elementary transformations of exchanging rows and
columns simultaneously. \hfill$\square$

\medskip
The following is called a vector of $W_{\lambda}$:
\begin{equation}
\begin{array}{rcccccc}
\underline{e}=(e_d\ \  & e_{d-1}& e_{d-2}&\cdots & e_3 & e_2 & e_1\\
               e_d\ \  & e_{d-1}& e_{d-2}&\cdots & e_3 & e_2 & \\
                   & \cdots &        & \cdots &  &  & \\
               e_d\ \  & e_{d-1}& e_{d-2}&        &  &  & \\
                e_d\ \  & e_{d-1}&       &        &  &  & \\
                 e_d). & & &        &  &  &
\end{array} \label{(5)}
\end{equation}

For example, given a Jordan form $J_4(\lambda)^2\oplus
J_2(\lambda)^3$, we have $m_1=2+0+3+0=5$, $m_2=2+0+3=5$,
$m_3=2+0=2$, $m_4=2$,
$$
\begin{array}{llccc}
\underline{e}=(\!\!\!\!\!&2\  & 0 & 3& 0\\
&2\ & 0& 3& \\ &2\ & 0& & \\
& 2),& & &
\end{array}\qquad
W_{\lambda}=\left(\begin{array}{c} { \unitlength 0.6mm
\begin{picture}(50,50)

 \put(0,35){\dashbox{1}(15,15){$\lz I_5$}}
 \put(15,35){\dashbox{1}(15,15){$ I_5$}}
 \put(15,20){\dashbox{1}(15,15){$\lz I_5$}}
\put(30,20){\dashbox{1}(10,15){$\begin{array}{c} I_2 \\
0 \end{array}$}}
 \put(30,10){\dashbox{1}(10,10){$\lz I_2$}}
 \put(40,10){\dashbox{1}(10,10){$ I_2$}}
 \put(40,0){\dashbox{1}(10,10){$\lz I_2$}}
\end{picture}}
\end{array}
\right).
$$

Let $\lambda^1, \lambda^2, \cdots, \lambda^{\alpha}$ be  a set of
eigenvalues which are pairwise different, then
$$
W=\mbox{diag} (W_{\lz^1}, W_{\lz^2}, \cdots, W_{\lz^{\alpha}})
$$
is called a {\it Weyr matrix} with a vector $$\e=(\e^1, \e^2,\cdots, \e^{\alpha}),$$ where $\e^l$ is a vector of
$W_{\lambda^l}$ defined in Formula (5). From now on we will fix an order on the base field $\rK$, such that
$$\lz^1<\lz^2<\cdots<\lz^{\alpha}.$$ For example, if $\rK$ is the field of complex numbers, we may use the
lexicographic order: $a+bi<c+di$ provided $a<c$ or $a=c$ but $b<d$.

{\bf Corollary 2.5.1} Any square matrix over $\rK$ is similar to a
unique Weyr matrix under a fixed order of $\rK$. \hfill$\square$

Let $W$ be a Weyr matrix, $X$ be a square matrix having the same
size as $W$. Then an equation $WX= XW$ gives
$$ X=diag(X_1, \cdots,X_l, \cdots, X_{\alpha}).$$
 Denote $X_l$ by
$Y$, and $d_l$ by $d$, then $Y=(Y_{ij})_{d\times d}$, where
$Y_{ij}=0$ when $i>j$, and
$$
Y_{ij}=\left(\begin{array}{cccc}
 U_{11}^h&U_{12}^h&\cdots&U_{1j'}^h\\
 &\cdots&&\cdots\\
U_{h1}^h&U_{h2}^h&\cdots&U_{hj'}^h\\
 &U_{h+1,2}^h&\cdots&U_{h+1,j'}^h\\
  &&\ddots&\vdots\\
&&&U_{i'j'}^h
\end{array}\right)_{i'\times j'}
$$
when $i\leqslant j$, where $h=j-i+1$, $i'=d-i+1$, $j'=d-j+1$,
 and $U_{gr}^h$ are $e_{g'}\times e_{r'}$ free matrices, where
$g'=d-g+1$, $r'=d-r+1$. For example, if  $W=W_{\lambda}$ is given
above, then
$$
X=\left(\begin{array}{c}
 {\unitlength 1.1mm
\begin{picture}(50,50)

 \put(0,35){\dashbox{1}(15,15){$\begin{array}{cc} U_{11}^1\!\! &\!\! U_{13}^1 \\
 &\!\! U_{33}^1\end{array}$}}
 \put(15,35){\dashbox{1}(15,15){$\begin{array}{cc} U_{11}^2\!\! &\!\! U_{13}^2 \\
 0&\!\! U_{33}^2\end{array} $}}
\put(30,35){\dashbox{1}(10,15){$\begin{array}{c} U_{11}^3 \\
 U_{31}^3\end{array} $}}
\put(40,35){\dashbox{1}(10,15){$\begin{array}{c} U_{11}^4 \\
 U_{31}^4\end{array} $}}

 \put(15,20){\dashbox{1}(15,15){$\begin{array}{cc} U_{11}^1\!\! &\!\! U_{13}^1 \\
 &\!\! U_{33}^1\end{array}$}}
\put(30,20){\dashbox{1}(10,15){$\begin{array}{c} U_{11}^2 \\
 0\end{array} $}}
\put(40,20){\dashbox{1}(10,15){$\begin{array}{c} U_{11}^3 \\
 U_{31}^3\end{array} $}}

 \put(30,10){\dashbox{1}(10,10){$U_{11}^1$}}
 \put(40,10){\dashbox{1}(10,10){$ U_{11}^2$}}
 \put(40,0){\dashbox{1}(10,10){$U_{11}^1$}}
\end{picture}}
\end{array}
\right)_{4\times 4}.
$$

\subsection{Reductions}

\kg Let $(\K,\M,H)$ be a bimodule problem, $(A,B)$ be the triangular basis, and $(A^*, B^*)$  the dual basis.
Define two matrices
$$ R_0=\sum_{\I}(E_\I)^*\otimes_k
E_{\I}+\sum_{(\I,\J)}\sum_{w=1}^{\bar{r}_{\I\J}}(\zeta^w_{\I\J})^*
\otimes_k \zeta_{\I\J}^w , $$
$$
 N_0=\sum_{(\I,\J)}\sum_{w=1}^{r_{\I\J}} (\rho_{\I\J}^w)^*\otimes
\rho_{\I\J}^w.
$$
 where  the tensor product is the usual tensor product of matrices,
  i.e if $$C=(c_{ij})_{p\times q},\ \  D=(d_{ij})_{t\times s},\ \ {\rm then}\ \ C\otimes_k D=
  (C\cdot d_{ij})_{t\times s}.$$ And we regard the coefficient functions $(E_{\I})^*$,
  $(\zeta_{\I\J}^w)^*$, $(\rho_{\I\J}^w)^*$ as indeterminates. Clearly,
  if the values of
  $$((E_{\I})^*, (\zeta_{\I\J}^w)^*\mid \forall\,\I, 1\leq w\leq \overline{r}_{\I\J},
  (\I,\J)\in T/\sim\times T/\sim)$$ range over $k^{\dim\K}$, we
  obtain the algebra $\K$;
and if the values of $$((\rho_{\I\J}^w)^*\mid 1\leq w\leq
r_{\I\J}, (\I,\J)\in T/\sim\times T/\sim)$$ run over $k^{\dim\M}$,
we obtain the bimodule $\M$. The $(p_{\I\J}^w, q_{\I\J}^w)$-th
elements of $N_0$ are said to be {\it free}, and the others are
said to be {\it dependent} (see Formula (1) of 2.3).

Next we fix a size vector $\n$ of $(T, \sim)$, consider the
partitioned upper triangular matrix algebra $\K_{\n\times\n}$,
bimodule $\M_{\n}$ and matrix $H_{\n}$. it is obvious that
$\K_{\n\times\n}$ is Morita equivalent to $\K$ if $n_{\I}\neq 0$,
$\forall\; \I\in T/\!\sim$; and $\K_{\n\times\n}$ is not basic if
there exists some $\I\in T/\!\sim$, such that $n_{\I}>1$. Let
\begin{equation}\begin{array}{l}\displaystyle
R=\sum_{\I}(\overline{E_\I})^*\otimes
E_{\I}+\sum_{(\I,\J)}\sum_{w=1}^{\bar{r}_{\I\J}}(\overline{\zeta^w_{\I\J}})^*
\otimes \zeta_{\I\J}^w ,\\[+2ex]
\displaystyle N=\sum_{(\I,\J)}\sum_{w=1}^{r_{\I\J}}
(\overline{\rho_{\I\J}^w})^*\otimes \rho_{\I\J}^w,\end{array}
\label{(6)}
\end{equation}
where $(\overline{E_\I})^*$, $(\overline{\zeta^w_{\I\J}})^*$ and
$(\overline{\rho_{\I\J}^w})^*$ are $n_{\I}\times n_{\I}$,
$n_{\I}\times n_{\J}$ and $n_{\I}\times n_{\J}$ matrices with all
the entries being indeterminates  respectively. It is clear that
if the values of $(\overline{E_\I})^*$ and
$(\overline{\zeta^w_{\I\J}})^*$ range over $$\prod_{\I\in
T/\!\sim} M_{n_{\I}}(k)\times \prod_{\I, \J\in T/\!\sim}
M_{n_{\I}\times n_{\J}}(k)^{\overline{r}_{\I\J}},$$ we obtain
$\K_{\n\times\n}$; and if the values of
$(\overline{\rho_{\I\J}^w})^*$ run over $$\prod_{\I, \J\in
T/\!\sim} M_{n_{\I}\times n_{\J}}(k)^{r_{\I\J}},$$ we obtain
$\M_{\n}$. And $(N,R,\n)$ is called a {\it triple of bimodule
problem } $(\K,\M,H)$. The $(p_{\I\J}^w, q_{\I\J}^w)$-th blocks of
$N$ are said to be {\it free}, and the others are said to be {\it
dependent} (see Formula (1) of 2.3).

 The matrix $N_0$ may be regarded as a ``universal matrix" of $\M$,
 and $R_0$ a ``universal matrix" of $\K$ or a
 ``universal endomorphism" of
$N_0$. Similarly, $N$ may be regarded as a ``universal matrix" of
$\M_{\n}$, and $R$ a ``universal matrix" of $\K_{\n}$ or a
``universal endomorphism" of $N$. Consider the matrix equation
$$(N+H_{\n})R=R({N}+H_{\n}).$$ Let $(p,q)$ be as in  Definition 2.4.1,
then $N_{pq}$ is the first
non-zero block of $N$.
 Then Corollary 2.4.1 ensures that
$ H_{\n}R\equiv_{\prec(p,q)} RH_{\n} $, i.e. the blocks before
$(p,q)$ are really equal on both sides.  And the $(p,q)$-block of
the equation is
$$\begin{aligned}
&&(H_{\n})_{p1}X_{1q}+\cdots+(H_{\n})_{p,q-1}X_{q-1,q}+
(H_{\n})_{pq}X_{qq} +N_{pq}X_{qq}\\ &=&X_{pp}N_{pq} +
X_{pp}(H_{\n})_{pq}
+X_{p,p+1}(H_{\n})_{p+1,q}+\cdots+X_{pt}(H_{\n})_{tq}.
\end{aligned}$$
 Since $(H_{\n})_{ij}=0$  when $i \nsim j$,
$(H_{\n})_{ij}=h_{ij}I_{n_i}$ when $i\sim j$ by \textrm{IV} of
Definition 2.2.1, and $(H_{\n})_{pq}X_{qq}= X_{pp}(H_{\n})_{pq}$
in both cases $p\nsim q$ and $p\sim q$,   the above equation is
equivalent to the following:
\begin{equation}
N_{pq}X_{qq}-X_{pp}N_{pq}=\sum_{l=p+1}^tX_{pl}h_{lq}-
\sum_{l=1}^{q-1}h_{pl}X_{lq}.
\end{equation}
Let us define a new equation
\begin{equation}
\sum_{l=p+1}^tx_{pl}h_{lq}- \sum_{l=1}^{q-1}h_{pl}x_{lq}=0.
\end{equation}
And define a constant matrix $\overline{N}_{pq}$ with the same
size of $N_{pq}$ in the following $3$ cases \cite{D, CB1, S}:

\medskip
{\bf Regularization.} If the equation system  II of Definition
2.2.1 at $(\P,\Q)$ do not imply equation (8), then let
$\overline{N}_{pq}=0 $. And  denote zero here by $\emptyset$, i.e.
$\overline{N}_{pq}=\emptyset$, to distinguish $0$'s coming from
regularization or from other cases.

If the equation system  II of Definition 2.2.1 at $(\P,\Q)$
implies equation (8), then
\begin{equation}
N_{pq}X_{qq}=X_{pp}N_{pq}. \label{8}
\end{equation}

{\bf Edge reduction.} When $ p\nsim q$, let
$$ \overline{N}_{pq}=\left(\begin{array}{cc} 0 &
I_r \\
0 & 0\end{array}\right)_{n_p\times n_q}$$ for some $r\leq
\min\{n_p,n_q\}$.

 {\bf Loop reduction.} When $ p\sim q$, then $X_{pp}=X_{qq}$, let
$$ \overline{N}_{pq}=W$$
 for some Weyr matrix $W$.

 \medskip
 As soon as $\overline{N}_{pq}$ has been fixed, we obtain a new
 triple $(N',R',\n')$ as follows.
We first define a new vector $\n'$ according to the $3$ cases
respectively. Regularization: $\n'=\n$; Edge reduction: $\n'$ is
obtained from $\n$ by splitting each component $n_i$, $\forall\;
i\in\P$,  into $(n_{i_1}, n_{i_2})$ when  $n_{i_1}=r$,
$n_{i_2}=n_{\P}-r$; and $n_j$, $\forall j\in\Q$, into $(n_{j_1},
n_{j_2})$ when $n_{j_1}=n_{\Q}-r$, $n_{j_2}=r$; Loop reduction:
$\n'$ is obtained from $\n$ by splitting each component $n_i$,
$\forall\; i\in\P$,  into $\underline{e}=(\underline{e}^1, \cdots
\underline{e}^l, \cdots \underline{e}^{\alpha})$ given in 2.5. Let
\begin{equation}
\left\{\begin{array}{l}
(H_{\n})'=H_{\n}+\overline{N}_{pq}\otimes \rho,\\
N'=N-N_{pq}\otimes \rho,\\
R'\mbox{ is a restriction of }R \mbox{ given by }
R'(H_{\n})'\equiv_{\preceq (p,q)}(H_{\n})'R'
\end{array}\right.\label{(10)}
\end{equation}
Then $R'$ satisfies equation system  \textrm{II} of Definition
2.2.1 and equation (7), if we replace $N_{pq}$ by
$\overline{N}_{pq}$. $R'$ determines a partitioned upper
triangular matrix algebra $(\K_{\n})'$, which is not basic in
general. Similarly $N'$ determines a
$(\K_{\n})'$-$(\K_{\n})'$-bimodule $(\M_{\n})'$. The procedure of
regularization, edge reduction and loop reduction is called a {\it
reduction for the triple $(N,R,\n)$}.

We will construct an induced bimodule problem $(\K',\M',H')$ from
$(\K,\M,H)$ given by each reduction for $(N, R, \n)$, such that
$\K'$ is basic in 3.1. we will also show in 3.1 that
$$(\K_{\n})'=\K'_{\n'},\ \ (\M_{\n})'=\M'_{\n'}\ \ {\rm and}\ \
(H_{\n})'=H'_{\n'}.$$ Thus $(N', R', \n')$ is a triple of $(\K',
\M', H')$.

 \newpage


\newpage
\bcen\section{Induced bimodule problems}\ecen

In this section we will describe the induced bimodule problem
$(\K',\M',H')$ given by one of three reductions respectively, and
give a nice functor $\vartheta: Mat(\K',\M') \rightarrow
Mat(\K,\M)$. In the end, we will present an explicit
correspondence between bimodule problems and their dual structure
--- bocses. In particular we show the differentials of the arrows of a
bocs in terms of the multiplication of matrices.

\medskip

\subsection{Induced bimodule problems}

\kg Now we  construct an induced bimodule problem $(\K',\M',H')$
from the bimodule problem $(\K,\M,H)$ after one of the three
reductions given in 2.6 respectively.

{\bf Regularization.} \textrm{I}. Let  $T'=T$, $\sim'=\sim$. Then
$\n'$ is a size vector of $(T',\sim')$.

\textrm{II}. $\K'=\{S'\in \K\mid S'\mbox{ satisfies the equation
(8) of 2.6}\}$.

\textrm{III}. $\M'=\{M'\in\M\mid M' \mbox{ satisfies the equation
} z_{pq}=0\}$, which is a $\K'$-$\K'$-bimodule by Proposition
2.4.1.

\textrm{IV}. $H'=H+\emptyset\otimes\rho$ which implies
$H'_{\n'}=(H_{\n})'$  and $(H'_{\n'})_{pq}=\overline{N}_{pq}$.

{\bf Edge reduction.} The equation at the $(p,q)$-block of $R
'(H_{\n})'\equiv_{\preceq(p,q)}(H_{\n})'R'$ is as follows:
 $$\left(\begin{array}{cc} 0 & I_r \\
0 & 0 \end{array}\right)X'_{qq} =X'_{pp}\left(\begin{array}{cc} 0 & I_r \\
0 & 0 \end{array}\right),$$ which yields
$$
X'_{pp}=\left(\begin{array}{cc}
X_{p_1p_1} & X_{p_1p_2} \\
0        & X_{p_2p_2}
\end{array}\right), \quad
X'_{qq}=\left(\begin{array}{cc}
X_{q_1q_1} & X_{q_1q_2} \\
0        & X_{q_2q_2}
\end{array}\right),
$$
and  $X_{p_1p_1}=X_{q_2q_2}$.

\textrm{I}. Define a size vector $\wt{\n}\in (T,\sim)$ by
$\wt{n}_{\I}=1$, $\forall\ \I\in T/\!\sim \setminus \{\P, \Q\}$;
and (i) $\wt{n}_{\P}=\wt{n}_{\Q}=2$ when $\min\{n_{\P},
n_{\Q}\}>r>0$; (ii) $\wt{n}_{\P}=2$, $\wt{n}_{\Q}=1$ when
$n_{\P}>r$, $n_{\Q}=r>0$; (iii) $\wt{n}_{\P}=1$, $\wt{n}_{\Q}=2$
when $n_{\P}=r>0$, $n_{\Q}>r$;  (iv) $\wt{n}_{\P}=\wt{n}_{\Q}=1$
when $n_{\P}=n_{\Q}=r>0$;  (v)  $\wt{n}_{\P}=\wt{n}_{\Q}=1$ when
$r=0$. Let $T'$ be obtained from $T$ by splitting
 $i$, $\forall i\in \P$, into $i_1,i_2$ in cases (i) and (ii), into
 $i_1$ in cases (iii) and (iv), into $i_2$ in case (v),  and splitting
  $j$, $\forall j\in \Q$, into $j_1,j_2$ in cases (i) and (iii), into
 $j_2$ in cases (ii) and (iv), into $j_1$ in case (v).
 Define
 $$
 (\P\cup \Q)'=\{i_1,j_2\mid \forall i\in\P,\ \forall j\in \Q\},$$
 $$\P'=\{i_2\mid \forall i\in \P\}\ \mbox{and}\  \Q'=\{j_1\mid \forall j\in\Q\}.$$
Then
$$
T'/\sim'=(T/\!\sim\setminus\{\P,\Q\})\cup\{\P',\Q',(\P\cup\Q)'\}.
$$
$\n'$ is a size vector of $(T',\sim')$.

\textrm{II}.  $\K'=\{S'\in \K_{\wt{\n}\times\wt{\n}}\mid S'\mbox{
satisfies the following additional equations}\}$:
$$
x_{i_1i_2}-x_{p_1p_2}=0,\  \forall i\in\P  \mbox{  and  }
x_{j_1j_2}-x_{q_1q_2}=0,\quad \forall j\in\Q.
$$

\textrm{III}. $\M'=\{M'\in\M_{\wt{\n}}\mid M' \mbox{ satisfies }
z_{pq}=0\} $, which is a $\K'$-$\K'$-bimodule by Proposition
2.4.1.

\textrm{IV}. $H'=H_{\wt{\n}}+ U\otimes\rho$, where $$U =\left(\begin{array}{cc} 0&1\\ 0&0 \end{array} \right),\
\ \left(\begin{array}{c} 1\\ 0
\end{array} \right),\ \ (0 \quad 1 ),\ \  (1),\ \ (0)$$ according to cases
(i)--(v) of $\wt{\n}$ respectively. $H'_{\n'}=(H_{\n})'$ and
$(H'_{\n'})_{pq}=\overline{N}_{pq}$.

\medskip

{\bf Loop reduction.}  The equation
$$WX'_{qq}=X'_{pp}W$$
gives a new index set. Namely  we define a  linearly ordered set
of indices according to $\underline{e}=(\underline{e}^1, \cdots
\underline{e}^l, \cdots \underline{e}^{\alpha})$ in order to
define $T'$:

$$
\begin{array}{llllllll}
\underline{\delta}^{il}=(\!\!&\delta_{1,d}^{il}&\delta_{1,d-1}^{il}&\delta_{1,d-2}^{il}
&\cdots&\delta_{1,3}^{il}&\delta_{1,2}^{il}&\delta_{1,1}^{il}\\
&\delta_{2,d}^{il}&\delta_{2,d-1}^{il}&\delta_{2,d-2}^{il}
&\cdots&\delta_{2,3}^{il}&\delta_{2,2}^{il}&\\
&&\cdots&&\cdots& & & \\
&\delta_{d-2,d}^{il}&\delta_ {d-2,d-1}^{il}&\delta_ {d-2,d-2}^{il}&&&&\\
&\delta_{d-1,d}^{il}&\delta_ {d-1,d-1}^{il}&&&&&\\
&\delta_{d,d}^{il}\quad)&&&&&&
\end{array}
$$
for each $i\in \P$, $1\leq l\leq \alpha$. For example, the row and
column indices of the blocks of $X$ given in example of 2.5 can be
expressed by the elements of $\delta$ in the table below:

 \ \
\begin{center}\unitlength=1mm

\begin{picture}(80,40)
\multiput(0,0)(0,5){7}{\line(1,0){80}} \put(0,0){\line(0,1){40}}
\multiput(20,0)(10,0){7}{\line(0,1){40}}
\put(0,40){\line(1,0){80}} \put(0,35){\line(4,-1){20}}
\put(20,30){\line(-4,3){13}} \put(1,31){\tiny row} \put(5,35){$h$}
\put(11,37){{\tiny column}} \put(8,26){$\delta_{14}$}
\put(8,21){$\delta_{12}$} \put(8,16){$\delta_{24}$}
\put(8,11){$\delta_{22}$} \put(8,6){$\delta_{34}$}
\put(8,1){$\delta_{44}$} \put(23,33){$\delta_{14}$}
\put(33,33){$\delta_{12}$} \put(43,33){$\delta_{24}$}
\put(53,33){$\delta_{22}$} \put(63,33){$\delta_{34}$}
\put(73,33){$\delta_{44}$} \put(24,26){$1$} \put(34,26){$1$}
\put(44,26){$2$} \put(54,26){$2$} \put(64,26){$3$}
\put(74,26){$4$} \put(34,21){$1$}\put(44,21){$2$} \put(54,21){$2$}
\put(64,21){$3$} \put(74,21){$4$} \put(44,16){$1$}
\put(54,16){$1$} \put(64,16){$2$} \put(74,16){$3$}
 \put(54,11){$1$} \put(64,11){$2$}\put(74,11){$3$}
\put(64,6){$1$}\put(74,6){$2$} \put(74,1){$1$} \thicklines
\put(20,20){\line(1,0){60}} \put(40,10){\line(1,0){40}}
\put(60,5){\line(1,0){20}} \put(40,10){\line(0,1){20}}
\put(60,5){\line(0,1){25}} \put(70,0){\line(0,1){30}} \thinlines
\end{picture}
\end{center}

 \textrm{I}. $T'$ is obtained from $T$ by splitting
$i\in \P$ into
$$
\{\delta_{gr}^{il}\mid e_r^l\ne 0, 0\leq g\leq r; r=d_l, d_l-1,
\cdots, 2,1; l=1,2,\cdots, \alpha\}.
$$
where $e_r^l$ is the $(d_l-r+1)$-th component of
$\underline{e}^l$. And
$$
T'/\!\!\!\sim'=(T/\!\!\!\sim \setminus \P)
 \cup \{ \P^l_r \mid 1\leq r\leq d_l, l=1,2,\cdots,\alpha\},
$$
 where $\P_r^l=\{\delta_{gr}^{il}\in T'\mid
 1\leq g\leq r, \forall\ i\in\P\}$. Then $\n'$ is a size vector of
 $(T', \sim')$.

 Moreover we define a size vector $\wt{\n}$ of $(T, \sim)$ with
$\wt{n}_{\I}=1$, $\forall \ \I\in T/\!\sim \setminus \{\P\}$, and
$$
\wt{n}_{\P}=\#\{\delta_{gr}^{pl}\in T'\}=\sum_{l=1}^\alpha
\sum_{r=1}^{d_l}\delta_r^l\cdot r,
$$
where $\delta_r^l=1$, if $e_r^l\neq 0$, and $\delta_r^l=0$, if
$e_r^l=0$.

 \textrm{II}.  $\K'=\{S'\in \K_{\wt{\n}\times\wt{\n}}\mid S'\mbox{
 satisfies
the following additional equations}\}$:
$$
x_{\delta_{g'_1r}^{il}\delta_{g'_2s}^{il}}-x_{\delta_{g_{\scriptscriptstyle
1}r}^{pl} \delta_{g_{\scriptscriptstyle 2}s}^{pl}}=0,
$$
 where  $g'_2-g'_1=g_2-g_1=h-1$, $s-r<h$,
 $g_1,g'_1\leq r$, $g_2, g'_2\leq s$, $1\leq r,s\leq d_l$,
 $l=1,2,\cdots,\alpha$, $\forall \ i\in \P$; and
$$
x_{\delta_{g_1r}^{il}\delta_{g_2s}^{il}}=0
$$
where $g_2-g_1=h-1$, $s-r\geq h$, $g_1\leq r$, $g_2\leq s$, $1\leq
r,s\leq d_l$, $l=1,2,\cdots,\alpha$, $\forall\ i\in\P$.

\textrm{III}. $\M'=\{M'\in\M_{\wt{\n}}\mid M' \mbox{ satisfies }
z_{pq}=0\} $ which is a $\K'$-$\K'$-bimodule by  Proposition
2.4.1.

 \textrm{IV}. $H'=H_{\wt{\n}}+\wt{W}\otimes \rho$ with $ \wt{W}$
 being a Weyr matrix similar to
$\bigoplus\limits_{l=1}^{\alpha}\bigoplus\limits_{r=d_l}^1
J_{r}(\lambda_l)^{\delta_{r}^l}$, where $\delta_{r}^l$ are as in
the end of  \textrm{I}. Thus $H'_{\n'}=(H_{\n})'$ and
$(H'_{\n'})_{pq}=\overline{N}_{pq}$.

{\bf Corollary 3.1.1}\ $n'_{\P^l_r}=e^l_{r}$  where $ 1\le r\le
d_l,\ l=1,\cdots,\alpha.$\hfill $\Box$

\medskip
{\bf Proposition 3.1.1}. In the above 3 procedures, $H'$ satisfies \textrm{IV} of Definition 2.2.1. Let
$d'(S')=S'H'-H'S'$ for any $S'\in \K'$, we obtain a derivation from $\K'$ to $\M'$. Thus $(\K', \M', H')$ is a
bimodule problem, and $(N', R', \underline{n}')$ is a triple of $(\K', \M', H' )$.

{\bf Proof.} $S'H'\equiv_{\preceq (p,q)}H'S'$ implies that
$d'(S')\in \M'$. \hfill$\square$

\medskip
The procedure of regularization, edge reduction or loop reduction
given above is called a {\it reduction of bimodule problem
$(\K,\M,H)$}, and $(\K', M',H')$ is called an {\it induced
bimodule problem} of $(\K,\M,H)$. There is a map
$$(\vartheta_1,\vartheta_0): (\K',\M')\rightarrow
(\K_{\wt{\n}\times\wt{\n}}, \M_{\wt{\n}}),$$ such that the algebra
homomorphism $\vartheta_1: \K'\rightarrow
\K_{\wt{\n}\times\wt{\n}}$, and the $\K'$-$\K'$-bimodule
homomorphism $\vartheta_0:\M'\rightarrow\M_{\wt{\n}}$ are both
natural embeddings.

\subsection{Deletions and size changes}

\kg In this  subsection we first define the 4-th reduction
algorithm, so called {\it deletion}. Let $(\K, \M, H)$ be a
bimodule problem,  $\wt{\n}$ a size vector of  $(T, \sim)$ with
$\wt{n}_{\I}=1$ or $ 0$.  Let
$$T'=\{i\mid \wt{n}_i=1\},\
\sim'=\sim|_{T'}, \ \K'=\K_{\wt{\n}\times\wt{\n}}, \
\M'=\M_{\wt{\n}}, \ H'=H_{\wt{\n}}.$$ Then we obtain an induced
bimodule problem $(\K',\M',H')$ from $(\K,\M,H)$ by {\it deleting
a subset of equivalent classes} $\{\I\mid \wt{n}_{\I}=0 \}$ from
 $T/\sim$. And we have a a map $$(\vartheta_1,\vartheta_0):
(\K',\M')\rightarrow (\K, \M),$$ such that the algebra homomorphism $\vartheta_1: \K'\rightarrow \K$, and the
$\K'$-$\K'$-bimodule homomorphism $\vartheta_0: \M'\rightarrow\M$ are both natural embeddings.

The following observations are straightforward.

(1) Given any size vector $\n$ of $(T, \sim)$, we define a size vector $\wt{\n}$ according to $\n$ with
$\wt{n}_{\I}=1$ when $n_{\I}\neq 0$, and $\wt{n}_{\I}=0$ when $n_{\I}= 0$. Then we obtain an induced bimodule
problem $(\K',\M',H')$. Let $\n'\in (T', \sim')$ with $n'_i=n_i$, $\forall\ i\in T'$, and a triple $(N', R',
\n')$ over $(\K',\M',H')$ be given by Formula (6) of 2.6. If $(N, R, \n)$ is a triple of $(\K,\M)$, then
 $$N'=N, \ R'=R, \ H'=H,$$
  i.e. the two triples are essentially the same.

(2) A size vector $\n$ of $(T,\sim)$ is said to be {\it sincere}
if $n_{\I}\neq 0$ for all $\I\in T/\!\!\sim$. In  most
considerations, a size vector is usually assumed to be sincere by
a deletion if necessary.

 (3) If $\n$ is sincere, then after doing anyone of the $3$ reductions given in 3.1, the induced size vector
$\n'$ is still sincere.

\medskip Next we make some size changes of triples over the induced
bimodule problem.

 {\bf  Proposition 3.2.1} Let $(\K,\M,H)$ be a bimodule
problem, $(N,R,\n)$ be a triple of $(\K,\M,H)$,  and $(N',R',\n')$
be induced from $(N,R,\n)$ by one of the $3$ reductions in 3.1.
Then for any size vector $\hat{\n}'$ of $(T',\sim')$, there exist
a size vector $\hat{\n}$ of $ (T,\sim)$ thus a triple $(\hat{N},
\hat{R}, \hat{\n})$ of $(\K,\M,H)$, such that

(1) If $\hat{n}_{\P}=0$ or $\hat{n}_{\Q}=0$, then $(\hat{N}',
\hat{R}', \hat{\n}')$ is obtained from $(\hat{N}, \hat{R},
\hat{\n})$ by  identity.

(2) If $\hat{n}_{\P}\neq 0$ and $\hat{n}_{\Q}\neq 0$, then the
equation system  \textrm{II} of Definition 2.2.1 and equation (8)
of 2.6 for $(\K,\M,H)$ at the pair $(\P,\Q)$  determine the same
type of reductions from $(N,R,\n)$ to $(N',R',\n')$ and from
$(\hat{N}, \hat{R}, \hat{\n})$ to $(\hat{N}', \hat{R}',
\hat{\n}')$, where
$$\overline{\widehat{N}}_{pq}= (H'_{\hat{\n}'})_{pq}=\left\{\begin{array}{lll}\emptyset,
\quad & \mbox{ in \  regularization};\\[1mm]
\left(\begin{array}{cc} 0&\hat I\\ 0&0 \end{array} \right),&\mbox{ in \  edge \ reduction};\\[4mm]
\hat W, & \mbox{ in \  loop\ reduction}.\end{array}\right. $$

{\bf Proof.} The size vector $\hat{\n}$ is defined from
$\hat{\n}'$ according to the $3$ reductions respectively.
\begin{equation}
\left\{\begin{array}{ll}
\mbox{regularization:} & \hat{\n}=\hat{\n}'; \\
\mbox{edge reduction:}&   \hat{n}_{\I}=\hat{n}'_{\I}, \forall\
\I\in T/\!\!\sim\setminus\{\P,\Q\},\\
& \hat{n}_{\P}=\hat{n}'_{\P'}+\hat{n}'_{(\P\cup \Q)'},
\hat{n}_{\Q}=\hat{n}'_{(\P\cup \Q)'}+\hat{n}'_{\Q'};\\
\mbox{ Loop reduction:}& \hat{n}_{\I}=\hat{n}'_{\I}, \forall\
\I\in T/\!\sim \setminus \{\P\},
\hat{n}_{\P}=\sum_{l,r}r\hat{n}'_{\P_r^l}.
\end{array}\right.
\end{equation}
Let  $C=\{\I\in T/\!\sim \mid \hat{n}_{\I}=0\}$.  When $\P\in C$,
or $\Q\in C$, we set $\hat{\n}=\hat{\n}'$, then
$$\hat{H}_{\hat{\n}}=\hat{H}'_{\hat{\n}'},\ \hat{N}=\hat{N}',\
\hat{R}=\hat{R}'.$$ Thus the procedure from $(\hat{N}, \hat{R},
\hat{\n})$ to $(\hat{N}', \hat{R}', \hat{\n}')$  is  an identity.
When  $\P\notin C$ and $\Q \notin C$, the equation system II of
Definition 2.2.1 and (8) of 2.6 for $(\K,\M,H)$ at the pair
$(\P,\Q)$ allow us to make a reduction from $(\hat{N}, \hat{R},
\hat{\n})$ to $(\hat{N}', \hat{R}', \hat{\n}')$ in the same way as
that from $(N,R,\n)$ to $(N', R',\n')$. \hfill$\square$

The procedure given in the proposition is called a {\it size
change}.

\subsection{Reduction functors and reduction sequences}

\kg Let $(\K,\M,H)$ be a bimodule problem and $(\K',\M',H')$ be an
induced bimodule problem given by one of the $4$ reductions. There
is a nice functor
$$\vartheta: Mat(\K',\M')\longrightarrow
Mat(\K,\M)$$ which is induced from $(\vartheta_1,\vartheta_0)$ given in 3.1 and 3.2 and is defined as follows.
For any size vector $\m'$ of $(T',\sim')$ and any matrix $M'\in \M'_{\m'}$, let
$$\vartheta(M')=M'+(H'_{\m'})_{pq}\otimes \rho$$
 in case of
regularization, edge or loop  reductions.   We write
$\hat{\n}'=\m'$, and let $\hat{\n}=\m$ be calculated by Formula
(11) of 3.2.  In case of deletion let $$\vartheta(M')=M',\quad
\m=\m'.$$
 Then  $\m$ is
 a size vector of $(T, \sim)$ in all the cases and $\vartheta(M')\in \M_{\m}$. On
the other hand, the action of $\vartheta$ on morphisms is always
identity.

{\bf Lemma 3.3.1} $M'+H'_{\m'}=\vartheta(M')+H_{\m}$.

{\bf Proof.} Because of $H'_{\m'}=(H'_{\m'})_{pq}\otimes \rho +H_{\m},$ and the definition of $\vartheta(M')$.
\hfill$\square$

{\bf Proposition 3.3.1 } $\vartheta$ is a functor, which yields an
equivalence between category $Mat(\K',\M')$ and a full subcategory
of $Mat(\K,\M)$.

More precisely, $\varphi': M'\rightarrow L'$ is a morphism in
$Mat(\K',\M')$ if and only if $\vartheta(\varphi')=\varphi':
M\rightarrow L$ is a morphism in $Mat(\K,\M)$, where
$M=\vartheta(M')$, $L=\vartheta(L')$.

{\bf Proof.} Suppose that the size vectors of $M', L'$ are $\m', \underline{l}'$, and the size vectors of $M,L$
are $\m, \l$, calculated by Formula (11) of 3.2 respectively. Then Lemma 3.3.1 tells us that
$$M'+H'_{\m'}=M+H_{\m},\ \  L'+H'_{\l'}=L+H_{\l},$$ therefore
$$(M'+H'_{\m'})\varphi'=\varphi'(L'+H'_{\underline{l}'})\ \
{\rm iff}\ \
  (M+ H_{\m})\varphi'
=\varphi'(L+H_{\underline{l}}).$$
 Thus if $\varphi'\in
\K'_{\m'\times\l'}$, then $\varphi'\in\K_{\m\times\l}$, since $\K'_{\m'\times\l'}\subseteq \K_{\m\times\l}$.
Conversely, if $\varphi'\in \K_{\m\times\l}$ is a morphism from $M$ to $N$, then the first formula ensures that
$H'_{\m'}\varphi'\equiv_{\preceq(p,q)}\varphi'H'_{\l'}$ i.e. $\varphi'\in\K'_{\m'\times\l'}.$ \hfill $\square$

 The functor $\vartheta$ defined above is called a {\it
reduction functor}.

Let us go back to 2.6 and 3.1. Suppose that a size vector $\n$ of
$(T,\sim)$ is sincere, and a triple $(N,R,\n)$ is given by Formula
(6) of 2.6, moreover $(N',R',\n')$ is obtained from $(N,R,\n)$ by
one of the $3$ reductions. Then if we change superscript
``$\prime$" to ``$1$", and continue the reduction for
$(N^1,R^1,\n^1)$, we will obtain again a triple $(N^2,R^2,\n^2)$.
Thus a {\it reduction sequence} of triples follows by induction:
\bcen\unitlength=1mm
\begin{picture}(140,30)
\put(0,25){\mbox{$(N,R,\n),$}}
\put(20,25){\mbox{$(N^1,R^1,\n^1),$}}
\put(45,25){\mbox{$\cdots$,}}
\put(55,25){\mbox{$(N^r,R^r,\n^r),$}}
\put(80,25){\mbox{$\cdots$,}}
\put(90,25){\mbox{$(N^{s-1},R^{s-1},\n^{s-1}),$}}
\put(125,25){\mbox{$(N^s,R^s,\n^s)$}}

\put(5,15){\mbox{$N_{pq}$}} \put(25,15){\mbox{$N^1_{p_1q_1}$}}
\put(45,15){\mbox{$\cdots$}} \put(60,15){\mbox{$N^r_{p_rq_r}$}}
\put(80,15){\mbox{$\cdots$}}
\put(100,15){\mbox{$N^{s-1}_{p_{s-1}q_{s-1}}$}}

\mput(12,13)(1,-1){10}{\circle*{0.5}}
\mput(34,13)(1,-1){10}{\circle*{0.5}}
\mput(69,13)(1,-1){10}{\circle*{0.5}}
\mput(118,13)(1,-1){10}{\circle*{0.5}}

\put(62,0){\mbox{$\cdots$}} \put(105,0){\mbox{$\cdots$}}

\put(25,0){\mbox{$\overline{N}_{pq}$}}
\put(45,0){\mbox{$\overline{N}^1_{p_1q_1}$}}
\put(80,0){\mbox{$\overline{N}^r_{p_rq_r}$}}
\put(130,0){\mbox{$\overline{N}^{s-1}_{p_{s-1}q_{s-1}}$}}

\put(-10, 10){\mbox{$(*)$}}
\end{picture}\ecen
 where $p_r\in \P_r$, $q_r\in \Q_r$, $\P_r,\Q_r\in T^r/\!\sim^r $
with $N_{p_rq_r}^r$ being the first non-zero block of $N^r$,
$\overline{N}^r_{p_rq_r}=(H^{r+1}_{\n^{r+1}})_{p_rq_r}$, $r=0, 1,
\cdots, s-1$.

 There is also a corresponding {\it reduction sequence of induced bimodule problems}:
$$
 (**) \quad (\K, \M,H), (\K^1, \M^1,H^1), \cdots,
(\K^r, \M^r, H^r), \cdots, (\K^{s-1}, \M^{s-1}, H^{s-1}), (\K^s,
\M^s, H^s)
$$

 Furthermore, the sequence $(**)$ yields a {\it sequence
of reduction functors} $$\vartheta_{r-1,r}:
 Mat(\K^{r}, \M^{r})
\rightarrow Mat(\K^{r-1}, \M^{r-1}),\ r=1,2,\cdots,s,$$ where
$(\K^0,\M^0,H^0)=(\K,\M,H)$. A composition of reduction functors
is also called a {\it reduction functor}. In particular, let
$$
\vartheta_{0,r}=\vartheta_{0,1}\vartheta_{1,2} \cdots
\vartheta_{r-1,r}: Mat(\K^r, \M^r) \rightarrow Mat(\K, \M),
$$
$r=1,\cdots, s$. Then $\vartheta_{0,r}$ are reduction functors.

{\bf Corollary 3.3.1} Let $(\K,\M,H)$ be a bimodule problem with a
reduction sequence ($*$), and let $\hat{\n}^s$ be a size vector of
$(T^s,\sim^s)$. Then
$$
(\hat{*})\qquad (\hat{N},\hat{R},\hat{\n}),
(\hat{N}^1,\hat{R}^1,\hat{\n}^1),\cdots,
(\hat{N}^r,\hat{R}^r,\hat{\n}^r), \cdots,
(\hat{N}^{s-1},\hat{R}^{s-1},\hat{\n}^{s-1}),
(\hat{N}^s,\hat{R}^s,\hat{\n}^s)
$$
is also a reduction sequence, provided that
$(\hat{N}^{r-1},\hat{R}^{r-1},\hat{\n}^{r-1})$ is obtained from
$(\hat{N}^r,\hat{R}^r,\hat{\n}^r)$ by size change given in
Proposition 3.2.1 inductively for $r=s,s-1, \cdots,1$.
\hfill$\square$

\subsection{Canonical forms}

\kg Our special interests in the sequence ($*$) lead to a case of
$N^s=(0)$, then $\M^s=\{(0)\}$. If that is the case, we shall
refer to $(N^s,R^s, \n^s)$ as $(N^{\infty},R^{\infty},
\n^{\infty})$, $(\K^s,\M^s, H^s)$ as $(\K^{\infty},\M^{\infty},
H^{\infty})$, and $\vartheta_{0,s}$ as $\vartheta_{0,\infty}$. The
canonical form of an individual matrix of $Mat(\K,\M)$ will emerge
from this special case in the present subsection.

{\bf Lemma 3.4.1} Let $(\K,\M,H)$ be a bimodule problem,
$M\in\M_{\m}$, and a triple $(N,R,\m)$ be given by Formula (6) of
2.6. Then there exists a unique choice of $\overline{N}_{pq}$
given by one of the $3$ reductions of 3.1, such that $M\simeq
\vartheta(M')$ for some $M'\in Mat(\K',\M')$.

{\bf Proof.} Quote the matrix equation (7) of 2.6, where we change
the first block $N_{pq}$  to the fixed block $M_{pq}$ of $M$, and
the second $N_{pq}$ to $\overline{N}_{pq}$:
\begin{equation}
M_{pq}X_{qq}-X_{pp}\overline{N}_{pq}=\sum_{l=p+1}^t X_{pl}h_{lq}
-\sum_{l=1}^{q-1} h_{pl}X_{lq}
\end{equation}
We will determine $\overline{N}_{pq}$ based on  the $3$ reductions
respectively.

Regularization. If the equation system \textrm{II} of Definition
2.2.1 does not imply equation (8) of 2.6, set
$\overline{N}_{pq}=\emptyset$. Define a matrix $S \in \K_{\m
\times \m}$ having the identities as the diagonal blocks, and
satisfying
$$
M_{pq}= \sum_{l=p+1}^t S_{pl}h_{lq} -\sum_{l=1}^{q-1}
h_{pl}S_{lq}.
$$
Suppose now that the equation system  \textrm{II}  implies (8).

{Edge reduction.} Set $\overline{N}_{pq}=\left(\begin{array}{cc}
0& I_r \\ 0 & 0\end{array}\right) $ with $r=rank(M_{pq})$. Define
a diagonal block matrix $S \in \K_{\m \times \m}$, where
$S_{\I}=I_{n_{\I}}$ for any $\I\in T/\!\!\sim \setminus \{\P,\Q\}$
and $S_{\P}, S_{\Q}$ satisfy $M_{pq}S_{qq}=S_{pp}\left(\begin{array}{cc} 0& I_r \\
0 & 0\end{array}\right) $.

{Loop reduction.} Set $\overline{N}_{pq}=W$, the Weyr matrix
similar to $M_{pq}$. Define a diagonal block matrix $S \in \K_{\m
\times \m}$ with $S_{\I}=I_{n_{\I}}$ for any $\I\in T/\!\!\sim
\setminus \{\P\}$ and $S_{\P}$ satisfying $M_{pq}S_{qq}=S_{pp} W$.

In all the $3$ cases, let
\begin{equation}
M'=S^{-1}(M+H_{\m})S- (\overline{N}_{pq}\otimes\rho+H_{\m}).
\end{equation}
 Then $M'\in \M_{\m'}$. Thus
  $$(M+H_{\m})S=S(M'+H'_{\m'}),$$
i.e. $M\simeq \vartheta(M')$ in $Mat(\K,\M)$. \hfill$\square$

{\bf Lemma 3.4.2} Let $(\K,\M,H)$ be a bimodule problem, $M,L\in
Mat(\K,\M)$. If $M\simeq L$, then there exists an induced bimodule
problem $(\K',\M',H')$ obtained by one of  $3$ reductions of 3.1,
such that $M\simeq \vartheta(M')$, $L\simeq \vartheta(L')$ for
some $M',L'\in Mat(\K',\M')$ and $M'\simeq L'$.

{\bf Proof.} Since $M\simeq L$, $M$ and $N$ have the same size vector $\m$ of $(T,\sim)$ by Corollary 2.2.1.
Suppose $S\in\K_{\m\times \m}$ is invertible such that $(M+H_{\m})S=S(L+H_{\m})$ with the $(p,q)$-block
$$
M_{pq}S_{qq}-S_{pp}L_{pq}=\sum_{l=p+1}^t S_{pl}h_{lq}
-\sum_{l=1}^{q-1} h_{pl}S_{lq}.
$$
If the equation system \textrm{II} of 2.2 does not imply (8) of
2.6, set $\overline{N}_{pq}=\emptyset$.

 Otherwise, the equality becomes
$M_{pq}S_{qq}=S_{pp}L_{pq}$. Set
$\overline{N}_{pq}=\left(\begin{array}{cc} 0& I_r \\ 0 &
0\end{array}\right) $ whenever $p\nsim q$, where
$r=rank(M_{pq})=rank(L_{pq})$; or $\overline{N}_{pq}=W$ whenever
$p\sim q$, where $W$ is a Weyr matrix similar to $M_{pq}$ and
$L_{pq}$. Therefore $M\simeq \vartheta(M')$, $L\simeq
\vartheta(L')$ for some $M',L'\in \M_{\m'}$ by Lemma 3.4.1, i.e.
$$(M+H_{\m})U=U(M'+H'_{\m'})\ \ {\rm and}\ \ (L+H_{\m})V=V(L'+H'_{\m'})$$ for
some invertible $U,V\in \K_{\m\times\m}$. Finally
 $$(M'+H'_{\m'})U^{-1}SV=U^{-1}SV(L'+H'_{\m'})$$
 with $U^{-1}SV\in
\K'_{\m'\times\m'}$, i.e. $M'\simeq L'$ in $Mat(\K',\M')$.
\hfill$\square$

{\bf  Theorem 3.4.1\cite{S}}.  Let $(\K,\M,H=0)$ be a bimodule
problem, $M\in \M_{\m}$. Then there exists a unique reduction
sequence $(*)$ starting from $(N,R, \m)$, ending at $(N^{\infty},
R^{\infty}, \m^{\infty})$, such that
$$M\simeq \vartheta_{0,\infty}(M^{\infty})$$ for some
$M^{\infty}\in Mat(\K^{\infty},\M^{\infty})$ under both the fixed
orders of matrices and of  the base field (see Definition 2.3.1
and 2.5).

{\bf Proof.}  Suppose we have already had a unique sequence $(*)$ up to $(N^r, R^r, \m^r)$, and some $M^r\in
\M_{\m^r}^r$, such that $\vartheta_{0,r}(M^r)\simeq M$. Lemma 3.4.1 shows a unique induced
 triple $(N^{r+1}, R^{r+1},\m^{r+1})$ with
$\vartheta_{r,r+1}(M^{r+1})\simeq M^r$ for some
 $M^{r+1}\in \M_{\m^{r+1}}^{r+1}$. On the other hand,
 Lemma 3.4.2 ensures that any choice of $M^r$ in an isomorphism class
 of $Mat(\K^r, \M^r)$ does not influence the determination of
 $\overline{N}_{p_rq_r}$. Therefore the induced triple $(N^{r+1},
 R^{r+1}, \m^{r+1})$ is uniquely determined.  Because of the finiteness of the
 size of $N$, we finally reach to  $N^s=(0)$ for some positive
 integer $s$. \hfill$\square$

{\bf Remarks.} From now on, we assume that the original bimodule
problem $(\K,\M,H)$ satisfies the condition $H=0$.

(1) If that is the case, $\vartheta_{0,\infty}(M^{\infty})\simeq M$ implies that $(M+0)S
=S(0+H_{\m^{\infty}}^{\infty})$ for some invertible $S\in \K_{\m\times\m}$, or equivalently,
$$H_{\m^{\infty}}^{\infty}=S^{-1}MS.$$ The matrix $H_{\m^{\infty}}^{\infty}$ is called  the {\it canonical form}
of $M$, which can be regarded as a generalized Jordan form. It is clear that
$\K^{\infty}_{\m^{\infty}\times\m^{\infty}}$ is the endomorphism ring of
 $H_{\m^{\infty}}^{\infty}$.

 Without  $H=0$
  the trouble  is that there is
 no guarantee such that $H_{\m^{\infty}}^{\infty}$ is similar to $M$, which is not
 convenient.

(2) The restriction $H=0$  is natural, since the most bimodule problems
 considered in practice, including our main subject $P_1(\Lambda)$,
  have such a property.

{\bf Corollary 3.4.1} Let $(\K,\M,H=0)$ be a bimodule problem, $
M, L\in Mat(\K,\M) $. Then $M\simeq L $ if and only if  $M$ and $
L$ have the same canonical forms. \hfill$\square$

\subsection{Krull-Schmidt property}

\kg Following Ringel \cite{R2}, an $k$-additive category with the
finite dimensional objects is called a {\it Krull-Schmidt
category}, provided it has finite direct sums and split
idempotents.

{\bf Proposition 3.5.1} \cite{CB2}.  Let $(\K,\M,H)$ be a bimodule
problem, then $Mat(\K,\M)$ is a Krull-Schmidt category.

{\bf Proof.} Clearly, $Mat(\K,\M)$ has finite direct sums.

Now suppose $\varphi: M\rightarrow M$ is an idempotent with $M\in\M_{\m}$ and $\varphi\in \K_{\m\times\m}$.
Since idempotents split in the finite dimensional algebra $\K$, there exists a size vector $\l$ and some
matrices $\pi\in\K_{\m\times\l}$, $\iota\in\K_{\l\times\m}$ such that $\pi\iota=\varphi$ and
$\iota\pi=id_{\K_{\l\times\l}}$. Let $L=\iota M\pi- \iota d(\pi) $, we will show that  $\pi: M\rightarrow L$ and
$\iota: L\rightarrow M$ are morphisms of $Mat(\K,\M)$.

In fact, since $\pi=\varphi\pi$, $d(\pi)=d(\varphi)\pi +\varphi
d(\pi)$, $M\pi- \pi L= M\pi - \varphi M\pi +\varphi d(\pi)=
M\varphi\pi-\varphi M\pi+ d(\pi)-d(\varphi)\pi= (M\varphi-\varphi
M- d(\varphi))\pi + d(\pi)=d(\pi)$. Similarly, since
$id_{L}=\iota\pi$, $\iota=\iota\varphi$, we have $d(\iota)\pi
+\iota d(\pi)=0$ and $d(\iota)=d(\iota) \varphi+ \iota
d(\varphi)$. Thus $L\iota-\iota M=\iota M\varphi-\iota
d(\pi)\iota-\iota M= \iota M\varphi + d(\iota)\varphi -\iota M=
\iota M\varphi -\iota \varphi M +d(\iota)- \iota d(\varphi)= \iota
(M\varphi-\varphi M-d(\varphi)) + d(\iota)=d(\iota)$.
\hfill$\square$

\medskip
{\bf Lemma 3.5.1.} Let $(\K,\M,H)$ be a bimodule problem, and $\n$
be a size vector of $(T,\sim)$. If $\varphi\in \K_{\n\times\n}$ is
an idempotent, there exists some transformation matrix
$\chi\in\K_{\n\times\n}$, such that
$\overline{\varphi}=\chi^{-1}\varphi\chi$ is a diagonal idempotent
with
$$
(\overline{\varphi})_{\L}=\left(\begin{array}{cc} I_{d_{\L}} & 0 \\
0&0 \end{array}\right)_{n_{\L}\times n_{\L}},
$$
for some $0\leq d_{\L}\leq n_{\L}$, and for each $\L\in T/\sim$.

{\bf Proof.} (Given by Hu Yongjian)  Suppose that
$\varphi=(S_{ij})_{t\times t}$ is an upper triangular partitioned
 matrix with $S_{ij}(i<j)$ being $n_i\times n_j$ matrices
satisfying the equation system II of Definition 2.2.1.

(1) It is clear that $\forall \L\in T/\sim$, $S_{\L}^2=S_{\L}$.
Thus there exists some invertible matrix $U_{\L}$, such that
$U_{\L}^{-1}S_{\L}U_{\L}=\left(\begin{array}{cc} I_{d_{\L}}&0\\
0&0
\end{array}\right)_{n_{\L}\times n_{\L}}$, since our base field
$k$ is algebraically closed. Let $\chi_0=\rm{diag}(U_{11}, U_{22},
\cdots, U_{tt})$, and $\varphi_1 =\chi_0^{-1}\varphi \chi_0$

(2) Denote $\varphi_1$ by $A$, then
$$A=\left( \begin{array}{c}
\unitlength 1mm\begin{picture}(96,72)

\put(0,60){\dashbox{1}(16,12){$\begin{array}{cc} I_{d_1}\! & 0  \\
0& 0  \end{array}$ }} \put(16,60){\dashbox{1}(8,12){$S_{12}^1$ }}
\put(24,60){\dashbox{1}(8,12){$S_{12}^2$ }}
\put(32,60){\dashbox{1}(16,12){$\cdots$}}
\put(48,60){\dashbox{1}(8,12){$S_{1l}^1$}}
\put(56,60){\dashbox{1}(8,12){$S_{1l}^2$ }}
\put(64,60){\dashbox{1}(16,12){$\cdots$}}
\put(80,60){\dashbox{1}(8,12){$S_{1t}^1$ }}
\put(88,60){\dashbox{1}(8,12){$S_{1t}^2$ }}

\put(16,48){\dashbox{1}(16,12){$\begin{array}{cc} I_{d_2}\! & 0  \\
0& 0  \end{array}$ }} \put(32,48){\dashbox{1}(16,12){$\cdots$ }}
\put(48,48){\dashbox{1}(8,12){$S_{2l}^1$}}
\put(56,48){\dashbox{1}(8,12){$S_{2l}^2$ }}
\put(64,48){\dashbox{1}(16,12){$\cdots$}}
\put(80,48){\dashbox{1}(8,12){$S_{2t}^1$ }}
\put(88,48){\dashbox{1}(8,12){$S_{2t}^2$ }}

\put(32,36){\dashbox{1}(16,12){$\ddots$ }}
\put(48,36){\dashbox{1}(8,12){$\vdots$}}
\put(56,36){\dashbox{1}(8,12){$\vdots$}}
 \put(64,36){\dashbox{1}(16,12){$\cdots$}}
\put(80,36){\dashbox{1}(8,12){$\vdots$}}
\put(88,36){\dashbox{1}(8,12){$\vdots$ }}

\put(48,24){\dashbox{1}(16,12){$\begin{array}{cc} I_{d_{l}}\! & 0  \\
0& 0  \end{array}$ }} \put(64,24){\dashbox{1}(16,12){$\cdots$}}
\put(80,24){\dashbox{1}(8,12){$\; S_{lt}^1$ }}
\put(88,24){\dashbox{1}(8,12){$\; S_{lt}^2$ }}

\put(64,12){\dashbox{1}(16,12){$\ddots$ }}
\put(80,12){\dashbox{1}(8,12){$\vdots$}}
\put(88,12){\dashbox{1}(8,12){$\vdots$ }}

\put(80,0){\dashbox{1}(16,12){$\begin{array}{cc} I_{d_t}\! & 0  \\
0& 0  \end{array}$ }}
\end{picture}
\end{array}\right),$$
Let $A_l$ be the patitioned submatrix of $A$ consisting of the
intersections of $1$--$l$ block-rows and block columns. Let
$$\xi_l=\left(\begin{array} {c} S_{1l}^1\\ S_{2l}^1\\ \vdots\\
S_{l-1,l}^1 \end{array} \right), \quad
\eta_l=\left(\begin{array}{c} S_{1l}^2\\ S_{2l}^2\\ \vdots\\
S_{l-1,l}^2\end{array}\right)$$

Then $$A_l=\left(\begin{array}{c|c|c} A_{l-1}&\xi_l&\eta_l\\
\hline & I_{d_l} & 0 \\  & 0 & 0\end{array}\right).$$

(3) We claim that
$$\begin{aligned} (\textrm{i}) &\  A_l^2=A_l \mbox{ is an
idempotent matrix.}\\
(\textrm{ii}) &\  A_{l-1}\xi_l=0.\\ (\textrm{iii})&\
(I-A_{l-1})\eta_l=0.
\end{aligned}$$
In fact $A_t^2=A^2=A=A_t$. If $A_l^2=A_l$, them
$$
\left(\begin{array}{c|cc} A_{l-1}^2 & A_{l-1}\xi_l+ \xi_l &
A_{l-1}\eta_l \\ \hline & I_d & 0\\ & 0 & 0\end{array}\right)
=\left(\begin{array}{c|cc} A_{l-1}&\xi_l&\eta_l\\ \hline & I_{d} &
0 \\  & 0 & 0\end{array}\right).$$ Our claim follows by  induction
on $l=t, t-1, \cdots, 1$

(4) We define a transformation matrix $P_l$ inductively on $l=1,
2, \cdots, t $. Let $P_1=I_{n_1}$, than
$P_1^{-1}A_1P_1=\left(\begin{array}{cc} I_{d_1}&0\\ 0&0
\end{array}\right)$. We assume that there exists an invertible
matrix $P_{l-1}$, such that
$P_{l-1}^{-1}A_{l-1}P_{l-1}=\rm{diag}\left(\left(\begin{array}{cc}
I_{d_1}&0\\ 0&0 \end{array}\right) \cdots, \left(\begin{array}{cc}
I_{d_{l-1}}&0\\ 0&0 \end{array}\right)\right) $. Let $P_l=
\left(\begin{array}{c|c|c} P_{l-1}&\xi_l&-\eta_l\\ \hline 0&
I_{d_l} & 0 \\  & 0 & I_{n_l-d_l}\end{array}\right)$, then
$$
\begin{aligned}
P_l^{-1} A_l P_l= & \left(\begin{array}{c|c|c}
P_{l-1}^{-1}&-P_{l-1}^{-1}\xi_l&P_{l-1}^{-1}\eta_l\\ \hline & I_{d_l} & 0 \\
& 0 & I_{n_l-d_l}\end{array}\right) \left(\begin{array}{c|c|c}
A_{l-1}&\xi_l&\eta_l\\ \hline & I_{d_l} & 0 \\  & 0 &
0\end{array}\right) \left(\begin{array}{c|c|c}
P_{l-1}&\xi_l&-\eta_l\\ \hline 0& I_{d_l} & 0 \\  & 0 &
I_{n_l-d_l}\end{array}\right) \\
=& \left(\begin{array}{c|c|c} P_{l-1}^{-1} A_{l-1}&0&P_{l-1}^{-1}\eta_l\\
\hline & I_{d_l} & 0 \\  & 0 & 0\end{array}\right)
\left(\begin{array}{c|c|c} P_{l-1}&\xi_l&-\eta_l\\ \hline &
I_{d_l} & 0 \\  & 0 & I_{n_l-d_l}\end{array}\right)\\
 =& \left(\begin{array}{c|c|c} P_{l-1}^{-1}A_{l-1}P_{l-1}&0& 0\\ \hline &
I_{d_l} &  \\  &  & I_{n_l-d_l}\end{array}\right)
\end{aligned}
$$
by (ii), (iii) of (3). Thus we may take
$$P=\left( \begin{array}{c} \unitlength 1mm\begin{picture}(80,60)

\put(0,48){\dashbox{1}(16,12){$I_{n_1}$ }}
\put(16,48){\dashbox{1}(8,12){$\xi_1$ }}
\put(24,48){\dashbox{1}(8,12){$-\eta_1$ }}

\put(16,36){\dashbox{1}(16,12){$I_{n_2}$ }}
\put(32,36){\dashbox{1}(8,24){$\xi_2$ }}
\put(40,36){\dashbox{1}(8,24){$-\eta_2$ }}

\put(48,24){\dashbox{1}(16,36){ }}
\put(32,24){\dashbox{1}(16,12){$I_{n_3}$ }}

\put(48,12){\dashbox{1}(16,12){$\ddots$ }}
\put(64,12){\dashbox{1}(8,48){$\xi_t$ }}
\put(72,12){\dashbox{1}(8,48){$-\eta_t$ }}

\put(64,0){\dashbox{1}(16,12){$I_{n_t}$ }}
\end{picture}
\end{array}\right).$$

(5) Denote $P$ by $\chi_1$, we prove that $\chi_1\in
\K_{\n\times\n}$. In fact, since $A\in \K_{\n\times\n}$, we have
$$
\sum_{\I\ni i<j\in \J} c_{ij}^l S_{ij}=0
$$
for $1\leq l\leq q_{\I\J}$ and each $(\I,\J)\in (T/\sim)\times
(T/\sim)$, where $S_{ij}=(S_{ij}^1\ |\ S_{ij}^2)$. Therefore $
\sum\limits_{\I\ni i<j\in \J} c_{ij}^l S_{ij}^1=0 $ and $
\sum\limits_{\I\ni i<j\in \J} c_{ij}^l (-S_{ij}^2)=0 $.
Consequently,
$$
\sum_{\I\ni i<j\in \J} c_{ij}^l P_{ij}=0,
$$
where $P_{ij}=(S_{ij}^1\ |\ -S_{ij}^2)$ for $i<j$. Thus
$P\in\K_{\n\times\n}$.

Let $\chi=\chi_0\chi_1\in \K_{\n\times\n}$, then $\chi^{-1}
\varphi\chi= \rm{diag}\left(\left(\begin{array}{cc} I_{d_1} & 0\\
0 & 0\end{array} \right), \cdots, \left(\begin{array}{cc} I_{d_l} & 0\\
0 & 0\end{array} \right) \right)$
 for some $0\leq d_l\leq n_{\L}$ and each $\L\in T/\sim$. The proof
 is completed. \hfill$\square$

\medskip

Lemma 3.5.1 suggests an alternative definition of direct sums in
$Mat(\K,\M)$, which is  simple and concrete. We first define the
direct sum of general matrices.

{\bf Definition 3.5.1} Let $(T,\sim)$ be a set of integers with a
relation given in
 I of Definition 2.2.1, and let $C=(C_{ij})_{t\times t}$, $D=(D_{ij})_{t\times
t}$ be two matrices partitioned by $(T,\sim)$. Then a matrix
$$ \left(\begin{array}{cc} C_{ij} & 0\\
0 & D_{ij} \end{array}\right)_{t\times t}$$
 is said to be a {\it
direct sum} of $C$ and $D$, and  is denoted by $C\oplus_T D $, or
$C\oplus D$ for simplicity.

\medskip

{\bf Lemma 3.5.2} Let $(\K, \M, H)$ be a bimodule problem. If
$\m_1, \m_2$ are size vectors of $(T, \sim) $ and $\m=\m_1 +\m_2$,
then $H_{\m}=H_{\m_1}\oplus_T H_{m_2}$.

{\bf Proof.}  Let $ H=(h_{ij})_{t\times t}$. Given any pair $(i,j)$ with $1\leq i, j\leq t$, if $i\nsim j$, then
$h_{ij} =0$, and $(H_{\m_l})_{ij}=0$;  if $i\sim j$ then $(H_{\m_l})_{ij}=h_{ij} I_{(m_l)_i}$ for $l=0,1,2$,
where $\m_0=\m$. Thus $$(H_{\m})_{ij}= (H_{\m_1})_{ij}\oplus (H_{\m_2})_{ij},$$ since $m_i=(m_1)_i+(m_2)_i,$
$m_j=(m_1)_j+(m_2)_j$. The assertion follows block-wise.
 \hfill $\square$

{\bf Lemma 3.5.3} Let $(\K, \M, H)$ be a bimodule problem.

(1) If $M_1\in \M_{\m_1}$, $M_2\in \M_{\m_2}$ and $M=M_1\oplus_T
M_2$, then $$M+H_{\m}= (M_1+H_{\m_1})\oplus_T (M_2+H_{\m_2}).$$
And $M=M_1\oplus M_2$ in the sense of Proposition 3.5.1.

(2) Conversely, $\forall\, M\in\M_{\m}$, if $\varphi: M\rightarrow M$ is an idempotent, then there exists some
$M_0\simeq M$, such that $M_0=M_1\oplus_T M_2$ in the sense of Definition 3.5.1.

{\bf Proof.} (1) The direct sum over $T$ is give by lemma 3.5.2.
And the idempotent can be obtained easily.

(2) $\forall\, M\in\M_{\m}$ and any  idempotent $\varphi: M\rightarrow M$, there exists an invertible matrix
$\chi\in \K_{\m\times\m}$, such that $\overline{\varphi}=\chi^{-1}\varphi\chi$ is a diagonal idempotent given by
Lemma 3.5.1. Define $$M_0=\chi^{-1}(M+H_{\m})\chi-H_{\m},$$ then $M\simeq M_0$, and $\overline{\varphi}:
M_0\rightarrow M_0$ is a morphism. Thus $(M_0+H_{\m})\overline{\varphi}=\overline{\varphi}(M_0+H_{\m})$,
$\forall\, i,j\in T$, we have
$$\left((M_0)_{ij}+(H_{\m})_{ij}\right) \left(\begin{array}{cc} I_{d_j} & 0\\
0 & 0\end{array} \right)= \left(\begin{array}{cc} I_{d_i} & 0\\
0 & 0\end{array} \right) \left((M_0)_{ij}+(H_{\m})_{ij}\right),$$
i.e. $$(M_0)_{ij}=\left(\begin{array}{cc} M_{ij}^1 & 0\\
0 & M_{ij}^2 \end{array}\right),$$ such that the size of $(M_1)_{ij}$ is $d_i\times d_j$, and that of
$(M_2)_{ij}$ is $(m_i-d_i) \times (m_j-d_j)$.\hfill$\square$

\medskip
$M\in Mat(\K,\M)$ is said to be {\it indecomposable} if $M\simeq
M_1\oplus_T M_2$ implies that the size of $M_1$ or $M_2$ equals
$0$.

\medskip

{\bf Theorem 3.5.1.} Let $(\K,\M,H=0)$ be a bimodule problem with
a reduction sequence ($*$) given by Theorem 3.4.1 for some $M\in
Mat(\K, \M) $.

(1) Define a size vector $\underline{m(\I)}$ of $(T^{\infty},
\sim^{\infty})$ for any $\I\in T^{\infty}/\sim^{\infty}$, such
that $m(\I)_{\I}=1$ and $m(\I)_{\J}=0$, $\forall \ \J\neq \I$.
Then the zero matrix $O(\I)\in Mat(\K^{\infty}, \M^{\infty})$ of
size $\underline{m(\I)}$ is indecomposable. And
$$\{O(\I)\mid \forall \
\I\in T^{\infty}/\sim^{\infty} \}$$ is a complete set of
isomorphism classes of indecomposables  of $Mat(\K^{\infty},
\M^{\infty})$. Moreover,
$$M^{\infty}=\oplus_{T^{\infty}}O(\I)^{m_{\I}^{\infty}}$$ in
$Mat(\K^{\infty}, \M^{\infty})$.

(2) Let $\vartheta: Mat(\K^{\infty}, \M^{\infty})\rightarrow
Mat(\K,\M)$ be the reduction functor. Then
$$\vartheta(O(\I))=H^{\infty}_{\underline{m(\I)}}\mbox{ and denote it by
}H^{\infty}(\I) \mbox{ for simplicity}.$$ Thus  $H^{\infty}(\I)$
is the canonical form of a direct summand of $M$ for any $\I\in
T^{\infty}/\!\sim^{\infty}$. And $M\simeq \oplus_{T^{\infty}}
(H^{\infty}(\I))^{m_{\I}^{\infty}}$.

(3) If $M\simeq \oplus M_l^{e_l}$ in the sense of Proposition
3.5.1, such that  $M_l$ are pairwise non-isomorphic, then there
exists some $\I_l\in T^{\infty}/\sim^{\infty}$, such that
$M_l\simeq H^{\infty}(\I_l)$ and $e_l=m_{\I_l}^{\infty}$.

{\bf Proof.} (1) Since $$O(\I) \in \M^{\infty}_{\underline
{m(\I)}}, \mbox{\ \ and  \ \ }  \K^{\infty}_{\underline {m(\I)}
\times \underline {m(\I)}},$$ the endomorphism ring of $O(\I)$, is
local, $O(\I)$ is indecomposable. The direct sum comes from
$\m^{\infty}=\sum_{\I} \underline {m(\I)}^{m^{\infty}_{\I}}$

 (2) $\forall \;\I\in T^{\infty}/\!\sim^{\infty}$, denote
 $\underline
{m(\I)}$ by $\hat \n^{\infty}$. Then corollary 3.3.1 gives a
reduction sequence with a minimal end term. Thus $H^{\infty}(\I)$
is a canonical form of the objects in an isomorphism classes. Let
$\iota: H^{\infty}(\I) \rightarrow H^{\infty}_{\underline
m^{\infty}}$ and $\pi: H^{\infty}_{\underline m^{\infty}}
\rightarrow H^{\infty}(\I)$ be the natural morphisms given by the
direct summand respectively, and $M\varphi=\varphi
H^{\infty}_{\m^{\infty}}$ be given by Remark (1) of 3.4 for some
morphism $\varphi\in \K_{\m \times \m}$. Then the morphisms
$\iota\varphi: H^{\infty}(\I) \rightarrow M$ and
${\varphi}^{-1}\pi: M \rightarrow H^{\infty}(\I)$ tell us that
$H^{\infty}(\I)$ is isomorphic to a direct summand of $M$. The
direct sum is given by lemma 3.5.3 (1).

(3) Since $\oplus M_l^{e_l}\simeq \oplus
(H^{\infty}(\I))^{m_{\I}^{\infty}}$, the assertion follows from
Krull-Schmidt theorem. \hfill $\square$

\medskip
For example, see the last part of 2.2, we have $M\simeq
M_1\oplus_TM_2$, if we exchange the first and second rows and
columns inside the $i_{10}$-row and column blocks simultaneously.
Where
\begin{center}\unitlength=1mm
\begin{picture}(80,30)
\put(0,0){\framebox(30,30){}} \mput(0,5)(0,5){5}{\line(1,0){30}}
\mput(5,0)(5,0){5}{\line(0,1){30}} \put(22,16){$1$}
\put(21,21){$\lambda$} \put(27,26){$1$} \put(-10,13){$M_1=$}
\put(40,13){$M_2=$} \put(36,13){,}

\put(32, 1){$i_{12}$} \put(32,6){$i_{10}$} \put(32, 11){$i_{5}$}
\put(32, 16){$i_{3}$} \put(32, 21){$i_{2}$} \put(32,26){$i_{1}$}

\put(50,2){\framebox(25,25)} \mput(50,7)(0,5){4}{\line(1,0){25}}
\mput(55,2)(5,0){4}{\line(0,1){25}} \put(72,13){$0$}
\put(72,18){$1$}

\put(77, 3){$i_{10}$} \put(77,8){$i_{5}$} \put(77, 13){$i_{3}$}
\put(77, 18){$i_{2}$} \put(77, 23){$i_{1}$} \put(83,13){,}

\end{picture}
\end{center}
where the endomorphism rings are
$${\K}_{\m_1^\infty\times\m_1^{\infty}}^\infty=\{{\rm diag}(s_1,
s_1, s_1, s_1, s_1, s_1)\mid \forall s_1\in k\},$$ and $${\K}_{\m_2^\infty\times\m_2^{\infty}}^\infty=\{{\rm
diag}(s_2, s_2, s_2, s_2, s_2)\mid \forall s_2\in k\}$$ respectively.

\subsection{Corresponding  bocses}

\kg The present subsection is devoted to describing  the dual
structure of a bimodule problem $(\K,\M,H)$, i.e. the
corresponding bocs $\mathfrak{A}$ of $(\K,\M,H)$. The notion of
bocs has been studied intensively in Kiev school since 70's of
last century, and it has been used successfully for the proof of
Drozd's theorem. We illustrate here an explicit relation between
$(\K,\M,H)$ and $\mathfrak{A}$ in order to simplify the
calculation  in the proof of our main theorem.

Let $(\K,\M,H)$ be a bimodule problem defined in 2.2.1 with a
triangular basis $(A,B)$. Then the algebra $(\K,\,\cdot\,,1)$
yields a coalgebra structure $(\Omega,\mu,\varepsilon)$ with the
dual basis $$\{E^*_{\I} \ | \  \forall \ \I   \in T/\sim \} \cup
B^*.$$ Namely if
 the multiplication of  algebra $\K$ is given by
\begin{equation}
\zeta_{\I\L}^u\zeta_{\L\J}^v=\sum_\omega \c\zeta_{\I\J}^w
\end{equation}
 with the
structure constants $c\in \rK$, and
$$\zeta_{\I\J}^wE_{\J}=\zeta_{\I\J}^w,  \ \ \  E_{\I}
\zeta_{\I\J}^w=\zeta_{\I\J}^w,$$
 then
\begin{equation}\mu(\zeta_{\I\J}^w)^*=
(\zeta_{\I\J}^w)^*\otimes(E_{\J})^* +(E_{\I})^*\otimes
(\zeta_{\I\J}^w)^* + \delta_2(\zeta_{\I\J}^w)^*,
\end{equation}
where
\begin{equation}
\delta_2(\zeta_{\I\J}^w)^*=\sum_{\L,u,v}\c
(\zeta_{\I\L}^u)^*\otimes (\zeta_{\L\J}^v)^*,
\end{equation}
and
\begin{equation}
\mu(E_{\I})^*=(E_{\I})^*\otimes (E_{\I})^* .\end{equation}
 Therefor we obtain a linear map $\mu: \Omega\rightarrow \Omega\otimes\Omega.$
 On the other hand, let $\Gamma'=\rK\times\rK\times\cdots\times
\rK$ with $s$ copies of $k$ be a trivial algebra, where the
equivalent classes $\I_1,\I_2, \cdots,\I_s$ may be regarded as the
vertices in the ordinary quiver of $\Gamma'$. Then
  we define a
linear map $\varepsilon: \Omega\rightarrow \Gamma'$ given by
$\varepsilon(E_{\I})^*= 1_{\I}$ and
$\varepsilon(\zeta_{\I\J}^w)^*=0$. Thus we have the counit
$\varepsilon$ and comultiplication $\mu$ of $\Omega$, which
satisfy the laws of coalgebra $(\mu \otimes id) \mu = (id \times
\mu) \mu$, and $(\varepsilon \otimes id) \mu = id, (id \otimes
\varepsilon)\mu = id$, where $id$ stands for the identity morphism
on $\Omega$.

On the other hand $\K$-$\K$ bimodule $\M$ yields a left and right
comodule structure $\Gamma_m$, which has the dual basis $A^*
=\{(\rho_{\I\J}^w)^*\}$, and the comodule actions are given
according to Proposition 2.4.1:
$$\gamma_R(\rho_{\I\J}^w)^* =\sum_{\L,u,v} ~_l\c(\zeta_{\I\L}^u)^*\otimes
(\rho_{\L\J}^v)^*,$$
$$\gamma_L(\rho_{\I\J}^w)^* = \sum_{\L,u,v} ~_r\c (\rho_{\I\L}^u)^*\otimes
(\zeta_{\L\J}^v)^*,$$ then $\gamma_L: \Gamma_m \rightarrow \Omega
\otimes \Gamma_m$ satisfies the left comodule low $(\mu\otimes
id_{\Gamma_m})\gamma=(id_{\Omega} \otimes \gamma)\gamma$ and
$\gamma_R: \Gamma_m \rightarrow \Gamma_m\otimes\Omega$ satisfies
the right comodule low
$(id_{\Gamma}\otimes\mu)\gamma=(\gamma\otimes id_{\Omega})\gamma$.

Next we reconstruct $\Gamma_m$ to an algebra $\Gamma$ which is
freely generated by $A^*=\{(\rho_{\I\J}^w)^*\}$ over $\Gamma'$.
Then $\Omega$ can be viewed as a $\Gamma$-$\Gamma$-bimodule as
follows. Write $\Omega=\Omega_0\oplus \bar{\Omega}$ as a direct
sum of $\rK$-vector spaces, then $\Omega_0\cong\Gamma$,
 $\bar{\Omega}$ is freely generated
by $B^*=\{(\zeta_{\I\J}^w)^*\}$ over  $\Gamma$. Furthermore
suppose that
\begin{equation} d(\zeta_{\I\J}^u)=\zeta_{\I\J}^uH-H\zeta_{\I\J}^u
=\sum_wc_{\I\J}^{uw}\rho_{\I\J}^w,\end{equation}
then the derivation  $d:\K\rightarrow \M$ has a duality $Dd:
D\M\rightarrow D\K$ given by
\begin{equation}
(Dd)(\rho_{\I\J}^w)^*=\sum_u c_{\I\J}^{uw}(\zeta_{\I\J}^u)^*.
\end{equation}
which is in fact from $D\M$ to $D(\rad\K)$, because of
$d(E_{\I})=0$ by  Corollary 2.2.1.
 According to formula (18), and the left and
right module actions given in Proposition 2.4.1, we define a map
$\delta_1$ on $A^*$ given by $$\delta_1(\rho_{\I\J}^w)^*=(
Dd)(\rho_{\I\J}^w)^*+\gamma_L(\rho_{\I\J}^w)^*-\gamma_R(\rho_{\I\J}^w)^*,$$
more precisely
\begin{equation}
\delta_1(\rho_{\I\J}^w)^*=\sum_u
c_{\I\J}^{uw}(\zeta_{\I\J}^u)^*+\sum_{\L,u,v}(~_l\c(\zeta_{\I\L}^u)^*\otimes
(\rho_{\L\J}^v)^*-~_r\c (\rho_{\I\L}^u)^*\otimes
(\zeta_{\L\J}^v)^*)
\end{equation}
 Then
 $\Omega$ satisfies the relations:
\begin{equation}
(\rho_{\I\J}^w)^*\otimes(E_{\J})^*-
(E_{\I})^*\otimes(\rho_{\I\J}^w)^*=\delta_1(\rho_{\I\J}^w)^*.
 \end{equation}
  The reason will be given in Lemma 3.6.1. We denote the structure defined
   above by $\mathfrak{A}$, and write
$\mathfrak{A}=(\Gamma, \Omega)$, which
 is called a {\it bocs}, see \cite{D, CB1}. Note that both $\mu$ and
 $\varepsilon$ can be viewed as $\Gamma$-$\Gamma$-bimodule morphisms.

Now we show that the maps $\delta_1:\Gamma \rightarrow \Omega$ and
$\delta_2:\bar{\Omega} \rightarrow \bar{\Omega} \otimes
\bar{\Omega}$ extended linearly are determined just by
multiplications of matrices.

If we write $x_{ij}$, $z_{ij}$ for the $(i,j)$-entries of $R_0$,
$N_0$ which satisfy 2.2.1.$\textrm{II}$ and $\textrm{III}$
respectively, then
\begin{equation}
R_0= \sum_i E^*_{ii}E_{ii}+ \sum_{i,j}\left(\sum_w
(\zeta_{\I\J}^w)^*b_{\I\J}^w(i,j)\right)E_{ij}=\sum_{i,j}
x_{ij}E_{ij},
\end{equation}
\begin{equation}
N_0= \sum_{i,j}\left(\sum_w
(\rho_{\I\J}^w)^*a_{\I\J}^w(i,j)\right)E_{ij}=\sum_{i,j}
z_{ij}E_{ij}.
\end{equation}
Recall Formulae (1) and (2) of 2.3, the entries
$z_{p^w_{\I\J}q^w_{\I\J}}=(\rho^w_{\I\J})^*$ of $N_0$ are said to
be free, and $z_{ij}=\sum_w a_{\I\J}^w(i,j)(\rho^w_{\I\J})^*$ for
$(i,j)\ne (p^w_{\I\J}, q^w_{\I\J})$, $\forall w, \forall (\I,\J)
\in T/\sim\times T/\sim$, are said to be dependent. Similarly, the
entries $x_{\bar{p}^w_{\I\J}\bar{q}^w_{\I\J}}=(\zeta^w_{\I\J})^*$
of $R_0$ are said to be {\it free} and $x_{ij}=\sum_w
b_{\I\J}^w(i,j)(\zeta^w_{\I\J})^*$ for $(i,j)\ne
(\bar{p}^w_{\I\J}, \bar{q}^w_{\I\J})$, $\forall w, \forall (\I,\J)
\in T/\sim\times T/\sim$, are said to be {\it dependent}. Now
consider the matrix equation $(N_0+H)R_0=R_0(N_0+H)$, the
$(p^w_{\I\J}, q^w_{\I\J})$-entry , which we denote  by $p,q$ for
simplicity, gives an equation:
\begin{equation}
z_{pq}x_{qq}-x_{pp}z_{pq}=\Big(\sum_{p<l}
x_{pl}h_{lq}-\sum_{l<q}h_{pl}x_{lq}\Big)
+\Big(\sum_{p<l}x_{pl}z_{lq}-\sum_{l<q}z_{pl}x_{lq}\Big)
\end{equation}
 By Formulae (1) and (2) of 2.3 the left and right module actions
given in  Proposition 2.4.1 are equivalent respectively to the
formulae
\begin{equation}
\sum_{(i,j)\in \I\times
\J}(\sum_{l\in\L}b_{\I\L}^u(i,l)a_{\L\J}^v(l,j))E_{ij}=
\sum_{(i,j)\in \I\times\J} (\sum_w~_l\c a_{\I\J}^w(i,j))E_{ij},
\end{equation}
and
\begin{equation}
 \sum_{(i,j)\in \I\times
\J}(\sum_{l\in\L} a_{\I\L}^u(i,l)b_{\L\J}^v(l,j))E_{ij}=
\sum_{(i,j)\in \I\times \J}(\sum_w~_r\c a_{\I\J}^w(i,j) )E_{ij},
\end{equation}
 Furthermore Formulae (18),(2) and (3) yield
\begin{equation}
\sum_{(i,j)\in\I\times\J}\Big(\sum_{l\in\L}(b_{\I\L}^u(i,l)
h_{lj}-h_{il}b_{\L\J}^u(l,j))\Big)E_{ij}=\sum_{(i,j)\in\I\times\J}\Big(\sum_w
c_{\I\J}^{uw}a_{\I\J}^w(i,j)\Big)E_{ij}.
\end{equation}
 Denote $(p_{\I\J}^w,
q_{\I\J}^w)$ still by $(p,q)$. Since $a^w_{\I\J}(p,q)=1$,
$a^{w'}_{\I\J}(p,q)=0$, $\forall w'\ne w$, we have
\begin{eqnarray*}
Dd((\rho_{\I\J}^w)^*) &\stackrel{(19),(27)}{=}&
\sum_u\left(\sum_{l\in\L}(b_{\I\L}^u(p,l)h_{lq}-
h_{pl}b_{\L\J}^u(l,q))\right)(\zeta_{\I\J}^u)^*\\
&=&\sum_{l\in\L}\left(\sum_ub_{\I\L}^u(p,l)(\zeta_{\I\J}^u)^*\right)h_{lq}
-\sum_{l\in\I}h_{pl}\left(\sum_ub_{\L\J}^u(l,q)(\zeta_{\I\J}^u)^*\right)\\
&\stackrel{(22)}{=}& \sum_l(s_{pl}h_{lq}-h_{pl}s_{lq}).
\end{eqnarray*}
Moreover,
\begin{align*}
\sum_{\L, u ,v} &\left( ~_l\c(\zeta_{\I\L}^u)^*\otimes
(\rho_{\L\J}^v)^*- ~_r\c(\rho_{\I\L}^u)^*\otimes(\zeta_{\L\J}^v)^*
\right) \\
\stackrel{(25),(26)}{=}& \sum_{\L,u,v}\sum_{l\in\L} \left(
b_{\I\L}^u(p,l)a_{\L\J}^v(l,q)
(\zeta_{\I\L}^u)^*\otimes(\rho_{\L\J}^v)^*-
a_{\I\L}^u(p,l)b_{\L\J}^v(l,q)(\rho_{\I\L}^u)^*\otimes(\zeta_{\L\J}^v)^*
\right)\\
\stackrel{(22),(23)}{=}& \sum_{l}(x_{pl}z_{lq}-z_{pl}x_{lq}).\\
\end{align*}

{\bf Lemma 3.6.1} The Formula (24) is equivalent to
$$(\rho_{\I\J}^w)^*\otimes E_{\J}^*-E_{\I}^*\otimes
(\rho_{\I\J}^w)^*=\delta_1(\rho_{\I\J}^w)^*$$. \hfill$\square$

On the other hand, The Formulae (14) and (2) show that
\begin{equation} \sum_{\I\ni i<j\in \J}(\sum_lb_{\I\L}^u(i,l)
b_{\L\J}^v(l,j))E_{ij}=\sum_{\I\ni i<j\in\J}(\sum_w \c
b_{\I\J}^w(i,j)) E_{ij}.
\end{equation}
If the basis element $(\zeta_{\I\J}^w)^*$ is given by Formula (2)
of 2.3, then the $({p}^w_{\I\J}, {q}^w_{\I\J})$-entry of $R_0\cdot
R_0$ equals
\begin{align*}
 \sum_{\bar{p}\leq l\leq \bar{q}}x_{\bar{p}l}x_{l\bar{q}}&\stackrel{(22)}{=}
(\zeta_{\I\J}^w)^*\otimes E_{\J}^*+E_{\I}^*\otimes
(\zeta_{\I\J}^w)^*+
\sum_{\L,u,v}\sum_{l\in\L}b_{\I\L}^u({p},l)b_{\L\J}^v(l,{q})
(\zeta_{\I\L}^u)^*\otimes (\zeta_{\L\J}^v)^*\\
&\stackrel{(28)}{=} (\zeta_{\I\J}^w)^*\otimes
E_{\J}^*+E_{\I}^*\otimes (\zeta_{\I\J}^w)^*+\sum_{\L,u,v} \c
(\zeta_{\I\L}^u)^*\otimes (\zeta_{\L\J}^v)^*
\end{align*}

{\bf Lemma 3.6.2} $\mu(\zeta_{\I\J}^w)^*)=\sum_{\bar{p}\leq l\leq
\bar{q}}x_{\bar{p}l}x_{l\bar{q}}.$  \hfill$\square$

{\bf Theorem 3.6.1} Let $(\K,\M,H)$ be a bimodule problem with a
triangular basis $(A,B)$. Then
$$
\delta_1(\rho_{\I\J}^w)^*=\Big(\sum_{p<l}
x_{pl}h_{lq}-\sum_{l<q}h_{pl}x_{lq}\Big)
+\Big(\sum_{p<l}x_{pl}z_{lq}-\sum_{l<q}z_{pl}x_{lq}\Big)
$$
where $p=p_{\I\J}^w$, $q=q_{\I\J}^w$. In particular, for the first
basis element $\rho\in A$,
$$
\delta(\rho^*)=\sum_{p<l} x_{pl}h_{lq}-\sum_{l<q}h_{pl}x_{lq}
$$
 On the other hand,
$$
\delta_2(\zeta_{\I\J}^w)^*=\sum_{\bar{p}<l<\bar{q}}x_{\bar{p}l}x_{l\bar{q}}
$$
where $\bar{p}=\bar{p}_{\I\J}^w$, $\bar{q}=\bar{q}_{\I\J}^w$.
\hfill$\square$

{\bf Corollary 3.6.1} The bocs $\mathfrak{A}$ is normal and
triangular, see \cite{Ro, BK, CB1}.

{\bf Proof.} The normality is given by Formula (17) in the sense
of \cite{Ro}, which follows from the algebra $\K$ being upper
triangular. The triangularity follows from the triangularity of
the basis $(A,B)$. In fact, $\delta_1$ on $A^*$ is triangular,
i.e. $(\rho_{\I\J}^w)^*\succ (\rho_{\I\J}^u)^* \mbox{ and }
(\rho_{\I\J}^v)^*$, by Proposition 2.4.1; $\delta_2$ on $B^*$ is
also triangular i.e. $(\zeta_{\I\J}^w)^*\succ (\zeta_{\I\J}^u)^*
\mbox{ and } (\zeta_{\I\J}^v)^*$, since $\zeta^w_{\I\J}\succ
\zeta^u_{\I\J}$ and $\zeta^v_{\I\J}$ in Formula (14), see
\cite{CB1, BK}. \hfill$\square$

\newpage

\bcen
\section{Freely parameterized bimodule problems}
\ecen

\subsection{Parameters}

 \kg Let $(\K,\M,H=0)$ be a bimodule problem, and $(*)$ be a reduction sequence.
  It has been shown in 2.5, that every Weyr matrix in the third
sequence of ($*$) possesses some fixed eigenvalues. From now on we
allow some of the eigenvalues of Weyr matrices to be {\it
parameters}, and $(\K^r, \M^r, H^r)$ is said to be a {\it
parameterized bimodule problem}. Suppose that
$$P^r=\{\lambda_1, \lambda_2, \cdots, \lambda_{i_r}\}$$
is the set of parameters appearing in $H_{\n^r}^r$, then $P^1\subseteq P^2\subseteq\cdots\subseteq P^s.$ Denote
the domain of $P^r$  by $D^r$, then $D^r\subseteq \rK^{i_r}$, $r=1,2,\cdots,s$.

{\bf Proposition 4.1.1 \cite{S}.}  $D^s$ is determined by some
polynomial equations and inequalities in $i_s$ indeterminates.
Therefore it is a locally closed subset of $\rK^{i_s}$.

{\bf Proof. } We use induction on $r$. The case of $r=0$ is trivial since $H_{\n}=0$ has no parameter. Suppose
the assertion is true for $r$. Now we perform the $(r+1)$-reduction. Denote by $C$ the coefficient matrix of
equation system II of Definition 2.2.1 at $(\P_r, \Q_r)$;  by $\overline{C}$ that of equation system II and
equation (8) of 2.6 at $(\P_{r},\Q_{r})$.

Case 1.  Equation system II of Definition 2.2.1 implies the equation 2.6(8), if and only if
$\mbox{rank}C=\mbox{rank}\bar{C},$ which is determined by some polynomial equations in $x_1, x_2, \cdots,
x_{i_r}$.

Case 2.  Equation system II of Definition 2.2.1.does not imply the equation 2.6(8), if and only if
$\mbox{rank}C+1=\mbox{rank}\bar{C},$ which is determined by a polynomial inequality in $x_1, x_2, \cdots,
x_{i_r}$.\hfill$ \square$

{\bf Remarks.} (1) The parameters in $P^r\setminus P^{r-1}$ are
free at the $r$-th reduction unless otherwise assumption, since
the restriction equations or inequalities for these parameters
only appear in the reductions after  the $r$-th step.

(2)  If $N^s=(0)$ in $(*)$, then the parameterized matrix $H_{\n^{\infty}}^{\infty}(\lambda_1,\lambda_2, \cdots,
\lambda_{i_s})$ can be regarded as a {\it representation variety} over $(\K, \M)$,  which is locally closed.

(3) Note that it is difficult to define a representation category over a parameterized bimodule problem in
general, (see  examples of 4.6). Thus we will focus on a very special situation of the parameterized bimodule
problems defined in the next subsection.

\subsection{Freely parameterized  bimodule Problems}

\kg {\bf Definition 4.2.1}  A  bimodule problem with $i$
parameters is said to be {\it freely parameterized}, provided

(1) for  any $\I\in T/\!\!\sim$, there exists at most one
parameter $\lambda$ attached to $\I$;

(2) all parameters are algebraically independent, and the domain
$$
D=\rK\setminus \{\mbox{the roots of }
g(\lambda_1,\cdots,\lambda_i)=g_1(\lambda_1) \cdots
g_i(\lambda_i)\} $$   \hfill$\square$

An equivalent class $\I\in T/\!\!\sim$ is said to be {\it
non-trivial} if there exists a parameter $\lambda_{\I}$ attached
to $\I$; otherwise, $\I$ is said to be {\it trivial}. Suppose that
 $\{\I_1,\cdots, \I_i\}$ is a complete set of non-trivial equivalent
classes of $(T,\sim)$, and $(\J_1,\cdots,\J_j)$ is a complete set
of trivial ones respectively.

Now we construct a bocs based on the freely parameterized bimodule
problem $(\K,\M,$ $H)$. Let $\Gamma'$ be a category with
indecomposable objects $\I_1,\cdots,\I_i; \J_1,\cdots,\J_j$ and
morphisms
$$
\Gamma'(\I_l,\I_l)=\rK[\lambda_l,g_l(\lambda_l)^{-1}], \quad
l=1,\cdots,i;\quad \Gamma'(\J_l,\J_l)=\rK,\quad l=1,\cdots, j.
$$
Then $\Gamma'$ is equivalent to the following algebra
$$\rK[\lambda_1,g_1(\lambda_1)^{-1}]\times\cdots\times
\rK[\lambda_i, g_i(\lambda_i)^{-1}]\times \rK_1\times \cdots
\times \rK_j,$$ where $\rK_1=\cdots=\rK_j=\rK$. $\Gamma'$ is
called a {\it minimal category}, and can be  shown in the
following diagram \cite{CB1}:
\begin{center}\unitlength=1mm
\begin{picture}(90,10) \put(9,2){\circle*{0.6}}
\put(7,5){\circle{4}} \put(8,3){\vector(1,0){0}}
\put(8,-2){$\I_1$} \put(-1,6){$\lambda_{1}$}

\put(24,2){\circle*{0.6}} \put(22,5){\circle{4}}
\put(23,3){\vector(1,0){0}} \put(23,-2){$\I_2$}
\put(14,6){$\lambda_{2}$}

\put(30,3){$\cdots$}

\put(48,2.5){\circle*{0.6}} \put(47,5){\circle{4}}
\put(48,3){\vector(1,0){0}} \put(48,-2){$\I_i$}
\put(39,6){$\lambda_{i}$}

\put(60,3){\circle*{1}} \put(70,3){\circle*{1}}
\put(78,2.5){$\cdots$} \put(90,3){\circle*{1}} \put(58,-1){$\J_1$}
\put(68,-1){$\J_2$} \put(88,-1){$\J_j$}

\end{picture}
\end{center}

Denote  $A^*$ by $\{a_1,\cdots,a_n\}$,  $B^*$ by $\{
v_1,\cdots,v_m\}$ according to the order of Definition 2.4.1, in
order to simplify the notations. Then $\mathfrak{A}$ is said to be
a {\it layered bocs} with a layer
$$
L=(\Gamma'; \omega; a_1, a_2, \cdots, a_n; v_1, v_2,\cdots,v_m),
$$
where $\omega: \Gamma'\rightarrow \Omega$ is a
$\Gamma'$-$\Gamma'$-bimodule map given by
$\omega(E_{\I})=(E_{\I})^*$ for any $\I\in T/\!\!\sim$. Note that
$\Gamma'$ is of infinite dimension if $i>0$.

A layered bocs  can be illustrated as a {\it differential
biquiver} as follows. The set of vertices of the biquiver is
$T/\!\!\sim$, the set of solid arrows is $\{\lambda_1, \cdots,
\lambda_i; a_1, \cdots, a_n \}$ and the set of dotted arrows is
$\{v_1,\cdots, v_m\}$ \cite{CB2}. The examples given in 2.1
correspond to the following differential biquivers, if the
$k$-basis of $rad(\Lambda)$, $a,b,c,d,f$ still stand for the
$k$-basis of $rad(\widetilde {\Lambda})$, as well as$ \left(
        \begin{array}{cc}
        0&\rad {\widetilde \Lambda}\\ 0&0
        \end{array}
  \right)$  with $a\prec b\prec
c\prec d\prec f$, and  $a^*, b^*, c^*, d^*, f^*$ are dual basis.
$${\unitlength=1mm
\begin{array}{c}
 \begin{picture}(80, 30) \put(-30,20){\mbox{\bf Example 1.}}
 \put(15,4){\circle*{1.00}}\put(15,21){\circle*{1.00}}
\put(15,20.00){\vector(0,-1){15}} \qbezier(13,20)(9,13)(13,6)
\put(13,6){\vector(1,-1){1}}\qbezier(17,20)(21,13)(17,6)
\put(17,6){\vector(-1,-1){1}} \put(13,13){$b^*$} \put(9,13){$a^*$}
\put(20,13){$c^*$}

\qbezier[10](14,21)(9,23)(11,25) \qbezier[10](14,21)(13,26)(11,25)
\qbezier[10](16,21)(21,23)(19,25)
\qbezier[10](16,21)(17,26)(19,25)
\qbezier[10](15,22)(13,24.5)(15,25.5)
\qbezier[10](15,22)(17,24.5)(15,25.5)

\qbezier[10](14,4)(9,2)(11,0) \qbezier[10](14,4)(13,-1)(11,0)
\qbezier[10](16,4)(21,2)(19,0) \qbezier[10](16,4)(17,-1)(19,0)
\qbezier[10](15,3)(13,1.5)(15,-0.5)
\qbezier[10](15,3)(17,1.5)(15,-0.5)

\put(22,20){$1'$} \put(22,3){$1$}

\put(35,10){\mbox{$ \begin{array}{l} \dz(a^*)=0,\\ \dz(b^*)=0, \\
\dz(c^*)=0. \end{array}$}}
\end{picture}  \\
  \begin{picture}(80, 30) \put(-30,20){\mbox{\bf Example 2.}}
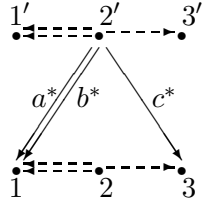
\put(0,1.50){\circle*{1.00}} \put(11,1.50){\circle*{1.00}}
\put(22,1.50){\circle*{1.00}} \put(2.00,1.50){\vector(-1,0){1}}
\mput(3.00,1.50)(2,0){4}{\line(1,0){1}}
\put(2.00,2.50){\vector(-1,0){1}}
\mput(3.00,2.50)(2,0){4}{\line(1,0){1}}
\put(20.00,2.0){\vector(1,0){1}}
\mput(12.00,2.0)(2,0){4}{\line(1,0){1}}

 \put(-1,-2){\mbox{$1$}}
\put(11,-2){\mbox{$2$}} \put(22,-2){\mbox{$3$}}

\put(0,19.50){\circle*{1.00}} \put(11,19.50){\circle*{1.00}}
\put(22,19.50){\circle*{1.00}} \put(2.00,19.50){\vector(-1,0){1}}
\mput(3.00,19.50)(2,0){4}{\line(1,0){1}}
\put(2.00,20.50){\vector(-1,0){1}}
\mput(3.00,20.50)(2,0){4}{\line(1,0){1}}
\put(20.00,20.0){\vector(1,0){1}}
\mput(12.00,20.0)(2,0){4}{\line(1,0){1}}

 \put(-1,21){\mbox{$1'$}}
\put(11,21){\mbox{$2'$}} \put(22,21){\mbox{$3'$}}

\put(11.00,18.00){\vector(-2,-3){10}}
\put(10.00,18.00){\vector(-2,-3){10}}
\put(12.00,18.00){\vector(2,-3){10}}

\put(2,10){\mbox{$a^*$}} \put(8,10){\mbox{$b^*$}}
\put(18,10){\mbox{$c^*$}}

\put(35,10){\mbox{$ \begin{array}{l} \dz(a^*)=0,\\ \dz(b^*)=0, \\
\dz(c^*)=0. \end{array}$}}
\end{picture}  \\
  \begin{picture}(80, 40) \put(-30,20){\mbox{\bf Example 3.}}
\put(10,0){\circle*{1.00}} \put(10,20){\circle*{1.00}}
\qbezier(7,19)(5,11)(7,3) \put(7,3){\vector(1,-1){1}}
\qbezier(13,19)(15,11)(13,3) \put(13,3){\vector(-1,-1){1}}

\qbezier[10](9,20)(4.5,19)(5,21) \qbezier[10](9,20)(5.5,22)(5,21)
\qbezier[10](9,21)(6,23)(7,24) \qbezier[10](9,21)(9,25)(7,24)
\qbezier[10](11,20)(15.5,19)(15,21)
\qbezier[10](11,20)(14.5,22)(15,21)
\qbezier[10](11,21)(14,23)(13,24)
\qbezier[10](11,21)(11,25)(13,24)

\qbezier[10](9,0)(4,1)(5,-1) \qbezier[10](9,0)(5,-2)(5,-1)
\qbezier[10](9,-1)(6,-3)(7,-4) \qbezier[10](9,-1)(9,-5)(7,-4)
\qbezier[10](11,0)(16,1)(15,-1) \qbezier[10](11,0)(15,-2)(15,-1)
\qbezier[10](11,-1)(14,-3)(13,-4)
\qbezier[10](11,-1)(11,-5)(13,-4)

\put(17,20){$1'$} \put(17,0){$1$}

\put(9.00,19.00){\vector(0,-1){18}} \put(11.00,19.00){\vector(0,-1){18}}

\put(4,10){\mbox{$a^*$}} \put(8,10){\mbox{$b^*$}}
\put(10.50,10){\mbox{$c^*$}} \put(14,10){\mbox{$d^*$}}

\put(35,8){\mbox{$ \begin{array}{l} \dz(a^*)=0,\\ \dz(b^*)=0, \\
\dz(c^*)=v'b^*-b^*v, \\ \dz(d^*)=\az v'a^*-a^*v+ u'b^*-\az b^*u.
\end{array}$}}
\end{picture}
\\
 \begin{picture}(80, 60) \put(-30,50){\mbox{\bf Example 4.}}
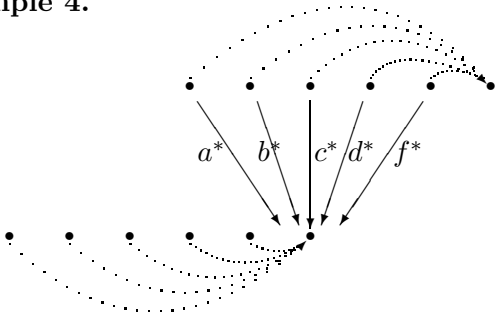
\mput(5,40)(8,0){6}{\circle*{1.00}}
\mput(21,20)(-8,0){6}{\circle*{1.00}}
\put(21,38){\vector(0,-1){17}} \put(14 ,38){\vector(1,-3){5.5}}
\put(6,38){\vector(2,-3){11}} \put(28 ,38){\vector(-1,-3){5.5}}
\put(36,38){\vector(-2,-3){11}}

\qbezier[25](5,41)(27,60)(44,41) \qbezier[20](13,41)(30,55)(44,41)
\qbezier[20](21,41)(31,51)(44,41)
\qbezier[20](29,41)(32,46)(44,41)
\qbezier[10](37,41)(39,43)(44,41)
\put(43.5,41.3){\vector(1,-1){1}}

\qbezier[10](20,19)(17,17)(13,19) \qbezier[20](20,19)(12,14)(5,19)
\qbezier[20](20,19)(8,10)(-3,19) \qbezier[20](20,19)(4,6)(-11,19)
\qbezier[25](20,19)(0,1)(-19,19)\put(19.5,18.5){\vector(1,1){1}}

\put(6,30){\mbox{$a^*$}} \put(14,30){\mbox{$b^*$}}
\put(21.50,30){\mbox{$c^*$}} \put(26,30){\mbox{$d^*$}}
\put(32,30){\mbox{$f^*$}}

\put(60,25){\mbox{$ \begin{array}{l} \dz(a^*)=0,\\ \dz(b^*)=0, \\
\dz(c^*)=0, \\ \dz(d^*)=0,  \\ \dz(f^*)=0. \end{array}$}}
\end{picture}
\end{array}}
$$

A freely parameterized bimodule problem is said to be {\it minimal}, if $\M=\{(0)\}$. Correspondingly a triple
$(N,R,\n)$ is said to be {\it minimal}, if $N=(0)$. Moreover, a layered bocs $\mathfrak{A}=(\Gamma, \Omega)$ is
called a {\it minimal bocs}, if $\Gamma=\Gamma'$, then the layer $L=(\Gamma'; \omega; v_1, \cdots, v_m)$ without
any solid arrows. The corresponding differential byquiver consists of isolated vertices,  vertices with single
solid loops, and dotted arrows.

From now on we will not worry about any differences between a
freely parameterized bimodule problem and  its layered bocs. We
will use freely the two concepts depending on convenience of
concerned problems.

{\bf Remark} Let $(\K,\M,H=0)$ be a bimodule problem, and $(*)$ be
a reduction sequence. We usually obtain a freely parameterized
bimodule problem as the following way: every parameter appears in
a loop reduction, i.e.
$\overline{N}_{p_{r-1}q_{r-1}}=W\oplus(\lambda)$ with $W$ being a
Weyr matrix of fixed eigenvalues and the domain
$D^r|_{\{\lambda\}}=\rK\ \setminus\{\mbox{the eigenvalues of
}W\},$ for some step $1\leq r< s $. Then the end term
$(N^s,R^s,\n^s)$ of $(*)$ is said to be a {\it freely
parameterized triple}. In this case, $(\K^s,\M^s, H^s)$, the end
term of the corresponding
 sequence $(**)$, is a freely parameterized bimodule
problem.

\subsection{The differential of the first arrow}

\kg We will illustrate all the possibilities of the differential
of the first arrow of a layered bocs in this subsection, which are
not mentioned before in any references.

Let $\mathfrak{A}=(\Gamma,\Omega)$ be a bocs with a layer
$L=(\Gamma';\omega;a_1,\cdots,a_n; v_1,\cdots,v_m)$. Then the
differential $\delta(a_1)$ of the first arrow $a_1:\P\rightarrow
\Q$ has the following possibilities:

{\bf A1.}   { \unitlength 1mm
\begin{picture}(24, 8)
\put(10.00,1.00){\circle*{1.00}} \put(5.00,2.00){\oval(5,5)[t]}
\put(5.00,2.00){\oval(5,5)[bl]} \put(4.00,-1.00){\vector(3,1){3}}
\put(15.00,2.00){\oval(5,5)[t]} \put(15.00,2.00){\oval(5,5)[br]}
\put(16.00,-1.00){\vector(-3,1){3}}
\put(10.00,5.00){\makebox{$\P$}}
\put(0.00,2.0){\makebox{$\lambda$}}
\put(20.00,2.00){\makebox{$a_1$}}
\end{picture},}  $\P =\Q$,
$\Gamma'(\P,\P)=\rK[\lambda,g_{\P}(\lambda)^{-1}]$ or

\vspace{2mm}

 \unitlength 1mm
 \hspace{.8cm}\begin{picture}(40, 6)
\put(10.00,0.00){\circle*{1.00}} \put(25.00,0.00){\circle*{1.00}} \put(11.00,3.00){$\P$} \put(24.00,3.00){$\Q$}
\put(5.00,1.00){\oval(5,5)[t]} \put(5.00,1.00){\oval(5,5)[bl]} \put(4.00,-2.00){\vector(3,1){3}}
\put(30.00,1.00){\oval(5,5)[t]} \put(30.00,1.00){\oval(5,5)[br]} \put(31.00,-2.00){\vector(-3,1){3}}
\put(10.50,0.00){\line(1,0){13}} \put(23.50,0.00){\vector(1,0){1}} \put(0.00,1.0){\makebox{$\lambda$}}
\put(15.50,2.00){\makebox{$a_1$}} \put(35.00,1.00){\makebox{$\mu$}}
\end{picture},  $\P\neq\Q$,
$\Gamma'(\P,\P)=\rK[\lambda, g_{\P}(\lambda)^{-1}]$ and
 $\Gamma'(\Q,\Q)=\rK[\mu,$ $ g_{\Q}(\mu)^{-1}]$. Then
\begin{equation}
\delta(a_1)=\sum_{j=1}^mf_j(\lambda, \mu)v_j
\end{equation}
where we may assume that $f_j(\lambda,\mu)\in k[\lambda,\mu]$ by dividing $v_j$ by some power of $g(\lambda)$;
and $\lambda v_ j$ stands for the left multiplication of $v_j$ by $\lambda$, and $\mu v_j$ for the right
multiplication  by $\lambda$ in the first case or by $\mu$ in the second case.

Recall from \cite{CB1}, let $f(\lambda, \mu)$ be the highest
common factor of the $f_j(\lambda,\mu)$ and let $q_j(\lambda,
\mu)=f_j(\lambda, \mu)/f(\lambda,\mu)$. Since $q_j(\lambda, \mu)$
are coprime, there are polynomials $s_j(\lambda, \mu)$ and a
non-zero polynomial $c(\lambda)$, such that
\begin{equation}
c(\lambda)=\sum_{j=1}^ms_j(\lambda, \mu)q_j(\lambda, \mu)
\end{equation}
 in ring $\rK[\lambda,\mu]$. Thus
$$
1=\sum_{j=1}^mc(\lambda)^{-1}s_j(\lambda, \mu)q_j(\lambda, \mu)
$$
 in a Hermite ring $S=\rK[\lambda, \mu, g_{\P}(\lambda)^{-1}c(\lambda)^{-1},
  g_{\Q}(\mu)^{-1}]$. So there is an invertible matrix $Q$ in $M_{m}(S)$ with
the first row $(q_j(\lambda, \mu))_{j=1, \cdots, m }$, and we can
make a change of basis of the form
$$
(w_1, w_2, \cdots, w_m)^T=Q(v_1, v_2, \cdots, v_m)^T
$$
so that \begin{equation} \delta (a_1)=f(\lambda, \mu)w_1.
\end{equation}

{\bf A2.}\hspace{0.7cm}
 {\unitlength 1mm
\begin{picture}(33, 6)
\put(10.00,-1.00){\circle*{1.00}}
\put(25.00,-1.00){\circle*{1.00}} \put(5.00,0.00){\oval(5,5)[t]}
\put(5.00,0.00){\oval(5,5)[bl]} \put(4.00,-3.00){\vector(3,1){3}}
\put(10.50,-1.00){\line(1,0){13}}
\put(23.50,-1.00){\vector(1,0){1}}
\put(10.00,1.00){\makebox{$\P$}} \put(27.00,0.00){\makebox{$\Q$}}
\put(0.00,0.0){\makebox{$\lambda$}}
\put(16.00,1.00){\makebox{$a_1$}}
\end{picture},}  $\P\neq \Q$, $\Gamma'(\P,\P)=\rK[\lambda,g_{\P}(\lambda)^{-1}]$,
$\Gamma'(\Q,\Q)=\rK$.

\noindent (Or dually,
 {\unitlength 1mm
\hspace{0cm}\begin{picture}(35, 6)
\put(5.00,-1.00){\circle*{1.00}} \put(20.00,-1.00){\circle*{1.00}}
\put(25.00,0.00){\oval(5,5)[t]} \put(25.00,0.00){\oval(5,5)[br]}
\put(26.00,-3.00){\vector(-3,1){3}}
\put(5.50,-1.00){\line(1,0){13}}
\put(18.50,-1.00){\vector(1,0){1}} \put(3.00,1.00){\makebox{$\P$}}
\put(18.00,1.00){\makebox{$\Q$}} \put(29.00,0.0){\makebox{$\mu$}}
\put(12.00,1.00){\makebox{$a_1$}}
\end{picture},}
$\Gamma'(\P,\P)=\rK$, $\Gamma'(\Q,\Q)=\rK[\mu,g_{\Q}(\mu)^{-1}]$).
Then
\begin{equation}
\delta(a_1)=\sum_{j=1}^mf_j(\lambda)v_j,
\end{equation}
where $f_j(\lambda)\in k[\lambda]$. If $f(\lambda)$ is the highest common factor of $f_j(\lambda)$, then $
q_j(\lambda)=f_j(\lambda)/f(\lambda)$, $j=1,\cdots,m$, are coprime, and there are polynomials $s_j(\lambda)$
such that $1=\sum_{j=1}^ms_j(\lambda)q_j(\lambda)$. After a basis transformation as in Case 1 given by an
invertible matrix $Q\in M_{m} (\rK[\lambda, g_{\P}(\lambda)^{-1}])$, we have
\begin{equation}
\delta(a_1)=f(\lambda)w_1.
\end{equation}

{\bf A3.}\ \  {\unitlength 1mm \ \begin{picture}(17, 6)
\put(10.00,-1.00){\circle*{1.00}}
 \put(5.00,0.00){\oval(5,5)[t]} \put(4.00,-3.00){\vector(3,1){3}}
\put(5.00,0.00){\oval(5,5)[bl]} \put(12.00,-2.00){\makebox{$\P$}}
\put(-2.00,0.0){\makebox{$a_1$}}
\end{picture}},\ $\P=\Q$, $\Gamma'(\P,\P)=\rK$; or

{\unitlength 1mm \hspace{0.5cm}\begin{picture}(23, 5)
\put(5.00,0.00){\circle*{1.00}} \put(20.00,0.00){\circle*{1.00}}
\put(18.50,0.00){\vector(1,0){1}} \put(5.50,0.00){\line(1,0){13}}
\put(3.00,-4.00){\makebox{$\P$}} \put(18.00,-4.00){\makebox{$\Q$}}
\put(11.00,2.00){\makebox{$a_1$}}
\end{picture}},
$\P\neq \Q$, $\Gamma'(\P,\P)=\rK$, $\Gamma'(\Q,\Q)=\rK$. Then
\begin{equation}
\delta(a_1)=\sum_{j=1}^m f_jv_j,
\end{equation}
where $f_j\in k$. After a basis transformation  given by an
invertible matrix $Q\in M_{m}(\rK)$, we have $ \delta(a_1)=w_1$,
or $\delta(a_1)=0$.

{\bf Proposition 4.3.1} Let $(\K,\M,H)$ be a freely parameterized bimodule problem, and $(\K',\M',H')$ be
induced from  $(\K,\M,H)$ by one of the $3$ reductions of 3.1.  Let $a_1: \P\longrightarrow\Q$ be the first
arrow. Then $(\K',\M',H')$ is still freely parameterized and preserves all the free parameters if and only if
$\delta(a_1)$ and the reduction for $a_1$ has the following forms according to A1,A2,A3 respectively.

{\bf A1} (1) $\delta(a_1)=0$, i.e. $f_j(\lambda,\mu)=0$ for all $j$ in Formula (29).

\hspace{1.2cm} When $\P=\Q$, we set $\overline{N}_{pq}= (\lambda^0)$ for a fixed eigenvalue $\lambda^0$. Then
$g'_{\P}(\lambda)=g_{\P}(\lambda)$.

\hspace{1.2cm} When $\P\ne\Q$, we set $\overline{N}_{pq}=(0)$, Then $g'_{\P}(\lambda)=g_{\P}(\lambda)$,
$g'_{\Q}(\mu)=g_{\Q}(\mu)$.

\hspace{6mm} (2) $\delta(a_1)=f(\lambda,\mu)w_1\ne 0$ in Formula (30), where $f(\lambda,\mu)$ must be equal to

\hspace{1.2cm} $f_{\P}(\lambda)f_{\Q}(\mu)$, and $c(\lambda)$ is given by Formula (30).

 \hspace{1.2cm} When $\P=\Q$, we set $\overline{N}_{pq}=\emptyset_{1\times
1}$. Then
$g'_{\P}(\lambda)=g_{\P}(\lambda)c(\lambda)f_{\P}(\lambda)$.

\hspace{1.2cm} When $\P\ne \Q$, we set $\overline{N}_{pq}=\emptyset_{1\times 1}$. Then
$g'_{\P}(\lambda)=g_{\P}(\lambda)c(\lambda)f_{\P}(\lambda)$,

\hspace{2.6cm} $g'_{\Q}(\mu)= g_{\Q}(\mu)f_{\Q}(\mu)$.

{\bf A2}. (1)  $\delta(a_1)=0$, i.e. $f_j(\lambda)=0$ for all $j$ in Formula (32), we set $$\overline{N}_{pq}=
(0\ \cdots \ 0\ 0 )_{1\times n_q}, \ \ {\rm then} \ \
  g'_{\P}(\lambda)=g_{\P}(\lambda);\quad \mbox{ or }$$
$$\overline{N}_{pq}= (0\
\cdots \ 0\ 1 )_{1\times n_q},\
  g'_{(\P\cup\Q)'}(\lambda)=g_{\P}(\lambda).\qquad \qquad$$

\hspace{6mm} (2)  $\delta(a_1)=f(\lambda)w_1\ne 0$ in Formula (33), we set
 $$\overline{N}_{pq}=\emptyset_{1\times n_q} \ \ {\rm then} \ \
  g'_{\P}(\lambda)=g_{\P}(\lambda)f(\lambda).$$

\hspace{6mm} Or dually,

\hspace{6mm}(1$^{\prime}$)  $\delta(a_1)=0$, we set
 $$\overline{N}_{pq}= (0\ \cdots \ 0\
0 )^T_{1\times n_p}, \ \ {\rm then} \  \
  g'_{\Q'}(\mu)=g_{\Q}(\mu);\quad \mbox{ or }$$
$$\overline{N}_{pq}= (1\ \cdots \ 0\ 0 )^T_{1\times n_p},\
  g'_{(\P\cup\Q)'}(\mu)=g_{\Q}(\mu).\qquad\qquad$$

\hspace{6mm} (2$^{\prime}$) $\delta(a_1)=f(\mu)w_1\ne 0$, we set $\overline{N}_{pq}= \emptyset_{n_p\times 1}$,
then $g'_{\Q}(\mu)=g_{\Q}(\mu)f(\mu)$.

{\bf A3}. (1) $\delta(a_1)=0$, we set
$$\overline{N}_{pq}=\left\{\begin{array}{ll} W,\quad &\mbox{when}\
\P=\Q,\\[2mm]
\left(\begin{array}{cc} 0 & I \\ 0 & 0 \end{array}\right),
&\mbox{when} \P\ne\Q,\end{array}\right.$$

\hspace {1.2cm} where $W$ is  a Weyr matrix with some fixed eigenvalues.

\hspace {6mm}(2)  $\delta(a_1)=w$, we set $\overline{N}_{pq}=\emptyset$.

\hspace{6mm}(3)  $\delta(a_1)=0$ and $\P=\Q$, then $\overline{N}_{pq}$ may be equal to $W\oplus (\lambda)$. Thus
a new

\hspace{1.2cm} parameter $\lambda$  appears at this step, and the  domain of $\lambda$ is $$\rK\setminus\{
\mbox{the eigenvalues of } W\}$$.

{\bf Proof.} In case A1, we must have $n_{\P}=1$, $n_{\Q}=1$, since there are free parameters $\lambda$ attached
to $\P$, and $\mu$ to $\Q$ respectively.

(1) If $\P=\Q$, $\overline{N}_{pq}$ can not be taken as a
parameter, and if $\P\ne \Q$, $\overline{N}_{pq}$ can not be taken
as $1\times 1$ identity matrix because of condition (2) of
Definition 4.2.1.

(2) If $f(\lambda,\mu)=0$, $\lambda,\mu$ would be algebraically
dependent, which is a contradiction to condition (3) of Definition
 4.2.1. Thus $f(\lambda,\mu)\ne 0$. And the condition (3) of 4.2.1
 forces that $f(\lambda,\mu)=f_{\P}(\lambda)f_{\Q}(\mu)$. Thus the
proof of A1 is completed.

Cases A2 and A3 are easy.\hfill$\square$

{\bf Remark.} The reduction sequences of freely parameterized bimodule problems is very special in the set of
general reduction sequences of parameterized bimodule problems. In fact, in Formula (31) of case A1, if
$f(\lambda,\mu)=f_{\P}(\lambda) f_{\Q}(\mu) f_0(\lambda,\mu)$, such that $f_0(\lambda,\mu)\in k[\lambda,\mu]$ is
not a constant, and does not contain any non-constant factor in $k[\lambda]$ or $k[\mu]$, then we are not  able
to continue the reduction, such that the induced triple is still freely parameterized and preserves the
parameters $\lambda$ and $\mu$. We will see in 5.2 that such   a situation can not occur in time case, but
really occurs in wild case see Examples in 7.4 and 7.6.

{\bf Corollary 4.3.1} Let $(\K,\M,H)$ be a bimodule problem having
a reduction sequence $(*)$ of freely parameterized triples given
in the Remark of 4.2, and $\hat{\n}^s$ be a size vector of
$(T^s,\sim^s)$ such that $\hat{n}^s_{\I}=1$ or $0$ whenever $\I\in
T^s/\!\!\sim^s$ is non-trivial. Then $(\hat{*})$ of Corollary
3.3.1 is still a reduction sequence of freely parameterized
triples.

{\bf Proof.} The only problem is, whether the domain at each step
of $(*)$ can be regarded as  a domain at the same step of
$(\hat{*})$. Suppose it is the case up to $(\hat{N}^r, \hat{R}^r,
\hat{\n}^r)$.
 Now consider the $(r+1)$-th reduction. Let
 $$C^r=\{\I\in
T^r/\sim^r\mid \hat{n}^r_{\I}=0 \}.$$

(1) If $\P_r\in C^r$ or $\Q_r\in C^r$, we have an identity functor
from $(\hat{N}^r, \hat{R}^r,$ $ \hat{\n}^r)$ to $(\hat{N}^{r+1},
\hat{R}^{r+1}, \hat{\n}^{r+1})$ by Proposition 3.2.1.  Thus
$\hat{D}^r=\hat {D}^{r+1}$. However there may be some restriction
from $D^r$ to $D^{r+1}$ or there may be a new parameter appearing
in $ D^{r+1}$.

(2) If $\P_r\notin C^r$, $\Q_r\notin C^r$, the equation system
\textrm{II} of Definition 2.2.1 and equation (8) of 2.6 at
$(\P_r,\Q_r)$ are the same both for $\K^r$ and $\hat{\K}^r$ by
Proposition 3.2.1 once again. Thus items (1) and (2) of A1,A2,A3
of Proposition 4.3.1 give the same restriction from $\hat D^r$ to
$\hat D^{r+1}$, as that from ${D}^r$ to ${D}^{r+1}$.

In item (3) of A3, the restriction on the new parameter $\lambda$
in
 $\hat{D}^{r+1}$ is equal to or less than that
in $D^{r+1}$, since the eigenvalues of $\hat{W}$ may be equal to
or less than those of $W$ by Proposition 3.2.1.

By induction the sequence of the domains of $(*)$ suits that of
$(\hat{*})$ well. Therefore $(\hat{*})$ is also a reduction
sequence of freely parameterized triples. \hfill$\square$

If 3.3 $(**)$ is  a sequence of freely parameterized bimodule
problems,  there exists a sequence of corresponding layered
bocses:
$$
(***)\qquad  \mathfrak{A},\  \mathfrak{A}^1,\ \cdots,\
\mathfrak{A}^r,\  \cdots,\  \mathfrak{A}^s.
$$

\subsection{Valuation matrices}

\kg Let $(\K,\M,H)$ be a parameterized bimodule problem with a set
of parameters $\{\lambda_1,\cdots,$ $\lambda_i\}$ of domain $D$.
It is natural to ask whether the parameters can be regarded as the
solid arrows. If the answer is yes, then $\lambda_{\I_1},
\lambda_{\I_2}, \cdots,\lambda_{\I_i}$ and $a_1, a_2,\cdots, a_n$
are free generators of $\Gamma$ over some trivial category
$\rK\times\cdots\times \rK$, and we are able to define a wider
representation category over $(\K,\M,H)$ such that
$\lambda_{\I_l}$ can be taken as Weyr matrices of fixed
eigenvalues. This idea will be realized in the rest of the section
without difficulty when $(\K,\M,H)$ is a freely parameterized
bimodule problem. And a more detailed observation shows that the
$2$ conditions of Definition 4.2.1 are necessary for our purpose:
 Condition (1) gives the
minimality of $\Gamma'$; Condition (2) is the most important
condition in order to define the morphisms in a wilder category
(see  Examples in 4.6).

{\bf Definition 4.4.1} Let $(\K,\M,H)$ be a bimodule problem with
$i$ free parameters.

 (1) By a {\it valuation of a parameter} $\lambda_\I$ we mean the assigning to
$\lambda_{\I}$  a Weyr matrix
$$W(\underline{\lambda}_{\I}^0)=W_{\lambda_{\I
}^{1}}\oplus\cdots\oplus W_{\lambda_{\I}^{\alpha_{\I }}}$$ (see
2.5), where $\lambda_{\I}^{1}, \cdots, \lambda_{\I}^{\alpha_\I}$
are pairwise different eigenvalues taken from the domain
$D\big|_{\{\lambda_\I\}}$. When $\I$ runs over the non-trivial
equivalent classes of $T/\sim$, we then obtain a set of
valuations: $$W(\underline{\lambda}^0_{\I_1}),
W(\underline{\lambda}^0_{\I_2}),
\cdots,W(\underline{\lambda}^0_{\I_i}).$$

(2)  Define a size vector $\m$ of $(T,\sim)$, such that $m_{\I}$
equals the size of $W(\underline{\lambda}_{\I}^0)$ whenever the
parameter $\lambda_{\I}$ is attached to $\I$, and $m_j$ can be
taken as any non-negative integer for any trivial equivalent class
$\J\in  T/\sim$. Then we have a matrix $H_{\underline{m}}, $ which
is not freely parameterized if there exists some non-trivial
vertex $\I$, such that $m_\I>1$. Write $\hat{\m}=\m$, then
$(\hat{N}, \hat{R}, \hat{\m})$ is a triple of $(\K, \M, H)$.

(3) Let $\underline{e}_{\I}$ be the size vector determined by
$W(\underline{\lambda}_{\I}^0)$ (see 2.5). We define a vector
$\tilde{\underline{m}}$ such that $\tilde{m}_\J=m_{\J}$ for any
trivial $\J\in T/\sim$, and $\tilde \m_\I=(\underline e_\I)$ for
any non-trivial $\I\in T/\sim$. A matrix
$\wt{H}_{\wt{\m}}(\underline{\lambda}^0)$ is said to be a {\it
valuated matrix} of $H$, provided that
$H_{\wt{\m}}(\underline{\lambda}^0)$ is obtained from $H_{{\m}}$
by substituting each $W(\underline{\lambda}_{\I}^0)$ for
$\lambda_{\I} I_{m_{\I}}$.

(4) A triple $(\wt{N},\wt{R},\wt{\m})$ is constructed based on
$(\hat{N},\hat{R},\hat{\m})$, such that $\wt{N}=\hat{N}$, but
$\wt{N}$ has a size vector $\wt{\m}$; and $\wt{R}$ satisfies the
matrix equation
$$\wt{R}\wt{H}_{\wt{\m}}(\underline{\lambda}^0)
\equiv_{\prec(p,q)} \wt{H}_{\wt{\m}}(\underline{\lambda}^0)
\wt{R}.$$

    It will be proved  in the next subsection that
the triple $(\wt{N},\wt{R},\wt{\m})$ is induced by a sequence of
reductions.

\subsection{The  reduction sequence towards a valuation matrix}

\kg {\bf Lemma 4.5.1} Let $f(\lambda,\mu)=\sum_{i,j\geq 0} a_{ij}\lambda^i\mu^j\in \rK[\lambda,\mu] $. Let
$X=(x(\alpha,\beta))_{m\times n}\in M_{m\times n}(k)$, and $L, R$ be upper triangular nilpotent matrices. Then
the $(\alpha,\beta)$-entry in the matrix
$$
\sum_{i,j\geq 0} a_{ij}(\lambda I_m+L)^iX(\mu I_n+R)^j
$$
 equals
$$\sum_{i,j\geq 0} \Big(a_{ij}\lambda^i\mu^jx(\alpha,\beta) +
\sum\limits_{(\alpha',\beta')\prec(\alpha,\beta)}d(\lambda,\mu,\alpha',\beta') x(\alpha',\beta')\Big),$$ where
$d(\lambda,\mu,\alpha',\beta')$ are polynomials in $\lambda, \mu $ as coefficients of $x(\alpha',\beta')$, and
$(\alpha',\beta')\prec (\alpha,\beta)$ are given under the matrix order of Definition  2.3.1. \hfill$\square$

Let $\sum\limits_{i,j} c_{ij}(\lambda,\mu)x_{ij}=0$ be a linear equation in the indeterminates $x_{ij}$, and the
coefficients $c_{ij}(\lambda,\mu)\in \rK[\lambda,\mu]$. Let $W(\underline{\lambda}^0)= W_{\lambda^{1}}\oplus
\cdots \oplus W_{\lambda^{g}}$, $W(\underline{\mu}^0)= W_{\mu^{1}}\oplus \cdots \oplus W_{\mu^{h}}$  be two Weyr
matrices with sizes $m$ and $n$ respectively. Let matrix $X_{ij}=(X_{ij}^{\gamma\eta})_{g\times h }$ be
partitioned according to $W(\underline{\lambda}^0)$ and $W(\underline{\mu}^0)$, i.e. $X_{ij}^{\gamma\eta}$ are
$m_{\gamma}\times n_{\eta}$ matrices whenever the size of $W_{\lambda^{\gamma}}$ is $m_{\gamma}$, and that of
$W_{\mu^{\eta }}$ is $n_{\eta}$, thus $X_{ij}^{\gamma\eta}=(x_{ij}^{\gamma\eta}(\alpha,\beta))_{m_{\gamma}\times
n_{\eta}}$. For example,
\begin{center}\unitlength=1mm
\begin{picture}(130,50)
\put(0,15){\framebox(20,20)}\put(0,25){\line(1,0){20}}
\put(10,15){\line(0,1){20}} \put(2,29){$W_{\lambda^1}$}
\put(12,19){$W_{\lambda^2}$}

\put(30,15){\framebox(45,20)}\put(30,25){\line(1,0){45}}
\put(40,15){\line(0,1){20}} \put(55,15){\line(0,1){20}}
\put(32,29){$X_{ij}^{11}$} \put(44,29){$X_{ij}^{12}$}
\put(62,29){$X_{ij}^{13}$} \put(32,19){$X_{ij}^{21}$}
\put(44,19){$X_{ij}^{22}$} \put(62,19){$X_{ij}^{23}$}

\put(85,0){\framebox(45,45)}\put(85,20){\line(1,0){45}}
\put(85,35){\line(1,0){45}}\put(95,0){\line(0,1){45}}
\put(110,0){\line(0,1){45}} \put(87,39){$W_{\mu^1}$}
\put(99,26){$W_{\mu^2}$} \put(118,9){$W_{\mu^3}$}
\end{picture}
\end{center}

Now suppose that we have a bimodule problem $(\K,\M,H=0)$ and a
reduction sequence $(*)$ of freely parameterized triples given in
the Remark of 4.2. Consider the $(r+1)$-th reduction of $(*)$. Let
$a_1: \P_r\rightarrow \Q_r$ be the first arrow of
$(\K^r,\M^r,H^r)$, where $\P_r$ (resp. $\Q_r$) is either
non-trivial with a free parameter $\lambda$ (resp. $\mu$) or
trivial. Assume that the equation system \textrm{II} of Definition
2.2.1 at $(\P_r,\Q_r)$ for $\K^r$ is given by
$$
(\textrm{II}):\qquad \sum_{\P_r\ni i<j\in \Q_r}
c_{ij}^l(\lambda,\mu) x_{ij}=0,
$$
and the equation (8) of 2.6 is given by
$$
(8):\qquad \sum_{\P_r\ni i<j\in \Q_r} c_{ij}(\lambda,\mu)
x_{ij}=0.
$$
We denote the coefficient matrix of (II) by $C(\lambda,\mu);$ that
of (II) and (8) by $C'(\lambda,\mu)$. Furthermore assume that we
also have the matrix equations:
$$
(\wt{\textrm{II}}):\qquad \sum_{\P_r\ni i<j\in \Q_r}
c_{ij}^l(W(\underline{\lambda}^0),W(\underline{\mu}^0))X_{ij}=0,
$$
$$
(\wt{8}):\qquad \sum_{\P_r\ni i<j\in \Q_r}
c_{ij}(W(\underline{\lambda}^0),W(\underline{\mu}^0))X_{ij}=0,
$$
or equivalently,
$$
(\wt{\textrm{II}}_{\gamma\eta}):\qquad \sum_{\P_r\ni i<j\in \Q_r}
c_{ij}^l(W_{\lambda^{\gamma}},W_{\mu^{\eta}})X_{ij}^{\gamma\eta}=0,
$$
$$
(\wt{8}_{\gamma\eta}):\qquad \sum_{\P_r\ni i<j\in \Q_r}
c_{ij}(W_{\lambda^{\gamma}},W_{\mu^{\eta}})X_{ij}^{\gamma\eta}=0
$$
for $1\leq \gamma\leq g$, $1\leq \eta\leq h$.   Which are
equivalent to a system of linear equations:
$$
(\wt{\textrm{II}}_{\gamma\eta}(\alpha,\beta)):\ \sum_{\P_r\ni
i<j\in\Q_r}\Big(c_{ij}^l(\lambda^{\gamma},\mu^{\eta
})x_{ij}^{\gamma\eta}(\alpha,\beta) +
\sum_{(\alpha',\beta')\prec(\alpha,\beta)}d_{ij}^l(\lambda^{\gamma
}, \mu^{\eta },\alpha',\beta')
x_{ij}^{\gamma\eta}(\alpha',\beta')\Big)=0,
$$
$$
(\wt{8}_{\gamma\eta}(\alpha,\beta)):\ \sum_{\P_r\ni
i<j\in\Q_r}\Big(c_{ij}(\lambda^{\gamma},\mu^{\eta})x_{ij}^{\gamma\eta}(\alpha,\beta)
+
\sum_{(\alpha',\beta')\prec(\alpha,\beta)}d_{ij}(\lambda^{\gamma},
\mu^{\eta},\alpha',\beta')
x_{ij}^{\gamma\eta}(\alpha',\beta')\Big)=0,
$$
for $1\leq \alpha\leq m_{\gamma}$, $1\leq \beta\leq n_{\eta}$, and each pair $(\gamma,\eta)$. Next we fix a pair
$(\gamma,\eta)$, and for each index pair $(\alpha,\beta)$, we list the equations at the $(\alpha,\beta)$-entry.
If we count the index pairs $(\alpha, \beta)$ according to the matrix order of Definition 2.3.1, then we obtain
an equation system $\wt{\textrm{II}}^0_{\gamma\eta}=\bigcup_{(\alpha,\beta)}
\wt{\textrm{II}}_{\gamma\eta}(\alpha,\beta)$ whose coefficient matrix is
$$\wt{C}(W_{\lambda^{\gamma}},W_{\mu^{\eta}})=\left(
\begin{array}{cccc}
C(\lambda^{\gamma}, \mu^{\eta}) & 0& \cdots & 0\\
* & C(\lambda^{\gamma}, \mu^{\eta}) & \cdots & 0\\
\vdots & \vdots & \ddots & \vdots \\
* & * & \cdots & C(\lambda^{\gamma}, \mu^{\eta})
\end{array}\right)_{(m_{\gamma}n_{\eta})\times (m_{\gamma}n_{\eta})}.
$$
And the coefficient matrix of $\wt{\textrm{II}}^1_{\gamma\eta}=\bigcup_{(\alpha,\beta)}
\wt{\textrm{II}}_{\gamma\eta}(\alpha,\beta)\cup (\wt{8}_{\gamma\eta}(\alpha,\beta))$  is
$$\wt{C}'(W_{\lambda^{\gamma}},W_{\mu^{\eta}})=\left(
\begin{array}{cccc}
C'(\lambda^{\gamma}, \mu^{\eta}) & 0& \cdots & 0\\
* & C'(\lambda^{\gamma}, \mu^{\eta}) & \cdots & 0\\
\vdots & \vdots & \ddots & \vdots \\
* & * & \cdots & C'(\lambda^{\gamma}, \mu^{\eta})
\end{array}\right)_{(m_{\gamma}n_{\eta})\times
(m_{\gamma}n_{\eta})}.
$$
We denote the coefficient matrix of
$\bigcup_{(\gamma,\eta)}(\wt{\textrm{II}}^0_{\gamma\eta})$ by
$\wt{C}(W(\underline{\lambda}^0),W(\underline{\mu}^0))$, and that
of $\bigcup_{(\gamma,\eta)}(\wt{\textrm{II}}^1_{\gamma\eta})$ by
$\wt{C}'(W(\underline{\lambda}^0),W(\underline{\mu}^0))$. Then
$$
\wt{C}(W(\underline{\lambda}^0),W(\underline{\mu}^0))=\mbox{diag}
(\wt{C}(W_{\lambda^{\gamma}},W_{\mu^{\eta}}))_{(\gamma,\eta)},
$$
$$
\wt{C}'(W(\underline{\lambda}^0),W(\underline{\mu}^0))=\mbox{diag}
(\wt{C}'(W_{\lambda^{\gamma}},W_{\mu^{\eta}}))_{(\gamma,\eta)},
$$

{\bf Proposition 4.5.1} (1) If $\rank C'(\lambda,\mu)=\rank C(\lambda,\mu)$ in the domain $D^{r+1}$, then
$$
\rank \wt{C}'(W_{\lambda^{\gamma}},W_{\mu^{\eta}})= \rank
\wt{C}(W_{\lambda^{\gamma}},W_{\mu^{\eta}}),
$$
when $(\lambda^{\gamma}, \mu^{\eta})\in D^{r+1}|_{(\lambda,\mu)}$. Moreover
$$
{\rm rank} \wt{C}'(W(\underline{\lambda}^0),W(\underline{\mu}^0))= {\rm rank}
\wt{C}(W(\underline{\lambda}^0),W(\underline{\mu}^0))$$ when $(\lambda^{\gamma}, \mu^{\eta})\in
D^{r+1}|_{(\lambda,\mu)}$ for $1\leq \gamma \leq g$, $1\leq \eta \leq h$.

(2) If $\rank C'(\lambda,\mu)=\rank C(\lambda,\mu)+1$ in the domain $D^{r+1}$, then
$$
\rank\wt{C}'(W_{\lambda^{\gamma}},W_{\mu^{\eta}})= \rank
\wt{C}(W_{\lambda^{\gamma}},W_{\mu^{\eta}})+ m_{\gamma}\times
n_{\eta},$$
  when $(\lambda^{\gamma}, \mu^{\eta})\in
D^{r+1}|_{(\lambda,\mu)}$. Moreover
$$
{\rm rank} \wt{C}'(W(\underline{\lambda}^0),W(\underline{\mu}^0))=
{\rm rank}
\wt{C}(W(\underline{\lambda}^0),W(\underline{\mu}^0))+m_{\P_r}m_{\Q_r},
$$
when $(\lambda^{\gamma},\mu^{\eta})\in D^{r+1}|_{(\lambda,\mu)}$
for $1\leq \gamma \leq g$, $1\leq \eta \leq h$.

{\bf Proof.} (1)   It is not difficult to see that $${\rm rank}
\wt{C}(W_{\lambda^{\gamma}},W_{\mu^{\eta}})= m_{\gamma}n_{\eta}
{\rm rank}(C(\lambda,\mu)), \ {\rm and} \
\rank\wt{C}'(W_{\lambda^{\gamma}},W_{\mu^{\eta}})=
m_{\gamma}n_{\eta}{\rm rank}(C^{\prime}(\lambda,\mu)).$$ In fact,
the  the left hand sides in two equalities have the nice
expressions above, which imply that $``\geq"$ hold in both
formulae. On the other hand, if some row of $C(\lambda,\mu)$ is a
linear combination of other rows, then the corresponding block row
of $\wt C(W_{\lambda^{\gamma}},W_{\mu^{\eta}})$ is the same linear
combination of the same block rows.  The similar assertion is true
for $C^{\prime}$ and $\wt {C}'$. Therefor $``\leq"$ hold in both
formulae.

(2) The reason is: $\sum_{(\gamma,\eta)} m_{\gamma}
n_{\eta}=m_{\P_r} n_{\Q_r}$. \hfill $\square$

 Let $(\K,\M,H=0)$ be a bimodule problem with a
reduction sequence $(*)$ of freely parameterized triples given in
the Remark of 4.2. Suppose that $H^s$ in the end term of $(**)$
possesses a set of parameters $\lambda_{\I_1}, \lambda_{\I_2},
\cdots, \lambda_{\I_i}$. Let a set of valuation matrices
$$ \mbox{(Vm)}: \quad
\{ W(\underline{\lambda}^0_{\I_1}),
W(\underline{\lambda}^0_{\I_2}), \cdots,
W(\underline{\lambda}^0_{\I_i})\}
$$
be given by item (1) of Definition 4.4.1, and let $H^s_{\underline
m^s}$, $\wt{H}^s_{\wt{\m}^s}(\underline{\lambda}^0)$ be given in
the items (2), (3) of Definition 4.4.1  respectively. Assume that
$(\underline {\hat m}, \underline {\hat m}^1,\cdots,\underline
{\hat m}^s)$ is a sequence of size vectors calculated by Formula
(11) of 3.2 inductively from $\underline {\hat m}^s = \underline
m^s $, and $(H_{\hat{\m}}, H^1_{\hat{\m}^1}, \cdots
H^s_{\hat{\m}^s})$ is a sequence of matrices given by $(H, H^1,
\cdots, H^s)$ and $(\hat{\m}, \hat{\m}^1, \cdots, \hat{\m}^s)$.
Then we are able to construct $\wt{\underline{m}}^r$ and
$\wt{H}^r_{\wt{\underline m}^r}(\underline \lambda^0(r))$
according to the valuation matrices (Vm) for $r=s,s-1,\cdots,1,0$,
by Definition 4.4.1.

Thus we have a sequence of triples $(\hat *)$ given in Corollary
4.3.1 (but the restriction of $m_{\I}=1$ or $0$ for any
non-trivial $\I$ is relaxed). And we have a sequence of valuation
matrices
$$
0=\wt{H}_{\wt{\m}}(\underline{\lambda}^{0}(0)), \wt{H}^1_{\wt{\m}^1}(\underline{\lambda}^{0}(1)), \cdots,
\wt{H}^r_{\wt{\m}^r}(\underline{\lambda}^{0}(r)), \cdots,
\wt{H}^{s-1}_{\wt{\m}^{s-1}}(\underline{\lambda}^{0}(s-1)), \wt{H}^s_{\wt{\m}^s}(\underline{\lambda}^{0}(s)),
$$
and a sequence of triples
$$
(\wt{*})\quad (\wt{N},\wt{R}, \wt{\m}), (\wt{N}^1,\wt{R}^1,
\wt{\m}^1), \cdots, (\wt{N}^r,\wt{R}^r, \wt{\m}^r), \cdots,
(\wt{N}^{s-1},\wt{R}^{s-1}, \wt{\m}^{s-1}), (\wt{N}^s,\wt{R}^s,
\wt{\m}^s),
$$
such that each triple $(\wt{N}^r,\wt{R}^r, \wt{\m}^r)$ is obtained
from $(\hat{N}^r, \hat{R}^r, \underline {\hat m}^r)$ by Definition
4.4.1.

{\bf Theorem 4.5.1} The sequence $(\wt{*})$ is a reduction sequence. And either the type of reductions from
$(\wt{N}^r,\wt{R}^r,\wt{\m}^r)$ to $(\wt{N}^{r+1},\wt{R}^{r+1}, \wt{\m}^{r+1})$ is  the same as that from
$({N}^r,{R}^r, {\m}^r)$ to $({N}^{r+1},{R}^{r+1}, {\m}^{r+1})$,  or the triple $(\wt{N}^{r+1},\wt{R}^{r+1},
\wt{\m}^{r+1})$ is obtained from $(\wt{N}^r,\wt{R}^r,\wt{\m}^r)$ by identity,  $0\leq r\leq s-1$.

(We stress that the reductions in $(\wt{*})$ are ordered according to the matrix order in $(\hat{*})$, which
differs from the usual order of Definition 2.3.1.)

{\bf Proof.}   For $r=0$, we set $\wt{R}=\hat{R}$,
$\wt{N}=\hat{N}$,\ $\wt{\m}=\hat{\m}$. Suppose we have already a
sequence of reductions up to $(\wt{N}^r, \wt{R}^r, \wt{\m}^r)$,
now consider the $(r+1)$-th reduction towards $(\wt{N}^{r+1},
\wt{R}^{r+1}, \wt{\m}^{r+1})$. If $(\hat{N}^{r+1}, \hat{R}^{r+1},
\hat{\m}^{r+1})$ is obtained from $(\hat{N}^{r}, \hat{R}^{r},
\hat{\m}^{r})$ by identity, so is $(\wt{N}^{r+1}, \wt{R}^{r+1},
\wt{\m}^{r+1})$ from $(\wt{N}^r, \wt{R}^r, \wt{\m}^r)$, we have
done. Otherwise, let us go back to Proposition 4.3.1, and let
$a_1: \P_r\rightarrow \Q_r$ be the first arrow of $(\K^r, \M^r,
H^r)$.

In item (1) of A1,A2,A3 (as well as ($1^{\prime}$) of A2) of
Proposition 4.3.1, $\delta(a_1)=0$, which means that the equation
system (II) implies $(8)$. Then item (1) of Proposition 4.5.1
tells that $(\wt{\textrm{II}})$ implies $(\wt{8})$. Thus we are
able to set
$$\overline{\wt{N}}^{\, r}_{p_rq_r}=\overline{\hat{N}}^{\,
r}_{p_rq_r}$$
 by an edge reduction, or a loop reduction with fixed
eigenvalues.

In item (2) of A1,A2,A3 (as well as ($2'$) of A2) of Proposition
4.3.1, $\delta(a_1)\neq 0$, which means that the equation system
(II) does not imply $(8)$. Then item (2) of Proposition 4.5.1
tells that
$$
\rank\wt{C}'(W_{\underline\lambda^0},W_{\underline\mu^0})= \rank
\wt{C}(W(\underline{\lambda}^0),W(\underline{\mu}^0))+ m_{\P_r}
n_{\Q_r}.$$ Thus we are able to set
$$\overline{\wt{N}}^{\,
r}_{p_rq_r}=\emptyset=\overline{\hat{N}}^{\, r}_{p_rq_r}$$ by
regularization.

In item (3) of A3 of Proposition 4.3.1, $\delta(a_1)=0$. We are
able to set
$$\overline{\wt{N}}^{\, r}_{p_rq_r}=\hat{W}\oplus
W(\underline{\lambda}^0_{\I}),$$
 whenever $\overline{\hat{N}}^{\,
r}_{p_rq_r}=\hat{W}\oplus \lambda I_{m_{\I}}$ by a loop reduction,
where $\I\in T^{r+1}/\sim^{r+1}$.

Our theorem follows by induction. \hfill$\square$

\subsection{Valuation representation categories}

\kg Let $(\K,\M,H)$ be a freely parameterized bimodule problem
given by Definition 4.2.1, and
 $\mathfrak{A}=(\Gamma,\Omega)$ be the layered bocs
 corresponding to
 $(\K,\M,H)$.  The
 present subsection is devoted to showing a valuation representation category
$\wt{Mat}(\K, \M)$, and equivalently $\wt{R}(\mathfrak{A})$.

 Given any size vector $\m$ and any valuation matrix
$\widetilde H_{\wt{\m}}(\underline{\lambda}^0)$ defined in 4.4.1,
an object of size vector $ \wt {\underline m}$ in $\wt{Mat}(\K,
\M)$ is defined as a matrix $M\in \M_{\wt{\m}}$ with an {\it
adjoint valuation matrix}
$\wt{H}_{\wt{\m}}(\underline{\lambda}^0)$. If $L$ is also a matrix
of size $\wt{\l}$ with an adjoint valuation matrix
$\wt{H}_{\wt{\l}}(\underline{\mu}^0)$, then a morphism $\varphi:
M\rightarrow L $ in $\wt{Mat}(\K, \M)$ is a matrix
$\varphi\in\K_{\wt{\m}\times\wt{\l}}$, such that
$$(M+\wt{H}_{\wt{\m}}(\underline{\lambda}^0))\varphi=\varphi(L+
\wt{H}_{\l}(\underline{\mu}^0)).$$
 Thus we obtain a category of
representations, which is called a {\it valuation representation
category} of $(\K,\M,H)$, and we denote it by $\wt{Mat}(\K, \M,
H)$.

Representation category $\wt{R}(\mathfrak{A})$ has been
well-defined in \cite{D} and \cite{CB1}. The objects of
$\wt{R}(\mathfrak{A})$ are $\Gamma$-modules. In fact, if $M\in \wt
{Mat}(\K,\M)$,  an object $\overline{M}\in \wt{R}(\mathfrak{A})$
corresponding to  $M$ is defined by a set of vector spaces:
$\overline{M}_{\I}=\rK^{m_{\I}}$, $\forall\, \I\in T/\!\!\sim$; a
set of linear maps:
$\overline{M}(\lambda_{\I})=W(\underline{\lambda}^0_{\I})$ for
$\lambda_{\I}$ attached to $\I$, and
$\overline{M}(\rho^w_{\I\J})^*=M_{p^w_{\I\J}q^w_{\I\J}}$ for any
basis element $\rho^w_{\I\J}\in A$ (see Formula (1) in 2.3). Thus
$\overline{M}$ is a $\Gamma$-module. If $\overline{L}\in
\wt{R}(\mathfrak{A})$ corresponds to  $L \in {\wt {Mat}(\K,\M)}$,
and $\varphi(=S): M \rightarrow L$ is a morphism in $\wt {Mat}
(\K,\M)$, a morphism $\overline{\varphi}: \overline{M}\rightarrow
\overline{L}$ corresponding to $ \varphi$ in
$\wt{R}(\mathfrak{A})$ is defined by a set of linear maps:
$(\varphi_{\I}, \varphi(\zeta^w_{\I\J})^*\mid \forall \I\in
T/\sim, \forall \quad  \zeta^w_{\I\J}\in B)$, where
$\varphi(\zeta^w_{\I\J})^*=S_{\overline p^w_{\I\J} \overline
q^w_{\I\J}}$. It is clear that $\varphi_{\I}:
\overline{M}_{\I}\rightarrow \overline{L}_{\I}$ satisfy
$$\overline{M}(\rho^w_{\I\J})^*\varphi_{\J}-\varphi_{\I}
\overline{L}(\rho^w_{\I\J})^* =\varphi(\delta(\rho^w_{\I\J})^*)$$
 (see the last part of 3.6). If $\overline{\psi}: \overline{L}\rightarrow \overline{E}$
is also a morphism,  then the composition
$\overline{\eta}=\overline{\varphi}\overline{\psi}:
\overline{L}\rightarrow \overline{E}$ is given by
$\overline{\eta}=(\eta_{\I}, \eta(\zeta^w_{\I\J})^* \mid \forall
\,\I\in T/\!\!\sim, \forall \zeta^w_{\I\J}\in B)$, such that
$\eta_{\I}=\varphi_{\I}\psi_{\I}$, and
$$
\eta(\zeta^w_{\I\J})^*=\varphi(\zeta^w_{\I\J})^*\psi_{\J}+\varphi_{\I}
\psi(\zeta^w_{\I\J})^*+\sum\limits_{\L,u,v}\c
\varphi(\zeta^u_{\I\L})^*\psi(\zeta^v_{\I\L})^*
$$
 (see also the last part of 3.6).

 Now suppose the  freely parameterized bimodule problem $(\K, \M,H)$
 (respectively $\mathfrak A$) is minimal.
 Given any non-trivial $\I\in T/\!{\sim}$ with an attached
parameter $\lambda_{\I}$ of domain $\rK\setminus\{\mbox{the roots
of } g_{\I}(\lambda_{\I})\}$, let $\underline{m(\I,d)}$ be a size
vector of $(T, \sim)$, such that $m(\I,d)_{\I}=d$ for some
$d\in\mathbb{N}$ and $m(\I,d)_{\L}=0$, $\forall\, \L\in
T/\!{\sim}\setminus\{\I\}$. For any $\lambda_{\I}^0\in\rK$ with
$g_{\I}(\lambda_{\I}^0)\ne 0$, we define a zero matrix
$O(\I,d,\lambda_{\I}^0)\in \wt{\M}_{\underline{\wt{m}(\I,d)}}$
with an adjoint valuation matrix
$\wt{H}_{\underline{\wt{m}(\I,d)}}(J_d(\lambda_{\I}^0))$, and we
denote it by $\wt{H}(J_d(\lambda_{\I}^0))$ for simplicity.

{\bf Corollary 4.6.1} (1) $O(\I,d,\lambda_{\I}^0)$ is
indecomposable.

 (2) $\{O(\J)\mid \forall \mbox{ trivial } \J\in T/\sim\}$ and
 $$ \{O(\I, d, \lambda^0_{\I})\mid
\forall \mbox{ non-trivial } \I\in T/\!\!\sim, \forall\;
d\in\mathbb{N}, \forall\, \lambda^0_{\I}\in k \mbox{ with }
g_{\I}(\lambda^0_{\I})\neq 0\}$$ form a complete set of
indecomposables of $\wt{Mat}(\K,\M)$.

 {\bf Proof.} (1) Let $O(\I,d,\lambda_{\I}^0)$ be an object of $\wt{Mat}(K,
 \M, H)$ with an adjoint matrix
 $\wt{H}(J_d(\lambda^0_{\I})).$
  We denote $O(\I,d,\lambda_{\I}^0)$ by  $O$ for
simplicity.
   Suppose  $\overline{O}\in
 \wt{R}(\mathfrak{A})$ corresponds to $O$. Thus
 $\overline{O}_{\I}=\rK^d$, $\overline{O}_{\J}=0$, and
 $\overline{O}(\lambda_{\I})=J_d(\lambda_{\I}^0)$. Given any
 morphism $\overline{\varphi}: \overline{O}\rightarrow
 \overline{O}$, then
\begin{equation}
\overline{\varphi}_{\I}=\left( \begin{array}{cccc} x_1 & x_2 &
\cdots & x_d \\ & x_1 & \cdots & x_{d-1} \\ & & \ddots & \vdots \\
& & & x_1 \end{array}\right)
\end{equation}
is given by $
J_d(\lambda_{\I}^0)\overline{\varphi}_{\I}=\overline{\varphi}_{\I}
J_d(\lambda_{\I}^0) $, and $\overline{\varphi}_{\J}=0$. Therefore
$End_{\wt{R}(\mathfrak{A}^{\infty})}(\overline{O})$ is local.
Consequently, $\overline{O}$ is indecomposable, so is $O$.

(2) See \cite[6.2]{CB1}. \hfill$\square$

From now on we will not worry about any difference between
$\wt{Mat}(\K,\M)$ and $\wt{R}(\mathfrak{A})$. And we denote them
still by $Mat(\K,\M)$ and $R(\mathfrak{A})$ respectively for
simplicity.

If $(**)$ is a reduction sequence of freely parameterized bimodule
problems given in 3.3, the {\it reduction functors}
$${\vartheta}_{r-1,r}: {Mat}(\K^{r},\M^{r})\rightarrow
{Mat}(\K^{r-1}, \M^{r-1})$$
 are defined as follows. Given any
matrix $M^{r}\in \M^{r}_{\wt{\m}^{r}}$ with an adjoint valuation
matrix $\wt{H}_{\wt{\m}^{r}}^{r}(\underline{\lambda}^0(r))$, then
$$
{\vartheta}_{r-1,r}(M^{r})=M^{r}+\wt{H}^{r}_{\wt{\m}^{r}}(\underline{\lambda}^0(r))_{p_{r-1}q_{r-1}}
\otimes \rho^{r-1}
$$
with an adjoint valuation matrix
$\wt{H}^{r}_{\wt{\m}^{r}}(\underline{\lambda}^0(r))$, where
$\rho^{r-1}$ is the first basis element of $\M^{r-1}$. Moreover,
the composition
$$
{\vartheta}_{0r}={\vartheta}_{01}\cdots {\vartheta}_{r-1,r}:
{Mat}(\K^r,\M^r)\rightarrow {Mat}(\K,\M)
$$
is obtained  by sending $M^r$ to $M^r+
\wt{H}^{r}_{\wt{\m}^{r}}(\underline{\lambda}^0(r))$. And the
action of the reduction functors on morphisms is an identity as
the same as that in 3.3. In particular if the end term of $(**)$
is minimal, then $${\vartheta}_{0,\infty}(O(\I,d,\lambda_{\I}^0))=
\wt{H}^{\infty}(J_d(\lambda_{\I}^0)),$$  which is indecomposable
by Corollary 4.6.1.

Finally we claim why the conditions  (2) of Definition 4.2.1 must
hold for a freely parameterized bimodule problem by some examples
to end the subsection.

{\bf Example 1.} Let $\mathfrak{A}=(\Gamma,\Omega)$ be a local
bocs with a layer $L=(\Gamma'; \omega; a; v)$

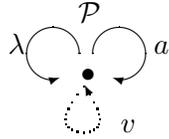
\vspace{3mm}
\begin{center}
\begin{picture}(50,20)\unitlength=0.5pt
\put(0,60){\oval(40,40)[t]} \put(0,60){\oval(40,40)[bl]}
\put(0,40){\vector(3,1){5}} \put(50,60){\oval(40,40)[t]}
\put(50,60){\oval(40,40)[br]} \put(50,40){\vector(-3,1){5}}

\put(21,39){$\bullet$} \qbezier[10](23,35)(-5,2)(23,0)
\qbezier[10](23,35)(47,2)(23,0) \put(27,30){\vector(-1,2){3}}

\put(51,0){\makebox{$v$}} \put(19,84){\makebox{$\P$}}
\put(-34,60){\makebox{$\lambda$}} \put(76,60){\makebox{$a$}}
\end{picture}
\end{center}
where $\Gamma'(\P,\P)=\rK[\lambda]$, $\delta(a)=v\lambda-\lambda
v$, i.e. $f(\lambda,\mu)=-\lambda+\mu$. If we require
$f(\lambda,\mu)\ne 0$, and set $M_{\P}=\rK$,
$M(\lambda)=(\lambda^0)$, $\forall\, \lambda^0\in\rK$, then there
is not any morphism $\varphi: M\rightarrow M $. In fact , if
$\varphi$ were a morphism, then $M(a)\varphi_{\P}-\varphi_{\P}
M(a)=v\lambda^0-\lambda^0 v=0$, which would lead to a
contradiction to the requirement of $f(\lambda^0,\lambda^0)\ne 0$.

{\bf Example 2.} We have  the same picture as Example 2, but
$\delta(a)=2v\lambda-\lambda v$, i.e.
$f(\lambda,\mu)=-\lambda+2\mu$. If we require $f(\lambda,\mu)=0$,
and set $M_{\P}=\rK$, $M(\lambda)=(\lambda^0)$, $\forall\
\lambda^0\in\rK$, and $L_{\P}=\rK$, $L(\lambda)=(\mu^0)$, then
there is not any morphism $\varphi: M\rightarrow L $ when
$\lambda^0\ne 2\mu^0 $. In fact , if $\varphi$ were a morphism,
then $M(a)\varphi_{\P}-\varphi_{\P} L(a)=2v\lambda^0-\mu^0 v\neq
0$, which would lead to a contradiction to the requirement of
$f(\lambda^0,\mu^0) =0$.
\newpage

\bcen
\section{Minimally wild bimodule problems}\ecen

\kg This section is devoted to classifying the minimally wild bimodule problems.

\medskip
\subsection{The wild theorem}
\kg We will give an alternative proof of the well-known Drozd's
wild theorem in this subsection. Our proof may not be an essential
improvement, since the idea of the original proof is perfect. But
the new proof lowers the dimension from the original 43 to 20.
This makes the matrices simpler.

{\bf Definition 5.1.1}(\cite{CB1}) We say a  layered bocs $\mathfrak{A}$ is {\it wild}, if there is a functor
$$F: \rK\langle x,y\rangle\mbox{-mod} \rightarrow  R(\mathfrak{A})$$ which preserves iso-classes and
indecomposability.

{\bf Theorem 5.1.1}(\cite{D}, \cite{CB1}) Let
$\mathfrak{A}=(\Gamma, \Omega)$ be a  bocs
 with layer $L=(\Gamma';\omega;
a_1,\cdots,a_n; v_1,$ $\cdots,$ $ v_m)$. Suppose that $a_1:
\P\rightarrow \Q$. Then $\mathfrak{A}$ must be wild in the
following two cases:

{\bf Case 1.} $\Gamma'(\P,\P)=\rK[\lambda, g_{\P}(\lambda)^{-1}]\ \mbox{and}\ \Gamma'(\Q,\Q)=\rK[\mu,
 g_{\Q}(\mu)^{-1}]$; $$\delta(a_1)=f(\lambda, \mu)v_1$$ where
 $f(\lambda, \mu)\in \rK[\lambda, \mu, g_{\P}(\lambda)^{-1}, g_{\Q}(\mu)^{-1}]$ is non-invertible.

 {\bf  Case 2.} $\Gamma'(\P,\P)=\rK[\lambda, g_{\P}(\lambda)^{-1}]\ and \ \Gamma'(\Q,\Q)=\rK$;
 $$\delta(a_1)=0,\ \ \ \mbox{ or dual}.$$

{\bf Proof. } Case 1. Without loss of generality we may assume
that $\mathfrak{A}$ is local when  $\P=\Q;$ or $\mathfrak{A}$ has
two vertices when $\P\neq\Q:$
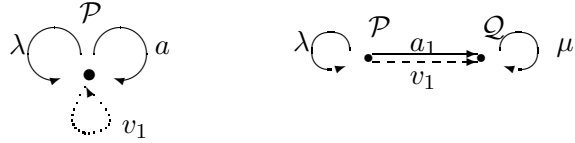
\begin{center}
\begin{picture}(60,20)\unitlength=0.5pt
\put(-10,60){\oval(40,40)[t]} \put(-10,60){\oval(40,40)[bl]}
\put(-10,40){\vector(3,1){5}} \put(40,60){\oval(40,40)[t]}
\put(40,60){\oval(40,40)[br]} \put(40,40){\vector(-3,1){5}}

\put(11,39){$\bullet$} \qbezier[10](13,35)(-8,2)(16,0)
\qbezier[10](13,35)(50,2)(13,0) \put(17,29){\vector(-1,2){3}}

\put(41,0){\makebox{$v_1$}} \put(9,84){\makebox{$\P$}}
\put(-44,60){\makebox{$\lambda$}} \put(66,60){\makebox{$a$}}

\unitlength 1mm 
\put(40.00,10.00){\circle*{1.00}}
\put(55.00,10.00){\circle*{1.00}} \put(35.00,11.00){\oval(5,5)[t]}
\put(35.00,11.00){\oval(5,5)[bl]}
\put(34.00,8.00){\vector(3,1){3}} \put(60.00,11.00){\oval(5,5)[t]}
\put(60.00,11.00){\oval(5,5)[br]}
\put(61.00,8.00){\vector(-3,1){3}}
\put(40.50,10.50){\line(1,0){13}}
\put(53.50,10.50){\vector(1,0){1}}
 \put(40.50,9.50){\line(1,0){1}}
\put(42.50,9.50){\line(1,0){1}} \put(44.50,9.50){\line(1,0){1}}
\put(46.50,9.50){\line(1,0){1}} \put(48.50,9.50){\line(1,0){1}}
\put(50.50,9.50){\line(1,0){1}} \put(52.50,9.50){\line(1,0){1}}
\put(53.50,9.50){\vector(1,0){1}}
\put(30.00,11.0){\makebox{$\lambda$}}
\put(40.00,13.00){\makebox{$\P$}}
\put(55.00,13.00){\makebox{$\Q$}}
 \put(45.50,11.00){\makebox{$a_1$}}
 \put(45.50,6.50){\makebox{$v_1$}}
\put(65.00,11.00){\makebox{$\mu$}}
\end{picture}
\end{center}
We may also assume that $f(\lambda, \mu)\in \rK[\lambda, \mu]$
after dividing $v_1$ by some power of
$g_{\P}(\lambda)g_{\Q}(\mu)$, write
\begin{equation}
f(\lambda, \mu)=\alpha (\lambda-\lambda^0)+\beta(\mu-\mu^0)
+\gamma_1(\lambda-\lambda^0)^2+\gamma_2(\lambda-\lambda^0)(\mu-\mu^0)+\gamma_3
(\mu-\mu^0)^2+\cdots
\end{equation}
where $f(\lambda^0,\mu^0)=0$, $g_{\P}(\lambda^0)\ne 0$,
$g_{\Q}(\mu^0)\ne 0$. If $\P\ne\Q$, and $\alpha\beta\ne 0$, we set
$m_{\P}=9$, $m_{\Q}=11$, and define an object $M\in
R(\mathfrak{A})$, such that $ M(\lambda)\simeq
J_1(\lambda^0)\oplus J_3(\lambda^0)\oplus J_5(\lambda^0), \quad
M(\mu)\simeq J_2(\mu^0)\oplus J_4(\mu^0)\oplus J_5(\mu^0)$. Then
$$
M(\lambda)=\left(\begin{array}{c} \unitlength 1mm
\begin{picture}(58,54)
\put(0,36){\dashbox{1}(20,18){$\begin{array}{ccc} \lambda^0 & 0 &
0
\\  & \lambda^0 & 0  \\  &  & \lambda^0  \end{array}$ }}
\put(20,36){\dashbox{1}(13,18){$\begin{array}{cc} 1 & 0  \\ 0 & 1
\\ 0 & 0   \end{array}$ }}
\put(20,24){\dashbox{1}(13,12){$\begin{array}{cc} \lambda^0 & 0
\\  & \lambda^0     \end{array}$ }}
\put(33,24){\dashbox{1}(13,12){$\begin{array}{cc} 1 & 0 \\ 0 & 1
\end{array}$ }}
\put(33,12){\dashbox{1}(13,12){$\begin{array}{cc} \lambda^0 & 0
\\  & \lambda^0     \end{array}$ }}
\put(46,12){\dashbox{1}(6,12){$\begin{array}{c} 1  \\ 0
\end{array}$ }}
\put(46,6){\dashbox{1}(6,6){$\lambda^0$}}
\put(52,6){\dashbox{1}(6,6){1}}
\put(52,0){\dashbox{1}(6,6){$\lambda^0$}}
\end{picture}
\end{array}\right)
$$
and
$$
S_{\P}=\left(\begin{array}{c} \unitlength 1mm
\begin{picture}(74,55)

\put(0,36){\dashbox{1}(24,18){$\begin{array}{ccc} s_1 & s_{13}^1 &
s_{15}^1 \\  & s_3 & s_{35}^1  \\  &  & s_5  \end{array}$ }}
\put(24,36){\dashbox{1}(16,18){$\begin{array}{cc} s_{11}^2 & s_{13}^2  \\
0 & s_{33}^2\\ 0 & 0   \end{array}$ }}
\put(40,36){\dashbox{1}(16,18){$\begin{array}{cc} s_{11}^3 & s_{13}^3  \\
s_{31}^3 & s_{33}^3\\ 0 & s_{53}^3   \end{array}$ }}
\put(56,36){\dashbox{1}(8,18){$\begin{array}{c} s_{11}^4   \\
 s_{31}^4\\  0   \end{array}$ }}
\put(64,36){\dashbox{1}(8,18){$\begin{array}{c}  s_{11}^5   \\
 s_{31}^5 \\  s_{51}^5 \end{array}$ }}

\put(24,24){\dashbox{1}(16,12){$\begin{array}{cc} s_1 & s_{13}^1
\\  & s_3     \end{array}$ }}
\put(40,24){\dashbox{1}(16,12){$\begin{array}{cc} s_{11}^2 & s_{13}^2 \\
0 & s_{33}^2 \end{array}$ }}
\put(56,24){\dashbox{1}(8,12){$\begin{array}{c} s_{11}^3   \\
 s_{31}^3   \end{array}$ }}
\put(64,24){\dashbox{1}(8,12){$\begin{array}{c}  s_{11}^4   \\
 s_{31}^4  \end{array}$ }}

\put(40,12){\dashbox{1}(16,12){$\begin{array}{cc} s_1 & s_{13}^1
\\  & s_3\end{array}$ }}
\put(56,12){\dashbox{1}(8,12){$\begin{array}{c} s_{11}^2\\ 0
\end{array}$ }}
\put(64,12){\dashbox{1}(8,12){$\begin{array}{c} s_{11}^3\\
s_{31}^3\end{array}$ }}

\put(56,6){\dashbox{1}(8,6){$s_1$}}
\put(64,6){\dashbox{1}(8,6){$s_{11}^2$}}
\put(64,0){\dashbox{1}(8,6){$s_1$}}
\end{picture}
\end{array}\right)
$$
 is given by
$M(\lambda)S_{\P}-S_{\P}M(\lambda)=0$. On the other hand
$$
M(\mu)=\left(\begin{array}{c} \unitlength 1mm
\begin{picture}(72,66)
\put(0,48){\dashbox{1}(20,18){$\begin{array}{ccc} \mu^0 & 0 & 0
\\  & \mu^0 & 0  \\  &  & \mu^0  \end{array}$ }}
\put(20,48){\dashbox{1}(20,18){$\begin{array}{ccc} 1 & 0 & 0 \\
0 & 1 & 0\\ 0 & 0 & 1   \end{array}$ }}
\put(20,30){\dashbox{1}(20,18){$\begin{array}{ccc} \mu^0 & 0 & 0
\\  & \mu^0 & 0\\  &  & \mu^0   \end{array}$ }}
\put(40,30){\dashbox{1}(13,18){$\begin{array}{cc} 1 & 0 \\ 0 & 1
\\0 & 0
\end{array}$ }}
\put(40,18){\dashbox{1}(13,12){$\begin{array}{cc} \mu^0 & 0
\\  & \mu^0     \end{array}$ }}
\put(53,18){\dashbox{1}(13,12){$\begin{array}{cc} 1 & 0  \\ 0 & 1
\end{array}$ }}
\put(53,6){\dashbox{1}(13,12){$\begin{array}{cc}\mu^0 & 0 \\ &
\mu^0 \end{array}$}}
\put(66,6){\dashbox{1}(6,12){$\begin{array}{c} 1 \\ 0
\end{array}$}} \put(66,0){\dashbox{1}(6,6){$\mu^0$}}
\end{picture}
\end{array}\right),
$$
and
$$
S_{\Q}=\left(\begin{array}{c} \unitlength=0.9mm
\begin{picture}(84,68)

\put(0,50){\dashbox{1}(25,18){$\begin{array}{ccc} \ t_1 & t_{12}^1
& t_{14}^1 \\  & t_2 & t_{24}^1  \\  &  & t_4
\end{array}$ }}
\put(25,50){\dashbox{1}(28,18){$\begin{array}{ccc} \ \ t_{11}^2 & t_{12}^2 & t_{14}^2  \\
\ \ t_{21}^2 & t_{22}^2 & t_{24}^2\\ \ \ 0 & 0 & t_{44}^2
\end{array}$}}

\put(25,32){\dashbox{1}(28,18){$\begin{array}{ccc} \quad\ t_1 &
t_{12}^1 &t_{14}^1
\\&t_2 &t_{24}^1 \\&&t_4     \end{array}$ }}
\put(53,32){\dashbox{1}(32,36){$\begin{array}{cccc} &&& \\
&&&\\&&*&\\&&&\\&&&\\&&& \end{array}$ }}
\put(53,0){\dashbox{1}(32,32){$\begin{array}{cccc} &&&  \\
 &&&\\&&*&\\ &&&   \end{array}$ }}
\end{picture}
\end{array}\right).
$$
is given by $M(\mu)S_{\Q}-S_{\Q}M(\mu)=0$.

The $9 \times 11$ matrix $S(\delta(a_1))$ equals
$$
\left(\begin{array}{cccccccccc} *&*&*&&&&&&\\
*&*&*&&&&&\\ 0&0&0&&&&*&\\ *&*&*&&&&&\\
{\alpha v_{71}}&{\alpha v_{72}}&*&&&&&\\ \s{\alpha v_{81}+\atop
\gamma_1 v_{91}} & \s{\alpha
v_{82}+\atop \gamma_1 v_{92}}&*&&&&&\\
0&0&0&{\beta v_{71}}&{\beta v_{72}}&*&&\\
{\alpha v_{91}}&{\alpha v_{92}}&{\alpha v_{93}}& \s{\alpha
v_{94}+\beta
v_{81}\atop +\gamma_2 v_{91}}& \s{\alpha v_{95}+\beta v_{82}\atop +\gamma_2 v_{92}}&*&&\\
 0&0&0&{\beta v_{91}}&{\beta v_{92}}&{\beta v_{93}}& \s{\beta
v_{94}+\atop \gamma_3 v_{91}}& \s{\beta v_{95}+\atop \gamma_3
v_{92}} &*& \cdots\end{array} \right)
$$
Let
$$
M(a_1)=\left(\begin{array}{ccccccccccc} \emptyset&\emptyset&\emptyset&&&&&&&&\\
\emptyset&\emptyset&\emptyset&&&&&&&&\\
1&\xi&\emptyset&&&&&&&&\\
\emptyset&\emptyset&\emptyset&&&&&&&&\\
\eta&1&\emptyset&&&&\hbox{\huge$\emptyset$}&&&&\\
1&\emptyset&\emptyset&&&&&&&&\\
0&0&1&&&&&&&&\\
0&0&1&&&&&&&&\\
0&0&0&&&&&&&&
\end{array} \right)_{9\times 11}
$$
where $\eta$ and $\xi$ are two algebraically independent parameters with a domain $k\times k$. Thus $\forall\,
S\in End_{\mathfrak{A}}(M)$, we have $S_{\P}=sI_9+L_1$
  and $S_{\Q}=sI_{11}+L_2$
for some upper triangular nilpotent matrices $L_1,L_2$, and $ s\in
\rK $.

If $\alpha\beta=0$, we can take some lower sizes to construct a
canonical matrix with two independent parameters $\eta,\xi$,
\cite{ZX}.

If $\P=\Q$, and there exists a pair $(\lambda^0, \mu^0)$ with
$\lambda^0=\mu^0$, $f(\lambda^0, \mu^0)=0$, (i.e. $\lambda^0$ is a
root of $f(\lambda, \lambda)$), $g_{\P}(\lambda^0)\neq 0$, we may
set $m_{\P}=4$, $M(\lambda)\simeq J_3(\lambda^0)\oplus
J_1(\lambda^0)$. Thus
$$
M(\lambda)=\left(\begin{array}{c} \unitlength 1mm
\begin{picture}(25,25)
\put(0,12){\dashbox{1}(12,12){$\begin{array}{cc} \lambda^0 & 0
\\  & \lambda^0     \end{array}$ }}
\put(12,12){\dashbox{1}(6,12){$\begin{array}{c} 1  \\ 0
\end{array}$ }}
\put(12,6){\dashbox{1}(6,6){$\lambda^0$}}
\put(18,6){\dashbox{1}(6,6){1}}
\put(18,0){\dashbox{1}(6,6){$\lambda^0$}}
\end{picture}
\end{array}\right),\qquad \mbox{and }\quad
S_{\P}=\left(\begin{array}{c} \unitlength 1mm
\begin{picture}(33,25)

\put(0,12){\dashbox{1}(16,12){$\begin{array}{cc} s_1 & s_{13}^1
\\  & s_3\end{array}$ }}
\put(16,12){\dashbox{1}(8,12){$\begin{array}{c} s_{11}^2\\ 0
\end{array}$ }}
\put(24,12){\dashbox{1}(8,12){$\begin{array}{c} s_{11}^3\\
s_{31}^3\end{array}$ }}

\put(16,6){\dashbox{1}(8,6){$s_1$}}
\put(24,6){\dashbox{1}(8,6){$s_{11}^2$}}
\put(24,0){\dashbox{1}(8,6){$s_1$}}
\end{picture}
\end{array}\right)
$$
is given by $M(\lambda)S_{\P}=S_{\P}M(\lambda)$. When
$\alpha\beta\neq 0$,
$$S(\delta(a_1))=\left(\begin{array}{cccc} *&*&*&*\\ 0&*&*&*\\
\alpha v_{41}&*&*&*\\ 0&0&\beta v_{41}&* \end{array} \right),
\quad M(a_1)=\left(\begin{array}{cccc} \emptyset&\emptyset&\emptyset&\emptyset\\
\xi&\emptyset&\emptyset&\emptyset\\
\eta&\emptyset&\emptyset&\emptyset\\ 0&1&\emptyset&\emptyset
\end{array} \right),$$
where $\eta$ and $\xi$ are independent parameters with a domain $k\times k$. Thus given any $S\in
End_{\mathfrak{A}}(M)$, $S_{\P}=sI_4+L$ for some upper triangular nilpotent matrix $L$ and $ s\in \rK$. When
$\alpha\beta= 0$ or $\alpha\neq 0, \beta=0$, the method is similar \cite{ZX}.

If $\P=\Q$, and all the pairs $(\lambda^0, \mu^0)$ with
$f(\lambda^0,\mu^0)=0$, $g_{\P}(\lambda^0)g_{\Q}(\mu^0)\neq 0$
have the property that $\lambda^0\neq \mu^0$, we may use a rolled
up version $\overline{M}$ of the construction given in  the case
$\P\neq \Q$, i.e.  let $m_{\P}=20$,
$\overline{M}(\lambda)=M(\lambda^0)\oplus M(\mu^0)$ and
$\overline{M}(a_1)=\left(\begin{array}{cc} 0&M(a_1)\\
0&0 \end {array}\right)$.

 Case 2. \unitlength 1mm
\begin{picture}(32, 8)
\put(10.00,-1.00){\circle*{1.00}}
\put(25.00,-1.00){\circle*{1.00}} \put(5.00,0.00){\oval(5,5)[t]}
\put(5.00,0.00){\oval(5,5)[bl]} \put(4.00,-3.00){\vector(3,1){3}}
\put(10.50,-1.00){\line(1,0){13}}
\put(23.50,-1.00){\vector(1,0){1}}
\put(10.00,1.00){\makebox{$\P$}} \put(27.00,0.00){\makebox{$\Q$}}
\put(0.00,0.0){\makebox{$\lambda$}}
\put(16.00,1.00){\makebox{$a_1$}}
\end{picture}, (or we have a dual diagram). Let $m_{\P}=5$, $ m_{\Q}=2$,
$M(\lambda)=J_5(\lambda^0)$, where  $g_{\P}(\lambda^0)\neq 0$,
then \vspace{2mm}
$$
S_{\P}=\left(\begin{array}{ccccc} s_1&s_2&s_3&s_4&s_5\\
& s_1&s_2&s_3&s_4\\ &&s_1&s_2&s_3\\&&&s_1&s_2\\ &&&&s_1
\end{array}
\right),  \quad S_{\Q}=\left(\begin{array}{cc} t_{11}& t_{12}\\
t_{21} & t_{22} \end{array}\right); \quad
M(a_1)=\left(\begin{array}{cc} \xi&\emptyset\\ \eta&\emptyset\\
\emptyset&\emptyset\\ 1&\emptyset  \\
0&1\end{array}\right),
$$
 where $\eta$ and $\xi$ are independent parameters with a domain $k\times k$. Thus $\forall S\in End_{\mathfrak{A}}(M)$,
  $S_{\P}=sI_5,
\quad S_{\Q}=sI_2$, for any $ s\in \rK$. The proof for the dual
case is similar.

It is clear from the construction of $M$, that we are able to
define a local layered bocs $\mathfrak{B}$ induced from
$\mathfrak{A}$ just before the appearing of $\xi$ in all the cases
(see also \cite{ZX}), i.e. $\mathfrak{B}$ has a layer
$L_{\mathfrak{B}} =(\Gamma'_{\mathfrak{B}}; \omega_{\mathfrak{B}};
b_0, b_1, \cdots, b_n; u_1, \cdots, u_m )$  with a unique
non-trivial vertex $\I$, and
$\Gamma'_{\mathfrak{B}}(\I,\I)=\rK[\eta]$. Moreover,
$\delta(b_0)=0$ still by the construction, and we are able to set
$b_0 =\xi$. We denote the reduction functor from $R(\mathfrak{B})$
to $R(\mathfrak{A})$ by $G$.
 And define a functor
$$F: \rK\langle \eta, \xi\rangle\mbox{-mod} \rightarrow R(\mathfrak{B})$$ by
sending any  $\rK\langle \eta, \xi\rangle$-module \\[10pt]
\begin{center}\unitlength=0.5pt
\begin{picture}(20,10)
 \put(0,60){\oval(40,40)[t]}
\put(0,60){\oval(40,40)[bl]} \put(0,40){\vector(3,1){5}}
\put(50,60){\oval(40,40)[t]} \put(50,60){\oval(40,40)[br]}
\put(50,40){\vector(-3,1){5}}
\put(21,39){$\bullet$} \put(19,84){\makebox{$\rK^m$}} \put(-55,80){\makebox{$L(\eta)$}} \put(-75,50){$L=$}
\put(76,80){\makebox{$L(\xi)$}}
\end{picture}
\end{center}
for any dimension $m$ and any linear maps $L(\eta), L(\xi)$ to $\overline{L}=F(L)$, such that
$$\overline{L}_{\P}=\rK^m,\ \overline{L}(\eta)=L(\eta),\ \overline{L}(b_0)=L(\xi),\ \overline{L}(b_i)=0\ {\rm
for}\  i=1,2,\cdots,n.$$ If $L'\in \rK\langle\eta,\xi\rangle\mbox{-mod}$  and $\overline{L}'=F(L')$, then a
morphism $\varphi: L\rightarrow L'$ is sent to a morphism $\overline{\varphi}: \overline{L}\rightarrow
\overline{L}'$, such that $\overline{\varphi}_{\P}=\varphi$, and $\overline{\varphi}(u_j)=0$ for
$j=1,2,\cdots,m$. Thus, if $L$ is indecomposable, so is $\overline{L}$; and if $\overline{L}\simeq
\overline{L}'$, then $L\simeq L'$ both by Corollary 2.2.2. Finally, the functor $GF: \rK\langle \eta,
\xi\rangle\mbox{-mod} \rightarrow R(\mathfrak{A})$
 preserves indecomposability and isomorphism classes, so that
 $\mathfrak{A}$ is of wild type. \hfill$\square$

\subsection{The tame theorem}

\kg In this subsection we will give a slightly different statement
of the well-known tame theorem \cite{D,CB1}  by constructing a
finite set of minimal local  bimodule problems whose valuated
representation categories cover all the canonical forms of the
indecomposable matrices of dimension at most $n$ for each
non-negative integer $n$.

{\bf Lemma 5.2.1} Let $(\K,\M,H)$ be a non-wild bimodule problem having a reduction sequence $(*)$. Suppose that
$\overline{N}^{r_0}_{p_{r_0},q_{r_0}}=W\oplus (\lambda)$ with $\lambda$ attached to $\I\in
T^{r_0+1}/\sim^{r_0+1}$. Given any $r> r_0$, if the first arrow $a_1$ in $(\K^r,\M^r,H^r)$ starts or ends at
$\I$,  then the reduction for $a_1$ must be a regularization. Consequently, the equivalent class $\I$  remains
until the end of the sequence$(*)$, i.e. $\I\in T^{\infty}/\sim^{\infty}$. In particular, if $R^{\infty}$ is
local, then $\overline{N}_{p_{r_0}q_{r_0}}^{r_0}=(\lambda)$.

{\bf Proof.} Item (1) of A1 and A2 in  Proposition 4.3.1 can not
occur because of non-wildness,  and item (2) of A1,A2 always
yields a regularization. \hfill$\square$

{\bf Lemma 5.2.2} Let $(\K,\M,H)$ be a non-wild bimodule problem. Then any  reduction sequence $(*)$ of freely
parameterized triples, which preserves all the free parameters, can reach a minimal triple.

{\bf Proof.} By the same reason as in the proof of Lemma 5.2.1.
\hfill $\square$

{\bf Lemma 5.2.3} Let $(\K,\M,H=0)$ be a non-wild bimodule problem having a reduction sequence $(*)$ and
$N^s=(0)$. Then $H^{\infty}(\J)$ for any trivial $\J\in T^{\infty}/\sim^{\infty}$ defined in 3.5, and
$\wt{H}^{\infty}(J_d(\lambda_{\I}^0))$ for any non-trivial $\I\in T^{\infty}/\sim^{\infty}$ defined in 4.6  are
both canonical forms under the fixed orders of matrix and the base field.

{\bf Proof.} Let $\hat{\n}^{\infty}=\underline{m(\J)}$ when $\J$
is trivial. We obtain a reduction sequence $(\hat{*})$ in
Corollary 3.3.1 under the
 orders mentioned in the lemma. Thus $H^{\infty}(\J)$ is a canonical
 form.

Let $\hat{\n}^{\infty}=\underline{m(\I,1)}$ when $\I$ is
non-trivial, which has dimension $1$ at $\I$ and $0$ at any
$\J\ne\I$. We obtain a reduction sequence $(\hat{*})$ in Corollary
4.3.1, which contains a unique parameter $\lambda_{\I}$ and
$\hat{R}^{\infty}$ is local. Thus
$\overline{\hat{N}}^{r}_{p_{r}q_{r}}=(\lambda_{\I})$ for some
$1\leq r \le s$ by Lemma 5.2.1. It is obvious that the order of
the base field is not destroyed in
$\wt{H}^{r+1}(J_d(\lambda_{\I}^0))$. On the other hand, all the
reductions after  the $(r+1)$-th step are regularization in
$(\wt{*})$ given by $J_d(\lambda_{\I}^0)$  by Lemma 5.2.1 and
Theorem 4.5.1, which are made according to the partitions of
$(\widehat{*})$. But they can also be made according to the order
of Definition 2.3.1, whenever we write the matrix equations (8) of
2.6 element-wise.  Therefore
$\wt{H}^{\infty}(J_d(\lambda_{\I}^0))$ is a canonical form under
the fixed orders of both matrix and base field. \hfill$\square$

{\bf Theorem 5.2.1} Let $(\K,\M,H=0)$ be a non-wild bimodule problem. Suppose we fix the  orders both on matrix
and the base field. Then for any fixed $n\in\mathbb{N}$, there exists a finite set of sequences of freely
parameterized triples $\{*_1,*_2,\cdots,*_h\}$ such that

(1) the end terms of corresponding $(**)$'s:
$$
(\K^{\infty}_1, \M^{\infty}_1, H^{\infty}_1), (\K^{\infty}_2,
\M^{\infty}_2, H^{\infty}_2), \cdots, (\K^{\infty}_h,
\M^{\infty}_h, H^{\infty}_h)
$$
are minimal and local;

 (2) for any indecomposable $M\in
Mat(\K,\M)$ of size at most $n$, there exists a unique sequence $(*_l)$ in the set and a matrix
$$M_l^{\infty}\in {Mat}(\K_l^{\infty}, \M_l^{\infty})\ \ {\rm with}\ \  M\simeq
{\vartheta}^l_{0,\infty}(M_l^{\infty}).$$

(3) ${\vartheta}^l_{0,\infty}(M^{\infty}_l)$ is the canonical form
of $M$.

{\bf Proof.} We use the induction firstly on  size $n$. $n=0$ is
trivial.  Suppose that the statement is true for a non-negative
integer $(n-1)$. We will start from a triple $(N,R,\n)$ of size
$n$. We use the induction for the second time  on the reduction
steps $r$.

 When $r=0$, there are only finitely many choices of size vector $\n$, since the equation
$$n=n_1+n_2+\cdots+n_t$$ has only finitely many solutions of non-negative integers. Suppose that for each fixed
triple $(N,R,\n)$, there are only finitely many possibilities to construct the reduction sequences up to  the
$r$-th step. Consider one of such sequences. We will show that the possibilities to construct
$\overline{N}^r_{p_rq_r}$ are also finite. Thus our conclusion holds, since a reduction sequence towards a
minimal bimodule problem can be taken at most $n^2$ steps.

 We treat the $(r+1)$-th reduction still according to Proposition 4.3.1. Item (1) of A1 and A2 can not occur
because of non-wildness.  Item (2) of A1,A2, A3 gives
$\overline{N}^r_{p_rq_r}=\emptyset$, so we have only one choice.

 In item (1) of A3 and $\P\neq\Q$,
$$\overline{N}^r_{p_rq_r}= \left(\begin{array}{cc} 0& I_d \\ 0& 0\end{array}\right)_{n^r_{p_r}\times
n^r_{q_r}},$$ then we have finitely many choices for $d=0,1,\cdots, \min\{n^r_{p_r}, n^r_{q_r}\}$.

 In item (1) of A3 and $\P=\Q$ or item (3) of A3. If $R^r$ is not local, then
either $T^r/\sim^r\setminus\{\P_r\}\ne \emptyset$, or
$T^r/\sim^r=\{\P_r\}$, and $n_{\P_r}>1$. Let $(\hat{*})$ be a
reduction sequence of Corollary 4.3.1 given by $\hat{\n}^r$, such
that $\hat{n}^r_{\P_r}=1$, $\hat{n}^r_{\J}=0$, $\forall \J\in
T^r/\sim^r$, $\J\ne \P_r$, then $\hat{n}^r<n^r$. The induction
hypothesis and  Lemma 5.2.1 tell us that the sequence $(\hat{*})$
with $\overline{\hat{N}}^r_{p_rq_r}=(\lambda)$ already exists in
our set and reaches a minimal triple with  domain
$\rK\setminus\{\mbox{the roots of } \hat{g}(\lambda)\}$. Thus
Lemma 5.3.3 ensures that $\wt{H}^{\infty}(J_d(\lambda^0))$, for
some $d\in \mathbb{N}$ and $\hat{g}(\lambda^0)\ne 0$, have covered
the canonical forms of indecomposables of size at most $n$ such
that $\overline{\hat{N}}^r_{p_rq_r}$ has eigenvalues $\lambda^0\in
(\rK \setminus \{ \text{the roots of}\ \hat{g}(\lambda)\}  )$.
Therefore
$$\overline{N}^r_{p_rq_r}\simeq \oplus_{l,d}
J_d(\lambda^l)^{e_{l,d}}$$
 with $\hat g(\lambda^l)=0$ such that
$$n_{p_r}=\sum\limits_{l,d} de_{l,d}.$$ But such a equation in
variables $e_{l,d}$ has only finitely many solutions of
non-negative integers, i.e. $\overline{N}_{p_rq_r}^r$ has only
finitely many choices.

If $R^r$ is local, then $T^r/\sim^r=\{\P_r\}$ and  $n_{p_r}=1$.
Item (3) of A3 gives $\overline{N}^r_{p_rq_r}=(\lambda)$ with a
domain $\rK\setminus \{\mbox{the roots of }g(\lambda)\}$, or
$\overline{N}^r_{p_rq_r}=(\lambda^l)$ with $g(\lambda^l)= 0$, thus
$\overline{N}^r_{p_rq_r}$ has only finitely many choices.

 Therefore we obtain finitely many reduction
sequences $(*)$  of size $n$. Finally all the sequences having
local end terms are retained, and others are excluded, since their
local direct summands have already been included in the case of
$(n-1)$. The proposition follows by double induction  on $n$ and
$r$. \hfill$\square$

\subsection{Local wild bimodule problems}

\kg This subsection is devoted to constructing a reduction sequence consisting of local parameterized  bimodule
problems (not necessarily freely!).  In the subsection we will not distinguish the multiplication of the
parameters whether they act from left or right. And we denote such kind of differentials by $\delta^c$, where
$c$ in $\delta^c$ means the commutativity of the multiplication from left and right.

{\bf Lemma 5.3.1.} Let $(\K,\M,H)$ be a bimodule problem with $T/\!\!\sim=\{\P\}$, and $H$ have a free parameter
$\lambda_0$ with the domain $\rK\setminus \{{\rm the\ roots\ of\ }f(\lambda_0)\}$. Suppose that $\n$ is a size
vector with $n_{\P}=1$ and  $(N,R,\n)$ is given by Formula (6) of 2.6. Then we are able to construct a reduction
sequence $(*)$ of 3.3 satisfying the following properties:

(1) $\n=\n^1=\cdots=\n^r=\cdots=\n^s$, $N^s=(0)$. $H^s$ contains
parameters $\lambda_0,\lambda_1,\cdots,\lambda_{\gamma}$ with
$\gamma\geq 0$, which are algebraically independent.

(2) When  $\delta^c(\rho^{r_l-1})^*=0$ in $(\K^{r_l-1},
\M^{r_l-1})$, we set
 $(H^{r_l})_{p_{r_l-1},q_{r_l-1}}=(\lambda_{l})$, where
$l=1,\cdots,\gamma$.

(3) When  $\delta^c(\rho^{r-1})^*=w^{r-1}\ne 0$ in $(\K^{r-1},\M^{r-1})$, we set
$(H^{r})_{p_{r-1},q_{r-1}}=\emptyset$ and $w^{r-1}=0$ after a localization given by $c^{r}(\lambda_0,
\lambda_1,\cdots,\lambda_l )$, where $ r_l <r<r_{l+1}$, $l=0,1,\cdots,\gamma$, $r_0=0$, $r_{\gamma+1}=s+1$.

(4) The domain of the parameters in $H^{r}$, for $ r_l < r< r_{l+1}$, is an open subset of $\rK$, determined by
$f^{r}(\lambda_0, \lambda_1,\cdots,\lambda_l)\neq 0$, where
\begin{equation}
f^{r}(\lambda_0, \lambda_1,\cdots,\lambda_l)=f(\lambda_0)\left(
\prod_{l'=0}^{l-1}\prod_{r_{l'} < r'<r_{l'+1}}c^{r'}(\lambda_0,
\lambda_1,\cdots,\lambda_{l'})\right) \prod_{r_l< r'\leq r}
c^{r'}(\lambda_0, \lambda_1,\cdots,\lambda_l),
\end{equation}
and $c^{r'}(\lambda_0, \lambda_1, \cdots, \lambda_{l'})$ are
polynomials in $\lambda, \lambda_1, \cdots, \lambda_{l'}$
appearing at $r'$-th reductions.

(5) $f^{r_{l}}(\lambda_0,\lambda_1, \cdots, \lambda_{l-1})=f^{r_{l}-1}(\lambda_0,\lambda_1, \cdots,
\lambda_{l-1})$, for $l=1,\cdots, \gamma$.

{\bf Proof.} (1) is obvious. We will prove (2), (3), (4) and (5)
by induction on $r$. $r=0$ is clear. Suppose that we have already
had a sequence of parameterized bimodule problems

$$(\K,\M,H), (\K^1,\M^1,H^1), \cdots,
(\K^{r}, \M^{r}, H^{r})$$
 and parameters $\lambda_0,
\lambda_1, \cdots, \lambda_{l}$ in $H^{r}$ with a polynomial
$$f^{r}(\lambda_0, \lambda_1,\cdots,\lambda_{l})\ne 0$$
 for some
$r\geq r_{l}$. Let $(\rad\K^{r})^*$ have a $\rK$-basis $\{v_1,
v_2,\cdots,v_m \}$.  If $a_1=(\rho^{r})^*$ is the first arrow of
$(\K^{r}, \M^{r}, H^{r})$, then
\begin{equation}
\delta^c(a_1)=\sum_{j=1}^m f_j(\lambda_0, \lambda_1,
\cdots,\lambda_l)v_j
\end{equation}

 If $\delta^c(a_1)\ne 0$, we make a basis change in the field of rational
functions $\rK(\lambda_0, \lambda_1,\cdots,\lambda_{l})$, say
$(w_1,\cdots,w_m)^T=Q(v_1,\cdots,v_m)^T$ with the first row of $Q$
being $(f_1, f_2, \cdots, f_m)$. Then $\delta(a_1)=w^{r+1}$ and
thus $(H^{r+1})_{p_rq_r}=\emptyset$. Denote by $c^{r+1}(\lambda_0,
\lambda_1,\cdots,\lambda_{l})$ the product of the denominators of
the entries of  $Q$ and $Q^{-1}$, let
$$f^{r+1}(\lambda_0,
\lambda_1, \cdots,\\ \lambda_{l})=f^{r}(\lambda_0, \lambda_1,
\cdots, \lambda_{l})c^{r+1}(\lambda_0, \lambda_1, \cdots,
\lambda_{l}).$$ In case of $\delta^c(a_1)=0$, we set
$r_{l+1}=r+1$, $(H^{r+1})_{p_rq_r}=(\lambda_{l+1})$, and
$$f^{r+1}(\lambda_0, \lambda_1, \cdots,
\lambda_{l})=f^{r}(\lambda_0, \lambda_1, \cdots, \lambda_{l}).$$
Thus we complete the proof from $(\K^{r}, \M^{r}, H^{r})$ to
$(\K^{r+1}, \M^{r+1}, H^{r+1})$. Our assertion follows by
induction. \hfill $\square$

\subsection{Triangular formulae}

\kg Based on  Lemma 5.3.1, we will  in this subsection construct
two kinds of freely parameterized  bimodule problems which are
established by Bautista. First we fix some notations.

Let $\mathfrak{A}=(\Gamma, \Omega)$ be a local bocs, with a layer $L=(\Gamma'; \omega; a_1, a_2, \cdots, a_n;
v_1, v_2, \cdots, v_m)$. Suppose that  the indecomposable object of $\Gamma'$ is $\P$, and
$\Gamma'(\P,\P)=\rK[\lambda, f(\lambda)^{-1}]$. Denote by $\Delta$ the $\Gamma'$-$\Gamma'$-bimodule freely
generated by $a_1, a_2, \cdots, a_n$, $$\Delta^{\otimes
p}=\underbrace{\Delta\otimes_{\Gamma'}\Delta\otimes\cdots\otimes_{\Gamma'}\Delta}_p$$ and $\Delta^{\otimes
0}=\Gamma'$ . Then
\begin{equation}
\left\{\begin{array}{l} \Gamma=\bigoplus\limits_{p\ge 0}\Delta^{\otimes p},\\
\overline{\Omega}=\bigoplus\limits_{p,q\ge
0}\bigoplus\limits_{j=1}^m \Delta^{\otimes p}v_j \Delta^{\otimes
q}.\end{array}\right.
\end{equation}
For any $v\in\overline{\Omega}$, denote by $v^0$ the projection of
$v$ into $\oplus_{j=1}^m \Gamma' v_j \Gamma'$.

{\bf Definition 5.4.1} Given any constants $\lambda_0^0, \lambda_1^0, \cdots, \lambda_{\gamma-1}^0\in\rK$ with
$$f^{r_{\gamma}}(\lambda_0^0, \cdots, \lambda_{\gamma-1}^0)\ne 0$$ given in Lemma 5.3.1, and denote
$\lambda_{\gamma}$ by $\nu$, we construct  a local bocs $\mathfrak{A}_{(\lambda_0^0, \lambda_1^0, \cdots,
\lambda_{\gamma-1}^0)}$ with a layer $$L=(\Gamma'; \omega; b_1,\cdots,b_j; u_1, \cdots, u_m),$$ such that if
$\gamma=0$, $$\mathfrak{A}_{(\lambda_0^0, \lambda_1^0, \cdots, \lambda^0_{\gamma-1})}=\mathfrak{A},$$ thus
$\Gamma'(\P,\P)=\rK[\nu,f(\nu)^{-1}]$;  If $\gamma>0$, $$\mathfrak{A}_{(\lambda_0^0, \lambda_1^0, \cdots,
\lambda_{\gamma-1}^0)}=\mathfrak{A}^{r_{\gamma}},$$ thus $\Gamma'(\P,\P)=k[\nu]$, and $\nu$ appears for the
first time in $\mathfrak{A}^{r_\gamma}$ .

{\bf Lemma 5.4.1} There exists a sequence of  localizations given by $$\sigma_1(\nu), \sigma_2(\nu), \cdots,
\sigma_j(\nu),$$ and a triangular formula:
\begin{equation}
 \left\{
\begin{array}{ccl}
\dz(b_1)^0 &=& g_{11}(\nu,\kappa)u_1 \\
\dz(b_2)^0 &=& g_{21}(\nu,\kappa)u_1 +g_{22}(\nu,\kappa)u_2 \\
\cdots   &  & \qquad\cdots \\
\dz(b_j)^0 &=& g_{j1}(\nu,\kappa)u_1
+g_{j2}(\nu,\kappa)u_2+\cdots+g_{j,j}(\nu,\kappa)u_j
\end{array} \right.
\end{equation}
where $\nu$ stands for the left multiplication by $\nu$, and $\kappa$ for the right multiplication by $\nu$;
$$g_{ll'}(\nu, \kappa)\in\rK[\nu,\kappa,\sigma_l(\nu)^{-1}\sigma_l(\kappa)^{-1}],\ g_{ll}(\nu, \nu)\ne 0,\
l'\leq l,\ l=1,2 \cdots, j.$$

{\bf Proof.} $j=0$ is trivial. If $j>0$, then $\delta^c(b_1)\ne 0$, since $\nu=\lambda_{\gamma}$ is the last
parameter. Therefore $\delta(b)\neq 0$. Without loss of generality, we may divide the dotted arrows by some
power of $f(\nu)f(\kappa)$ when $\gamma=0$, such that the coefficients in $\delta(b_1)$ are all polynomials.
After a localization given by $c_1(\nu)$ in Formula (30) of 4.3, $\delta(b_1)=g_{11}(\nu,\kappa)u_1$ with
$g_{11}(\nu,\nu)\ne 0$. Let $\sigma_1(\nu)=f(\nu)c_1(\nu)$ if $\gamma=0$, and $\sigma_1(\nu)=c_1(\nu)$ if
$\gamma>0$. In case of $j=1$, the proof is completed. Otherwise,
$$\delta(b_2)^0=g_{21}(\nu,\kappa)u_1+ \sum_{l=2}^m
g'_{2l}(\nu,\kappa)u_{2l},$$
 and $g'_{2l}(\nu,\kappa)\in
 k[\nu,\kappa,\sigma_1(\nu)^{-1}\sigma_1(\kappa)^{-1}]$.
Without loss of generality we may divide $u_{2l}$ by some power of
 $\sigma_1(\nu)\sigma_1(\kappa)$, such that
$g'_{2l}(\nu,\kappa)\in \rK[\nu,\kappa]$. Suppose the highest common factor of $g'_{2l}(\nu,\kappa)$,
$l=2,\cdots,m$, is $g_{22}(\nu,\kappa)$, and $c_2(\nu)$ is obtained by (30) to make a basis change from $u_{22},
\cdots, u_{2m}$ to $u_2, u_{33}, \cdots, u_{3m}$ in a Hermit ring $\rK[\nu,\kappa,
\sigma_2(\nu)^{-1}\sigma_2(\nu)^{-1}]$ with $\sigma_2 (\mu)=\sigma_1(\mu)c_2(\mu)$. Then $g_{22}(\nu,\nu)\neq
0$. By induction we finally reach to the integer $j$ and Formula (40). Which depends on a sequence of
localizations given by
\begin{equation}
\sigma_l(\nu)=\left\{\begin{array}{ll} f(\nu)\prod_{p=1}^l
c_p(\nu),
&\mbox{ when }\gamma=0;\\
\prod_{p=1}^l c_p(\nu), &\mbox{ when }\gamma>0,
\end{array}\right. \quad\mbox{ write }
\sigma(\nu)=\sigma_j(\nu).
\end{equation}

The proof is completed.\hfill$\square$

\medskip
{\bf Definition 5.4.2} If for any $\lambda_0^0, \lambda_1^0,
\cdots, \lambda_{\gamma-1}^0$, $g_{ll}(\nu,\kappa)\in
\rK[\nu,\kappa, \sigma_l(\nu)^{-1} \sigma_l(\kappa)^{-1}]$  are
invertible, we may fix some constants $\lambda_0^0, \lambda_1^0,
\cdots, \lambda_{\gamma-2}^0$ with
 $$f^{r_{\gamma-1}}(\lambda_0^0,
\lambda_1^0, \cdots, \lambda_{\gamma-2}^0)\ne 0.$$
  Denote
$\lambda_{\gamma-1}$ by $\lambda$, we then construct a local bocs
$\mathfrak{A}_{(\lambda_0^0, \lambda_1^0, \cdots,
\lambda_{\gamma-2}^0)}=(\Gamma,\Omega)$ with a layer
 $$L=(\Gamma'; \omega'; a_1,\cdots,
a_i, b_1, \cdots, b_j; v_1,\cdots,v_m),$$ such that if $\gamma=1$,
$$\mathfrak{A}_{(\lambda_0^0, \lambda_1^0, \cdots,
\lambda_{\gamma-2}^0)}=\mathfrak{A},$$  thus
$\Gamma'(\P,\P)=\rK[\lambda,f(\lambda)^{-1}]$; if $\gamma>1$,
$$\mathfrak{A}_{(\lambda_0^0, \lambda_1^0, \cdots,
\lambda_{\gamma-2}^0)}=\mathfrak{A}^{r_{\gamma-1}},$$ thus
$\Gamma'(\P,\P)=\rK[\lambda]$, and $\lambda$ appears for the first
time in $\mathfrak{A}^{r_{\gamma-1}}$ .

 {\bf Lemma 5.4.2} There exists a localization given by $\tau(\lambda)$
  and a triangular formula:
\begin{equation}
 \left\{
\begin{array}{ccl}
\dz(a_1)^0 &=& h_{11}(\lambda,\mu)w_1 \\
\dz(a_2)^0 &=& h_{21}(\lambda,\mu)w_1\quad +h_{22}(\lambda,\mu)w_2 \\
\cdots   &  & \qquad\cdots \\ \dz(a_{i-1})^0 &=&
h_{i-1,1}(\lambda,\mu)w_1
+h_{i-1,2}(\lambda,\mu)w_2+\cdots+h_{i-1,i-1}(\lambda,\mu)w_{i-1}\\
\dz(a_i)^0 &=& h_{i1}(\lambda,\mu)w_1\quad\
+h_{i2}(\lambda,\mu)w_2\quad+\cdots+h_{i,i-1}(\lambda,\mu)w_{i-1}
\quad+\wt{h}(\lambda,\mu)\wt{w}
\end{array} \right.
\end{equation}
where $$h_{ll'}(\lambda,\mu)\in k[\lambda,\mu,
\tau(\lambda)^{-1}\tau(\mu)^{-1}],\ h_{ll}(\lambda,\lambda)\ne 0\
{\rm for}\ l'\leq l\  {\rm and}\ l=1,2,\cdots, i;$$  $w_1,w_2,
\cdots,w_{i-1}$ are free generators of $\overline{\Omega}$,
$\wt{w}$ is a linear form of $w_1, \cdots, w_{i-1}$,  or $\wt{w}$
is linearly independent of $\{w_1, w_2, \cdots, w_{i-1}\}$, and
$\wt{h}(\lambda,\mu)\ne 0$ but $\wt{h}(\lambda,\lambda)=0$.

{\bf Proof.} If $\dz(a_1)=0$ in Formula (29) of case A1 of 4.3, then $i=1$, $\wt{w}=0$. Otherwise, after a
localization given by $c_1(\lambda)$ in  Formula (30) of 4.3, $\dz(a_1)^0=h_{11}(\lambda, \mu)w_1$. If
$h_{11}(\lambda, \lambda)=0$, then $i=1$, $\wt{w}\ne 0$. If $h_{11}(\lambda, \lambda)\ne 0$,  we use the similar
procedure given in the proof of Lemma 5.4.1. Then we reach to   Formula (42), since $a_i=\lambda_{\gamma}$ must
appear at some stage $i$.

We use the following notations
\begin{equation}
\tau(\lambda)=\left\{\begin{array}{ll} f(\lambda)\prod_{l=1}^i
c_l(\lambda), &\mbox{ when }\gamma=1;\\ \prod_{l=1}^i
c_l(\lambda), &\mbox{ when }\gamma>1, \end{array}\right.
\quad\mbox{ and } h(\lambda)=\tau(\lambda)\prod_{l=1}^{i-1}
h_{ll}(\lambda,\lambda).
\end{equation}
The proof is completed.\hfill $\square$

\subsection{The minimal size assumption}

\kg {\bf Lemma 5.5.1.} Let $(\K, \M, H=0)$ be a bimodule problem,
and $(*)$ be a reduction sequence consisting of freely
parameterized triples, such that the end term $\mathfrak{A}^s$ of
$(***)$ has two vertices, and is in one of the following
configurations of Theorem 5.1.1:
 \begin{center}\unitlength=1mm
\begin{picture}(70,10) \put(9.00,0.00){\circle*{1.00}}
\put(26.00,0.00){\circle*{1.00}} \put(5.00,1.00){\oval(5,5)[t]}
\put(5.00,1.00){\oval(5,5)[bl]} \put(4.00,-2.00){\vector(3,1){3}}
\put(30.00,1.00){\oval(5,5)[t]} \put(30.00,1.00){\oval(5,5)[br]}
\put(31.00,-2.00){\vector(-3,1){3}}
\put(10.50,0.50){\line(1,0){13}} \put(23.50,0.50){\vector(1,0){1}}
 \put(10.50,-0.50){\line(1,0){1}}
\put(12.50,-0.50){\line(1,0){1}} \put(14.50,-0.50){\line(1,0){1}}
\put(16.50,-0.50){\line(1,0){1}} \put(18.50,-0.50){\line(1,0){1}}
\put(20.50,-0.50){\line(1,0){1}} \put(22.50,-0.50){\line(1,0){1}}
\put(23.50,-0.50){\vector(1,0){1}}
\put(0.00,1.0){\makebox{$\lambda$}}
\put(10.00,3.00){\makebox{$\P$}} \put(25.00,3.00){\makebox{$\Q$}}
 \put(15.50,1.00){\makebox{$a_1$}}
 \put(15.50,-3.50){\makebox{$w_1$}}
\put(35.00,1.00){\makebox{$\mu,$}}

\put(-16.00,0.00){\makebox{Case (1)}}
\put(55.00,0.00){\makebox{Case (2)}}
\end{picture}
\begin{picture}(32, 8)
\put(9.00,-1.00){\circle*{1.00}} \put(26.00,-1.00){\circle*{1.00}}
\put(5.00,0.00){\oval(5,5)[t]} \put(5.00,0.00){\oval(5,5)[bl]}
\put(4.00,-3.00){\vector(3,1){3}}
\put(10.50,-1.00){\line(1,0){13}}
\put(23.50,-1.00){\vector(1,0){1}}
\put(10.00,1.00){\makebox{$\P$}} \put(27.00,0.00){\makebox{$\Q$}}
\put(0.00,0.0){\makebox{$\lambda$}}
\put(16.00,1.00){\makebox{$a_1$}} \put(30.00,0.00){\makebox{ (or
dual).}}
\end{picture}
\end{center}
If we assume in addition that the size of the triples of $(*)$ is minimal such that $(*)$ meets any
configuration of 5.1.1. Then the local bocs $\mathfrak{A}^s_{\P}$ at $\P$ satisfies Formula (40) of 5.4 with all
the $$g_{ll}(\lambda,\lambda')\in k[\lambda, \lambda', \sigma_l(\lambda)^{-1}\sigma_l(\lambda')^{-1}]$$ being
invertible, where $\lambda'$ stands for the right multiplication of $\lambda$. And so does $\mathfrak{A}^s_{\Q}$
at $\Q$.

{\bf Proof.} Suppose the contrary, i.e. either we have Formula (42) of 5.4, or there exists some $1\leq e\leq j$
with $g_{ee}(\lambda,\lambda')\in k[\lambda, \lambda', \sigma_e(\lambda)^{-1}\sigma_e(\lambda')^{-1}]$ being
non-invertible in Formula (40). Then we would meet case (1) of 5.1.1 in sequence $(\hat{*})$, which is obtained
from $(*)$ by deletion of $(T\setminus\{\P\})$ according to Corollary 4.3.1. But the size of the triples of
$(\hat{*})$ is less than those of $(*)$, which is a contradiction to the minimal size assumption. \hfill
$\square$

{\bf Lemma 5.5.2.} Let $(\K, \M, H=0)$ be a bimodule problem, and
$(*)$ be a reduction sequence consisting  of freely parameterized
triples, such that the end term $\mathfrak{A}^s$ of $(***)$ is
local and in the configuration of case (1) of Theorem 5.1.1:
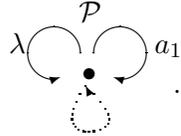
\begin{center}\unitlength=0.5pt
\begin{picture}(80,60)
 \put(-10,30){\oval(40,40)[t]} \put(-10,30){\oval(40,40)[bl]}
\put(-10,10){\vector(3,1){5}} \put(40,30){\oval(40,40)[t]}
\put(40,30){\oval(40,40)[br]} \put(40,10){\vector(-3,1){5}}

\put(11,9){$\bullet$} \qbezier[10](13,5)(-8,-28)(16,-30)
\qbezier[10](13,5)(50,-28)(13,-30) \put(17,-1){\vector(-1,2){3}}

 \put(9,54){\makebox{$\P$}}
\put(-44,30){\makebox{$\lambda$}} \put(66,30){\makebox{$a_1$}}

\end{picture}.
\end{center}\vskip 2mm
Suppose that $\lambda$ appears for the first time in
$\mathfrak{A}^r$ attached to a vertex $\P$ for some $1\leq r< s$.
If we assume in addition that the size of the triples of $(*)$ is
minimal,  such that $(*)$ meets any configuration of 5.1.1. Then
either

(1) $\mathfrak{A}^r$ is local, such that after performing the
procedure of Lemma 5.3.1 starting from $\mathfrak{A}^r$, we have
Formula (42) of 5.4, or Formula (40) with
$$g_{ee}(\lambda,\lambda')\in k[\lambda, \lambda',
\sigma_e(\lambda)^{-1}\sigma_e(\lambda')^{-1}]$$ being
non-invertible for some $1\leq e\leq j$; or

 (2) there exists some $r\leq l<s$, such that
 $\mathfrak{A}^l$ is in case (2) of Lemma 5.5.1.

{\bf Proof.} (1) If $\gamma>0$ in Lemma 5.3.1, we have done. If $\gamma=0$, such index $e$ must exist, since
$\delta(a_1)$ is not invertible.

(2) If $\mathfrak{A}^r$ is not local, then the induced local bocs $\mathfrak{A}^r_{\P}$ at $\P$ satisfies
Formula (40) of 5.4 with all the $g_{ll}(\lambda,\lambda')\in k[\lambda, \lambda',
\sigma_l(\lambda)^{-1}\sigma_l(\lambda')^{-1}]$ being invertible by the minimal size assumption. There must
exist some edges connecting $\P$ and other vertices. Otherwise $\P$ would be an isolated vertex, which
contradicts to $\mathfrak{A}^s$ being local. Continue the reduction from $\mathfrak{A}^r$, there must exist some
first edge $a_1$ connecting $\P$ and other vertex in some induced bocs $\mathfrak{A}^l$, such that
$\delta(a_1)=0$. Otherwise all such $a_1$ were deleted by regularizations, then $\P$ would be again an isolated
vertex. Thus $\mathfrak{A}^l$ is in  case (2) of Lemma 5.5.1 by the minimal size assumption. \hfill $\square$

\subsection{The classification of minimally wild bimodule
problems}

\kg {\bf Theorem 5.6.1.} Let $(\K,\M,H)$ be a bimodule problem,
and $(*)$ be a reduction sequence consisting of freely
parameterized bimodule problems. Suppose that
\begin{itemize}
\item[(1)] the size of the triples of $(*)$ is minimal; \item[(2)] there exists some minimal integer $r$ such
that $\mathfrak{A}^r$ is in one of the configurations of Theorem 5.1.1.
\end{itemize}
 Then there must exist a minimal integer $s$, such that the end term
  $\mathfrak{A}^s$ of $(***)$ is in one of the following five cases:
\begin{itemize}
\item[MW1.] \hskip 10mm \unitlength=1mm
\begin{picture}(30,10)
\put(0.00,0.00){\circle*{1.00}} \put(15.00,0.00){\circle*{1.00}}
\put(-5.00,1.00){\oval(5,5)[t]} \put(-5.00,1.00){\oval(5,5)[bl]}
\put(-6.00,-2.00){\vector(3,1){3}} \put(20.00,1.00){\oval(5,5)[t]}
\put(20.00,1.00){\oval(5,5)[br]}
\put(21.00,-2.00){\vector(-3,1){3}}
\put(0.50,0.50){\line(1,0){13}} \put(13.50,0.50){\vector(1,0){1}}
 \put(0.50,-0.50){\line(1,0){1}}
\put(2.50,-0.50){\line(1,0){1}} \put(4.50,-0.50){\line(1,0){1}}
\put(6.50,-0.50){\line(1,0){1}} \put(8.50,-0.50){\line(1,0){1}}
\put(10.50,-0.50){\line(1,0){1}} \put(12.50,-0.50){\line(1,0){1}}
\put(13.50,-0.50){\vector(1,0){1}}
\put(-10.00,1.0){\makebox{$\nu$}} \put(0.00,3.00){\makebox{$\P$}}
\put(15.00,3.00){\makebox{$\Q$}}
 \put(5.50,1.00){\makebox{$a_1$}}
 \put(5.50,-3.50){\makebox{$w_1$}}
\put(25.00,1.00){\makebox{$\kappa$}}
\end{picture} Bocs $\mathfrak{B}_1$ has two vertices $\P, \Q$,
and a layer  $L_1=(\Gamma'_1; \omega_1; a_1,$ \vskip 1mm \noindent$ \cdots, a_n;  v_1, \cdots, v_m )$, where
$\Gamma'_1(\P,\P)=k[\nu, g_{\P}(\nu)^{-1}]$, $\Gamma'_1(\Q, \Q)$ $=k[\kappa,g_{\Q}(\kappa)^{-1}]$, $\delta(a_1)$
$=f(\nu, \kappa)w_1$, such that $f(\nu, \kappa)\in k[\nu, \kappa, g_{\P}(\nu)^{-1}g_{\Q}(\kappa)^{-1}]$ is
non-invertible. Moreover the local  bocs $(\mathfrak{B}_1)_{\P}$ at $\P$ satisfies Formula (40) of 5.4, with all
the $g_{ll}(\nu,\nu')\in k[\nu, \nu', \sigma_l(\nu)^{-1}\sigma_l(\nu')^{-1}]$ given by Formula (41) being
invertible, and so does $(\mathfrak{B}_1)_{\Q}$ at $\Q$.

\item[MW2.] \begin{picture}(50, 8)
\put(10.00,-1.00){\circle*{1.00}}
\put(25.00,-1.00){\circle*{1.00}} \put(5.00,0.00){\oval(5,5)[t]}
\put(5.00,0.00){\oval(5,5)[bl]} \put(4.00,-3.00){\vector(3,1){3}}
\put(10.50,-1.00){\line(1,0){13}}
\put(23.50,-1.00){\vector(1,0){1}}
\put(10.00,1.00){\makebox{$\P$}} \put(27.00,0.00){\makebox{$\Q$}}
\put(0.00,0.0){\makebox{$\nu$}} \put(16.00,1.00){\makebox{$a_1$}}
\put(30.00,0.00){\makebox{ (or dual)}}
\end{picture} Bocs $\mathfrak{B}_2$ has two vertices $\P,\Q$,
and a layer  $L_2=(\Gamma'_2;$ \vskip 1mm  $ \omega_2; a_1,$ $
\cdots, a_n; v_1, \cdots, v_m )$, where $\Gamma'_2(\P,\P)=k[\nu,
g_{\P}(\nu)^{-1}]$, $\Gamma'_2(\Q,\Q)=k$,
 $\delta(a_1)=0$. Moreover the induced local bocs
$(\mathfrak{B}_2)_{\P}$ at $\P$ satisfies Formula (40) of 5.4, with all the $g_{ll}(\nu,\nu')\in k[\nu, \nu',
\sigma_l(\nu)^{-1}\sigma_l(\nu')^{-1}]$ given by Formula (41) being invertible.

\item[MW3.] \hskip 10mm \unitlength=0.5pt
\begin{picture}(60,60)
\put(-10,30){\oval(40,40)[t]} \put(-10,30){\oval(40,40)[bl]}
\put(-10,10){\vector(3,1){5}} \put(40,30){\oval(40,40)[t]}
\put(40,30){\oval(40,40)[br]} \put(40,10){\vector(-3,1){5}}

\put(11,9){$\bullet$} \qbezier[10](13,5)(-8,-28)(16,-30)
\qbezier[10](13,5)(50,-28)(13,-30) \put(17,-1){\vector(-1,2){3}}

 \put(9,54){\makebox{$\P$}}
\put(-44,30){\makebox{$\nu$}} \put(66,30){\makebox{$b$}}
\end{picture} \hskip 3mm
Bocs $\mathfrak{B}_3$ is local, and has a layer  $L_3=(\Gamma'_3;$ $ \omega_3; b_1,$ $ \cdots, b_j; v_1, \cdots,
v_m )$, \vskip 3mm where $\Gamma'_3(\P,\P)=k[\nu]$, and the differentials $b_1, b_2, \cdots,b_j$ are given by
Formula (40) of 5.4, such that $g_{ll}(\nu, \kappa)\in k[\nu, \kappa,\sigma_l(\nu)^{-1}, \sigma_l(\kappa)^{-1}]
$ for $1\leq l\leq e-1$ are invertible, but $g_{e,e}(\nu, \kappa)\in k[\nu, \kappa,\sigma_e(\nu)^{-1},
\sigma_e(\kappa)^{-1}]$ is not invertible.

 \item[MW4.]
\unitlength=1mm
\begin{picture}(20,8) \put(10,1){\circle*{0.6}}
\put(7,1){\circle{4}} \put(13,1){\circle{4}}
\put(10,4){\circle{4}} \put(9,2){\vector(0,1){0}}
\put(11,2){\vector(1,0){0}} \put(11,2){\vector(0,1){0}}
\put(2,0){$\lambda$} \put(13 ,4){$b$} \put(16,0){$a$}
\put(8.5,-4){$\P$}
\end{picture} Bocs $\mathfrak{B}_4$ is local, and has a layer  $L_4=(\Gamma'_4;$ $
\omega_4; a_1, \cdots, a_i;  b_1,$ $ \cdots, b_j; v_1,$ \vskip 2mm
$ \cdots, v_m )$, where $\Gamma'_4(\P,\P)=k[\lambda]$,  the
differentials $a_1,a_2, \cdots,a_i$ are given by Formula (42) of
5.4, such that
\begin{itemize}
  \item[(i)]  for any fixed $\lambda^0\in k$, $h(\lambda^0)\ne 0$
   in Formula (43), set
$a_i=\nu$, then the differentials of $b_1, \cdots, b_j$ given by
Formulae(40) have the property that  $ g_{ll}(\nu, \kappa)\in
k[\nu, \kappa, \sigma_l(\nu)^{-1}\sigma_l(\kappa)^{-1}]$ are all
invertible.

\item[(ii)] $\wt{w}$ is a linear combination of $w_1, \cdots,
w_{i-1}$ or $\wt{w}$ is linearly independent of  $w_1, \cdots,
w_{i-1}$, but $(\lambda-\mu)^2\mid \wt{h}(\lambda, \mu)$.
\end{itemize}
\item[MW5.] \unitlength=1mm
\begin{picture}(20,8) \put(10,1){\circle*{0.6}}
\put(7,1){\circle{4}} \put(13,1){\circle{4}}
\put(10,4){\circle{4}} \put(9,2){\vector(0,1){0}}
\put(11,2){\vector(1,0){0}} \put(11,2){\vector(0,1){0}}
\put(2,0){$\lambda$} \put(13 ,4){$b$} \put(16,0){$a$}
\put(8.5,-4){$\P$}
\end{picture} Bocs $\mathfrak{B}_5$ is the same as $\mathfrak{B}_4$, except
the hypothesis (ii) is changed to
\begin{itemize}
\item[(iii)] $\wt{w}$ is  linearly independent of $w_1, \cdots, w_{i-1}$, and $(\lambda-\mu)^2 \nmid
\wt{h}(\lambda, \mu)$.
\end{itemize}
\end{itemize}

{\bf Proof.} Suppose that $\mathfrak{A}^r$ is in one case of the
configurations of Theorem 5.1.1 with two vertices, then $r=s$ and
we obtain  MW1 or MW2 by Lemma 5.5.1. Now suppose $\mathfrak{A}^r$
is local, then $\lambda$  appears for the first time in
$\mathfrak{A}^r$ by Lemma 5.5.2. Denote $\lambda$ by $\lambda_0$,
it ensures that the domain of $\lambda_0$ in $\mathfrak{A}^r $ is
$k$. We perform the procedure given by Lemma 5.4.1. If there exist
some $\lambda_0^0, \lambda_1^0, \cdots, \lambda_{\gamma-1}^0$, and
$1\leq e\leq j$ in Lemma 5.4.1, such that $g_{ll}(\nu, \kappa)\in
k[\nu, \kappa,\sigma_l(\nu)^{-1}, \sigma_l(\kappa)^{-1}] $ for
$1\leq l\leq e-1$ are invertible, but $g_{e,e}(\nu, \kappa)\in
k[\nu, \kappa,\sigma_e(\nu)^{-1}, \sigma_e(\kappa)^{-1}]$ is not
invertible, we have MW3. Otherwise, $\gamma$ must be greater than
zero, since  $\mathfrak{A}^r $ is in case (1) of Theorem 5.1.1.
Thus we may fix some $\lambda_0^0, \lambda_1^0, \cdots,
\lambda_{\gamma-2}^0$ and obtain MW4 or MW5. \hfill $\square$

{\bf Remark}  We always assume that the reductions are given
according to the order of Definition 2.3.1, and the localizations
are given to the parameters at the left with respect to the
multiplication in case of the first arrows $a_1: \I \rightarrow
\J$ with $\I, \J$ both non-trivial. In particular Theorem 5.1.1
and 5.2.1 are proved by the above fixed procedure. Consequently if
$\mathfrak{A}^s$ given in Theorem 5.6.1 satisfy MW3 or MW4 or MW5,
then the parameter $\nu$ in $\mathfrak{B}_3$, or the parameter
$\lambda$ in $\mathfrak{B}_4$ or $\mathfrak{B}_5$ appears for the
first time in $\mathfrak{A}^s$, i.e. $\Gamma'_3(\P,\P)=k[\nu]$,
and $\Gamma'_4(\P,\P)=k[\lambda]$, $\Gamma'_5(\P,\P)=k[\lambda]$.

\newpage

 \bcen{\section{Exact categories and almost split conflations}}\ecen

\subsection{Some basic concepts}

\kg Let $\A$ be an additive category with Krull-Schmidt property.
We recall from \cite{GR}  and \cite{DRSS}  the following notions.
A pair $(\iota, \pi)$ of composable morphisms
$M\stackrel{\iota}{\longrightarrow}E\stackrel{\pi}{\rightarrow}L$
in $\A$ is called {\it exact} if $\iota$ is a kernel of $\pi$ and
$\pi$ is a cokernel of $\iota$.  Recall that $\iota$ is a {\it
kernel} of $\pi$, provided that (1) $\iota$ is monic; (2)
$\iota\pi=0$; (3) $\forall \varphi: N\rightarrow E$ with
$\varphi\pi=0$, there exists a $\varphi': N\rightarrow M$ such
that $\varphi'\iota=\varphi$. Dually we have the definition of
$\pi$ being a {\it cokernel} of $\iota$.

Let ${\cal E}$ be a class of exact pairs
$M\stackrel{\iota}{\longrightarrow}E\stackrel{\pi}{\longrightarrow}L$
which is closed under isomorphisms. The morphisms $\iota$ and
$\pi$ appearing in a pair $(\iota, \pi)$ of ${\cal E}$ are called
an {\it inflation} and a {\it deflation} of ${\cal E}$
respectively, and the pair itself is called a {\it conflation}.
The class ${\cal E}$ is said to be an {\it exact structure} on
$\A$, and $(\A, {\cal E})$ an {\it exact category} if the
following axioms are satisfied:

{\bf E1.} The composition of two deflations is a deflation.

{\bf E2.} For each $\varphi$ in $\A(L', L)$ and each deflation
$\pi$ in $\A(E,L)$, there are some $E'$ in $\A$, an $\varphi'$ in
$\A(E',E)$ and a deflation  $\pi': E'\rightarrow L' $ such that
$\pi'\varphi=\varphi' \pi$.

{\bf E3.} Identities are deflations. If $\varphi\pi$ is a
deflation, then so is $\pi$.

{\bf E3$~^{op}$.} Identities are inflations. If $\iota\varphi$ is
an inflation, then so is $\iota$.

An object $L$ in $\A$ is said to be {\it ${\cal E}$-projective}
(or {\it projective} for short ) if any conflation ending at $L$
is split. Dually an object $M$ in $\A$ is said to be {\it ${\cal
E}$-injective} (or {\it injective} for short ) if any conflation
starting at $M$ is split.

Let $\A$ be a Krull-Schmidt category. A morphism $\pi:
E\rightarrow L $ in $\A$ is called {\it right almost split} if it
is not a retraction and for any non-retraction $\varphi:
N\rightarrow L$, there exists a morphism $\psi: N\rightarrow E$
such that $\varphi=\psi\pi$. We say that $\A$ has {\it right
almost split morphisms} if for all indecomposable $L$ there exist
 right almost split morphisms ending at $L$. Dually we define
{\it left almost split morphisms}. We say that {\it $\A$ has
almost split morphisms} if $\A$ has right and left almost split
morphisms.

A morphism $\pi: E\rightarrow L$ is called {\it right minimal} if
every endomorphism $\eta: E\rightarrow E$ with the property that
$\pi=\eta\pi$ is an isomorphism. A {\it left minimal morphism}
$\iota: M\rightarrow E$ is defined dually.

Suppose that the Krull-Schmidt category $\A$ carries an exact
structure ${\cal E}$. Let
$$ (e) \qquad
M\stackrel{\iota}{\longrightarrow}E\stackrel{\pi}{\longrightarrow}L$$
be a conflation. Then the following assertions are equivalent.

(i) $\iota$ is minimal left almost split;

(ii) $\pi$ is minimal right almost split;

(iii) $\iota$ is left almost split and $\pi$ is right almost
split.

The conflation $(e)$ in the above condition is said to be an {\it
almost split conflation}. The exact category $(\A,{\cal E})$ is
said to {\it have almost split conflations} if

(i) $\A$ has almost split morphisms;

(ii) for any indecomposable non-projective $L$, there exists an
almost split conflation $(e)$ ending at $L$;

(iii)
 for any indecomposable non-injective $M$, there exists an almost
split conflation $(e)$ starting at $M$.

An almost split conflation
$M\stackrel{\iota}{\longrightarrow}E\stackrel{\pi}{\longrightarrow}L$
is uniquely determined by $M$ or $L$. We introduce the notation
$M=\tau(L)$ or $L=\tau^{-1}(M)$.

Let us fix an exact category $(\A,{\cal E})$. For given objects
$M$ and $L$ of $\A$, two conflations $(\iota,\pi)$ and
$(\iota',\pi')$ starting from $M$ and ending at $L$ respectively
are said to be {\it equivalent}, if there exists a commutative
diagram as follows:
$$
\begin{CD} M@>\iota>> E@>\pi>> L\\
@Vid VV @VV\eta V @VVid V\\
M@>>\iota'> E' @>>\pi'> L
\end{CD}
$$

\subsection{Exact categories $P(\Lambda)$ and $P_1(\Lambda)$}

\kg The results of this subsection are mainly quoted from
\cite{B2} and \cite{ZZ}. Let $\Lambda$ be a finite-dimensional
$\rK$-algebra, and $P(\Lambda)$ be defined in 2.1.  A sequence of
$P(\Lambda)$,
\begin{equation}
\begin{array}{c}
   (e) \qquad\quad  \begin{CD}
0@>>>P_1 @>f_1>> W_1 @>g_1>> Q_1 @>>>0\\
&& @V\alpha VV @VV\gamma V @VV\beta V\\
0@>>>P_0 @>>f_0> W_0 @>>g_0>Q_0 @>>>0
\end{CD}
\end{array}
\end{equation}
 is an {\it exact pair} if the two rows are exact and therefore split
in $\Lambda$-mod. If ${\cal E}$ stands for the set of such exact
pairs, then $(e)$ is a conflation.

{\bf Theorem 6.2.1} ${\cal E}$ is an exact structure and
$(P(\Lambda), {\cal E})$ is an exact category.

{\bf Proof.} E1, E3 and E3$^{op}$ are obvious. We only check E2.
If $(g_1, g_0)$ given in (44) is a deflation, and $(\varphi_1,
\varphi_0):
(Q'_1\stackrel{\beta'}{\longrightarrow}Q'_0)\rightarrow
(Q_1\stackrel{\beta}{\longrightarrow}Q_0)$ is a morphism in
$P(\Lambda)$. Then the following diagram gives a deflation $(g'_1,
g'_0)$  and a morphism $(\varphi'_1, \varphi'_0)$, which satisfy
E2.
\begin{center}\unitlength=0.8mm
\begin{picture}(50,40)
\put(0,10){$W_1$} \put(6,11){\vector(1,0){30}} \put(37,10){$Q_1$}
\put(0,35){$W'_1$} \put(6,36){\vector(1,0){30}}
\put(37,35){$Q'_1$} \put(3,34){\vector(0,-1){20}}
\put(39,34){\vector(0,-1){20}}

\put(6,35){\vector(1,-1){8}} \put(42,35){\vector(1,-1){8}}
\put(6,10){\vector(1,-1){8}} \put(42,10){\vector(1,-1){8}}

\put(15,24){$W'_0$} \put(51,24){$Q'_0$} \put(15,-1){$W_0$}
\put(51,-1){$Q_0$}\put(17,23){\vector(0,-1){20}}
\put(53,23){\vector(0,-1){20}} \put(21,25){\vector(1,0){30}}
\put(21,0){\vector(1,0){30}}

\put(17,38){$g'_1$} \put(30,27){$g'_0$} \put(22,7){$g_1$}
\put(33,2){$g_0$} \put(-2,20){$\varphi'_1$}
\put(12,15){$\varphi'_0$} \put(40,20){$\varphi_1$}
\put(54,12){$\varphi_0$} \put(8,26){$\gamma'$} \put(8,3){$\gamma$}
\put(46,31.5){$\beta'$} \put(46,6.5){$\beta$}
\end{picture}
\end{center}
where $W'_0=\{(w_0,q_0)\in W_0\oplus Q'_0\mid
g_0(w_0)=\varphi_0(q_0)\}$, $W'_1=\{(w_1,q_1)\in W_1\oplus
Q'_1\mid g_1(w_1)=\varphi_1(q_1)\}$,
$\gamma'=\left(\begin{array}{cc} \gamma & 0\\ 0&
\beta'\end{array}\right)$, $g'_0=\left(\begin{array}{c}  0\\
1\end{array}\right)$, $g'_1=\left(\begin{array}{c}  0\\
1\end{array}\right)$, $\varphi'_0=\left(\begin{array}{c}  1\\
0\end{array}\right)$,
$\varphi'_1=\left(\begin{array}{c}  1\\
0\end{array}\right)$.\hfill$\square$

\medskip
It is well known that  $\rad P\stackrel{\iota}{\hookrightarrow}P$
is a right almost split morphism in $\Lambda$-mod for any
indecomposable $P\in \Lambda$-proj. We denote by $r(P)$ a
projective cover of $\rad P$, and by $\rho$ the composition of
$r(P)\stackrel{\pi}{\longrightarrow}\rad
P\stackrel{\iota}{\hookrightarrow}P$, then
$r(P)\stackrel{\rho}{\longrightarrow}P$ is a minimal right almost
split morphism in $\Lambda$-proj. On the other hand, we denote by
$\Lambda$-inj the full subcategory of $\Lambda$-mod consisting of
injective $\Lambda$-modules.  It is also well known that
$I\stackrel{p}{\longrightarrow} I/\soc I$ is a left almost split
morphism in $\Lambda$-mod for any indecomposable $I\in
\Lambda$-inj. There are two nice equivalent functors
$F=D\Lambda\otimes_{\Lambda}-: \Lambda\mbox{-proj}\rightarrow
\Lambda\mbox{-inj}$ and $G=Hom_{\Lambda}(D(\Lambda), -):
\Lambda\mbox{-inj}\rightarrow \Lambda\mbox{-proj} $.  We denote
 $G(I)$ by $P$, and if $e(I)$ is an injective envelope of $I/soc(I)$,
 we denote $G(e(I))$ by $l(P)$. Let $\lambda$ be the image
of the composition $I\stackrel{p}{\longrightarrow} I/\soc\, I
\stackrel{i} {\hookrightarrow} e(I)$ under the functor $G$, then
$P\stackrel{\lambda}{\longrightarrow} l(P)$ is a  left almost
split morphism in $\Lambda$-proj. Now consider the following
objects of $P(\Lambda)$:
$$
J(P)=(P\stackrel{id}{\longrightarrow}P),\;
T(P)=(0\stackrel{0}{\longrightarrow}P),\;
Z(P)=(P\stackrel{0}{\longrightarrow}0) $$
 and
$$
R(P)=(r(P)\stackrel{\rho}{\longrightarrow}P),\;
L(P)=(P\stackrel{\lambda}{\longrightarrow}l(P)),\;
B(P)=(r(P)\stackrel{\rho\lambda}{\longrightarrow}l(P))
$$
We refer $U(P)$ to a minimal projective presentation of $\rad P$,
and $V(P)$ to the image of a minimal injective copresentation of
$I/\soc I$ under the functor $G$.

{\bf Lemma 6.2.1} (1) $\{J(P),T(P)\mid \forall\
\mbox{indecomposable } P\in\Lambda\mbox{-proj}\}$ is a complete
set of iso-classes of indecomposable ${\cal E}$-projectives.
Moreover, $R(P)\stackrel{(\rho\ id)}{\longrightarrow}J(P)$ and
$U(P)\stackrel{(0\ \rho)}{\longrightarrow}T(P)$ are right minimal
almost split morphisms in $P(\Lambda)$.

(2) $\{J(P),Z(P)\mid \forall\ \mbox{indecomposable }
P\in\Lambda\mbox{-proj}\}$ is a complete set of iso-classes of
indecomposable ${\cal E}$-injectives. Moreover, $J(P)\stackrel{(
id\ \lambda)}{\longrightarrow}L(P)$ and $Z(P)\stackrel{(\lambda\
0)}{\longrightarrow} V(P)$ are left minimal almost split morphisms
in $P(\Lambda)$.\hfill $\Box$

{\bf Proposition 6.2.1} The exact category $(P(\Lambda), {\cal
E})$ has almost split conflations.\hfill$\Box$

 We will pay most attention to $P_1(\Lambda)$ whose objects have no
direct summands $J(P)$. If ${\cal E}_1$ stands for the induced exact structure of ${\cal E}$ restricted to
$P_1(\Lambda)$, then $(P_1(\Lambda), {\cal E}_1)$ is also an exact category.

{\bf Lemma 6.2.2} (1) $\{T(P),L(P)\mid \forall\
\mbox{indecomposable } P\in\Lambda\mbox{-proj}\}$ is a complete
set of iso-classes of indecomposable ${\cal E}_1$-projectives.
Moreover, $U(P)\stackrel{(0\ \rho)}{\longrightarrow}T(P)$ and
$B(P)\stackrel{(\rho\ id)}{\longrightarrow}L(P)$ are right minimal
almost split morphisms in $P_1(\Lambda)$.

(2) $\{Z(P),R(P)\mid \forall\ \mbox{indecomposable }
P\in\Lambda\mbox{-proj}\}$ is a complete set of iso-classes of
indecomposable ${\cal E}_1$-injectives. Moreover,
$Z(P)\stackrel{(\lambda\ 0)}{\longrightarrow} V(P)$ and
$R(P)\stackrel{( id\ \lambda)}{\longrightarrow}B(P)$ are left
minimal almost split morphisms in $P_1(\Lambda)$.\hfill$\Box$

{\bf Proposition 6.2.2} The exact category $(P_1(\Lambda), {\cal
E}_1)$ has almost split conflations.\hfill$\Box$

{\bf Theorem 6.2.2} If $(e)$ given by Formula (44) is a conflation
in $P_1(\Lambda)$, and
$$
\begin{CD}  Cok(e):\qquad
 Cok(\alpha)@>Cok(f_1,f_0)>> Cok(\gamma) @>Cok(g_1,g_0)>>
 Cok(\beta) \rightarrow 0
\end{CD}
$$
is the corresponding exact sequence in $\Lambda$-mod, such that
$Cok(\alpha)$ is non-zero and neither injective nor simple. Then
$Coker(e)$ is an almost split sequence in $\Lambda$-mod, if and
only if $(e)$ is an almost split conflation in $P_1(\Lambda)$.

{\bf Proof.} By 2.3, 2.4 and 2.6 in \cite{ZZ}, Proposition 5.6 in
\cite{B2}. \hfill$\square$

\subsection{Morphisms in $R(\mathfrak{A})$}

\kg Let $\mathfrak{A}$ be a layered bocs.  The purpose of this
subsection is to transfer some morphisms of $R(\mathfrak{A})$ to
the morphisms of a module category over an algebra in order to use
some nice properties of module categories.

{\bf Lemma 6.3.1}\cite{Ro}. Let $\mathfrak{A}=(\Gamma,\Omega)$ be
a bocs with a layer $L=(\Gamma';\omega; a_1,\cdots,a_n; v_1,$
$\cdots,$ $ v_m)$, and  $M\in R(\mathfrak{A})$. Suppose that
 $$\{
M'_{\I}=M_{\I}\mid \forall\ \I\in \ind\Gamma'\}$$
 is a
set of $\rK$-vector spaces,
 $$\{\varphi_{\I}=id: M_{\I}\rightarrow
M'_{\I};\; \varphi(v_j): M_{\I_j}\rightarrow M'_{\J_j}\mid
\forall\ \I\in \ind\Gamma', \ v_j:\I_j\rightarrow \J_j\}$$ is a
set of linear maps. Then there exists a unique $M'\in
R(\mathfrak{A})$ such that $\varphi: M\rightarrow M'$ is an
isomorphism in $R(\mathfrak{A})$.

{\bf Proof.} We present a brief proof here for reader's
convenience. Let $M'(\lambda)=M(\lambda)$ for any free parameter
$\lambda$, and we set $M'(a_1)=M(a_1)-\varphi(\delta(a_1))$.
Suppose that $M'(a_1), \cdots,$ $ M'(a_{i-1})$ have been fixed,
then we set $M'(a_i)=M(a_i)-\varphi(\delta(a_i))$, since
$\delta(a_i)$ involves only $a_1, \cdots, a_{i-1}$.  Thus our
assertion follows by induction. \hfill$\square$

Let $\mathfrak{A}=(\Gamma,\Omega)$ be a  bocs with a layer
$L=(\Gamma';\omega;a_1,\cdots,a_n;v_1,\cdots,v_m)$. Then the bocs
 $\mathfrak{A}_0=(\Gamma,\Gamma)$
is called the {\it principal bocs} of $\mathfrak{A}$. It is clear
that the representation category $R(\mathfrak{A}_0)$ is just the
module category $\Gamma$-mod. If $\varphi: M\rightarrow N $ is a
morphism in $R(\mathfrak{A})$, then $\varphi=(\varphi_{\I},
\varphi(v_j)\mid \forall\, \I\in \ind \Gamma', j=1,2,\cdots,m)$
(see 4.6). We denote $(\varphi_{\I}\mid \forall\ \I\in
\ind\Gamma')$ by $\varphi_0$. From now on we assume that the bocs
is triangular also on $\{v_1, \cdots, v_m\}$, i.e. $\delta(v_j)$
involves only $v_1, \cdots, v_{j-1}$. If $\mathfrak{A}$
corresponds a bimodule problem, then it is obviously  triangular
on $\{v_1, \cdots, v_m\}$.

{\bf Proposition 6.3.1} \cite{O}. (1) If $\iota: M\rightarrow E$
is a morphism of $R(\mathfrak{A})$ with $\iota_0$ injective, then
there exists an isomorphism $\eta$ and a commutative diagram in
$R(\mathfrak{A})$:
$$
\begin{CD}
 &  & M @>\iota>> E\\
 && @VidVV @VV\eta V\\
0@>>>M @>>\iota^{\prime}>E^{\prime}
\end{CD}
$$
such that the second sequence is exact in $R(\mathfrak{A}_0)$.

(2) If $\pi: E\rightarrow L$ is a morphism of $R(\mathfrak{A})$
with  $\pi_0$ surjective, then there exists an isomorphism $\eta$
and a commutative diagram in $R(\mathfrak{A})$:
$$
\begin{CD}
  E @>\pi>> L\\
  @V\eta VV @VVidV\\
E^{\prime}@>>\pi^{\prime}>L @>>>0
\end{CD}
$$
such that  the second sequence is exact in $R(\mathfrak{A}_0)$.

{\bf Proof.} We present a complete proof here, since Ovsienko has
not published his paper.

(1) Our proof will be shown by induction on the dotted arrows.
First we fix a set of vector spaces $\{E^1_{\I}=E_{\I}\mid
\forall\ \I\in \ind \Gamma^{\prime} \}$ and a set of linear maps
$\{\varphi_{\I}=id: E_{\I}\rightarrow E^1_{\I}; \varphi(v_j)=0\mid
\forall\ \I\in \ind \Gamma^{\prime}, \mbox{ and } j\neq 1\}$.
$\varphi(v_1)$ will be defined below such that
$\iota^1=\iota\varphi: M\rightarrow E^1$ with $\iota^1(v_1)=0$.
Suppose $v_1:\I\rightarrow \J$, then $0=
\iota^1(v_1)=\iota_{\I}\varphi(v_1)+\iota(v_1)\varphi_{\J}$ (see
4.6), i.e. $\iota_{\I}\varphi(v_1)+\iota(v_1)=0$, which yields a
commutative diagram since $\iota_{\I}$ is injective.:
 \bcen \unitlength=0.8mm
\begin{picture}(40,25)

\put(0,0){$E_{\J}$} \put(26,0){$E^1_{\J}$} \put(0,20){$M_{\I}$}
\put(26,20){$E_{\I}$}

\put(7,1){\vector(1,0){18}} \put(7,21){\vector(1,0){18}}
\put(3,18){\vector(0,-1){13}}\put(28,18){\vector(0,-1){13}}

\put(-11,10.50){$-\iota(v_1)$} \put(29,10.50){$\varphi(v_1)$}
\put(15,22){$\iota_{\I}$} \put(10,-3){$\scriptstyle
\varphi_{\J}=id$}

\end{picture}
\ecen  Thus $E^1\in R(\mathfrak{A})$ is obtained by Lemma 6.3.1.
And a commutative diagram \bcen \unitlength=0.8mm
\begin{picture}(40,25)

\put(0,0){$M$} \put(26,0){$E^1$} \put(0,20){$M$} \put(26,20){$E$}

\put(7,1){\vector(1,0){18}} \put(7,21){\vector(1,0){18}}
\put(3,18){\vector(0,-1){13}}\put(28,18){\vector(0,-1){13}}

\put(-3,10.50){$id$} \put(29,10.50){$\eta^1=\varphi$}
\put(15,22){$\iota$} \put(14,-3){$\iota^1$}

\end{picture}
\ecen follows from  the structure with $\eta^1=\varphi$,
$\iota^1(v_1)=0$. The assertion is obtained by the triangularity
of $\{v_1,\cdots,v_m\}$ and induction.

(2) is obtained dually. \hfill$\square$

{\bf Proposition 6.3.2} \cite{O}. Let $\mathfrak{A}$ be a layered
bocs, and $(e):  M\stackrel{\iota}{\longrightarrow}
E\stackrel{\pi}{\longrightarrow} L$ with $\iota\pi=0$ be   a pair
of composable morphisms  in $R(\mathfrak{A})$. If
$$(e_0)\quad
0\longrightarrow M\stackrel{\iota_0}{\longrightarrow}
E\stackrel{\pi_0}{\longrightarrow} L\longrightarrow 0$$ is exact
in the category of vector spaces, then there exists an isomorphism
$\eta$ and a commutative diagram in $R(\mathfrak{A})$:
$$\begin{CD} (e) \qquad &     & M@>\iota>> E @>\pi>> L  \\
& &@VidVV @VV\eta V@VVidV \\
(e') \qquad   0 @>>> M @>>\iota^{\prime}>
E^{\prime}@>>\pi^{\prime}> L @>>>0
\end{CD}
$$
such that $(e')$ is an exact sequence in $R(\mathfrak{A}_0)$.

{\bf Proof. }  First we fix a set of vector spaces
$\{E^1_{\I}=E_{\I}\mid \forall\ \I\in \ind \Gamma^{\prime} \}$,
and a set of linear maps $\{\varphi_{\I}=id: E_{\I}\rightarrow
E^1_{\I}; \varphi(v_j)=0\mid \forall\ \I\in \ind \Gamma^{\prime},
j\neq 1\}$.  $\varphi(v_1)$ will be defined for the purpose of
constructing some $E^1\in R(\mathfrak{A})$ and some
$\iota^1=\iota\varphi: M\rightarrow E^1$, $\pi^1=\varphi^{-1}\pi:
E^1\rightarrow L $ such that $\iota^1(v_1)=0$, $\pi^1(v_1)=0$.
Suppose $v_1:\I\rightarrow \J$, then  $0=
\iota^1(v_1)=\iota_{\I}\varphi(v_1)+\iota(v_1)\varphi_{\J}$, i.e.
$\iota_{\I}\varphi(v_1)+\iota(v_1)=0$;
$0=\pi^1(v_1)=\varphi_{\I}^{-1}\pi(v_1)+\varphi^{-1}(v_1)\pi_{\J}$,
i.e. $\pi(v_1)+\varphi^{-1}(v_1)\pi_{\J}=0$ (see 4.6). Since
$\varphi\varphi^{-1}=id$,
$0=(\varphi\varphi^{-1})(v_1)=\varphi_{\I}\varphi^{-1}(v_1)+
\varphi(v_1)\varphi^{-1}_{\J}$, i.e.
$\varphi^{-1}(v_1)=-\varphi(v_1)$. Therefore
$\pi(v_1)-\varphi(v_1)\pi_{\J}=0$. The hypothesis
$(\iota\pi)(v_1)=0$ given by $\iota\pi=0$  yields a commutative
diagram \bcen \unitlength=0.8mm
\begin{picture}(65,30)
\put(0,25){$0$} \put(4,26){\vector(1,0){15}} \put(20,
25){$M_{\I}$} \put(27, 26){\vector(1,0){20}} \put(48,25){$E_{\I}$}

\put(22,24){\vector(0,-1){20}} \put(50,24){\vector(0,-1){20}}
\put(26.5,5){\vector(-1,-1){1}} \multiput(27,5)(1,1){19}{$\cdot$}

\put(20, 0){$E_{\J}$} \put(27, 1.7){\vector(1,0){20}}
\put(28,-0.2){\vector(-1,0){1}}
\multiput(29,-0.2)(2,0){9}{\line(1,0){1}} \put(48,0){$L_{\J}$}
\put(54,1){\vector(1,0){15}} \put(70, 0){$0$}

\put(35,-2.5){$\xi$} \put(35,3){$\pi_{\J}$}
\put(35,28){$\iota_{\I}$} \put(7,12){$-\iota(v_1)$} \put(27, 16){$
\varphi(v_1)$} \put(51,12){$\pi(v_1)$}

\end{picture}
\ecen We claim that there exists a linear map $\varphi(v_1):
E_{\I}\rightarrow E_{\J}$ such that
$\iota_{\I}\varphi(v_1)=-\iota(v_1)$ and
$\varphi(v_1)\pi_{\J=}\pi(v_1)$. In fact, if
$E_{\I}=\iota_{\I}(M_{\I})\oplus N$ is taken in the category of
vector spaces, we can define
$((x)\iota_{\I}+y)\varphi(v_1)=(x)(-\iota(v_1))+ (y)\pi(v_1)\xi$
with $\xi: L_{\J}\longrightarrow E_{\J}$ being a retraction of
$\pi_{\J}$. Thus $E^1\in R(\mathfrak{A})$ is obtained by Lemma
6.3.1, and a commutative diagram
$$
\begin{CD}
M @>\iota>> E @>\pi>> L\\
@VidVV @VV\eta^1=\varphi V @VVidV\\
M @>\iota^1>> E^1 @>\pi^1>> L
\end{CD}
$$ follows from the  structure, where $\eta^1=\varphi$,
$\iota^1(v_1)=0$, $\pi^1(v_1)=0$. The assertion is obtained by the
triangularity of $\{v_1,\cdots,v_m\}$ and induction.
\hfill$\square$

\subsection{An exact structure on $R(\mathfrak{A})$}

\kg There is a natural exact structure on $R(\mathfrak{A})$, which
we will illustrate in this subsection.

{\bf Lemma 6.4.1} (1) $\iota: M\rightarrow E$ is monic if
$\iota_0: M\rightarrow E$ is injective.

(2) $\pi: E\rightarrow L$ is epic if $\pi_0: E\rightarrow L$ is
surjective.

{\bf Proof.} (1) If $\iota_0$ is injective, Proposition 6.3.1 (1)
 gives a commutative diagram with $\iota': M\rightarrow E'$ in
 $R(\mathfrak{A}_0)$. Given any
morphism $\varphi: N\rightarrow M$ with $\varphi\iota=0$, we have
$\varphi\iota\eta=\varphi\iota'=0$. Then $\varphi_0\iota'_0=0$
yields $\varphi_0=0$; and for any $v:\I\rightarrow \J$,
$0=(\varphi_0\iota')(v)=\varphi(v)\iota'_J$ yields $\varphi(v)=0$.
Thus $\varphi=0$ and $\iota$ is monic.

(2) is obtained dually. \hfill $\square$

{\bf Lemma 6.4.2} A pair of composable morphisms $(e):
M\stackrel{\iota}{\longrightarrow}
E\stackrel{\pi}{\longrightarrow} L$ with $\iota\pi=0$ is exact in
$R(\mathfrak{A})$, if $(e_0): 0\longrightarrow
M\stackrel{\iota_0}{\longrightarrow}
E\stackrel{\pi_0}{\longrightarrow} L\longrightarrow 0$ is exact as
a sequence  of  vector spaces.

{\bf Proof.} (1) If $(e_0)$ is exact, Lemma 6.4.1 (1) tells us
that $\iota$ is monic.

(2) $\iota\pi=0$.

(3) Proposition 6.3.2 gives a commutative diagram.
 If $\varphi: N\rightarrow E $ with
$\varphi\pi=0$, then $\varphi(\eta\pi^{\prime})=0$. Let
$\xi=\varphi\eta$, then $\xi\pi^{\prime}=0$ implies that
$\xi_{\I}\pi_{\I}^{\prime}=0$ and $\xi(v)\pi_{\J}^{\prime}=0$ for
any vertex $\I$ and any dotted arrow $v:\I\rightarrow \J$. Let
$\varphi^{\prime}:N\rightarrow M$ be given by
$\varphi^{\prime}_{\I}\iota^{\prime}_{\I}=\xi_{\I}$,
$\varphi^{\prime}(v)\iota^{\prime}_{\J}=\xi(v),$ then we obtain
$\varphi^{\prime}\iota^{\prime}=\xi$. Thus
$\varphi^{\prime}\iota^{\prime}\eta^{-1}=\varphi$, i.e.
$\varphi^{\prime}\iota=\varphi$.

Therefore $\iota$ is a kernel of $\pi$. It can be proved dually
that $\pi$ is a cokernel of $\iota$.
 \hfill$\square$

\medskip
{\bf Definition 6.4.1} Let $\mathfrak{A}$ be a layered bocs, we
define a set ${\cal E}$ consisting of exact pairs
$M\stackrel{\iota}{\longrightarrow}
E\stackrel{\pi}{\longrightarrow} L$ in $R(\mathfrak{A})$, such
that $\iota\pi=0$ and
$$0\longrightarrow M\stackrel{\iota_0}{\longrightarrow}
E\stackrel{\pi_0}{\longrightarrow} L\longrightarrow 0$$ are exact
as vector spaces. \hfill$\square$

It is clear that ${\cal E}$ is closed under isomorphisms.

{\bf Lemma 6.4.3} (1) $\iota: M\rightarrow E$ is an inflation, if
and only if  $\iota_0$ is injective.

(2) $\pi: E\rightarrow L$ is a deflation, if and only if  $\pi_0$
is surjective.

{\bf Proof.} (1) If $\iota_0$ is injective,   Proposition 6.3.1
(1) gives an isomorphism $\eta:E\rightarrow E'$ in $R(\mathfrak
{A})$, and a morphism $\iota'$ in $R(\mathfrak{A}_0)$. Moreover
$\iota_0\eta_0=\iota'$ ensures that $\iota'$ is injective. Let
$L=\Coker \iota'$, and $\pi':E'\rightarrow E$ be the natural
projection. Define $\pi=\eta\pi'$, we have the following
commutative diagram:
\begin{center}
\unitlength=0.8mm
\begin{picture}(60,40)
\put(2,20){\makebox{$M$}} \put(11,19){\vector(2,-1){20}}
\put(11,22){\vector(2,1){20}} \put(33.50,30.50){\makebox{$E$}}
\put(33,6){\makebox{$E'$}}

 \put(41,31.50){\vector(2,-1){20}}
\put(42,8){\vector(2,1){20}} \put(62,19.50){\makebox{$L$}}
\put(34.50,29){\vector(0,-1){18}}

\put(35,20){\makebox{$\eta$}} \put(20,27){\makebox{$\iota$}}
\put(18,11){\makebox{$\iota'$}} \put(50,27){\makebox{$\pi$}}
\put(49,9){\makebox{$\pi'$}}

\end{picture}
\end{center}
Then $\iota\pi=\iota\eta\cdot \eta^{-1}\pi=\iota'\pi'=0$. And
$\eta_0^{-1}\pi'_0=\pi_0$ ensures that $\pi_0$ is surjective. Thus
$M\stackrel{\iota}{\longrightarrow}
E\stackrel{\pi}{\longrightarrow} L$ belongs to ${\cal E}$, $\iota$
is an inflation.

(2) can be proved dually. \hfill$\square$

\medskip
{\bf Theorem 6.4.1} ${\cal E}$ is an exact structure on
$R(\mathfrak{A})$, and $(R(\mathfrak{A}),{\cal E})$ is an exact
category.

{\bf Proof.} E1, E3 and E3$~^{op}$ are obvious. It  suffices to
prove E2. Suppose that $\pi: E\rightarrow L$ is a deflation, and
$\varphi: L'\rightarrow L$ is a morphism in $R(\mathfrak{A})$.
  Consider the following diagram
\begin{center}\unitlength=1mm
\begin{picture}(100, 40)
  \put(26,0){$E$} \put(32,1){\vector(1,0){30}} \put(63,0){$L$}

 \put(28,34){\vector(0,-1){30}}
\put(64,34){\vector(0,-1){30}}

 \put(26,35){$E'$}
\put(32,36){\vector(1,0){30}} \put(63,35){$L'$}

\put(30,34){\vector(1,-1){6}} \put(39,24){\vector(1,-1){6}} \put(35,24){$N$} \put(40,14){$E\oplus L'$}
\put(51,12){\vector(1,-1){10}} \put(57,25){${0 \choose 1}$} \put(55,8){$ {\s \pi \choose \s -\varphi}$}
\put(40,12){\vector(-1,-1){10}} \put(37,5){${1 \choose 0}$} \put(47,18){\vector(1,1){15}}

 \put(44,-1){$\pi$} \put(33, 31){$\kappa$} \put(44,37){$\pi'$} \put(23,17){$\varphi'$}
  \put(64.5,17){$\varphi$} \put(43,21){$\eta$}
 \end{picture}
\end{center}
 where
$({\pi\atop -\varphi})_0$ is surjective, since $\pi_0$ is already
surjective. Thus there exists some object $N$ and an isomorphism
$\eta: N\rightarrow E\oplus L'$ in $R(\mathfrak{A})$ such that
$\eta({\pi\atop -\varphi}): N\rightarrow L$ is surjective in
$R(\mathfrak{A}_0)$ by item (2) of Proposition 6.3.1. Let
$E'=\Ker(\eta({\pi\atop -\varphi}))$ and $\kappa: E'\rightarrow N$
be the natural embedding. We set $\pi'=\kappa\eta{0\choose 1}$.
Since for any $y\in L'$, $(y)\varphi_0\in L$,  there exists some
$x\in E$ with $(x)\pi_0=(y)\varphi_0$. Consequently
$(x,y)\eta_0^{-1}(\eta({\pi\atop -\varphi}))=0$, i.e. $(x,\
y)\eta_0^{-1}\in \kappa(E')$. And
$(x,y)\eta_0^{-1}\kappa^{-1}\pi'_0=y.$ Therefore $\pi'_0$ is
surjective and $\pi$ is a deflation by Lemma 6.4.3. Let
$\varphi'=\kappa\eta{1\choose 0}$. Since $\kappa\eta{\pi\choose
-\varphi}=0$, we have $\pi'\varphi=\varphi'\pi$. The proof of E2
is finished.
 \hfill$\square$

We learnt recently that \cite{BBP} proved the theorem similarly.

{\bf Corollary 6.4.1} Let $\frak A=(\Gamma, \Omega)$ be a layered
bocs and $(e): M\stackrel{\iota}\longrightarrow E\stackrel{\pi}
\longrightarrow L$ be a conflation of $R(\frak A)$. Then $(e) $ is
equivalent to an exact sequence
$$ (e')\quad 0\longrightarrow M\stackrel{\iota'}\longrightarrow E'\stackrel{\pi'}
\longrightarrow L \longrightarrow 0 $$
 in $R(\mathfrak{A}_0).$ Moreover we can choose some suitable basis such
 that
 $$ \iota'_{\L}=(0,I),\quad
 \pi'_{\L}=\left(\begin{array}{ll}I\\0\end{array}\right)$$
 for any $\L\in \ind\Gamma'$.

 {\bf Proof.} By Proposition 6.3.2.\hfill$\Box$

   We apply the  notions
given in 6.3 and 6.4 to a freely parameterized bimodule problem to
end this subsection.

{\bf Remark.} Let $(\K,\M,H)$ be a freely parameterized bimodule
problem.

(1) Let $\K_0=\{ S\in\K\mid S\mbox{ are diagonal matrices}\} $, then $(\K_0,\M)$ is said to be the {\it
principal bimodule problem} of $(\K,\M, H) $.

(2) Let $M,L\in Mat(\K,\M)$ with size vectors $\m,\l$
respectively. Let $S: M\rightarrow L$ be a morphism, and $S_0$ the
partitioned diagonal part of $S$ (see 2.2). Then $S_0$ is
injective (resp. surjective) if and only if $\mbox{rank} (S_0)=m$
(resp. $l$), if and only if $\mbox{rank} (S_{\I})=m_{\I}$ (resp.
$l_{\I}$) for any $\I\in T/\!\!\sim$.

(3) Let $\cal E$ be a class consisting of composable pairs
$M\stackrel{S}{\longrightarrow} E\stackrel{R}{\longrightarrow} L $
in $Mat(\K,\M)$ such that $SR=0$, and $0 \rightarrow
M\stackrel{S_0}{\longrightarrow} E\stackrel{R_0}{\longrightarrow}
L \rightarrow 0$ are exact as vector spaces. Then $\cal E$ is an
exact structure and $(Mat(\K,\M), \cal E)$ is an exact category.

\subsection{Homogeneous conflations}

\kg Let $(\K,\M,H)$ be a freely parameterized bimodule problem corresponding to a layered bocs $\mathfrak{A}$.
We have seen in 6.4 that $Mat(\K,\M)$ has an exact structure.

{\bf Definition 6.5.1} Let $(\K,\M,H)$ be a freely parameterized
bimodule problem.

(1) An indecomposable object $M$ in $Mat(\K,\M)$ is said to be
{\it homogeneous} or {\it $\tau$-invariant} if there exists an
almost split conflation $M\stackrel{\iota}{\longrightarrow} E
\stackrel{\pi}{\longrightarrow}M$.

(2) The category $Mat(\K,\M)$ is said to be {\it homogeneous} if
for any positive integer $n$, almost  all (except finitely many)
iso-classes of indecomposable objects with size at most $n$ are
homogeneous.

(3) The category $Mat(\K,\M)$ is said to be {\it strongly
homogeneous} if there exists neither projectives nor injectives,
and  all indecomposables are homogeneous.

Sometimes we say that a bimodule problem $(\K,\M, H)$ (or a bocs
$\mathfrak{A}$) is homogeneous, instead of saying $Mat(\K,\M)$ (or
$R(\mathfrak{A})$) for short.

{\bf Remark.} We {\bf do NOT require} that $Mat(\K,\M)$ has almost
split conflations globally in the sense of \cite{DRSS} when we
talk about homogeneous property. We only concentrate on the
individual almost split conflations with the homogeneous starting
and end terms. Furthermore we will see in subsection 7.2 that the
homogeneous property of a bimodule problem $(\K,\M,H)$ can be
inherited to any induced bimodule problem $(\K',\M',H')$. Thus we
are able to use reduction sequences and restrict our attention to
some minimally wild cases.

{\bf Examples.}  (1)\  Let $$
\K =\left\{\left(\begin{array}{cc} s & 0\\
0 & s\end{array}\right)\mid \forall s\in\rK\right\}, \M=\left\{
\left(\begin{array}{cc} a & 0\\
 0 &b\end{array}\right)\ | \forall\  a, b\in\rK\right\}, H=0
 , T=\{1,2\}, 1\sim 2.$$ Then $Mat(\K,\M)$
  is  equivalent to  $\rK\langle x,y\rangle$-mod. It
has been shown by  Vossieck in \cite{ZL} that $Mat(\K,\M)$ has no any almost split conflations.

(2) There are several sufficient conditions on a layered bocs to
have almost split conflations given in \cite{CB1, BK,BB} etc. In
particular, the Drozd bocs of a finite-dimensional algebra
$\Lambda$, or equivalently $P_1(\Lambda)$, has almost split
conflations.

(3) If an algebra $\Lambda$ is of tame representation type, then
$P_1(\Lambda)$ and also $\Lambda$-mod are homogeneous \cite{CB1}.

(4) Define an index set
  $T=\{1,2\},$\;  with $1\sim
2$. Let
$$\K=\left\{\left(\begin{array}{cc} s_1 & s_2 \\ 0 & s_1
\end{array}\right) \mid s_1,s_2\in\rK\right\},\
\M=\left\{\left(\begin{array}{cc}a & b \\ 0 & a
\end{array}\right) \mid a,b\in\rK\right\}.$$
 and $ H=0.$
 Then $Mat(\K,\M)$ has  almost split conflations, and it is
strongly homogeneous \cite{BCLZ}.

\subsection{Minimal bocses}

\kg Corollary 4.6.1 has shown the complete set of iso-classes of
indecomposables of a minimal bocs, we will show in the sebsection
the projectives, injectives and almost split conflations of a
minimal bocs.

\kg {\bf Proposition 6.6.1} Let $\mathfrak{A}=(\Gamma,\Omega)$ be
a minimal bocs with a layer $L=(\Gamma';\omega; v_1,\cdots,v_m)$.

(1) For any trivial $\J\in \ind \Gamma^{\prime}$, $O(\J)$ is
projective and injective in $R(\mathfrak{A})$.

(2) For any non-trivial $\I\in \ind \Gamma^{\prime}$ with
$\Gamma'(\I,\I)=k[\lambda, g_{\I}(\lambda)^{-1}]$, and for any
$\lambda^0\in k$, $g_{\I}(\lambda^0)\ne 0$.
$$
\begin{CD}
O(\I,1,\lambda^0)@>\iota_{1}>> O(\I,2,\lambda^0) @>\pi_{1}>>
O(\I,1,\lambda^0)
\end{CD}
$$
and
$$
\begin{CD}
O(\I,d,\lambda^0)@>(\pi_{d-1}, \iota_{d})>>
O(\I,d-1,\lambda^0)\oplus O(\I,d+1,\lambda^0)
@>{\iota_{d-1}\choose \pi_{d}}>> O(\I,d,\lambda^0)
\end{CD}
$$
are almost split conflations in $R(\mathfrak{A})$, where $d\in
\mathbb{N}, d\geq 2$.

{\bf Proof.} (1) We claim first that for any indecomposable
$M,N\in R(\mathfrak{A})$ with $M(\I)\ne 0, N(\J)\ne 0$ (see 4.6),
if $\I\ne\J$, then $\Hom_{\Gamma}(M,N)=0$, since
$$\Gamma=\Gamma^{\prime}=\rK[\lambda_1,g_1(\lambda_1)^{-1}]\times\cdots\times
\rK[\lambda_i,g_i(\lambda_i)^{-1}] \times
\rK_1\times\cdots\times\rK_j$$
 (see 4.2). If $\pi: E\rightarrow
O(\J)$ is a deflation in $R(\mathfrak{A})$, then $\pi_0$ is
surjective by the definition of ${\cal E}$ in 6.4. Thus item (2)
of Proposition 6.3.1 shows a commutative diagram with
$E'\stackrel{\pi'}{\longrightarrow} O(\J)\longrightarrow 0$ being
exact in $\Gamma$-mod. Therefore $E'=M\oplus O(\J)^e$ with $e\geq
1$ and $\pi'={0\choose \varphi}$, where $M$ does not contain any
direct summand $O(\J)$, and $\varphi$ is a surjective linear map
between vector spaces. Consequently $\varphi$ splits, so do $\pi'$
and $\pi$. Hence $O(\J)$ is projective.

The injectivity of $O(\J)$ can be proved dually.

(2) See Lemma 6.6 of  \cite{CB1}. \hfill$\square$

{\bf Corollary 6.6.1} Let $\mathfrak{A}$ be a minimal bocs, then
$R(\mathfrak{A})$
 has almost split conflations in the sense of
\cite{DRSS}. And $R(\mathfrak{A})$ is homogeneous.

 {\bf Proof.} There is a left almost split morphism:\
$O(\J)\to 0$, as well as a right almost split morphism: \ $0 \to
O(\J)$ for any trivial vertex $\J$. Therefore $R(\mathfrak{A})$
has almost split morphisms. The conclusion follows from
Proposition 6.6.1. \hfill$\square$

Next we show by an example that the exact structure $\mathcal{E}$ on $R(\mathfrak{A})$ defined in 6.4 collects
only part of exact pairs of composable morphisms.

{\bf Example.} Let $\mathfrak{A}=(\Gamma,\Omega)$ be a minimal bocs with three trivial vertices $\P,\Q,\L$, and
two dotted arrows $u:\P\dashrightarrow \Q$, $v: \Q\dashrightarrow\L$, such that $\delta(u)=0$, $\delta(v)=0$.
Define a pair of composable morphisms
$$
(e)  M\stackrel{\iota}{\longrightarrow} N \stackrel{\pi}{\longrightarrow} L
$$
where $M_{\P}=k$, $N_{\Q}=k$, $L_{\L}=k$, and $M,N,L$ are all indecomposable; $\iota(u)=id$, $\pi(v)=id$. Then
$\iota$ is a kernel of $\pi$, and $\pi$ is a cokernel of $\iota$. Thus $(e)$ is an exact pair, which does not
belong to $\mathcal{E}$. Moreover it is not difficult to see that $(e)$ is even an almost split pair in
$R(\mathfrak{A})$.

\newpage

\bcen{\section{ Almost split conflations in reductions}}\ecen

\subsection{Projective and injective objects}

\kg This subsection is devoted to showing the finiteness of projectives (resp. injectives) in a layered bocs and
also behavior of the projectives (resp. injectives) in reductions.

{\bf Lemma 7.1.1 \cite{B1}} Let $\mathfrak{A}=(\Gamma,\Omega)$ be a bocs with a layer $L=(\Gamma'; \omega; a_1,
\cdots, a_n; v_1, \cdots,$ $ v_m )$,  and $M\in R(\mathfrak{A})$  having  dimension vector $\m$. Suppose that
$\P\in $ $\ind\Gamma'$ is non-trivial. If $m_{\P}\neq 0$, then $M$ is neither projective nor injective.

{\bf Proof.} Let  $\lambda$ be a free parameter  attached to $\P$ and $M(\lambda)\simeq J_d(\lambda^0)\oplus W$,
where $d>0$, and $W$ is a Weyr matrix.  Let us define a  dimension vector $\n$ such that $n_{\P}=m_{\P}+1$,
$n_{\I}=m_{\I}, \, \forall\ \I\in $ $\ind\Gamma'$,  $\I \neq \P$. Let
$$
\varphi_{\P}=\left(\begin{array}{cc|c}I_d& 0_{d\times 1}& \\
\hline & & I \end{array}\right)_{m_{\P}\times n_{\P}}, \qquad
\psi_{\P}=\left(\begin{array}{c|c} I_d\\ 0_{1\times d}  & \\
\hline & I
\end{array}\right)_{n_{\P}\times m_{\P}}.
$$
We now define an object $N$ in $R(\mathfrak{A})$, such that $N(\lambda)=J_{d+1}(\lambda^0)\oplus W$, and
$N(\mu)=M(\mu)$ for any parameter $\mu\neq \lambda$,
$$
N(a_i)=\left\{\begin{array}{ll} M(a_i), & \mbox{when } \I\ne \P, \J\ne\P; \\ M(a_i)\varphi_{\P}, & \mbox{when }
\I\ne \P,
\J=\P; \\
\psi_{\P}M(a_i), & \mbox{when } \I= \P, \J\ne\P; \\
\psi_{\P}M(a_i)\varphi_{\P}, & \mbox{when } \I= \P=\J
\end{array}\right.
$$
 for any solid arrow $a_i: \I\rightarrow \J$.  Then we
 obtain a morphism   $\psi: N\rightarrow M$, such that
$\psi_{\I}=I_{m_{\I}}$,  for any $ \I\ne \P$, and $\psi(v_j)=0$, $j=1,\cdots, m$. Then
$N(\lambda)\psi_{\P}=\psi_{\P}M(\lambda)$ and $N(a_i)\psi_{\J}=\psi_{\I}M(a_i)$, i.e. $\psi$ is a morphism in
$R(\mathfrak{A}_0)$. It is obvious that  $\psi$ is a deflation but not split.
 Therefore $M$ is  not projective.

 The proof of $M$ being
non-injective is dual. \hfill$\square$

{\bf Proposition 7.1.1} (1) Let $(\K,\M,H)$  be a freely parameterized bimodule problem and $(\K',\M',H')$ be an
induced bimodule problem
 given by one of the 3 reductions of 3.1. If
 $M'\in Mat(\K',\M')$ is
non-projective (resp. non-injective ) then $\vartheta(M')$ is non-projective (resp. non-injective) in
$Mat(\K,\M)$.

(2) Let $(\K,\M,H=0)$ be a bimodule problem with a reduction sequences $(*)$  of freely parameterized triples.
If $M^s\in Mat(\K^s, \M^s)$ is non-projective (resp. non-injective), then $\vartheta_{r,s}(M^s)$ is
non-projective (resp. non-injective ) in $Mat(\K^r,\M^r)$ for $r=s-1, \cdots, 1, 0$.

\def\xxx{\end{document}}
{\bf Proof.} (1) Since $M'\in Mat(\K',\M')$ is non-projective, there must exists a non-split deflation $\pi':
N'\rightarrow M'$ in $Mat(\K', \M')$, then $\pi=\pi': \vartheta(N')\rightarrow\vartheta(M')$ is also a non-split
deflation since the surjectivity of $\pi'_0$ implies that of $\pi_0$. In fact, $\pi'_0=\pi_0$ in case of
deletion and regularization, and $\pi'_0$ is the diagonal part of $\pi_0$ according to the partition
$(T',\sim')$ in case edge and loop reductions. Therefore $M\in Mat(\K,\M)$ is not projective.

The proof for non-injective case is dual.

(2) By  (1) and  induction on $r=s, s-1, \cdots, 1,0$. \hfill $\square$

{\bf Corollary 7.1.1} Let $(\K,\M,H)$ be a freely parameterized bimodule problem. Then for each positive integer
$n$, there are only finitely many iso-classes  of projectives and injectives of size at most $n$ in
$Mat(\K,\M)$.

{\bf Proof.} It
 is clear that there are only finitely many choices of  size
 vectors $\m$ over $(T,\sim)$ with  $m\leq n$. Let $M\in Mat(\K,\M)$
 be a projective matrix of size vector
$\m$. Then $m_{\P}=0$, for any non-trivial $ \P\in T/\!\!\sim$
 by Lemma 7.1.1. Therefore $H_{\m}$ is a constant matrix.
Now we start a sequence of reductions according to Theorem 3.4.1 for $M$. We conclude that there does not  exist
any loop reduction in the sequence.  Otherwise, if the $(r+1)$-th reduction were a loop reduction, $M^r$ would
be non-projective in $Mat(\K^r, \M^r)$ still by Lemma 7.1.1. Therefore $M$ would be non-projective in
$Mat(\K,\M)$ by  proposition 7.1.1, which contradicts to our assumption. It follows obviously  that there are
only finitely many choices of such a sequence consisting of regularizations and edge reductions with size vector
$\m$. Summarizing the discussion, there are only finitely many iso-classes of projectives of size at most $n$ in
$Mat(\K,\M)$.

The proof for injectives is dual. \hfill $\square$

\subsection{Homogeneous property in reductions}

\kg We will give the key idea  for proving the main theorem  in this subsection, which was originally presented
by R. Bautista.

{\bf Lemma 7.2.1\cite{B1}.}  Let $(\K,\M,H)$ be a freely parameterized bimodule problem, and $(\K',\M',H')$ be
an induced bimodule problem
 given by one of the $4$ reductions in 3.1 or 3.2.

(1) Suppose that $M',L'\in Mat(\K',\M')$ have size vectors
$\m',\l'$, and $\vartheta(M'), \vartheta(L')\in Mat(\K,\M)$ have
size vectors $\m,\l$ respectively. If $\m=\l$, and $\m'\leqslant
\l'$, then $\m'=\l'$.

(2) If $\iota': M'\rightarrow E'$ is a morphism in $Mat(\K', \M' )$ with $\vartheta(\iota'):
\vartheta(M')\rightarrow \vartheta(E')$ being a left minimal almost split morphism of $Mat(\K,\M)$, then so is
$\iota'$. Dually if $\pi': E'\rightarrow L'$ is a morphism in $Mat(\K', \M')$ with $\vartheta(\pi'):
\vartheta(E')\rightarrow \vartheta(L')$ being a right minimal almost split morphism of $Mat(\K,\M)$,  then so is
$\pi'$ in $Mat(\K', \M')$.

(3) Let $(e') \ M'\stackrel{\iota'}{\longrightarrow} E'\stackrel{\pi'}{\longrightarrow} L'$ be a composable pair
of $Mat(\K',\M')$ with $M'$ non-injective and $L'$ non-projective. If $\vartheta(e')$ is an almost split
conflation of $Mat(\K,\M)$, then so is $(e')$ in $Mat(\K', \M')$.

{\bf Proof.} (1)   It is obvious in case of deletion. The other cases given by reductions  in 3.1 follow from
Formula (11) of 3.2.

(2) The functor $\vartheta: Mat(\K',\M')\rightarrow Mat(\K,\M)$ is fully faithful.

(3) Since $L'$ is non-projective, there exist some $N'$ and non-split morphism $\varphi': N'\rightarrow L'$ with
$\varphi'_0$ being surjective. But $\pi'$ is a right almost split morphism by (2), we have a morphism $\psi':
N'\rightarrow E'$ such that $\psi'\pi'=\varphi'$. Consequently $\psi'_0\pi'_0=\varphi'_0$ and $\pi'_0$ is
surjective. Similarly $\iota'_0$ is injective. Since $\vartheta(\iota')\vartheta(\pi')=0$, we have
$\iota'\pi'=0$ and
 $\iota'_0\pi'_0=0$. Therefore $\e'\geqslant \m'+\l' $. But $\e=\m+\l$, so
$\e'=\m'+\l'$ by (1). Hence $M'\stackrel{\iota'_0}{\longrightarrow} E'\stackrel{\pi'_0}{\longrightarrow} L'$ is
exact as a sequence of vector spaces, i.e. $(\iota',\pi')$ is a conflation. And $(e')$ is an almost split
conflation of $Mat(\K', \M')$ by (2). \hfill$\square$

{\bf Theorem 7.2.1} Let $\mathfrak{A}=(\Gamma, \Omega)$ be a layered bocs, and $\mathfrak{A'}=(\Gamma',
\Omega')$ be induced by one of the 4 reductions of 3.1 or 3.2 with a reduction functor $\vartheta:
R(\mathfrak{A'})\rightarrow R(\mathfrak{A})$. Suppose that an indecomposable $M'\in R(\mathfrak{A'})$  is
neither projective nor injective. If
 $$(e) \qquad  \vartheta(M')
\stackrel{\iota}{\longrightarrow} E\stackrel{\pi}{\longrightarrow} \vartheta(M')$$ is an almost split conflation
of $R(\mathfrak{A})$, then

(1) there exists some $E'\in R(\mathfrak{A'})$ such that $\vartheta(E')\simeq E$;

(2) there exists an almost split conflation

 $$(e') \qquad
M'\stackrel{\iota'}{\longrightarrow} E'\stackrel{\pi'}{\longrightarrow}M'$$ in $R(\mathfrak{A'})$ such that
$\vartheta(e')$ is equivalent to $(e)$.

{\bf Proof.} We assume that $\mathfrak{A}$ is the corresponding bocs of  bimodule problem $(\K,\M,H)$ and use
the notions of both structures freely. Denote $\vartheta(M')$ by $M$ for short, then we first set up two claims.

Claim 1. In  case of deletion, both conclusions (1) and (2) hold.

In fact, since
  $$(e_0)\quad 0\longrightarrow\vartheta(M')
\stackrel{\iota_0}{\longrightarrow} E
\stackrel{\pi_0}{\longrightarrow} \vartheta(M')\longrightarrow 0$$
is exact as a sequence of vector spaces, and $M=\vartheta(M')=M'$,
we conclude that $E_{\I}=0$, $\forall\ \I\notin T'/\!\!\sim'$.
Hence $E$ belongs to $R(\mathfrak{A'})$. Therefore $E'=E$,
$(e')=(e)$.

Claim 2. In other cases,  (2) is true if (1) holds. Indeed, suppose the commutative diagram below comes from
(1):
\begin{center}
\unitlength=1mm
\begin{picture}(60,40)
\put(0,20){\makebox{$\vartheta(M')$}} \put(11,19){\vector(2,-1){20}} \put(11,22){\vector(2,1){20}}
\put(33.50,30.50){\makebox{$E$}} \put(32,7){\makebox{$\vartheta(E')$}}

 \put(41,31.50){\vector(2,-1){20}}
\put(42,8){\vector(2,1){20}} \put(62,19.50){\makebox{$\vartheta(M')$}} \put(34.50,29){\vector(0,-1){18}}

\put(35,20){\makebox{$\eta$}} \put(20,27){\makebox{$\iota$}} \put(18,11){\makebox{$\iota\eta$}}
\put(50,27){\makebox{$\pi$}} \put(49,9){\makebox{$\eta^{-1}\pi$}}

\end{picture}
\end{center}
Since $\vartheta$ is fully faithful, there exists some $\iota':
M'\rightarrow E'$ with $\vartheta({\iota'})=\iota\eta$ and  $\pi':
E'\rightarrow M'$ with $\vartheta({\pi'})=\eta^{-1}\pi$. We obtain
a pair of morphisms $(e')$ in $R(\mathfrak{A'})$ satisfying the
hypothesis of item (3) of  Lemma 7.2.1. Therefore, $(e')$ is an
almost split conflation in $R(\mathfrak{A'})$ with $\vartheta(e')$
equivalent to $(e)$.

Then it suffices to prove  assertion (1) according to the $3$ reductions in 3.1 respectively.

 Regularization. (1) is valid because of $\vartheta$ being an
equivalence.

Edge reduction. We consider first an edge reduction given by $\left(\begin{array}{cc} 0 & 1\\ 0 &
0\end{array}\right)$. Let $\mathfrak{A''}$ be the induced bocs  and $\vartheta_1: R(\mathfrak{A}'')\rightarrow
R(\mathfrak{A}) $ be the  equivalent functor. Thus $\mathfrak{A}'$ is induced from $\mathfrak{A}''$ by a
deletion of $(T''/\sim'')\setminus (T'/\sim')$. And the following diagram commutes:
\begin{equation} \unitlength=1mm
\begin{array}{c}
\begin{picture}(40,30)
\put(13,5){$R(\mathfrak{A})$} \put(-1,20){$R(\mathfrak{A}')$} \put(26,20){$R(\mathfrak{A}'')$}
\put(10,21){\vector(1,0){15}} \put(5,19){\vector(1,-1){10}} \put(30,19){\vector(-1,-1){10}}
\put(15,22){$\vartheta_2$} \put(29,12){$\vartheta_1$} \put(5,12){$\vartheta$}
\end{picture}
\end{array}
\end{equation}
 Since $M'\in
R(\mathfrak{A}')$ is neither projective nor injective. $\vartheta_2(M')\in R(\mathfrak{A}'')$ is so by item (1)
of Proposition 7.1.1. And since $\vartheta_1$ is an equivalence, there exists some $E''\in R(\mathfrak{A}'')$
with $E\simeq \vartheta_1(E'')$. Therefore we have an almost split conflation
$$(e''): \vartheta_2(M') \stackrel{\iota''}{\longrightarrow}
E''\stackrel{\pi''}{\longrightarrow} \vartheta_2(M') $$
 in
$R(\mathfrak{A}'')$ with $\vartheta_1(e'')$ equivalent to $(e)$ by Claim 2. Furthermore, $M'=\vartheta_2(M')$
and $E'=E''$, $(e')=(e'')$ follow from deletion given by Claim 1. Thus, $\vartheta(E')\simeq E$ and
$\vartheta(e')$ is equivalent to $(e)$ as desired.

Loop reduction. Suppose the reduction is given by $W\oplus(\lambda)$ (or just $W$),  and the domain of $\lambda$
equals $\rK\setminus\{\mbox{the roots of } g(\lambda)\}$. If $E_{pq}\simeq W_1\oplus W_2$, where the eigenvalues
of $W_2$ belong to the domain of $\lambda$, but those of $W_1$ do not. Then we make a new loop reduction given
by
 $$ W\oplus
W_1\oplus (\lambda) \ ({\rm or}\ W\oplus W_1).$$  Consequently we obtain a new induced bocs $\mathfrak{A}''$
from $\mathfrak{A}$ by loop reduction, and $\mathfrak{A}'$ is induced from $\mathfrak{A}''$ by deletion of
$(T''/\sim'')\setminus (T'/\sim')$. Then there exist some $E''\in R(\mathfrak{A'})$ such that $E\cong
\vartheta_1(E'')$ and an almost split conflation $(e'')$ in $R(\mathfrak{A''})$ such that $\vartheta_1(e'')$ is
equivalent to $(e)$. The diagram (45) shows a similar proof for the case of loop reduction. \hfill $\square $

{\bf Corollary 7.2.1}  (1) Let $(\K,\M,H)$ be a bimodule problem
and $(\K',\M',H')$ induced from $(\K,\M,H)$ by one of the four
reductions given by 3.1. If $(\K,\M,H)$ is homogeneous, then so is
$(\K',\M',H')$.

(2) For a given bimodule problem $(\K,\M,H)$, if there exists a
reduction sequence  $(*)$  of freely parameterized triples such
that $Mat(\K^s, \M^s)$ is non-homogeneous, then $Mat(\K,\M)$ must
be non-homogeneous.

{\bf Proof.} (1) Given any fixed positive integer  $n$, we define a set of iso-classes of objects
\begin{align*}
\mathscr{S}=\{M_{\alpha}\in Mat(\K',\M')&\Big|\  M_{\alpha}\
\mbox {is indecomposable of size at most $n$,}\\
&\mbox{ and is neither projective nor injective}\}.
 \end{align*}
 Because of Corollary 7.1.1, $\mathscr{S}$ is a cofinite subset. If
$\mathscr{S}$ is a finite set, we are done. If $\mathscr{S}$ is an
infinite set, then $\{\vartheta(M_{\alpha})\mid M_{\alpha} \in
\mathscr{S}\}$ is an infinite set of objects in $Mat(\K,\M)$,
where $\vartheta(M_{\alpha})$ is neither projective nor injective
by Proposition 7.1.1, and the size of $\vartheta(M_{\alpha})$ is
bounded by $n$. Since  $Mat(\K,\M)$ is homogeneous, there exists a
cofinite subset ${\mathscr{S}}_0$ of $\mathscr{S}$, such that
$\forall \, M_{\alpha}\in {\mathscr{S}}_0$, there are almost split
conflations of $Mat(\K,\M)$:
$$ (e_{\alpha})\qquad \begin{CD}
\theta(M_{\alpha})@>{\iota}>> E_{\alpha}@>{\pi}>> \theta(M_{\alpha}) \end{CD}$$ Thus there are almost split
conflations of $R(\mathfrak{A'})$
$$(e'_{\alpha}) \qquad \begin{CD}
M_{\alpha}@>{\iota'}>>E'_{\alpha}@>{\pi'}>> M_{\alpha} \end{CD}$$
with $\vartheta(e'_{\alpha}) $ being equivalent to $(e_{\alpha}) $
by Theorem 7.2.1. Consequently $Mat(\K',\M')$ is homogeneous.

(2) If $Mat(\K,\M)$ were homogeneous, then  $Mat(\K^r, \M^r)$ would be homogeneous for $r=1, 2, \cdots, s$ by
(1) and induction on $r$, which is a contradiction to the hypothesis. \hfill $\square $

\subsection{Critically non-homogeneous bocses  $\mathfrak{B_1}$ and $ \mathfrak{B_2}$}

\kg {\bf Lemma 7.3.1} Let $f(\nu, \kappa)$ and $g(\nu, \kappa)$ be polynomials in two indeterminates $\nu,
\kappa$. If $f$ and $g$ are coprime, then there exists an infinite list of pairs $(\nu^0, \kappa^0)$, such that
$f(\nu^0, \kappa^0)=0$ but $g(\nu^0, \kappa^0)\neq 0$. Thus $\nu^0$ or $\kappa^0$ must range over a cofinite
subset  $D$ of $\rK$.

{\bf Proof.} By means of  Bezout's theorem. \hfill $\square$

{\bf Corollary 7.3.1} Let $g_{\P}(\nu),\sigma_{\P}(\nu)\in k[\nu]$ with $g_{\P}(\nu)\mid \sigma_{\P}(\nu)$, and
$g_{\Q}(\kappa),\sigma_{\Q}(\kappa)\in k[\kappa]$ with $g_{\Q}(\kappa)\mid \sigma_{\Q}(\kappa)$. If
$$f(\nu,\kappa)\in k[\nu,\kappa,g_{\P}(\nu)^{-1}, g_{\Q}(\kappa)^{-1}]$$ is not invertible, then there exists an
infinite list of pairs $(\nu^0, \kappa^0)$ with $f(\nu^0,\kappa^0)=0$, such that either

(1)  $\sigma_{\P}(\nu^0)g_{\Q}(\kappa^0)\neq 0$, and $\nu^0$ ranges over a cofinite subset $D_{\P}$ of $k$; or

(2)  $g_{\P}(\nu^0)\sigma_{\Q}(\kappa^0)\neq 0$, and $\kappa^0$ ranges over a cofinite subset $D_{\Q}$ of $k$.

{\bf Proof.} If $f(\nu, \kappa)\ne 0$, let $$f(\nu,\kappa)=
\overline{f}(\nu,\kappa)\overline{g}_{\P}(\nu)\overline{g}_{\Q}(\kappa),$$ such that $\overline{f}(\nu,\kappa)$
is a polynomial and $(\overline{f}(\nu,\kappa), g_{\P}(\nu)g_{\Q}(\kappa))=1$; and $\overline{g}_{\P}(\nu)$ is a
product of some positive or negative powers of the factors of $g_{\P}(\nu)$, $\overline{g}_{\Q}(\kappa)$ is that
of $g_{\Q}(\kappa)$. Suppose that $$\overline{f}(\nu,\kappa)= f_0(\nu,\kappa)f_{\P}(\nu)f_{\Q}(\kappa),$$ such
that $f_0(\nu,\kappa)$ has no
 non-constant factors taken from $k[\nu]$
or $k[\kappa]$.

(i) If $f_0(\nu,\kappa)$ is not a constant, then $(f_0(\nu,\kappa), \sigma_{\P}(\nu)\sigma_{\Q}(\kappa))=1$,
 then there exists an infinite list of pairs $(\nu^0,
\kappa^0)$, such that $f_0(\nu^0,\kappa^0)=0$, thus $f(\nu^0,\kappa^0)=0$, but
$\sigma_{\P}(\nu^0)\sigma_{\Q}(\kappa^0)\ne 0$. When $\nu^0$ ranges over a cofinite subset $D_{\P}$ of $k$, we
have (1); when $\kappa^0$ ranges over a cofinite subset $D_{\Q}$ of $k$, we have (2).

(ii) If $f_0(\nu,\kappa)$ is a constant, but $f_{\Q}(\kappa)$ is not,  then $(f_{\Q}(\kappa),
\sigma_{\P}(\nu)g_{\Q}(\kappa))=1$. $\forall \nu^0\in k$ with $\sigma_{\P}(\nu^0)\ne 0$, there exists some
$\kappa^0$ with $f_{\Q}(\kappa^0)=0$ but $g_{\Q}(\kappa^0)\neq 0$, we obtain item (1).

(iii) If $f_0(\nu,\kappa)$ and $f_{\Q}(\kappa)$ are both constants, then $f_{\P}(\nu)$ must be not,  thus
$(f_{\P}(\nu), g_{\P}(\nu)\cdot $ $\sigma_{\Q}(\kappa))=1$. $\forall \kappa^0\in k$ with
$\sigma_{\Q}(\kappa^0)\ne 0$, there exists some $\nu^0$ with $f_{\P}(\nu^0)=0$ but $g_{\P}(\nu^0)\neq 0$, we
obtain item (2).

Finally, if $f(\nu,\kappa)=0$, it is clear that we have (1) and (2). \hfill$\square$

{\bf Proposition 7.3.1} The bocs $\mathfrak{B_1}$  given in MW1 of Theorem 5.6.1 is not homogeneous.

{\bf Proof.} Consider the local bocses $(\mathfrak{B}_1)_{\P}$ at $\P$ and $(\mathfrak{B}_1)_{\Q}$ at $\Q$,
suppose the triangular formula (40) of 5.4 gives the polynomials $\sigma_{\P}(\nu)$ and $\sigma_{\Q}(\kappa)$ in
(41) respectively. Assume that we have item (1) of Corollary 7.3.1, i.e. there is an infinite list
$\{(\nu^0,\kappa^0)\}$, such that  $f(\nu^0,\kappa^0)=0$, $\sigma_{\P}(\nu^0)g_{\Q}(\kappa^0)\neq 0$ and $\nu^0$
ranges over a cofinite subset $D=D_{\P}$ of $k$. If $R(\mathfrak{B}_1)$ is homogeneous, we present an infinite
list of one dimensional objects $S_{\nu^0}$, $\forall \, \nu^0\in D$, with $$(S_{\nu^0})_{\P}=\rK, \,
(S_{\nu^0})_{\Q}=0, \, S_{\nu^0}(\nu)=(\nu^0), \, S_{\nu^0}(a_l)=0\ {\rm or} \ \emptyset \ {\rm for}\ l=1,
\cdots, i,$$ when we calculate $S_{\nu^0}$ according to Theorem 3.4.1.  If $R(\mathfrak{B_1})$ is homogeneous,
there must be a cofinite subset $D_0\subset D$, with $S_{\nu^0} $ being homogeneous for any $\nu^0\in D_0$. We
fix such a $\nu^0$, and denote $S_{\nu^0}$ by $S$
 for simplicity. Let $(e): \  S
\stackrel{\iota}{\longrightarrow} E \stackrel{\pi}{\longrightarrow} S $ be the almost split conflation
 starting and ending at $S$ in $R(\mathfrak{B_1})$, then $(e) $ is also an almost split
conflation of $R((\mathfrak{B}_1)_{\P})$  according to Theorem 7.2.1. Thus
$$E_{\P}=\rK^2,\ \  E_{\Q}=0,\ \  E(\nu)=\left(\begin{array}{cc} \nu^0&1\\
0&\nu^0
 \end{array}\right)$$ by Proposition 6.6.1  and $E(a_l)=\emptyset$ for all the solid
 loops $a_l$ at $\P$. We define an object $L\in R(\mathfrak{B}_1)$ with
 $$L_{\P}=\rK,\ L_{\Q}=\rK,\  L(\nu)=(\nu^0), \ L(\kappa)=(\kappa^0), \ L(a_1)=(1), \ L(a_l)=0,\
 \forall\,  l\geqslant2,$$ and we also define a morphism $\varphi: L\rightarrow S$
 given by $$\varphi_{\P}=(1),\ \ \varphi_{\Q}=(0),\ \ \varphi(v)=0$$ for any dotted arrow $v$.

\begin{center} \unitlength 0.7mm
\begin{picture}(100, 55)

\put(-5,25){\oval(5,5)[t]} \put(-5,25){\oval(5,5)[bl]} \put(-6.00,22.00){\vector(3,1){3}}

\put(-1.00,23.00){\makebox{$\rK^2$}} \put(0,22.00){\vector(0,-1){15}} \put(-0.50,4.0){\makebox{$0$}}

\put(93,25){\oval(5,5)[t]} \put(93,25){\oval(5,5)[br]} \put(94.00,22.00){\vector(-3,1){3}}
\put(85.00,23.00){\makebox{$\rK$}} \put(86,22.00){\vector(0,-1){15}} \put(85.50,4.0){\makebox{$0$}}

\put(2, 24.00){\vector(1,0){81}} \put(2, 5.00){\vector(1,0){81}}

\put(39,52){\oval(5,5)[t]} \put(39,52){\oval(5,5)[bl]} \put(38.00,49.00){\vector(3,1){3}}
\put(44.00,49.00){\makebox{$\rK$}} \put(45,48.00){\vector(0,-1){15}} \put(44.00,29.0){\makebox{$\rK$}}

\put(39,32){\oval(5,5)[t]} \put(39,32){\oval(5,5)[bl]} \put(38.00,29.00){\vector(3,1){3}}

 \put(43,47.00){\vector(-2,-1){35}} \put(48,47.00){\vector(2,-1){35}}

\put(43, 27.00){\vector(-2,-1){36}} \put(48,27.00){\vector(2,-1){36}}

\put(-18.50,28.0){\makebox{$\s J_2(\nu^0)$}} \put(30.50,55.0){\makebox{$\s (\nu^0)$}}
\put(33.50,35.4){\makebox{$\s (\kappa^0)$}} \put(95.50,27.0){\makebox{$\s (\nu^0)$}}

\put(18.50,41.0){\makebox{$\varphi'_{\P}$}} \put(67.50,40.0){\makebox{$\varphi_{\P}=id$}}
\put(46,38.0){\makebox{$id$}}

\put(-10,13.0){\makebox{$E:$}} \put(90.50,13.0){\makebox{$:S$}} \put(50,50.0){\makebox{$:L$}}
\put(43.50,19.00){\makebox{$\pi_{\P}$}}
\end{picture}
\end{center}
 We first claim that $\varphi$ is not a  split epimorphism. In fact
 if there were some morphism $\psi: S\rightarrow L$ with
 $\psi\varphi=id$, then we would have $$\varphi_{\P}=(1),\
 \psi_{\Q}=(0);\ \
 S(a_1)\psi_{\Q}-\psi_{\P}L(a_1)=\psi(\delta(a_1))=f(\nu^0, \kappa^0)\psi(w_1),$$ i.e.
 $-1=0$, which leads to a contradiction. Thus there exists some $\varphi': L\rightarrow E$
  such that  $\varphi'\pi=\varphi$, since $(e)$ is an almost
 split conflation. Then $\varphi'_{\P}\pi_{P}=\varphi_{\P}$, i.e.
 $(a, \  b){1\choose 0}=1$, $\varphi'_{\P}=(1,\, b)$. On the other hand,
$\varphi'_{\P} E(\nu)= L(\nu)\varphi'_{\P}$ implies that
 $$(1\
b)\left(\begin{array}{cc} \nu^0 & 1\\ 0 & \nu^0
\end{array}\right)=\nu_0(1,\ b),$$
   i.e.
 $ (\nu^0, \ 1+b\nu^0)=(\nu^0, \, \nu^0 b)$, which results in a contradiction.

The proof in the case of item (2) of Corollary 7.3.1  is dual. Therefore, as desired, $\mathfrak{B_1}$ is not
homogeneous.\hfill $\square$

\medskip
{\bf Proposition 7.3.2} The bocs $\mathfrak{B_2}$  given in MW2 of
Theorem 5.6.1 is not homogeneous.

{\bf Proof.} Consider the local bocs $(\mathfrak{B}_1)_{\P}$ at $\P$, suppose we have  the triangular formula
(40) and the polynomial $\sigma_{\P}(\nu)$  of (41) in 5.4.  Then the proof is similar to that of Proposition
7.3.1. \hfill $\square$

{\bf Example.} Consider Example 2 of 2.1 and 4.2.  A sequence of reductions for $a=(1), b=(\nu)$ yields an
induced  bocs $\mathfrak{B_2} $:  \unitlength 1mm
\begin{picture}(30, 7) \put(10.00,-1.00){\circle*{1.00}}
\put(25.00,-1.00){\circle*{1.00}} \put(5.00,0.00){\oval(5,5)[t]} \put(5.00,0.00){\oval(5,5)[bl]}
\put(4.00,-3.00){\vector(3,1){3}} \put(11,-1){\line(1,0){12}} \put(23,-1){\vector(1,0){1}}
\put(10.00,2.40){\makebox{$\P$}} \put(25.00,2.40){\makebox{$\Q$}} \put(0.00,0.0){\makebox{$\nu$}}
\put(15.50,1.00){\makebox{$c$}}
\end{picture}
 with two vertices $\P,\Q$ and $\Gamma'(\P,\P)=\rK[\nu]$,
 $\Gamma'(\Q,\Q)=k$, $\delta(c)=0$.

\subsection{A critically  non-homogeneous local bocs $\mathfrak{B_3}$}

\kg {\bf Proposition 7.4.1} The bocs $\mathfrak{B_3}$  given in
MW3 of Theorem 5.6.1 is not homogeneous.

{\bf Proof.} Consider $$g_{ee}(\nu,\kappa)\in k[\nu,\kappa, \sigma_e(\nu)^{-1}, \sigma_e(\kappa)^{-1}]$$ in
Formula (40) of 5.4, and $\sigma_e(\nu), \sigma(\nu)$ given in Formula (41). Suppose that we have item (1) of
Corollary 7.3.1, i.e. $g_{ee}(\nu^0,\kappa^0)=0$, $\sigma(\nu^0)\sigma_e(\kappa^0)\ne 0$ and  $\nu^0$ ranges
over a cofinite subset $D$ of $k$. Then we are able to present an infinite list of iso-classes of one
dimensional objects $S_{\nu^0}, \forall \ \nu^0\in D\setminus {the roots of \sigma(\nu)}$, such that
$$(S_{\nu^0})_{\P}=\rK,\  S_{\nu^0}(\nu)=(\nu^0).\  {\rm thus}\  S_{\nu^0}(b_l)=\emptyset, \ l=1, \cdots, j,$$ when we
calculate $S_{\nu^0}$ according to Theorem 3.4.1. If $R(\mathfrak{B_3})$ is homogeneous, there must be a
cofinite subset $D_0\subset D$ such that $S_{\nu^0}$ is homogeneous for any $ \nu^0\in D_0$. We fix such a
$\nu^0$, denote $S_{\nu^0}$ by $S$ for simplicity. Let $(e):\ S\stackrel{\iota}{\longrightarrow} E
\stackrel{\pi}{\longrightarrow} S$ be an almost split conflation of $R(\mathfrak{B_3})$, then $E_{\P}=\rK^2$,
and the eigenvalue of $E(\nu)$ is $\nu^0$, since the sequence
 $$(e_0):\
0\longrightarrow S \stackrel{\iota_0}{\longrightarrow} E \stackrel{\pi_0}{\longrightarrow} S\longrightarrow 0 $$
  is
exact over $\rK[\nu, \sigma(\nu)^{-1}]$. Thus $E(b_l)=\emptyset$, $l=1, \cdots, j$, according to Theorem 3.4.1
and
Proposition 4.5.1. We conclude that $E(\nu)=\left(\begin{array}{cc}\nu^0&1\\
0&\nu^0
\end{array}\right)$, otherwise $\left(\begin{array}{cc} \nu^0&0\\ 0&\nu^0 \end{array}
\right)$ would lead to $(e)$ to be  split.

We define an object $L\in R(\mathfrak{B_3})$, such that $L_{\P}=\rK^2$,
 $$L(\nu)=\left(\begin{array}{cc}\nu^0 &0\\ 0&\kappa^0
\end{array}\right),\quad L(b_e)=\left(\begin{array}{cc}0&1\\ 0&0
\end{array}\right),\quad
L(b_l)=\left(\begin{array}{cc}0&0\\ 0&0 \end{array}\right) $$
 for
$l\neq e$. And we also define a morphism $\varphi: L \rightarrow S $ given by
 $\varphi_{\P}={1\choose 0}$, $\varphi(v)={0\choose 0}$
 for any dotted arrow $v$. We first claim
that $\varphi$ is not a split epimorphism. In fact, if we have a morphism $\psi: S \rightarrow L$ with
$\psi\varphi=id$, then $\psi_{\P}=(1,\ a)$, and
$$ S(b_e)\psi_{\P}-\psi_{\P}L(b_e)=\psi(\delta(b_e))=g_{ee}(S(\nu),
L(\nu))\psi(u_e),$$
(see Formula (40) of 5.4). Thus,  $-(1,\ a)\left(\begin{array}{cc} 0&1\\
0&0
\end{array}\right)=(g_{ee}(\nu^0, \nu^0)u_{e1}, g_{ee}(\nu^0,
\kappa^0)u_{e2})$ or $(0,\ -1)=(*,\ 0)$, which results in a contradiction. Thus there exists some $\varphi': L
\rightarrow E$ such that $\varphi'\pi=\varphi$,  whenever $(e)$ is an almost split conflation. Then
$\varphi'_{\P}{1\choose 0} = {1\choose 0}$ yields
 $\varphi'_{\P}=\left(\begin{array}{cc} 1&b\\ 0&c
\end{array}\right)$ and $L(\nu)\varphi'_{\P}=\varphi'_{\P}E(\nu)$
yields
 $$\left(\begin{array}{cc} \nu^0&0\\ 0&\kappa^0
\end{array}\right)\left(\begin{array}{cc} 1&b\\ 0&c \end{array}\right)
=\left(\begin{array}{cc} 1&b\\ 0&c
\end{array}\right)\left(\begin{array}{cc} \nu^0&1\\ 0&\nu^0
\end{array}\right),$$
 i.e.
 $$\left(\begin{array}{cc} *&\nu^0b\\
0&*\end{array}\right)=\left(\begin{array}{cc} *&1+b\nu^0\\ 0&*
\end{array}\right)$$ which is a contradiction.

The proof of $\kappa^0$ running over a cofinite subset of $\rK$ is dual. Therefore, as desired, $\mathfrak{B_3}$
is not homogeneous. \hfill$\square$

{\bf Example.} Consider  Example 3 of 2.1 and 4.2, let $\alpha\neq 0, 1$. A sequence of reductions given by
 $$A=\left(\begin{array}{cc} 1&0\\
0&1\end{array}\right),\ \ B=\left(\begin{array}{cc} 0&1\\ 0&0
\end{array}\right),\ \ C=\left(\begin{array}{cc} \emptyset &\emptyset\\ \nu &\emptyset
\end{array}\right),\ \ D=\left(\begin{array}{cc} d_3& d_4\\ d_1& d_2
\end{array}\right)$$
  yields an induced local bocs $\mathfrak{B}_3$,
with a layer $L=(\Gamma; \omega; d_1, d_2, d_3, d_4; v,\cdots)$, where $\Gamma'(\P,\P)=\rK[\nu]$,
$\delta(d_1)=v \nu-(\alpha\nu) v$, i.e. $e=1$,  $g_{11}(\nu, \kappa)=-\alpha\nu+\kappa$, which is not invertible
in $k[\nu,\kappa]$.
 \bcen
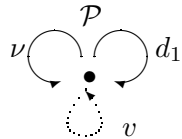
\begin{picture}(40,20)\unitlength=0.5pt
\put(0,60){\oval(40,40)[t]} \put(0,60){\oval(40,40)[bl]} \put(0,40){\vector(3,1){5}}
\put(50,60){\oval(40,40)[t]} \put(50,60){\oval(40,40)[br]} \put(50,40){\vector(-3,1){5}}
\put(25,32){\vector(-1,1){2}}

\put(21,39){$\bullet$} \qbezier[10](23,35)(-5,2)(23,0) \qbezier[10](23,35)(47,2)(23,0)

\put(51,0){\makebox{$v$}} \put(19,84){\makebox{$\P$}} \put(-34,60){\makebox{$\nu$}} \put(76,60){\makebox{$d_1$}}
\end{picture}
\ecen

\subsection{A  critically non-homogeneous local bocs $\mathfrak{B}_4$}

\kg {\bf Proposition 7.5.1}  The bocs $\mathfrak{B_4}$  given in
MW4 of Theorem 5.6.1 is not homogeneous.

{\bf Proof.} Suppose we have the contrary. We fix a $\lambda^0\in \rK$ with $h(\lambda^0)\neq 0$ (see Formula
(43) in 5.4), and continue the reductions such that $a_1=\emptyset, \cdots, a_{i-1}=\emptyset$. Let $a_i=(\nu)$
be a parameter, then we obtain an induced bocs $\mathfrak{B}'_4$ with a layer
 $$L'=((\Gamma')'; \omega'; b_1,\cdots,b_j;
u_1,\cdots,u_m),$$
  where the differentials $\delta(b_l)^0$ are given by Formula (40)
  and a polynomial $\sigma(\nu)$ is given by Formula (41) in
  Lemma 5.4.1.  It is clear that $\mathfrak{B}'_4$ is minimal if we set
   $(\Gamma')'(\P,\P) =\rK[\nu,\sigma(\nu)^{-1}] $.
    Denote by $\vartheta$ the reduction functor from
$R(\mathfrak{B}'_4)$ to $R(\mathfrak{B}_4)$. Define an infinite list of iso-classes of objects $\{S'_{\nu^0}|
\forall \, \nu^0\in \rK, \sigma(\nu^0)\neq 0\}$ in $R(\mathfrak{B}'_4)$, such that $$(S'_{\nu_0})_{\P}=\rK,\
S'_{\nu_0}(\nu)=(\nu_0)£¬ \ S'_{\nu_0}(b_l)=\emptyset \ {\rm for}\  l=1,2,\cdots,j.$$ Then we obtain an infinite
list $\{S_{\nu^0}=\vartheta(S'_{\nu^0})\}$ in $ R(\mathfrak{B}_4)$. If $R(\mathfrak{B_4})$ is homogeneous, there
must exist a cofinite subset $D_0\subset \rK \setminus \{ \mbox{the roots of } \sigma(\nu) \}$ such that for any
$\nu^0\in D_0$, $S_{\nu^0}$ is homogeneous. We fix such a $\nu^0\in D_0$, and denote $S_{\nu^0}$ by $S$ for
simplicity. If $(e ):\ S\stackrel{\iota}{\longrightarrow} E \stackrel{\pi}{\longrightarrow} S$ is an almost
split conflation of $R(\mathfrak{B}_4)$, there exists an almost split conflation $(e' ):\
S'\stackrel{\iota'}{\longrightarrow} E' \stackrel{\pi'}{\longrightarrow} S'$ of $R(\mathfrak{B}'_4)$, such that
$\vartheta(e')$ is equivalent to $(e)$ by Theorem 7.2.1.
Thus $$E'_{\P}=\rK^2,\ \  E'(\nu)=\left(\begin{array}{cc} \nu^0 & 1\\
0 & \nu^0 \end{array}\right),\ \ E'(b_l)=\emptyset,\ \ l=1,\cdots,j$$ by
 Proposition 6.6.1. Let us go back to $R(\mathfrak{B}_4)$, then we have
$S(\lambda)=(\lambda^0)$, $S(a_l)=\emptyset$;
 $E(\lambda)=\lambda^0 I_2$, $E(a_l)=\emptyset$,
$l=1, \cdots, i-1$.  Construct an object $L$ in $R(\mathfrak{B_4})$, such that
 \begin{align*}
& L_{\P}=\rK^2,\ \ L(\lambda)=\left(\begin{array}{cc} \lambda^0&1\\ 0&\lambda^0
\end{array}\right),\ \ L(a_l)=0  \mbox{ for }  l=1, \cdots,
i-1,\\
& L(a_i)=\left(\begin{array}{cc} \nu^0&0\\ 0&\nu^0
\end{array}\right),\ \ L(b_l)=0\ \mbox{ for}\ \ l=1, \cdots, j,
\end{align*}
  and a morphism
$\varphi: L\rightarrow S$ with $\varphi_{\P}={1\choose 0}$, $\varphi(v)={0\choose 0}$ for any dotted arrow $v$.
We first claim that $\varphi$ is not a split epimorphism. In fact if $\psi: S\rightarrow L$ with
$\psi\varphi=id$, then $\psi_{\P}=(1, \ a)$. But  $S_{\nu^0}(\lambda)\psi_{\P}=\psi_{\P}L(\lambda)$, i.e.
$\lambda^0(1, \ a )=(1,\ a) \left(\begin{array}{cc} \lambda^0&1\\
0&\lambda^0
\end{array}\right)$, or $(\lambda^0, \ \lambda^0 a)=(\lambda^0,\
1+a\lambda^0)$, which is a contradiction. Thus, there is some morphism $\varphi': L\rightarrow E$ such that
$\varphi'_{\P}\pi=\varphi$, whenever $(e)$ is an almost split conflation of $R(\mathfrak{B_4})$. Then
$\varphi'_{\P}{1\choose 0}={1\choose 0}$ leads to $\varphi'_{\P}=\left(\begin{array}{cc} 1&b\\ 0&c
\end{array}\right)$, and
$L(a_i)\varphi'_{\P}-\varphi'_{\P}E(a_i)=\varphi'(\delta(a_i))=0$ yields
 $$\left(\begin{array}{cc} \nu^0&0\\ 0&\nu^0
\end{array}\right) \left(\begin{array}{cc} 1&b \\0&c \end{array}\right)=
 \left(\begin{array}{cc} 1&b \\0&c \end{array}\right)
\left(\begin{array}{cc} \nu^0&1\\ 0&\nu^0
\end{array}\right),\ \rm{i.e.}\  \left(\begin{array}{cc} \nu^0&\nu^0 b\\ 0&\nu^0c\end{array}\right)
=\left(\begin{array}{cc} \nu^0&1+b\nu^0\\ 0&c\nu^0 \end{array} \right).$$
 It is a contradiction. Therefore, as desired,
  $\mathfrak{B_4}$ is not homogeneous. \hfill $\square$

{\bf Remark.} The assumption $(\lambda-\mu)^2\mid \wt{h}(\lambda, \mu )$  is necessary in the  proof. For
instance, if $\wt{h}(\lambda, \mu)=\lambda-\mu$,
$\varphi'(\wt{w})=\left(\begin{array}{cc}w_{11}&w_{12}\\w_{21}&w_{22}\end{array}\right)$, then
$$\begin{aligned}
\varphi'(\delta(a_i))&= \left(\begin{array}{cc}w_{11}&w_{12}\\w_{21}&w_{22}\end{array}\right)
\left(\begin{array}{cc}\lambda^0&0\\0&\lambda^0\end{array}\right)-
\left(\begin{array}{cc}\lambda^0&1\\0&\lambda^0\end{array}\right)
\left(\begin{array}{cc}w_{11}&w_{12}\\w_{21}&w_{22}\end{array}\right)\\
&=\left(\begin{array}{cc}-w_{21}&-w_{22}\\0&0\end{array}\right);
\end{aligned}$$ On the other hand,
$$\begin{aligned} \varphi'(\delta(a_i)) &=L(a_i)\varphi'_{\P}-\varphi'_{\P}E(a_i)\\
&=\left(\begin{array}{cc}\nu^0&0\\ 0&\nu^0
\end{array}\right)\left(\begin{array}{cc} 1&b\\
0&c\end{array}\right)
-\left(\begin{array}{cc} 1&b\\
0&c\end{array}\right)\left(\begin{array}{cc}\nu^0&1\\
0&\nu^0
\end{array}\right)\\ &=\left(\begin{array}{cc}  0& -1\\ 0&0 \end{array}\right)
\end{aligned}$$
which does not lead to any contradiction.

{\bf Example.} Consider Example 1 of 2.1 and 4.2. A sequence of reductions given by $a=1$, $b=\lambda$, yields
an induced local bocs $\mathfrak{B}_4$ with a layer $L=(\Gamma;\omega; c; v,\cdots)$,  where $\Gamma'(\P,
\P)=\rK[\lambda]$, $\delta(c)=0$. \unitlength=1mm
 \bcen
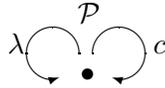
\begin{picture}(40,20)\unitlength=0.5pt
\put(0,60){\oval(40,40)[t]} \put(0,60){\oval(40,40)[bl]} \put(0,40){\vector(3,1){5}}
\put(50,60){\oval(40,40)[t]} \put(50,60){\oval(40,40)[br]} \put(50,40){\vector(-3,1){5}}

\put(21,39){$\bullet$}

 \put(19,84){\makebox{$\P$}}
\put(-34,60){\makebox{$\lambda$}} \put(76,60){\makebox{$c$}}
\end{picture}
\ecen

\subsection{An example in case of MW5}

\kg The minimally wild local bocs given in MW5 of Theorem 5.6.1 may be homogeneous, even strongly homogeneous
(see Example 4 of 6.5). Therefore a general wild bimodule problem does not necessarily have the non-homogeneous
property. We must concentrate on some special set of bimodule problems, so called bipartite bimodule problems
(see section 10), in which $P_1(\Lambda)$ are included. In this subsection,
 we will show an algebra, whose corresponding bimodule problem
 possess a reduction sequence with minimal size in case of MW5 of
 Theorem 5.6.1. We first prove two lemmas in order to treat the
 example and also the general case given in section 9.
  Following
Auslander and Reiten \cite{AR} we  recall the following lemma.

{\bf Lemma 7.6.1}\  Let  $\Gamma$ be an algebra (possibly infinite-dimensional). Denote by $\Gamma$-mod the
category of finite dimensional $\Gamma$-modules. Suppose that
$$(e)\qquad
0\longrightarrow M\stackrel{\iota}{\longrightarrow} E \stackrel{\pi}{\longrightarrow} L\longrightarrow 0 $$
 is an almost
split sequence in $\Gamma$-mod.

(1) If $L_1$ is a proper submodule of $L$, $E_1=(L_1)\pi^{-1}$, then the sequence
 $$0\longrightarrow M\stackrel{\iota_1}{\longrightarrow} E_1
\stackrel{\pi_1}{\longrightarrow} L_1\longrightarrow 0
 $$
  is split, where $\pi_1=\pi\mid_{E_1}$, $\iota_1=\iota$.

(2) If a submodule $ 0\neq K\subseteq M$, $M_2=M/K$ is a quotient module of $M$, and  $E_2=E/\iota(K)$, then the
sequence
$$
 0\longrightarrow M_2\stackrel{\iota_2}{\longrightarrow} E_2
\stackrel{\pi_2}{\longrightarrow} L\longrightarrow 0
$$ is split,
where $\iota_2, \pi_2$ are induced morphisms from $\iota$ and $\pi$ respectively.

{\bf Proof.} (2) (given by Pan Jun) Consider the following commutative diagram having two exact rows in
$\Gamma-$mod:
 \bcen \unitlength=0.8mm
\begin{picture}(105,30)
\put(8,0){$0$}\put(11,1){\vector(1,0){10}} \put(22,0){$M/K$}\put(34,1){\vector(1,0){12}} \put(47,0){$E/\iota
(K)$}\put(64,1){\vector(1,0){13}} \put(78,0){$L$}\put(84,1){\vector(1,0){12}} \put(97,0){$0$}
\put(-5,0){$(e_2):$}

\put(8,25){$0$}\put(11,26){\vector(1,0){12}} \put(24,25){$M$}\put(30,26){\vector(1,0){20}}
\put(51,25){$E$}\put(57,26){\vector(1,0){20}} \put(78,25){$L$}\put(84,26){\vector(1,0){12}} \put(97,25){$0$}
\put(-5,25){$(e):$}

\put(26,24){\vector(0,-1){20}} \put(53,24){\vector(0,-1){20}} \put(80,24){\vector(0,-1){20}}

\put(20,12){$\varphi_1$} \put(37,27){$\iota$} \put(65,27){$\pi$} \put(55,12){$\varphi_2$} \put(39,-2){$\iota_2$}
\put(69,-2){$\pi_2$} \put(81,12){$id$}

\put(31,16){$\s M/K$} \put(28,5){\vector(1,2){4}} \put(49,23){\vector(-3,-1){9}}
\multiput(45,2)(-1,1){9}{$\cdot$} \put(39,10){\vector(-1,1){3}}\put(31,7){$\s  id$} \put(43,8){$\s \xi$}
\put(45,20){$\s \psi$}
\end{picture}\ecen
where $\varphi_1$, $\varphi_2$ are natural projections respectively. Since $\varphi_1$ is not a split
epimorphism and $\iota$ is a left almost split morphism, there exists a morphism $\psi:\ E \rightarrow M/K$,
such that $\iota \psi=\varphi_1 id$. On the other hand, since the left square of the diagram determines a
pushout,  we have a morphism $\xi:\ E/\iota(K) \rightarrow M/K$, such that $\iota_2\xi=id$. Thus $(e_2)$ is
split as required.

(1) can be proved dually using pull back. \hfill$\square$

\medskip
{\bf Lemma 7.6.2}\ Let \unitlength 1mm
 \begin{picture}(30, 6) \put(10.00,1.00){\circle*{1.00}}
\put(25.00,1.00){\circle*{1.00}} \put(5.00,2.00){\oval(5,5)[t]} \put(5.00,2.00){\oval(5,5)[bl]}
\put(4.00,-1.00){\vector(3,1){3}} \put(11,1){\line(1,0){12}} \put(23,1){\vector(1,0){1}}
\put(10.00,4.00){\makebox{$\P$}} \put(25.00,4.00){\makebox{$\Q$}} \put(0.00,2.0){\makebox{$\lambda$}}
\put(15.50,3.00){\makebox{$c$}}
\end{picture} be a quiver, and
$\Gamma=kQ$ be the path algebra. We define an exact sequence
$$
(e) \qquad 0\longrightarrow S\stackrel{\iota}{\longrightarrow} E \stackrel{\pi}{\longrightarrow} S
\longrightarrow 0
$$
of $\Gamma$-mod, such that
$$
\begin{aligned}
&S_{\P}=\rK, S_{\Q}=\rK, S(\lambda)=(\lambda^0)\mbox{ for some }
\lambda^0\in\rK, S(c)=(1);\\
 &E_{\P}=\rK^2, E_{\Q}=\rK^2,
E(\lambda)=J_2(\lambda^0),
E(c)=I_2;\\
& \iota=(\iota_{\P}, \iota_{\Q}) \mbox{ with } \iota_{\P}=(0\ 1),
\iota_{\Q}=(0\ 1);\\
& \pi=(\pi_{\P}, \pi_{\Q}) \mbox{ with } \pi_{\P}={1\choose 0}, \pi_{\Q}={1\choose 0}.
\end{aligned}
$$
Then $(e)$ is not an almost split sequence of $\Gamma$-mod.

{\bf Proof.}\ Suppose that we have the contrary.  Let us define three $\Gamma$-modules $K, \overline{S},
\overline{E}$ respectively as follows: $K_{\P}=(0)$, $K_{\Q}=\rK$; $\overline{S}_{\P}=\rK$,
$\overline{S}_{\Q}=0$, $\overline{S}(\lambda)=(\lambda^0)$; $\overline{E}_{\P}=\rK^2$, $\overline{E}_{\Q}=\rK$,
$$\overline{E}(\lambda)=\left(\begin{array}{cc} \lambda^0&1\\
0&\lambda^0 \end{array}\right),\ \ \overline{E}(a_1)=\left(\begin{array}{c} 1\\
0 \end{array} \right).$$
 Then $\overline{E}$ is indecomposable, and the
commutative diagram below shows a contradiction to item (2) of Lemma 7.6.1.
\begin{center} \unitlength=0.7mm
\begin{picture}(100,90)
\put(-4,60){\oval(5,5)[t]} \put(-4,60){\oval(5,5)[bl]} \put(-3.5,57.5){\vector(3,1){2}}
\put(-4,30){\oval(5,5)[t]} \put(-4,30){\oval(5,5)[bl]} \put(-3.5,27.5){\vector(3,1){2}}

\put(1,23){\vector(1,-1){13}} \put(0,24){\makebox{$\rK$}} \put(1,53){\vector(1,-1){13}} \put(0,
54){\makebox{$\rK$}} \put(15,35){\makebox{$\rK$}} \put(16,5){\makebox{$0$}} \put(0,53){\vector(0,-1){23}}
\put(16,34){\vector(0,-1){23}}

\put(36,30){\oval(5,5)[t]} \put(36,30){\oval(5,5)[bl]} \put(36,27){\vector(3,1){3}}
\put(36,60){\oval(5,5)[t]} \put(36,60){\oval(5,5)[bl]} \put(36.5,57.5){\vector(3,1){2}}

 \put(41,23){\vector(1,-1){13}}
\put(39,24){\makebox{$\rK^2$}} \put(41,53){\vector(1,-1){13}} \put(39, 54){\makebox{$\rK^2$}}
\put(54,35){\makebox{$\rK^2$}} \put(55,5){\makebox{$\rK$}} \put(40,53){\vector(0,-1){23}}
\put(56,34){\vector(0,-1){23}}

\put(76,30){\oval(5,5)[t]} \put(76,30){\oval(5,5)[bl]} \put(76.5,27.5){\vector(3,1){2}}
\put(76,60){\oval(5,5)[t]} \put(76,60){\oval(5,5)[bl]} \put(76.5,57.5){\vector(3,1){2}}

 \put(81,23){\vector(1,-1){13}}
\put(80,24){\makebox{$\rK$}} \put(81,53){\vector(1,-1){13}} \put(80, 54){\makebox{$\rK$}}
\put(95,35){\makebox{$\rK$}} \put(95,5){\makebox{$\rK$}} \put(80,53){\vector(0,-1){23}}
\put(96,34){\vector(0,-1){23}}

\put(19,7){\vector(1,0){34}} \put(19,37){\vector(1,0){34}} \put(4,25.50){\vector(1,0){34}}
\put(4,55.50){\vector(1,0){34}} \put(59,7){\vector(1,0){34}} \put(59,37){\vector(1,0){34}}
\put(15,65){\makebox{$\rK$}} \put(44,25.50){\vector(1,0){34}} \put(44,55.50){\vector(1,0){34}}

\put(4,85.50){\vector(1,0){34}} \put(19,67){\vector(1,0){34}}

\put(0,83){\vector(0,-1){23}} \put(16,64){\vector(0,-1){23}} \put(40,83){\vector(0,-1){23}}
\put(56,64){\vector(0,-1){23}} \put(1,83){\vector(1,-1){13}} \put(0,84){\makebox{$0$}}
\put(15,65){\makebox{$\rK$}} \put(55,65){\makebox{$\rK$}} \put(41,83){\vector(1,-1){13}} \put(40,
84){\makebox{$0$}} \put(19,7){\vector(1,0){34}} \put(19,37){\vector(1,0){34}}

\put(-19,30){\makebox{$\s J_1(\lambda^0)$}} \put(-19,60){\makebox{$\s J_1(\lambda^0)$}}
\put(22,32){\makebox{$\scriptscriptstyle{J_2(\lambda^0)}$}}
\put(22,62){\makebox{$\scriptscriptstyle{J_2(\lambda^0)}$}} \put(62,32){\makebox{$\s J_1(\lambda^0)$}}
\put(62,62){\makebox{$\s J_1(\lambda^0)$}}

\put(-4,45){\makebox{$1$}} \put(37,45){\makebox{$I_2$}} \put(81,45){\makebox{$1$}}

\put(5,15){\makebox{$0$}} \put(45,15){\makebox{$\s {1\choose 0}$}} \put(85,15){\makebox{$1$}}

\put(16,21){\makebox{$0$}} \put(56,21){\makebox{$\s {1\choose 0}$}} \put(96,21){\makebox{$1$}}

\put(8,45){\makebox{$1$}} \put(48,45){\makebox{$I_2$}} \put(88,45){\makebox{$1$}}

\put(33,4){\makebox{$0$}} \put(73,4){\makebox{$1$}} \put(20,25){\makebox{$\s (0\ 1)$}} \put(62,25){\makebox{$\s
{1\choose 0}$}}

\put(28,38){\makebox{$\s (0\ 1)$}} \put(70,40){\makebox{$\s {1\choose 0 }$}} \put(20,55){\makebox{$\s (0\ 1)$}}
\put(62,55){\makebox{$\s {1\choose 0}$}}

\put(30,67.5){\makebox{$1$}} \put(20,86){\makebox{$0$}}

\put(-20,80){\makebox{$K:$}} \put(-30,50){\makebox{$(e)$}} \put(-30,20){\makebox{$(e_2)$}}
\end{picture}
\end{center}
Hence the assumption fails, $(e)$ is not an almost split sequence of $\Gamma$-mod. \hfill $\square$

\medskip
{\bf Proposition 7.6.1}\ Let $\mathfrak{B}=(\Gamma,\Omega)$  be a bocs with a layer
$L=(\Gamma;\omega;a_1,\cdots,a_n;v_1,\cdots,v_m)$, where $\Gamma(\P,\P)=\rK[\lambda, g_\P(\lambda)^{-1}],
\delta(a_1)=0$.
\begin{center}
\unitlength 1mm
 \begin{picture}(30, 6) \put(10.00,1.00){\circle*{1.00}}
\put(25.00,1.00){\circle*{1.00}} \put(5.00,2.00){\oval(5,5)[t]} \put(5.00,2.00){\oval(5,5)[bl]}
\put(4.00,-1.00){\vector(3,1){3}} \put(11,1){\line(1,0){12}} \put(23,1){\vector(1,0){1}}
\put(10.00,4.00){\makebox{$\P$}} \put(25.00,4.00){\makebox{$\Q$}} \put(0.00,2.0){\makebox{$\lambda$}}
\put(15.50,3.00){\makebox{$a_1$}}
\end{picture}.
\end{center}
 Then an edge reduction of $a_1=(1)$ yields an induced local
bocs $\mathfrak{B}'$. If $\mathfrak{B}'$ is minimal and tame, then $R(\mathfrak{B})$ is not homogeneous.

{\bf Proof.} Suppose we have the contrary. Let $D=\rK\setminus\{\mbox{the roots of } \sigma(\lambda)\}$  be the
domain of the minimal bocs
 $\mathfrak{B}'$ where
$\sigma(\lambda)$ is determined by  Formula (41) of 5.4, and we use $\lambda$ instead of $\nu$ . Then we define
an infinite set of iso-classes of objects of dimension $1$:
 $$ \{ S'_{\lambda^0}\mid S'_{\lambda^0}(\lambda)=(\lambda^0),
S'_{\lambda^0}(a_l)=\emptyset, l\geq 2, \forall \lambda^0\in D\}$$
 in $ R(\mathfrak{B}'_5)$. Let $\vartheta:
R(\mathfrak{B'})\rightarrow R(\mathfrak{B})$ be the reduction functor. We obtain also a set $ \{
S_{\lambda^0}=\vartheta(S'_{\lambda^0})\}$ in $ R(\mathfrak{B})$, where $$(S_{\lambda_0})_{\P}=\rK,\
(S_{\lambda_0})_{\Q}=\rK,\ (S_{\lambda^0})(\lambda)=(\lambda^0),\ (S_{\lambda_0})(a_1)=(1).$$ If
$R(\mathfrak{B})$ is homogeneous, there must be a cofinite subset $D_0\subset D$ such that for any $\lambda^0\in
D_0$, $S_{\lambda^0}$ is homogeneous. We fix such a $\lambda^0$,  and denote $S_{\lambda^0}$ by $S$ for
simplicity. If $(e): S\stackrel{\iota}{\longrightarrow} E \stackrel{\pi}{\longrightarrow}S $ is an almost split
conflation of $R(\mathfrak{B})$, then there exists an almost split conflation $(e'):\
S'\stackrel{\iota'}{\longrightarrow} E' \stackrel{\pi'}{\longrightarrow}S' $  of  $R(\mathfrak{B'})$ such that
$\vartheta(e')$ is equivalent to $(e)$ by  Theorem 7.2.1. Thus
 $$E'(\lambda)=\left(\begin{array}{cc} \lambda^0&1\\
0&\lambda^0\end{array}\right),\ \ E'(a_l)=\emptyset,\ l\geq  2$$ by Proposition  6.6.1.  Let us go back to
$R(\mathfrak{B})$, then
$$E_{\P}=\rK^2,\ \ E_{\Q}=\rK^2,\ \ E(a_1)=\left(\begin{array}{cc}1&0\\
0&1\end{array}\right).$$ It is clear that $\iota_{\P}=(0\ 1)$, $\iota_{\Q}=(0\ 1)$, $\pi_{\P}={1\choose 0}$,
$\pi_{\Q}={1\choose 0}$.

Let $\Gamma_0$ be the path algebra of quiver \unitlength 1mm
 \begin{picture}(30, 6) \put(10.00,1.00){\circle*{1.00}}
\put(25.00,1.00){\circle*{1.00}} \put(5.00,2.00){\oval(5,5)[t]} \put(5.00,2.00){\oval(5,5)[bl]}
\put(4.00,-1.00){\vector(3,1){3}} \put(11,1){\line(1,0){12}} \put(23,1){\vector(1,0){1}}
\put(10.00,4.00){\makebox{$\P$}} \put(25.00,4.00){\makebox{$\Q$}} \put(0.00,2.0){\makebox{$\lambda$}}
\put(15.50,3.00){\makebox{$a_1$}}
\end{picture}. We claim that $(e)$ is also
an almost split sequence in $\Gamma_0$-mod. In fact

(1) $(e)$ is a non-split exact sequence in $\Gamma_0$-mod. And $S$ is neither projective, nor injective.

(2) For any module $L$ and morphism $\varphi: L\rightarrow S$ in $\Gamma_0$-mod, if we set $L(a_l)=0$ for $l\geq
2$ and $\varphi(v_l)=0$ for $l=1,2,\cdots,m$, then $\varphi: L\rightarrow S$ can be regarded as a morphism of
$R(\mathfrak{B})$.

(3) If $\varphi$ is not a split epimorphism in $\Gamma_0$-mod, then it is not a split epimorphism in
$R(\mathfrak{B})$.

(4) There exists some $\varphi': L\rightarrow E$ with $\varphi'\pi=\varphi$ in $R(\mathfrak{B})$, whenever $(e)$
is an almost split conflation of $R(\mathfrak{B})$. Thus $$L(\lambda)\varphi'_{\P}=\varphi'_{\P} E(\lambda)\ \
{\rm and}\ \ L(a_1)\varphi'_{\Q}=\varphi'_{\P} E(a_1),$$ i.e. $\varphi'_0=(\varphi'_{\P},\varphi'_{\Q})$ is a
morphism in $\Gamma_0$-mod with $\varphi'_0\pi=\varphi$. Therefore $E\stackrel{\pi}{\longrightarrow} S$ is a
right almost split morphism in $\Gamma_0$-mod. The claim is proved.

Therefore $(e)$ is also an almost split sequence of $\Gamma_0$-mod, which leads to a contradiction to Lemma
7.6.2.
 \hfill$\square$

{\bf Example.} Consider Example 3 of 2.1 and 4.3. Let
$$\n=(2,2,2,2,2;2,2,2,2,2),$$ $(N,R,\n)$ be given in Formula (6)
of 2.6 and we denote $\overline {a^*}, \overline {b^*}, \overline
{c^*}, \overline {d^*}$ by $A, B, C, D$ respectively for
simplicity. Then a sequence of reductions given by
$$A=\left(\begin{array}{cc} 1&0\\ 0&1 \end{array} \right),\ \
  B=\left(\begin{array}{cc} 0&1\\ 0&0 \end{array} \right),\  \
  C=\left(\begin{array}{cc} c_3& c_4\\ \lambda& c_2\end{array}\right), \ \
  D=\left(\begin{array}{cc} d_3&d_4\\ d_1&d_2\end{array} \right)
$$
yields a local layered bocs. $\delta(c_2)^0$, $\delta(c_3)^0$, $\delta(c_4)^0$ and $\delta(d_1)^0$ satisfy
Formula (42) of 5.4. \bcen
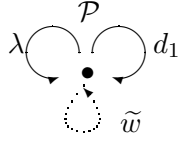
\begin{picture}(40,20)\unitlength=0.5pt
\put(0,60){\oval(40,40)[t]} \put(0,60){\oval(40,40)[bl]} \put(0,40){\vector(3,1){5}}
\put(50,60){\oval(40,40)[t]} \put(50,60){\oval(40,40)[br]} \put(50,40){\vector(-3,1){5}}
\put(25,32){\vector(-1,1){2}}

\put(21,39){$\bullet$} \qbezier[10](23,35)(-5,2)(23,0) \qbezier[10](23,35)(47,2)(23,0)

\put(51,0){\makebox{$\wt{w}$}} \put(19,84){\makebox{$\P$}} \put(-34,60){\makebox{$\lambda$}}
\put(76,60){\makebox{$d_1$}}
\end{picture}
\ecen If we continue the reductions given by $c_2=\emptyset$, $c_3=\emptyset$, $c_4=\emptyset$, then
$\delta(d_1)=\wt{w}\lambda-\lambda\wt{w}$ with $\wt{h}(\lambda, \mu)=-(\lambda-\mu)$ in Formula (42) and
$h(\lambda)=1$ in Formula (43). Furthermore, fix any $ \lambda^0\in k$, and denote $a_1$ by $\nu$, then
$\delta(d_2)^0$, $\delta(d_3)^0$, $\delta(d_4)^0$
 satisfy Formula (40) with $g_{ll}(\nu,\kappa)\in k[\nu,\kappa]$ being invertible for
 $l=2,3,4.$  Thus we obtain a sequence of
 parameterized bimodule problems with the end term in case of MW5.

Define a new size vector $$\wt{\n}=(2,2,2,2,2;3,3,3,3,3),$$ let $(\wt{N},\wt{R},\wt{\n})$ be given in Formula
(6) of 2.6. Then a sequence of reductions given by
$$\widetilde A=\left(\begin{array}{ccc} 0&1&0\\0& 0&1 \end{array} \right),\ \
\widetilde B=\left(\begin{array}{ccc} 0&0&1\\0& 0&0 \end{array}
\right),\ \ \widetilde C=\left(\begin{array}{ccc}
\emptyset&\emptyset&\emptyset\\ 0&
\lambda&\emptyset\end{array}\right), \ \
\widetilde D=\left(\begin{array}{ccc} d'_0&d_3&d_4\\
d_0&d_1&d_2\end{array} \right)
$$
yields an induced bocs $\mathfrak{B}$ with a layer $L=(\Gamma'; \omega; d_0, d_1,d_2,d'_0,d_3,d_4; v, \cdots)$,
where $\Gamma'(\P,\P)=\rK[\lambda]$, $\delta(d_0)=0$.
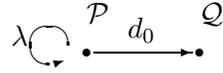
\begin{center}
\unitlength 1mm
 \begin{picture}(30, 6) \put(10.00,-1.00){\circle*{1.00}}
\put(25.00,-1.00){\circle*{1.00}} \put(5.00,0.00){\oval(5,5)[t]} \put(5.00,0.00){\oval(5,5)[bl]}
\put(4.00,-3.00){\vector(3,1){3}} \put(11,-1){\line(1,0){12}} \put(23,-1){\vector(1,0){1}}
\put(10.00,3.00){\makebox{$\P$}} \put(25.00,3.00){\makebox{$\Q$}} \put(0.00,0.0){\makebox{$\lambda$}}
\put(15.50,1.00){\makebox{$d_0$}}
\end{picture}
\end{center}
If we set $d_0=1$, then the induced bocs $\mathcal{B}'$ is minimal. Thus $\mathcal{B}$ and $\mathcal{B}'$ in our
example
 satisfy the hypothesis of Proposition 7.6.1, and
 the structure given in   Lemma 7.6.2 is as follows:
  $$\begin{array}{ccccc} &A&B&C&D\\ S&\left(\begin{array}{ccc} 0&1&0\\
  0&0&1 \end{array}\right)&\left(\begin{array}{ccc} 0&0&1\\ 0&0&0 \end{array}\right)&
  \left(\begin{array}{ccc} \emptyset&\emptyset&\emptyset\\ 0&\lambda^0&\emptyset \end{array}\right)&
  \left(\begin{array}{ccc} \emptyset&\emptyset&\emptyset\\
  1&\emptyset&\emptyset
  \end{array}\right)\\ &&&&\\
E&\left(\begin{array}{ccc}0_{2\times 2}&I_2&0\\ 0_{2\times 2}&0&I_2
\end{array}\right)&\left(\begin{array}{ccc}0_{2\times 2}&0&I_2\\
0_{2\times 2}&0&0
\end{array}\right)&\left(\begin{array}{ccc} \emptyset_{2\times
2}&\emptyset&\emptyset\\ 0_{2\times 2}&J_2(\lambda^0)&\emptyset
\end{array}\right)&\left(\begin{array}{ccc} \emptyset&\emptyset&\emptyset\\
I_2&\emptyset&\emptyset\end{array}\right)\\ &&&&\\
\overline{S}&\left(\begin{array}{cc} 1&0\\ 0&1
\end{array}\right)&\left(\begin{array}{cc} 0&1\\0&0
\end{array}\right)&\left(\begin{array}{cc} \emptyset&\emptyset\\
\lambda^0&\emptyset \end{array}\right)&\left(\begin{array}{cc}
\emptyset&\emptyset\\ 0&\emptyset \end{array}\right)\\ &&&&\\
\overline{E}&\left(\begin{array}{ccc} 0_{2\times 1}&I_2&0\\
0_{2\times 1}&0&I_2
\end{array}\right)&\left(\begin{array}{ccc} 0_{2\times 1}&0&I_2\\
0_{2\times 1}&0&0
\end{array}\right)&\left(\begin{array}{ccc} \emptyset_{2\times
1}&\emptyset&\emptyset\\ 0_{2\times 1}&J_2(\lambda^0)&\emptyset
\end{array}\right)&\left(\begin{array}
{ccc} \emptyset&\emptyset&\emptyset\\
({1\atop 0})&\emptyset&\emptyset
\end{array}\right)

\end{array}
$$
where $\overline{E}$ is indecomposable.

\newpage

\bcen
\section{One-sided differentials}
\ecen

\bigskip

\subsection{Matrix problems}

\kg Let $(\K,\M, H)$ be a bimodule problem of Definition 2.2.1 with a triangular basis $(A,B)$ of Definition
2.4.1. Suppose that $\rho_{\I_1\J_1}^{w_1}\in A$, $\M'$ is a $\K$-$\K$-subbimodule of $\M$ generated by
$$\{\rho_{\I\J}^w\in A\mid \rho_{\I\J}^{w}\succ \rho_{\I_1\J_1}^{w_1}\}$$ (see Proposition 2.4.1). Then $\M/\M'$
is a $\K$-$\K$-quotient bimodule of $\M$. And $(\K,\M/\M', d)$  is called a {\it partial bimodule problem} of
$(\K,\M,H)$, where $d: \K\rightarrow \M/\M'$ is a derivation given by $d(S)=(SH-HS)+\M'$ for any $S\in \K$. Thus
$(\K,\M/\M', d)$ belongs to the second set of 1.4, but may not satisfy Definition 2.2.1. Moreover, $(\K,\M/\M',
d)$ corresponds to a bocs $\mathfrak{B}$ given in \cite{CB2}.  Let $\mathfrak{A}=(\Gamma,\Omega)$ be the
corresponding bocs of $(\K,\M,H)$ with a layer $$L=(\Gamma'; \omega; a_1,\cdots,a_n; v_1,\cdots,v_m),$$
$\overline{\Gamma}$ be freely generated by $a_1, \cdots, a_{n_1}$ over $\Gamma'$, when
$a_{n_1}=(\rho^{w_1}_{\I_1\J_1})^*$. Since $$\delta_2(v_i)=\sum_{j,l<i} c_{jl}^i v_j\otimes v_l$$ with the
coefficients $c\in k$, which are independent of $a_1, \cdots, a_n$ by Formula (25) of 3.6, $\Omega$ is also a
$\overline{\Gamma}$-$\overline{\Gamma}$-coalgebra. Then $\mathfrak{B}=(\overline{\Gamma},\Omega)$ has a layer
$$\overline{L}=(\Gamma'; \omega; a_1, \cdots, a_{n_1}; v_1, \cdots, v_m),$$ which is called a {\it partial bocs}
of $\mathfrak{A}$. For their representation categories, we have
 $$Mat(\K,\M/\M')= \{M\in
Mat(\K,\M)\mid M_{p^w_{\I\J}q^w_{\I\J}}=0, \forall
(p^w_{\I\J},q^w_{\I\J})\succ(p^{w_1}_{\I_1\J_1},q^{w_1}_{\I_1\J_1}) \},$$ and $$ R(\mathfrak{B})= \{
\overline{M}\in R(\mathfrak{A})\mid M(a_j)=0, \forall\, j>n_1\}.$$ Both of them have exact structures inherited
from $Mat(\K,\M)$ and $R(\mathfrak{A})$ respectively. In order to treat with their reductions, we define another
algebraic structure, so called matrix problem.

\medskip
{\bf Definition 8.1.1} A {\it matrix problem} $(\K_1\times \K_2, \M, \mathcal{V})$ consists of the following
datum:
\begin{itemize}
\item [I$'$.\ \ ]\ A set of integers $T=\{1,2,\cdots,t\}$ and an equivalent relation $\sim$ on $T$.  Two subsets
$T_1, T_2\subset T$, which have $t_1, t_2$ elements respectively, and $\sim_1=\sim|_{T_1}$,
$\sim_2=\sim|_{T_2}$.

\item [II$'$.\ ]\ Two upper triangular matrix algebras $$\K_1=\{(s^1_{ij})_{t_1\times t_1}\} \ \ {\rm and}\ \
\K_2=\{(s^2_{ij})_{t_2\times t_2}\},$$ such that $$s^1_{ii}=s^1_{jj}, \ s^2_{ii}=s^2_{jj},\  s^1_{ii}=s^2_{jj}$$
when $i\sim j$, and a set of matrices $\mathcal{V}=\{(v_{ij})_{t_1\times t_2}\}$, $s_{ij}^1,\ s_{ij}^2$ for
$i<j$ and $v_{ij}$ satisfy the following equation system:
$$
\sum_{\I\ni i<j\in \J}(c_{ij}^{1l}x_{ij}^{1}+c_{ij}^{2l}x_{ij}^2)+\sum_{(i',j')\in
\I\times\J}c_{i'j'}^{3l}x_{i'j'}^3=0,
$$
which are indexed by $1\le l\le q_{\I\J}$ for some $q_{\I\J}\in \mathbb{N}$, and for each pair $$(\I,\J)\in
(T_1/\!\sim_1) \times (T_2/\!\sim_2).$$

\item [III$'$.] A $\K_1$-$\K_2$-bimodule $\M=\{(m_{ij})_{t_1\times t_2}\}$, where $m_{ij}$ satisfy the equation
system:
 $$\displaystyle
\sum_{(i,j)\in \I\times\J}d_{ij}^l z_{ij}=0$$
 which are indexed by $1\le l\le
q'_{\I\J}$ for some $q'_{\I\J}\in \mathbb{N}$, and for each pair $$ (\I,\J)\in (T_1/\!\sim_1) \times
(T_2/\!\sim_2).$$ \hfill $\Box$
\end{itemize}

{\bf Remark} We stress that the matrix problem given in Definition
8.1.1 is not necessarily a bimodule problem, since the derivation
is not defined. But it corresponds to a bocs. \hfill $\Box$

We also have a basis $A'_{\I\J}=\{\chi_{\I\J}^w\}$ and $A'$ of the solution space of equation system
\textrm{III}$^{\prime}$ similar to Formula (1) of 2.3. A basis of the solution space of equation system
\textrm{II}$^{\prime}$ for each pair $(\I,\J)$ is  denoted by $B'_{\I\J}$. And $B'$ is given similarly to
Formula (2) of 2.3. It is clear that $\rad\K_1$, $\rad\K_2$ and $\mathcal{V}$ are subspaces of the solution
space of the equation system \textrm{II}$^{\prime}$. We may choose basis of $\rad\K_1$, $\rad\K_2$, and
$\mathcal{V}$ respectively as follows.  First let $$W_0=\rad\K_1\cap\rad\K_2\cap\mathcal{V},$$ take a basis
$\{\xi_{\I\J}^{v_{123}}\}$ of $W_0$. Secondly, let $W_{12}$, $W_{23}$, $W_{31}$ be the complement spaces of
$W_0$ in $$\rad\K_1\cap\rad\K_2, \ \rad\K_2\cap\mathcal{V}, \ \mathcal{V}\cap\rad\K_1,$$ and take basis
$\{\xi_{\I\J}^{v_{12}}\}$, $\{\xi_{\I\J}^{v_{23}}\}$, $\{\xi_{\I\J}^{v_{31}}\}$ of $W_{12}$, $W_{23}$, $W_{31}$
respectively. Thirdly, let $W_1, W_2, W_3$ be the complement spaces of $$W_0\oplus W_{12}\oplus W_{31}, \
W_0\oplus W_{23}\oplus W_{12}, \ W_0\oplus W_{31}\oplus W_{23}$$ in $\rad\K_1$, $\rad\K_2$, $\mathcal{V}$, and
take basis $\{\xi_{\I\J}^{v_1}\}$, $\{\xi_{\I\J}^{v_2}\}$, $\{\xi_{\I\J}^{v_3}\}$ of $W_1$, $W_2$, $W_3$
respectively. Thus the dual basis
$$\{(\xi_{\I\J}^{v_{123}})^*,(\xi_{\I\J}^{v_{12}})^*,(\xi_{\I\J}^{v_{23}})^*,
(\xi_{\I\J}^{v_{31}})^*, (\xi_{\I\J}^{v_1})^*, (\xi_{\I\J}^{v_2})^*, (\xi_{\I\J}^{v_3})^*\}$$
 can be regarded as the
coefficient functions. Write
  $$\begin{aligned}
\{\xi_{\I\J}^{v'_1}\} &=\{\xi_{\I\J}^{v_{123}}, \xi_{\I\J}^{v_{12}}, \xi_{\I\J}^{v_{31}}, \xi_{\I\J}^{v_1}\},\\
\{\xi_{\I\J}^{v'_2}\} &=\{\xi_{\I\J}^{v_{123}}, \xi_{\I\J}^{v_{12}}, \xi_{\I\J}^{v_{23}}, \xi_{\I\J}^{v_2}\},\\
\{\xi_{\I\J}^{v'_3}\} &=\{\xi_{\I\J}^{v_{123}}, \xi_{\I\J}^{v_{23}}, \xi_{\I\J}^{v_{31}}, \xi_{\I\J}^{v_3}\}.
\end{aligned}$$

We can define the representation category $Mat(\K_1\times \K_2, \M)$ of a matrix problem $(\K_1\times \K_2, \M,
\mathcal{V})$ similarly to 2.2. Given any size vector $\n_1$ of $(T_1,\sim_1)$, $\n_2$ of $(T_2,\sim_2)$, such
that $(n_1)_i=(n_2)_i$ for any $i\in T_1\cap T_2$,  a triple $(N,R_1\times R_2, \n_1\times \n_2)$ is defined
similarly to Formula (6) of 2.6, such that $N, R_1, R_2$ and $V$ are $\n_1\times\n_2$, $\n_1\times\n_1$,
$\n_2\times\n_2$ and $\n_1\times\n_2$ partitioned matrices respectively. More precisely,
  $$\begin{aligned}
  N &=\sum_{(\I,\J)\in
(T_1/\!\sim_1 )\times ( T_2/\!\sim_2)} \sum_{w=1}^{r_{\I\J}}
\overline{(\chi_{\I\J}^w)^*}\otimes \chi_{\I\J}^w\\
R_1 &=\sum_{\I\in T_1/\!\sim_1}\overline{E^*_{\I}}\otimes E_{\I}+\sum_{(\I,\J)\in (T_1/\!\sim_1) \times
(T_1/\!\sim_1)}\sum_{v'_1=1}^{r^1_{\I\J}}\overline{(\xi_{\I\J}^{v'_1})^*}\otimes
\xi_{\I\J}^{v'_1}\\
R_2 &=\sum_{\I\in T_2/\!\sim_2}\overline{E_{\I}^*}\otimes E_{\I}+ \sum_{(\I,\J)\in (T_2/\!\sim_2)\times
(T_2/\!\sim_2)} \sum_{v'_2=1}^{r_{\I\J}^2}\overline{(\xi_{\I\J}^{v'_2})^*}\otimes
\xi_{\I\J}^{v'_2} \\
V &=\sum_{(\I,\J)\in (T_1/\!\sim_1) \times
(T_2/\!\sim_2)}\sum_{v'_3=1}^{r_{\I\J}^3}\overline{(\xi_{\I\J}^{v'_3})^*} \otimes \xi_{\I\J}^{v'_3}
\end{aligned}$$

 Then we start a reduction according to the matrix
equation
$$
R_1 N=V+NR_2, \quad {\rm or} \quad X_{pp}N_{pq}=V_{pq}+N_{pq}X_{qq},
$$
in the same way as in 2.6, where the second equation is the $(p,q)$-block of the first one.

Regularization.  If the equation $x_{pq}^3=0$ is not a linear combination of the equations in equation system
II$^{\prime}$, let $\overline{N}_{pq}=\emptyset$. Then we add $x_{pq}^3=0$ into the equation system
II$^{\prime}$. Let $R'_1$, $R'_2$, $V'$ be the restrictions of $R_1$, $R_2$, $V$ given by $x_{pq}^3=0$
respectively.

Suppose now that $x_{pq}^3=0$ is a linear combination of the equations in equation system II$^{\prime}$.

Edge reduction.  If $\P\ne\Q$, let
 $\overline{N}_{pq}=\left(\begin{array}{cc}0&I\\
0&0\end{array}\right)$.

Loop reduction. If $\P=\Q$, let $\overline{N}_{pq}=W$ for some Weyr matrix $W$.

In the last two cases, $R_1$ and $R_2$ are restricted by the equation
$$
X_{pp}\overline{N}_{pq}=\overline{N}_{pq}X_{qq}.$$ Denote by $R'_1$ and $R'_2$ the restricted algebras of $R_1$
and $R_2$ respectively, and
$$
V'= V-R'_1(\overline{N}_{pq}\otimes \chi)+(\overline{N}_{pq}\otimes \chi)R'_2.
$$ where $\chi$ is the first basis element of $A'$.

Let $$ N'= N-N_{pq}\otimes \rho $$ in all the $3$ cases, then we obtain a new  triple $(N', R'_1\times R'_2,
\n'_1\times \n'_2)$ similarly to 2.6. Consequently we also obtain a new matrix problem $(\K'_1\times \K'_2, \M',
\mathcal{V}')$ with a new index set $(T', \sim')$ similarly to 3.1. Finally we have a matrix equation over
$(\K'_1\times \K'_2, \M', \mathcal {V}')$:
$$
R'_1N'=V'+N'R'_2.
$$

Inductively we have a reduction sequence:
$$
(\triangle)\quad (N,R_1\times R_2,\n_1\times \n_2),\
 \cdots,\
(N^r,R_1^r\times R_2^r,\n_1^r\times \n_2^r),\ \cdots,\ (N^s,R_1^s\times R_2^s,\n_1^s\times \n_2^s).
$$
The notions of parametrization, and  free parametrization are
still valid for matrix problems.

\medskip

{\bf Proposition 8.1.1} Let $(\K,\M,H)$ be a bimodule problem. Then $(\K\times\K,\M,\mathcal{V})$  is a matrix
problem, where $T_1=T$, $T_2=T$, $\mathcal{V}=\{V=SH-HS\mid \forall S\in\K\}$.

{\bf Proof.} If $V=(v_{ij})_{t\times t}$, then $$v_{ij}=\sum_{l=i+1}^t s_{il}h_{lj}-\sum_{l=1}^{j-1}
h_{il}s_{lj}.$$ In fact $s_{ii}h_{ij}-h_{ij}s_{jj}=0$ in both cases
 $i\sim j$ and $i\nsim j$, and $v_{ij}=0$ when
$(i,j)\prec(p,q)$.  Then we obtain the required equation system II$^{\prime}$. \hfill $\square$

Let $(\K,\M,H)$ be a bimodule problem, and $\rho^{w_1}_{\I_1\J_1}$
be given in the beginning of the subsection, $(p,q)$ be defined in
2.4. Suppose in addition that if $$(p^w_{\I\J}, q^w_{\I\J})\succ
(p^{w_1}_{\I_1\J_1}, q^{w_1}_{\I_1\J_1}), \ \ {\rm then}\ \
p^w_{\I\J}< p^{w_1}_{\I_1\J_1}.$$ Let $\K_1$ be a  subalgebra of
$\K$ consisting of the submatrices of the matrices of $\K$ with
the $(i,j)$-th entries for $$p\geq i,j\geq p^{w_1}_{\I_1\J_1};$$
$\K_2$ be that for $$q\leq i,j\leq q^{w_1}_{\I_1\J_1};$$ let
$\overline{\M}$ consist of the submatrices of the matrices of $\M$
with the $(i,j)$-entries for $$p\geq i\geq p^{w_1}_{\I_1\J_1} \ \
{\rm and}\ \  q\leq j\leq q^{w_1}_{\I_1\J_1};$$ and $\mathcal{V}$
consist of those of
 $\{SH-HS\mid \forall\, S\in rad\K\}$.

{\bf Proposition 8.1.2} $(\K_1\times\K_2, \overline{\M}, \mathcal{V})$  defined above is a matrix problem, whose
representation category $Mat(\K_1\times \K_2, \overline{\M})$ is equivalent to $Mat(\K,\M/\M', d)$, as well as
$R(\mathfrak{B})$.

{\bf Proof.} \textrm{I$'$}. $T_1=\{p^{w_1}_{\I_1\J_1}, \cdots, p\}$, $T_2=\{q, \cdots, q^{w_1}_{\I_1\J_1}\}$.

\textrm{II$'$}. The equation system is given by Formula (18) of 3.6:
$$
\left\{\begin{array}{ll} x^1_{ij}=\sum_{w}(\zeta^w_{\I\J})^* b^w_{\I\J}(i,j) & \mbox{for } p\geq i,j\geq
p^{w_1}_{\I_1\J_1};\\[2ex]
x^2_{ij}=\sum_{w}(\zeta^w_{\I\J})^* b^w_{\I\J}(i,j) & \mbox{for } q\leq i,j\leq
q^{w_1}_{\I_1\J_1};\\[2ex]
{x^3_{ij}=\sum_l\Big(\sum_{w}(\zeta^w_{\I\L})^* b^w_{\I\L}(i,l)h_{lj}-\atop \qquad\quad
h_{il}\sum_w(\zeta^w_{\L\J})^*b_{\L\J}^w(l,j)\Big)} & \mbox{for } p\geq i\geq p^{w_1}_{\I_1\J_1}, q\leq j\leq
q^{w_1}_{\I_1\J_1}.
\end{array}\right.
$$

\textrm{III$'$}. $\overline{\M}$ is a $\K_1$-$\K_2$-bimodule, and the equation system is given by Formula (19)
of 3.6:
$$
z_{ij}=\sum_{w}(\rho^w_{\I\J})^* a^w_{\I\J}(i,j)  \mbox{  for } p\geq i\geq p^{w_1}_{\I_1\J_1}, q\leq j\leq
q^{w_1}_{\I_1\J_1}.\eqno\Box$$

{\bf Corollary 8.1.1} Let $(\K,\M,H)$ be a bimodule problem, $(\K_1\times\K_2, \overline{\M}, \mathcal{V})$ be
the matrix problem given by a partial bimodule problem $(\K,\M/\M', d)$. Suppose $(N_0,R,\n)$ is a triple of
$(\K,\M,H)$, and $(N, R_1\times R_2, V)$ is the corresponding triple of $(\K_1\times\K_2, \overline{\M},
\mathcal{V})$, such that $\forall\, i\in T_1$, $(n_1)_i=n_i$, and $\forall\, i\in T_2$, $(n_2)_i=n_i$. Then the
reduction sequence $(\triangle)$ of $(\K_1\times\K_2, \overline{\M}, \mathcal{V})$ starting from $(N, R_1\times
R_2, V)$ determines a reduction sequence $(*)$ of $(\K,\M,H)$ starting from $(N_0,R,\n)$.

{\bf Proof.} We use induction on the reduction step $r$. The case of $r=0$ is clear.

 Suppose that we have already a reduction sequence $(*)$ up to
$(N_0^r, R^r, \n^r)$ determined by $(\triangle)$, and $\chi^r$ is the first basis element of
$(\K_1^r\times\K_2^r, \overline{\M}^r, \mathcal{V}^r)$, $\rho^r$ is that of $(\K^r, \M^r, H^r)$, then
$\chi^r=\rho^r\mid_{\overline{\M}^r}$.

For step $(r+1)$, let $\overline{N}^r_{p_rq_r}$ be obtained at the $(r+1)$-th reduction of $(\triangle)$. Then
the size vector $\n^{r+1}$ is given by $\n^r$ and $\overline{N}^r_{p_rq_r}$ according to 3.1;
$H^{r+1}_{\n^{r+1}}= H^r_{\n^r}+ \overline{N}^r_{p_rq_r}\otimes \rho^r$;
$N_0^{r+1}=N_0^r-N^r_{p_rq_r}\otimes\rho^r$; and $R^{r+1}$ is a restriction of $R^r$ given in Formula $(10)$ of
2.6. \hfill$\square$

\subsection{One-sided differentials I}

\kg Let $(\K,\M,H)$ be a bimodule problem with a triangular basis $(A,B)$. Let $\rho^w_{\I_1\J_1}\in A$ be given
in the beginning of 8.1, and $\rho\in A$ be the first basis element given in 2.4.  We assume that
$$p^w_{\I_1\J_1}=p,$$ and if $$(p^w_{\I\J}, q^w_{\I\J})\succ (p,q^{w_1}_{\I_1\J_1}), \ \ {\rm then}\ \  p^w_{\I\J}<p.$$
 Thus the
partial bimodule problem $(\K,\M/\M', d)$ corresponds to a partial bocs $\mathfrak{B}$,
 and determines a matrix problem
$(\K_1\times\K_2, \overline{\M}, \mathcal{V})$, such that $$T_1=\{1\}, \ \ T_2=\{1,2,\cdots,t\};$$ and
$$\K_1=\{(s)\}, \ \ \K_2=\left\{S=(s_{ij})_{t\times t}\right\},$$ where $ s_{ij}=0 \mbox{ if }i>j, s_{ii}=s_{jj}
\mbox{ if }i\sim j $, and $s_{ii}=s$ if $i\sim 1$, $$\mathcal{V}=\{V=(\begin{array}{cccc}v_1& v_2& \cdots &
v_t\end{array})\},$$ moreover $s_{ij}\, (i<j)$ and $v_j$ satisfy equation system II$^{\prime}$ of Definition
8.1.1;
$$\overline{\M}=\{M=(\begin{array}{ccccccccc}\cdots& a& \cdots & b& \cdots & a & \cdots &
b&\cdots\end{array})_{1\times t}\},$$ and the entries of $M$ satisfy equation system III$^{\prime}$.

Therefore the partial bocs   $\mathfrak{B}$ has  a partial  layer
 $$L^p=(\Gamma'; \omega; a_1, \cdots, a_n;  b_1, \cdots, b_m; u, v, w,
q),$$ where the vertices $\P, \I_1, \cdots, \I_r\ $ are all trivial;  solid arrows  $$a: \P\rightarrow \I,\ \ b:
\P\rightarrow \P;$$ and dotted arrows $$u: \P\rightarrow \P, \ \ v: \P\rightarrow \I, \ \ w: \I\rightarrow \I, \
\ q: \I\rightarrow \P.$$ Note that $a, b$ stand for general solid arrows, and $u,v,w,q$ stand for general dotted
arrows.

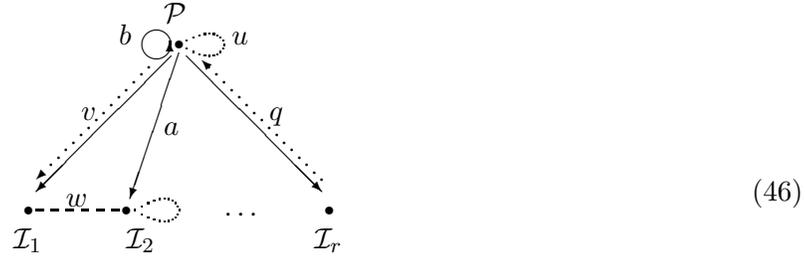
\begin{equation}{\unitlength=1mm
\begin{array}{c}
\begin{picture}(80,40)
\put(29,39){\circle{4}} \put(30.8,39){\vector(0,1){0.5}}
\put(32,39){\circle*{1}} \qbezier[10](32,39)(37,42)(38,39)
\qbezier[10](32,39)(37,36)(38,39)

\put(32,38){\vector(-1,-3){6.5}} \put(31,37.5){\vector(-1,-1){18}}
\put(33,37.5){\vector(1,-1){18}} \put(25,17){\circle*{1}}
 \put(12,17){\circle*{1}}  \put(52,17){\circle*{1}}
 \put(10,12){$\I_1$}  \put(50,12){$\I_r$}
 \put(25,12){$\I_2$}  \put(30,42){$\P$}

\qbezier[10](26,17)(31,20)(32,17)
\qbezier[10](26,17)(31,14)(32,17)

\mput(13,17)(2,0){6}{\line(1,0){1}}
\mput(28,36)(-1,-1){14}{\circle*{0.5}}
\put(15,23){\vector(-1,-1){2}}

\mput(36,36)(1,-1){16}{\circle*{0.5}}
\put(37,35){\vector(-1,1){2}}

\put(39,39){\makebox{$u$}} \put(19,29){\makebox{$v$}}
\put(44,29){\makebox{$q$}} \put(17,17.5){\makebox{$w$}}
\put(24,39){\makebox{$b$}} \put(30,27){\makebox{$a$}}
\put(38,15.5){\makebox{$\cdots$}}

\end{picture}
\end{array}}
\end{equation}
It is obvious, that the summands of $\delta(a_i)$ involve only the
monomials: (1) $v$, (2) $b^ev$, (3) $aw$; and those of $\delta(b)$
involve: (1) $u$, (2) $b^eu$, (3) $aq$ for some positive integer
$e$. Such differentials of solid arrows are said to be {\it
one-sided}.

Now we first consider the local partial bocs $\mathfrak{B}_{\P}$ at $\P$ obtained from $\mathfrak{B}$ by
deleting the vertices $\I_1, \I_2, \cdots, \I_r $. Thus $\mathfrak{B}_{\P}=(\Gamma_{\P}, \Omega_{\P})$ has a
partial layer $$L_{\P}^p=(\Gamma'_{\P}; \omega_{\P}; b_1, \cdots, b_m; u)$$ and the summands of $\delta(b_j)$
involve only the monomials of $u$ or $b^eu$. Then we have the following possibilities.

 $$ \mbox{\bf P1}\hspace{2cm} \left\{\begin{array}{lcl}
 \delta(b_1)^0 &= & g_{11}u_1 \\
 \delta(b_2)^0 &= & g_{21}u_1 \hspace{0.5cm}+ g_{22}u_2 \\
   &\cdots & \quad\cdots \\
\delta(b_m)^0 &= & g_{m1}u_1 \hspace{0.33cm}+ g_{m2}u_2\ \ \,
                 +\cdots+ g_{m,m}u_{m}
\end{array}\right.
$$
where $u_1, u_2, \cdots, u_{m}$ are linearly independent,
$g_{ll'}\in \rK, \ l\le l'$,\ $g_{ll}\neq 0$. Then
$\mathfrak{B}_{\P}$ is of finite type.

Otherwise we have the following formula.
 \begin{equation} \left\{\begin{array}{lcl}
 \delta(b_1)^0 &= & g_{11}u_1 \\
 \delta(b_2)^0 &= & g_{21}u_1 \hspace{0.5cm}+ g_{22}u_2 \\
  &\cdots & \quad\cdots \\
 \delta(b_{j-1})^0 &=& g_{j-1,1}u_1 + g_{j-1,2}u_2+ \cdots +
 g_{j-1,j-1}u_{j-1}\\
\delta(b_j)^0 &= & g_{j1}u_1 \hspace{0.5cm}+ g_{j2}u_2\quad
+\cdots + g_{j,j-1}u_{j-1}
\end{array}\right.
\end{equation}
where $u_1, u_2, \cdots, u_{j-1}$ are linearly independent,
$g_{ll'}\in \rK,\ l\le l'$, $g_{ll}\neq 0$. If we set
$b_l=\emptyset$, $u_l=0, l=1, \cdots, j-1$,  $b_j=\lambda$, then
\begin{equation}
\left\{\begin{array}{lcl}
\delta(b_{j+1})^0&=&h_{j+1,j+1}(\lambda)u_{j+1}\\
\delta(b_{j+2})^0&=&h_{j+2, j+1}(\lambda)u_{j+1}+ h_{j+2,
j+2}(\lambda)u_{j+2}\\
\cdots&& \quad \cdots\\
\delta(b_{m'})^0&=&h_{m', j+1}(\lambda)u_{j+1}\quad+h_{m',
j+2}(\lambda)u_{j+2}+\cdots +h_{m'm'}(\lambda)u_{m'}
\end{array}\right.
\end{equation}
where $u_{j+1},\cdots,u_{m'}$ are linearly independent,
$h_{ll}(\lambda)\neq 0$.

{\bf P2}\ $m'=m, h_{ll}(\lambda)$ are all non-zero constants. Then
 $\mathfrak{B}_{\P}$ is of tame type.

 {\bf P3}\   $m'=m,\ h_{11}(\lambda),\cdots,h_{e-1,e-1}(\lambda)$
 are non-zero constants  but $h_{e,e}(\lambda)$ is not a constant.
 Then
$\mathfrak{B}_{\P}$ is of wild type.

 {\bf P4}\   $m'<m,$ \ $ \mathfrak{B}_{\P}$ is of wild type.

\subsection{One-sided differentials II }

\kg Let us go back to the partial bocs $\mathfrak{B}$ given at the beginning of 8.2. Suppose that the local
partial bocs $\mathfrak{B}_{\P}$ is in case P2 of 8.2. If $$a_1\prec\cdots \prec a_{n_0}\prec b_j\prec
a_{n_0+1},$$ then we have either a triangular formula:
\begin{equation*} \left\{\begin{array}{lcl}
\delta(a_1)^0&=&f_{11}v_1\\
\delta(a_2)^0&=&f_{21}v_1+f_{22}v_2\\
\cdots&&\quad \cdots\\
\delta(a_{n_0})^0&=&f_{n_0 1}v_1+f_{n_0 2 }v_2+\cdots
+f_{n_0n_0}v_{n_0}
\end{array}\right.
\end{equation*}
where $v_1,v_2,\cdots,v_{n_0}$ are linearly independent, and
$f_{ll'}\in \rK$, $f_{ll}\neq 0$, or we have a triangular formula:
\begin{equation}
\left\{\begin{array}{lcl}
\delta(a_1)^0&=&f_{11}v_1\\
\cdots &&\quad \cdots\\
\delta(a_{n(1)-1})^0&=&f_{n(1)-1, 1}v_1+\cdots+f_{n(1)-1,
n(1)-1}v_{n(1)-1}\\
\delta(a_{n(1)})^0&=&f_{n(1), 1} v_1\quad +\cdots+f_{n(1),
n(1)-1}v_{n(1)-1}
 \end{array}\right.
 \end{equation}
 where $v_1,v_2,\cdots,v_{n(1)-1}$ are linearly independent, and $f_{ll'}\in \rK$, $f_{ll }\neq 0
 $. Set $a_l=\emptyset$, $v_l=0$, then we obtain an induced bocs by
  a series of regularizations. We continue  the
procedure, then obtain a sequence of positive integers $n(1), n(2),\cdots,n(i)$, and finally reach to the first
formula  by induction. Let
\begin{equation}I=\{1, 2, \cdots, n_0\}\setminus \{n(1), n(2), \cdots, n(i)\}\end{equation}

 Let us
write $\overline{a}_l=a_{n(l)}, l=1, \cdots, i$,
$\overline{b}=b_j$ for simplicity. The differentials  of
$\overline{b}$ with respect to $\overline{a}_1, \cdots,
\overline{a}_i $ is either
\begin{equation}
\begin{array}{ccl}
\delta(\overline{b})^0&=&\overline{a}_1(h_{11}q_1)\\
&&+\overline{a}_2(h_{21}q_1+h_{22}q_2)\\
&&+\cdots \quad \cdots\\
&&+\overline{a}_i(h_{i 1}q_1+h_{i 2}q_2+\cdots+h_{ii}q_i)
\end{array}
\end{equation}
where $q_1, q_2, \cdots, q_i$ are linearly independent and
$h_{11}h_{22}\cdots h_{ii}\ne 0$, or there exists some $p<i$, such
that
\begin{equation}
\delta(\overline{b})^0=\sum_{l=1}^i\overline{a}_l\Big(
\sum_{j=1}^ph_{lj}q_j\Big).
\end{equation}

{\bf Proposition 8.3.1} Suppose we are given an  one-sided differential partial bocs $\mathfrak{B}$,  in diagram
(46) of 8.2, such that $\mathfrak{B}_{\P}$ satisfies P2 of 8.2,  and $\delta(\overline{b})$ satisfies Formula
(52). Then $R(\mathfrak{B})$ is not homogeneous.

{\bf Proof.} Since $p<i$, there exists at least one vertex $\I\in
\{\I_1, \cdots, \I_r\} $ such that the number of the solid arrows:
$\P\rightarrow \I$ is more than that of the dotted arrows:
$\I\rightarrow \P$ given in (52). By a suitable reordering, we may
assume that the solid arrows  are $\overline{a}_1, \cdots,
\overline{a}_n$; and the dotted arrows are
 $q_1, \cdots, q_m$ respectively with $n>m$.  Let $\mathfrak{B'}$
be the induced partial bocs with two vertices $\P,\I$ obtained
from $\mathfrak{B}$ by deletion  and
 $$\delta(\overline{b})^0=\sum_{l=1}^n\overline{a}_l\Big(
 \sum_{j=1}^mh_{lj}q_j\Big)=\sum_{j=1}^m
\Big(\sum_{l=1}^n h_{lj}\overline{a}_l\Big) q_j$$
  in $\mathfrak{B'}$.
Suppose  that $R(\mathfrak{B})$ is homogeneous, then
$R(\mathfrak{B'})$ is homogeneous  by Corollary 7.2.1. Consider a
system of equations in $n$ indeterminates:
$$
\sum_{l=1}^nh_{lj}x_l=0, \quad j=1,\cdots, m
$$
Because the rank of the coefficient matrix is smaller than $n$,
there is a non-zero solution of the equation system, say
$x_1=\mu^1, x_2=\mu^2, \cdots, x_n=\mu^n$, where $\mu^1=0, \cdots,
\mu^{e-1}=0$, but $\mu^e\neq 0$ for some $1\leqslant e\leqslant
n$.

Define an object $S_{\nu^0}\in R(\mathfrak{B'})$ for any $\nu^0\in \rK$, such that $$(S_{\nu^0})_{\P}=\rK,\ \
(S_{\nu^0})_{\I}=0,\ \
 S_{\nu^0}(\overline{b})=(\nu^0),\ \  S_{\nu^0}(b)=\emptyset\ \  \text{for any loop}\ \
 b\ne \overline{b}.$$ If $R(\mathfrak{B'})$
is homogeneous, there exists a cofinite subset $D_0\subseteq \rK$ such that $S_{\nu^0}$ is homogeneous for any
$\nu^0\in D_0$. We fix an $\nu^0\in D_0$, denote $S_{\nu^0}$ by $S$ for simplicity, and suppose that $(e): S
\stackrel{\iota}{\longrightarrow} E \stackrel{\pi}{\longrightarrow} S $ is an almost split conflation in
$R(\mathfrak{B'})$. Then $(e)$ is also an almost split conflation in $R(\mathfrak{B'_{\P}})$ by Theorem 7.2.1.
Thus
 $$E_{\P}=\rK^2,\ \ E_{\I}=0,\ \  E(\overline b)=\left(\begin{array}{cc} \nu^0&1\\
0&\nu^0\end{array} \right), \mbox{ and } \ E(b)=\emptyset\ \mbox{ for}\ \  b\ne \overline{b}$$
 by  Proposition 6.6.1.  Now we construct
an object $L\in R(\mathfrak{B'})$, such that  $$L_{\P}=\rK,\  L_{\I}=\rK,\ L(\overline{a}_l)=(\mu^l), \
L(a_j)=\emptyset\  {\rm for}\  j \in I,  \ L(\overline b)=(\nu^0), \ L(b)=\emptyset \ {\rm for}\
b\ne\overline{b};$$
  and we also construct a morphism $\varphi:
L\rightarrow S $, such that $$\varphi_{\P}=1,\ \ \varphi_{\I}=0, \ \ \varphi(u)=0,\ \ \varphi(v)=0,\ \
\varphi(q)=0.$$ Then $\varphi$ is not a split epimorphism. In fact, if $\psi:  S\rightarrow L $ with $\psi
\varphi=id$, then $\psi_{\P}\varphi_{\P}=1$, $\psi_{\P}=1$, and
$$
-\mu^e=S(\overline{a}_e)\psi_{\I}-\psi_{\P}L(\overline{a}_e)
=\psi(\delta(\overline{a}_e))=0.
$$
We have  a contradiction. Thus there exists some $\varphi':
L\rightarrow E$, such that  $\varphi'\pi=\varphi$, whenever $(e)$
is an almost split conflation in $R(\mathfrak{B'})$. Then
$\varphi'_{\P}=(1,\ a)$ for some $a \in \rK$. On the other hand,
$$L(\overline{b})\varphi'_{\P}-\varphi'_{\P}E(\overline{b})
=\varphi'(\delta(\overline{b}))=0,\quad {\rm i.e.}\quad
 \nu^0(1, \ a )=(1,\ a) \left(\begin{array}{cc} \nu^0&1\\
0&\nu^0
\end{array}\right).$$
  We have  again a contradiction. Therefore the assumption
fails, and $R(\mathfrak{B})$ is not homogeneous. \hfill $\square$

The possibilities of the differentials of $ a_{n_0+1},\cdots, a_n$
with respect to $\overline b$ are given in the following two
formulae.
\begin{equation}
\left\{\begin{array}{lclll} \delta(a_{n_0+1})^0&=&\displaystyle
\sum_{l \in I}f_{n_0+1,l}(\overline{b})v_l
&+\ f_{n_0+1,n_0+1}(\overline{b})v_{n_0+1}&\\
\cdots&&\quad \cdots&&\\
\delta(a_n)^0&=&\displaystyle \sum_{l \in
I}f_{n,l}(\overline{b})v_l&+\ f_{n,n_0+1}(\overline{b})v_{n_0+1}
&+\cdots+f_{n,n}(\overline{b})v_n\\
\end{array}\right.
\end{equation}
where $\{v_l\ |\ \forall\ l \in I\}$ and
$\{v_{n_0+1},\cdots,v_n\}$ are linearly independent,
$f_{ll}(\overline b)\neq 0$. Or

\begin{equation*} \left\{\begin{array}{lclll}
\delta(a_{n_0+1})^0&=&\displaystyle \sum_{l \in I}f_{n_0+1,l}
(\overline{b})v_l&+\ f_{n_0+1,n_0+1}(\overline{b})v_{n_0+1}&\\
\cdots&&\quad \cdots&&\\
\delta(a_{n_1-1})^0&=&\displaystyle \sum_{l \in
I}f_{n_1-1,l}(\overline{b})v_l&+\
f_{n_1-1,n_0+1}(\overline{b})v_{n_0+1}
&+\ \cdots+f_{n_1-1,n_1-1}(\overline{b})v_{n_1-1}\\[+3ex]
\delta(a_{n_1})^0&=&\displaystyle \sum_{l \in
I}f_{n_1,l}(\overline{b})v_l&+\
f_{n_1,n_0+1}(\overline{b})v_{n_0+1}
&+\ \cdots+f_{n_1,n_1-1}(\overline{b})v_{n_1-1}\\
\end{array}\right.
\end{equation*}
where $\{v_l\ |\ \forall\ l \in I\}$ and
$\{v_{n_0+1},\cdots,v_{n_1-1}\}$ are linearly independent,
$f_{ll}(\overline{b})\neq 0$.

{\bf Proposition 8.3.2} Suppose we are given an one-sided differential partial bocs $\mathfrak{B}$ in diagram
(46) of 8.2, such that $\mathfrak{B}_\P$ satisfies P2 of 8.2; $\delta(\overline{b})$ satisfies Formula (51); and
$\delta(a)$'s with respect to $\overline b$ satisfy the preceding second formula.  Then $\mathfrak{B}$ is not
homogeneous.

{\bf Proof.} We set $\overline{a}_l=0$ for $l=1,\cdots,i$, $a_l=\emptyset$ for $l \in I$,
$\overline{b}=\lambda$, $b_l=\emptyset$ for $l\ne j$, and $a_{n_0+1}=\emptyset,\cdots,a_{n_1-1}=\emptyset$. If
$a_{n_1}:\ \P\rightarrow \I$, then we delete $T/\sim \setminus \{\P,\I\}$ and obtain an induced partial bocs
  {\unitlength 1mm \begin{picture}(30, 6) \put(10.00,1.00){\circle*{1.00}}
\put(25.00,1.00){\circle*{1.00}} \put(5.00,2.00){\oval(5,5)[t]}
\put(5.00,2.00){\oval(5,5)[bl]} \put(4.00,-1.00){\vector(3,1){3}}
\put(11,1){\line(1,0){12}} \put(23,1){\vector(1,0){1}}

\put(0.00,2.0){\makebox{$\lambda$}}
\put(15.50,3.00){\makebox{$a_{n_1}$}}
\end{picture}},
 where $\delta(a_{n_1})=0$.
The non-homogeneous property  of $\mathfrak{B}$ follows from
Proposition  7.3.2 and Corollary 7.2.1. \hfill$\square$

\subsection{The structure of  loop $\overline{b}$}

\kg Let the partial bocs $\mathfrak{B}$ be given as in 8.2 such
that the induced local partial bocs $\mathfrak{B}_{\P}$ is in case
P2 of 8.2,  $\delta(\overline{b})$ satisfies Formula (51) of 8.3,
 and $\delta(a)$'s with respect to $\overline b$  satisfy Formula (53) of 8.3.
We will prove in this subsection, that the block determined by
loop $\overline{b}$ will not contain  any parameter  after any
series of reductions, whenever the parameters appear in a local
bocs.

We first assume that the index set $I=\emptyset$ in Formula (50) of 8.3  and $m=1$, $\overline{b}=b_1$ given in
8.2. For any size vector $\n$ of $(T,\sim)$, we have a size vector $$\n_1\times\n_2 \ \ {\rm of}\ \ (T_1,
\sim_1)\times (T_2,\sim_2)$$ according to Corollary 8.1.1. Then  we start a sequence of reductions from triple
$(N, X\times R, \n_1\times\n_2 )$ given by  the matrix equation:
\begin{eqnarray}
&&X(D_1 \, \cdots\, D_i\, B\, A_{i+1}\,\cdots\, A_{i+j})\nonumber \\
&=&(0\cdots 0 \quad 0 \quad V_1\,\cdots \, V_j)+(D_1\, \cdots \, D_i\, B\, A_{i+1}\, \cdots\, A_{i+j})R
\end{eqnarray}
where matrices $D_1,\cdots,D_i$ and $B$ correspond to the solid arrows $\overline{a}_1,\cdots,\overline{a}_i$
and $\overline{b}$ respectively,  $A_l$ correspond to the solid arrows $a_l: \P\rightarrow \I$ for $l\geq i+1$;

$$ R= \left(\begin{array}{c} \unitlength=0.6mm
\begin{picture}(70,70)
\mput(0,40)(2,0){35}{\line(1,0){1}}
\mput(30,30)(2,0){20}{\line(1,0){1}}
\mput(30,30)(0,2){20}{\line(0,-1){1}}
\mput(40,0)(0,2){35}{\line(0,-1){1}} \put(2,63){\makebox{$Y_1$}}
\put(12,53){\makebox{$\ddots$}} \put(22,43){\makebox{$Y_i$}}
\put(18,59){\makebox{$W$}} \put(32,63){\makebox{$Q_1$}}
\put(34.5,53){\makebox{$\vdots$}} \put(32,43){\makebox{$Q_i$}}
\put(52,53){\makebox{$W$}} \put(32.5,33){\makebox{$X$}}
\put(42,33){\makebox{$V_1$}} \put(52,34.5){\makebox{$\cdots$}}
\put(62,33){\makebox{$V_j$}} \put(42,23){\makebox{$Y_{i+1}$}}
\put(52,13){\makebox{$\ddots$}} \put(61,3){\makebox{$Y_{i+j}$}}
\put(60,19){\makebox{$W$}}

\end{picture}
\end{array} \right)
$$
 where $X$, $Y_l$ are  $n_{\P}\times n_{\P}$,  $n_{\I}\times
n_{\I}$ invertible matrices, $V,W,Q$ correspond to the dotted
arrows $v, w, q$ respectively. In particular, $Q_1, \cdots, Q_i$
are blocks with entries being  linearly independent of all the
other entries according to Formula (51). Denote $(D_1\ \cdots\
D_i\ B\ A_{i+1}\ \cdots\ A_{i+j})$ by $N$.

{\bf Lemma 8.4.1}  All the reductions sequence ($\Delta$) starting
from the triple $(N, X\times R, \underline{n}_1\times\n_2)$
defined above must reach a step $r$, such that either

(1) the $r$-th step is an edge reduction before $A_{i+1}$, and
$(X_{p_{r-1},p_{r-1}})'=X_{(p_{r-1})_1,(p_{r-1})_1}$ (see the edge
reduction given in 3.1); or

(2) the $r$-th step is a
 loop reduction with a  Weyr matrix
$W^0$ inside $B$  before $A_{i+1}$. Where all the reductions
before $r$-th step are edge reductions and regularizations.

Moreover, the induced triple $(N^r, X^r\times R^r, \underline{n}^r_1\times\n_2^r)$
 has the property that $Q^r_l=Q_l$, and the
entries of the blocks $Q^r_l$ and the nilpotent blocks of $X^r$
are linearly independent and linearly independent of all the other
entries of $X^r,R^r,V^r$.

{\bf Proof.}  Assume that the number of the columns of matrix $ (D_1\cdots D_i\  B)$ equals  $m$. Our assertion
is shown by induction on $m$. If $m_B$ stands for the size of $B$ and $m=m_B$, then the first reduction is a
loop reduction for $B$, we obtain the case (2) for $r=1$. Otherwise $m>m_B$, $\overline{a}_1$ is an edge and
after an edge reduction, $X$ is restricted to
 $$X'=\left(\begin{array}{cc} X_{11}&X_{12}\\
0&X_{22} \end{array}\right).$$
 If $X_{22}= 0$, then $X'=X_{11},$ we obtain the case (1) for $r=1$.

 If $X_{22}\ne 0$, then $X_{22}$ is
independent of $X_{11}$ and $Y$'s,  if $X_{11}\ne 0$ and
$X_{22}\ne 0$, then $X_{12}$ is independent of $Y, Q, W, V$
respectively. Write
$$m_1=m-(\mbox{the number of columns of } D_1),$$ then $m_1 \le m$. Note that the first reduction for
$\overline{a}_1$ must be an edge reduction, but we write here also
a regularization for a unified statement in the further
reductions. If it is a regularization, then we still have $m_1<m$.

{} \vspace{2mm}

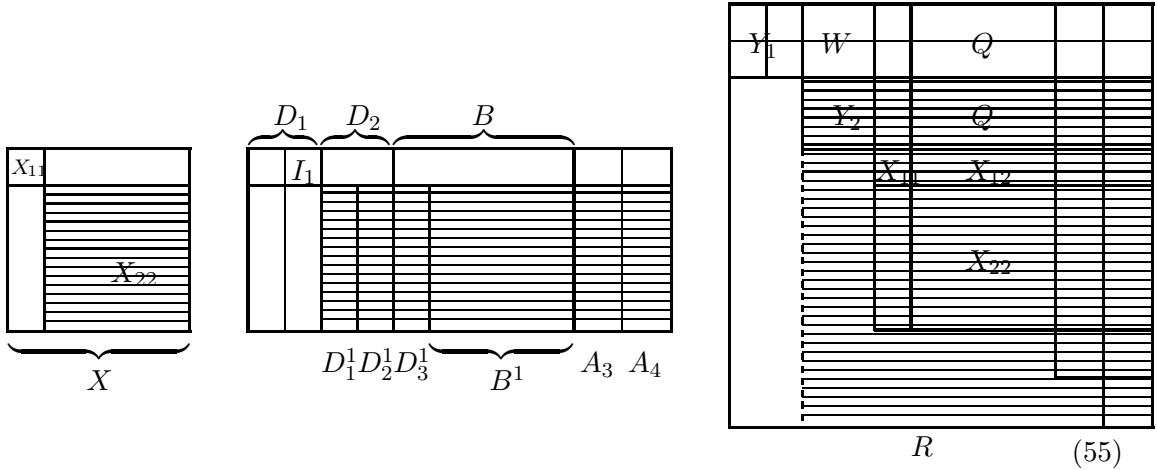
\begin{equation}\begin{array}{c} {\unitlength=0.8mm
\begin{picture}(185,60)\thicklines
\put(0,16){\framebox(30,30){}} \put(0,40){\line(1,0){30}}
\put(6,16){\line(0,1){30}}

\put(0.5,42){\makebox{$\s X\!_{11}$}}
\put(17,24){\makebox{$X_{22}$}}

\thinlines\mput(6,17.5)(0,1.5){16}{\line(1,0){24}}

\thicklines \put(40,16){\framebox(70,30){}}
\put(40,40){\line(1,0){70}} \put(52,16){\line(0,1){30}}
\put(64,16){\line(0,1){30}} \put(70,16){\line(0,1){24}}
\put(94,16){\line(0,1){30}} \put(58,16){\line(0,1){24}}

\put(47,41){\makebox{$I_1$}}

\thinlines\mput(52,17.8)(0,1.5){15}{\line(1,0){58}}
\put(102,16){\line(0,1){30}}\put(46,16){\line(0,1){30}}

\thicklines \put(120,0){\framebox(70,70){}}
\put(120,58){\line(1,0){70}} \put(132,46){\line(1,0){58}}
\put(144,40){\line(1,0){46} }\put(144,70){\line(0,-1){54}}
\put(150,70){\line(0,-1){54}} \put(132,70){\line(0,-1){24}}
\put(126,70){\line(0,-1){12}} \put(174,70){\line(0,-1){62}}
\mput(132,46)(0,-2){23}{\line(0,-1){1}}
\put(144,16){\line(1,0){46}}
 \put(174,8){\line(1,0){16}}

\thinlines\mput(132,1.8)(0,1.5){38}{\line(1,0){58}}
\put(182,70){\line(0,-1){70}} \put(120,64){\line(1,0){70}}

\put(123,62){\makebox{$Y_1$}} \put(135,62){\makebox{$W$}}
\put(160,62){\makebox{$Q$}} \put(137,50){\makebox{$Y_2$}}
\put(160,50){\makebox{$Q$}} \put(144,41){\makebox{$X_{11}$}}
\put(159,41){\makebox{$X_{12}$}} \put(159,26){\makebox{$X_{22}$}}

\put(40,46.5){\makebox{$\overbrace{\hskip 9mm }$}}
 \put(0,14){\makebox{$\underbrace{\hskip 24mm}$}}
\put(52,46.5){\makebox{$\overbrace{\hskip 9mm}$}}
\put(70,14){\makebox{$\underbrace{\hskip 19mm}$}}
\put(64,46.5){\makebox{$\overbrace{\hskip 24mm}$}}
\put(13,6){\makebox{$X$}} \put(52,9){\makebox{$D_1^1$}}
\put(58,9){\makebox{$D_2^1$}} \put(44,50){\makebox{$D_1$}}
\put(56,50){\makebox{$D_2$}}\put(64,9){\makebox{$D_3^1$}}\put(80,6){\makebox{$B^1$}}
\put(77,50){\makebox{$B$}} \put(95,9){\makebox{$A_3$}}
\put(103,9){\makebox{$A_4$}} \put(150,-5){\makebox{$R$}}
\end{picture}}
\end{array}
\end{equation}
\bigskip

In case of edge reduction and $X_{22}\ne 0$, we remove the block-row and the block-column where $X_{11}$ sits
from $X^1$; and those where $Y_1$ ( according to the partition of $(T_2,\sim_2))$ sits from $R^1$; as well as
$\overline D_1,$ and the upper rows parallel to $X_{11}$ from $N^1$. Then we obtain 3 matrices shown in the
shadowed part of diagram (55) for example, which still satisfy the hypothesis given in (54), if we use $X_{22}$
instead of $X$, $B^1$ instead of $B$, $D_1^1, \cdots D_{i_1}^1$ partitioned after the edge reduction  instead of
$D_1, \cdots, D_i$ respectively; the block $X_{12}$ goes to part $Q^1$, and $X_{11}$ goes to part $Y^1$.  In
case of regularization, we remove the block-row and the block-column where $Y_1$ sits from $R^1$, as well as
$\overline D_1$, and we also obtain 3 matrices satisfying the hypothesis given in (54).

Note that the above procedure is an illustration in order to help to understand the induction. Thus by  a
sequence of edge reductions and regularizations, or possibly a loop reduction at last,  we finally reach case
(1) or (2), since $m$ is a finite number. Consequently, we obtain an induced triple of matrices $(N^r, X^r\times
R^r, \underline{n}_1^r\times \underline n_2^r)$ before $A_{i+1}$,
 (see the diagram (56) below for example),
such that $Q_l^r$ and the nilpotent blocks of $X^r$ satisfy the
conditions of the lemma.\hfill$\square$

\medskip
{\bf Lemma 8.4.2} With the notations above. Let $N^r=(D^r, B^r,
A^r)$ be obtained from $N=(D,B,A)$ after all the reductions given
in Lemma 8.4.1. Then $Q^r$ and the nilpotent blocks of $X^r$ force
the rest blocks of $B^r$ equal to $\emptyset$ in any sequence of
further reductions under the order of Definition 2.3.1.

\vspace{12mm}

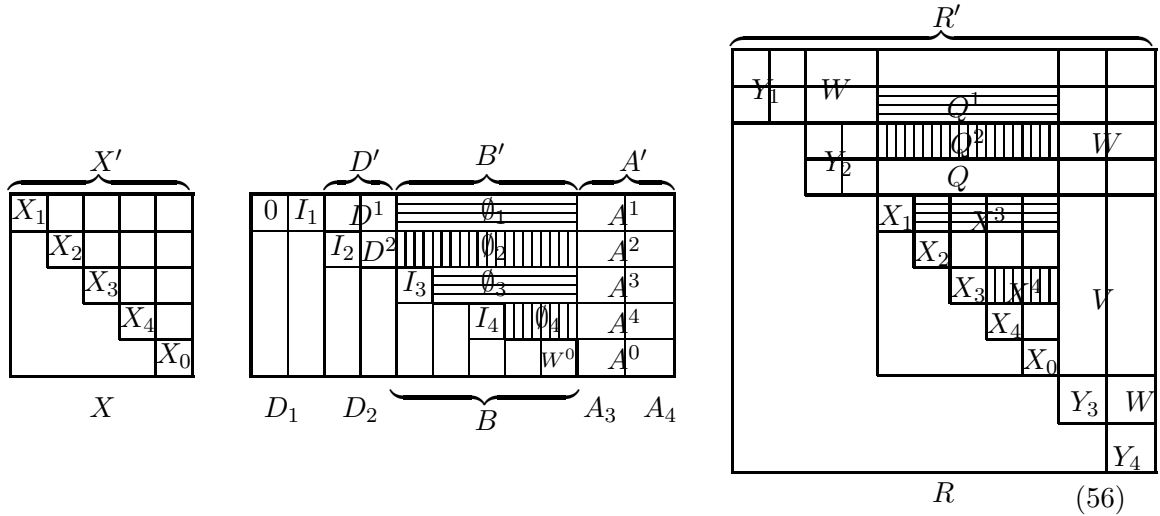
\begin{equation}\begin{array}{c} {\unitlength=0.8mm
\begin{picture}(185,60)\thicklines
\put(0,16){\framebox(30,30){}} \put(0,40){\line(1,0){30}} \put(6,34){\line(1,0){24}} \put(12,28){\line(1,0){18}}
\put(18,22){\line(1,0){12}}

 \put(6,34){\line(0,1){12}} \put(12,28){\line(0,1){18}}
\put(18,22){\line(0,1){24}} \put(24,16){\line(0,1){30}}

\put(0,42){\makebox{$X_1$}} \put(6,36){\makebox{$X_2$}}
\put(12,30){\makebox{$X_3$}} \put(18,24){\makebox{$X_4$}}
\put(24,18){\makebox{$X_0$}}

\thicklines \put(40,16){\framebox(70,30){}} \put(52,40){\line(0,1){6}} \put(58,34){\line(0,1){6}}
\put(70,28){\line(0,1){6}} \put(82,22){\line(0,1){6}} \put(94,16){\line(0,1){6}} \put(52,40){\line(1,0){6}}
\put(58,34){\line(1,0){12}}\put(70,28){\line(1,0){12}} \put(82,22){\line(1,0){12}}

\thinlines \put(52,16){\line(0,1){30}} \put(64,16){\line(0,1){30}}
\put(94,16){\line(0,1){30}} \put(58,16){\line(0,1){30}}
\put(46,16){\line(0,1){30}}

\put(40,40){\line(1,0){70}} \put(52,34){\line(1,0){58}}
\put(64,28){\line(1,0){46}} \put(76,22){\line(1,0){34}}

\put(70,16){\line(0,1){18}} \put(76,16){\line(0,1){12}}
\put(82,16){\line(0,1){12}} \put(88,16){\line(0,1){6}}

\put(47,42){\makebox{$I_1$}}\put(42,42){\makebox{$0$}}
\put(53,36){\makebox{$I_2$}} \put(65,30){\makebox{$I_3$}}
\put(77,24){\makebox{$I_4$}}
\put(88,17.5){\makebox{$\scriptstyle{W^0}$}}
\put(56,41){\makebox{$D^1$}} \put(58,35){\makebox{$D^2$}}
\put(78,42){\makebox{$\emptyset_1$}}
\put(78,36){\makebox{$\emptyset_2$}}
\put(78,30){\makebox{$\emptyset_3$}}
\put(87,24){\makebox{$\emptyset_4$}} \put(99,41){\makebox{$A^1$}}
\put(99,35){\makebox{$A^2$}} \put(99,29){\makebox{$A^3$}}
\put(99,23){\makebox{$A^4$}} \put(99,17){\makebox{$A^0$}}

\mput(64,41.5)(0,1.5){4}{\line(1,0){30}}
\mput(65.5,40)(1.5,0){20}{\line(0,-1){6}}
\mput(70,29.5)(0,1.5){3}{\line(1,0){24}}
\mput(83.5,28)(1.5,0){8}{\line(0,-1){6}}
\put(102,16){\line(0,1){30}}

\thicklines \put(120,0){\framebox(70,70){}} \mput(120,58)(0,6){2}{\line(1,0){70}} \put(132,52){\line(1,0){58}}
\put(132,46){\line(1,0){58}} \put(144,16){\line(1,0){46}} \put(174,8){\line(1,0){16}}
\put(132,70){\line(0,-1){24}} \put(144,70){\line(0,-1){54}} \put(174,70){\line(0,-1){62}}
\put(182,70){\line(0,-1){70}}

\put(144,40){\line(1,0){30}} \put(150,34){\line(1,0){24}}
\put(156,28){\line(1,0){18}} \put(162,22){\line(1,0){12}}

\put(150,34){\line(0,1){12}} \put(156,28){\line(0,1){18}}
\put(162,22){\line(0,1){24}} \put(168,16){\line(0,1){30}}

\thinlines\mput(144,59.5)(0,1.5){4}{\line(1,0){30}}
\mput(145.5,52)(1.5,0){20}{\line(0,1){6}}
\mput(150,41.5)(0,1.5){4}{\line(1,0){24}}
\mput(163.5,28)(1.5,0){7}{\line(0,1){6}}

\put(126,70){\line(0,-1){12}} \put(138,58){\line(0,-1){12}}

\put(144,41){\makebox{$X_1$}} \put(150,35){\makebox{$X_2$}}
\put(156,29){\makebox{$X_3$}} \put(162,23){\makebox{$X_4$}}
\put(168,17){\makebox{$X_0$}} \put(159,40){\makebox{$X^3$}}
\put(165,28){\makebox{$X^4$}}

\put(123,62){\makebox{$Y_1$}} \put(135,50){\makebox{$Y_2$}}
\put(176,10){\makebox{$Y_3$}} \put(183,1){\makebox{$Y_4$}}
\put(155.5,59){\makebox{$Q^1$}} \put(155.5,47.5){\makebox{$Q$}}
\put(156,53){\makebox{$Q^2$}} \put(134.5,62){\makebox{$W$}}
\put(179.5,53){\makebox{$W$}} \put(179.5,27){\makebox{$V$}}
\put(185,10){\makebox{$W$}}

\put(63,14){\makebox{$\underbrace{\hskip 25mm}$}}
\put(77,7){\makebox{$B$}} \put(42,9){\makebox{$D_1$}}
\put(55,9){\makebox{$D_2$}} \put(95,9){\makebox{$A_3$}}
\put(105,9){\makebox{$A_4$}} \put(13,9){\makebox{$X$}}
\put(153,-5){\makebox{$R$}}
\put(-0.5,46.5){\makebox{$\overbrace{\hskip 24mm}$}}
\put(52,46.5){\makebox{$\overbrace{\hskip 9mm}$}}
\put(64,46.5){\makebox{$\overbrace{\hskip 24mm}$}}
\put(94.5,46.5){\makebox{$\overbrace{\hskip 12.5mm}$}}
\put(119.5,70.5){\makebox{$\overbrace{\hskip 56mm}$}}

\put(13,50){\makebox{$X'$}} \put(56,50){\makebox{$D'$}}
\put(77,51){\makebox{$B'$}} \put(101,50){\makebox{$A'$}}
\put(153,74){\makebox{$R'$}}
\end{picture}}
\end{array}
\end{equation}

\bigskip

{\bf Proof.}\ \ For the sake of simplicity, we denote $X^r, N^r,
R^r$ by $X', N', R'$ respectively.

We first prove case (2) of Lemma 8.4.1 according to the partition
of
 $$ (T^r_1\times T^r_2, \, \sim^r_1\times \sim_2^r).$$
Assume that
  $I_1, I_2, \cdots, I_{\beta}$ are identity matrices appearing in
  t he preceding edge reductions. The assertion will be shown by
  induction on $\beta$. $\beta=0$ is trivial since $B'$ is then
  obtained by a loop reduction and there is no block left in $B'$.
  Denote by $X_1,X_2,\cdots, X_\beta$ the diagonal blocks of $X'$
  parallel to $I_1, I_2, \cdots, I_{\beta}$
  and $X_0$ that
at the lower right corner of $X'$ parallel to the Weyr matrix
$W^0$. If $I_1, \cdots, I_{\alpha}$ are sitting in part $D'$, and
$I_{\alpha+1}, \cdots, I_{\beta}$ are sitting in $B'$, then we
denote by $D^1, \cdots, D^{\alpha}$ the blocks of $D'$ parallel to
$I_1, I_2, \cdots, I_{\alpha}$, which remain in $N'$ after the
reductions in Lemma 8.4.1;   by $A^1, \cdots, A^{\beta}, A^0$ the
blocks of $A'$ parallel to $I_1, \cdots, I_{\beta}, W^0$ which
remain in $N'$ after the reductions in Lemma 8.4.1 respectively.
See the middle matrix of diagram (56) for example, where $\beta=4,
\alpha=2$.   On the other hand, we denote by $Q^1, \cdots,
Q^{\alpha}$ the row-blocks in part $Q$ parallel to $I_1,\cdots,
I_{\alpha}$ of $(N')^T$, the transpose of $N'$, and by
$X^{\alpha+1}, \cdots, X^{\beta}$ those of $X'$ parallel to
$I_{\alpha+1 }, \cdots, I_{\beta}$ of $(N')^T$ respectively in
$R'$. See the shadowed blocks $Q^1, Q^2, X^3, X^4 $ of diagram
(56) for example. Denote $Y',V',W', Q'$ of $R'$ still by $Y,V,W,Q$
respectively for simplicity. It is obvious that $X^l$ are higher
than $X_l$ in $X'$.

The differentials in part $A^0$ involve only the entries of $X_0,
Y, V, W $.
 They do not involve $Q$ and the nilpotent part of $X'$
at all.
 Therefore those blocks are still linearly independent of all the other blocks.
  The differential of
the block of $B'$ parallel to $I_{\beta}$ has a term
$I_{\beta}X^{\beta}$, which leads to $\emptyset_{\beta}$ by
regularization (see  $X^4$ and $\emptyset_4$ in (56) for example).
 Suppose that we have already $\emptyset_{\beta}, \emptyset_{\beta-1},
\cdots, \emptyset_{l+1}$ parallel to $I_{\beta}, I_{\beta-1}, \cdots, I_{l+1}$ in $B'$ given by regularizations
involving $X^{\beta}, X^{\beta-1}, \cdots,X^{l+1}$ for some $l\ge \alpha$.
 Now we consider the block parallel to $I_l$ in
$B'$.
 Since the differentials in parts $$A^0;\ \emptyset_{\beta},
 A^{\beta}; \ \emptyset_{\beta-1},
A^{\beta-1}; \ \cdots\ ;\  \emptyset_{l+1}, A^{l+1}$$ involve only the entries of $Y,V,W$, and some blocks
parallel to or lower than $X_{l+1}$ in $X'$;  on the other hand, $X^l$ is higher than $X_{l+1}$ (see diagram
(56) for example), we conclude that the entries of $X^l$ are still linearly independent  after the reduction for
$ A^{l+1}$. The differential of the block of $B'$ parallel to $I_l$ has a term $I_lX^l$, which leads to
$\emptyset_l$ by regularization (see $X^3$ and $\emptyset_3$ in (55) for example). Then we obtain a sequence of
$\emptyset$'s up to the $(\alpha+1)$-th block row in $B'$ by induction. Moreover all the reductions before
$D^{\alpha}$ involve only $X', Y, V, W$, and the entries of  $Q$ remain linearly independent.

Next the differential in parts $A^{\alpha+1}$ and $D^{\alpha}$ involve the entries of $X', Y, V, W$. The
differential of the row-block of $B'$ parallel to $I_{\alpha}$ has a term $I_{\alpha}Q^{\alpha}$, which leads to
$\emptyset_{\alpha}$ by regularization (see $Q^2$ and $\emptyset_2$ in (56) for example). Suppose we already
have $\emptyset_{\alpha}, \emptyset_{\alpha-1}, \cdots, \emptyset_{l+1}$ parallel to $I_{\alpha}, I_{\alpha-1},
\cdots, I_{l+1}$ given by regularizations involving $Q^{\alpha}, Q^{\alpha-1}, \cdots, Q^{l+1}$ for some integer
$l>0$. Now we consider the block parallel to $I^l$ in $B'$. Since the differentials for $$A^{\alpha+1},
D^{\alpha};\ \emptyset_{\alpha}, A^{\alpha}, D^{\alpha-1};\ \cdots\ ;\  \emptyset_{l+1}, A^{l+1},D^l$$ involve
only the entries of $X', Y, V, W $, and some entries of $Q$ lower than $Q^l$, we conclude that the entries of
$Q^l$ are still linearly independent after  the reduction for $D^l$. The differential of the block of $B'$
parallel to $I_l$ has a term $I_lQ^l$, which leads to $\emptyset_l$ by regularization (see $Q^1$ and
$\emptyset_1$ in (56) for example). Thus we obtain a sequence of $\emptyset$'s up to the  block parallel to
$I_1$ in $B'$ by induction.

Therefore the rest blocks of $B'$ are all equal to $\emptyset$ in any sequence of reductions under the order
given by Definition 2.3.1 according to the partition $(T^r_1 \times T^r_2,{\sim}^r_1 \times {\sim}^r_2)$.

The proof for the case (1) of Lemma 8.4.1 is as follows: if the $(r+1)$-th reduction is for a block in
$A_{i+1}$, we start the procedure from $A^\beta$, where $\beta\ne \alpha$; if the $(r+1)$-th is for a block in
$B'$, we start the procedure from $\emptyset_\beta$, where $\beta\ne\alpha$ or  $\beta=\alpha$; if the
$(r+1)$-th reduction is for a block before $B'$, we start the procedure from $D^\beta$, where$\beta=\alpha$.

The lemma still holds, if we use the usual partition step by step
in the reduction sequence. In fact we only need to make some
refinement for the block-rows.\hfill $\square$

{\bf Proposition 8.4.1} Let $(\K_1\times\K_2, \overline{\M},\mathcal{V})$ and $\mathfrak{B}$ be given in diagram
(46) of 8.2. Suppose that the partial bocs $\mathfrak{B}$ is wild and the induced local partial bocs
$\mathfrak{B}_{\P}$ is in case P2 of 8.2; $\delta(\overline{b})$ satisfies Formula (51) of 8.3; and the partial
differentials of $a$'s after $\overline{b}$ satisfy Formula (53) of  8.3. Given any size vector $\n_1\times\n_2$
we have a matrix equation:
$$
X(\cdots A \cdots B\cdots A \cdots B \cdots) =(V_1\ V_2\ \cdots\ V_n)+(\cdots A\cdots B\cdots A\cdots
B\cdots)R_0$$ Then the reductions given by this  equation are completely determined by the reductions given by
Formula (54), and $A_l$ goes to $\emptyset$ for any $l \in I$ in Formula (50) of 8.2, $B_l$ goes to $\emptyset$
for any $l \ne j$ in Formula (47) and (48) of 8.2.

Moreover, assume that $\mathfrak{B}$ is a partial bocs of $\mathfrak{A}$ appearing in a sequence $(*)$ with the
end term in case of MW5. Then block $\overline{B}$ corresponding to $\overline{b}$ will not contain any
parameter after any reductions.

 {\bf Proof.} (1) Consider the differentials of $b_l$
for any $l\ne j$ given by Formula (47) and (48) of 8.2 with $m'=m$, $h_{ll}(\lambda)$ being non-zero constants
ensure that the reduction for any block in  $B_l$ (corresponding to $b_l$), must be a regularization determined
by $U_l$ (corresponding to $u_l$), which does not influence the reductions of any other blocks in $(\cdots
A\cdots B\cdots A\cdots B\cdots)$.

(2) Consider the differentials of $a_l$, $l\in I$, given by Formula (49) of 8.3, and the partial differentials
of $a_{l'}$ given by Formula (53). Since $$\{v_l\mid l\in I\}\ \ {\rm and}\ \ \{v_{n_0+1}, \cdots, v_n\}$$ are
linearly independent, any reduction (in particular, regularization) of any part of block $A_{l'}$ for $l'>n_0$
(see the beginning of 8.3) does not influence the regularization of any part of block $A_l$ for $l\in I$. On the
other hand, the regularization of any block of $A_l$   does not influence the reductions of any other blocks in
$(\cdots A\cdots B\cdots A\cdots B\cdots)$.

(1) and (2) tell us that the reductions given by the matrix
equation in the proposition are completely determined by those
given by Formula (54).

(3) Finally we prove the last assertion of the proposition. If we
are in case (1) of Lemma 8.4.1,   then\ there\ is\ no\  loop\
reduction\
 (even\ no\ loop\ and\ edge\ reductions)\ in
 $\overline B$,
the assertion holds obviously. Otherwise the only possible loop reduction inside $\overline{B}$ yields a Weyr
matrix $W^0$ in case (2) of Lemma 8.4.1.  We claim that the induced triple $(N^{r-1},X^{r-1}\times R^{r-1},
\underline n_1^{r-1}\times\underline n_2^{r-1})$  before $W_0$ appearing is not local. In fact, if there exist
some diagonal blocks $X_i$ of $X^{r-1}$ with $i>0$, then $X_0$ is independent of all the $X_i$'s. Otherwise,
$X_0=X^{r-1}$, there must exist some $Y_j$ as a diagonal block of $R^{r-1}$ by wildness of $\mathfrak{B}$ and
tameness of $\mathfrak{B}_{\P}$. Then $X_0$ is independent of all the $Y_j$'s. In both cases the induced partial
bocs at this step is not local. We conclude that $\overline{B}$ will not contain any parameter after any
reductions, since MW5 requires that any parameter appears for the first time in a local bocs by Corollary 5.6.1.
The proof is finished. \hfill$\square$

\subsection{Some lemmas in algebra}

\kg This subsection is devoted to giving several lemmas in algebra
which are needed  in Section 9.

Let $$\{x_1, \cdots, x_m\}\ \ {\rm and}\ \ \{y_1,\cdots, y_m\}$$  be two sets of indeterminates. We say that
$$a_1x_1+a_2x_2+\cdots+a_mx_m$$ and
$$a_1y_1+a_2y_2+\cdots+a_my_m$$ are the same linear forms, where $a_1,\cdots, a_m\in k$ are  constants.

{\bf Lemma 8.5.1}\ Let $(\K_1\times \K_2, \M, {\cal V})$ be a matrix problem, and $(N, R_1\times R_2,
\n_1\times\n_2)$ be a triple. Suppose that ${\cal N}$ is a matrix with constant entries of size
$\n_1\times\n_2$. If we fix a row-index $\alpha$ and an index $l\in \mathcal{L}\in T_2/\!\!\sim_2$, then any
linear forms at the $(\alpha,\beta)$-entry of $V+{\mathcal{N}}R_2$ in the indeterminates
$$(\overline{\xi_{\I\L}^{v_{123}}})^*,\ (\overline{\xi_{\I\L}^{v_{ij}}})^*,\ (i,j)=\{(12), (23), (31)\}\ {\rm
and}\ (\overline{\xi_{\I\L}^{v_{i}}})^*,\ i=1,2,3,\ \forall\, \I\in T/\!\sim$$ given in 8.1 at $\beta$-th column
of $V$ and $R_2$ are the same for all $\beta \in l$, (see the figure below). \hfill$\square$

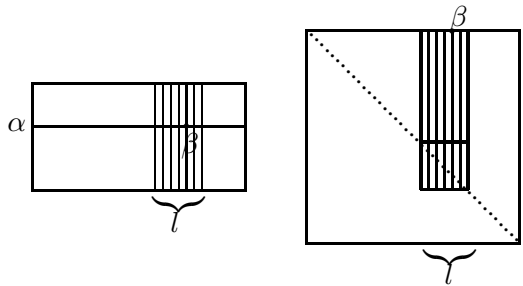
\begin{center}\unitlength=0.7mm

\begin{picture}(110, 40)
\thicklines\put(8,10){\framebox(40,20)} \put(8,22){\line(1,0){40}}
 \put(37, 22){\circle*{1}} \put(36, 17){$\beta$} \put(3,
21){$\alpha$}\put(30.5, 9){$\underbrace{\hskip 3mm}$}
\put(34,2){$l$}

\thinlines\mput(31,10)(1.5,0){7}{\line(0,1){20}}

\thicklines\put(60,0){\framebox(40,40)}
\thicklines\mput(81.5,40)(1.5,0){7}{\line(0,-1){30}}
\put(81.5,19){\line(1,0){9}} \put(81.5,10){\line(1,0){9}}

 \put(87.5,40){\circle*{1}} \put(87,41){$\beta$} \put(82,
-1){$\underbrace{\hskip 3mm}$} \put(86,-8){$l$}

\mput(99,-1)(-1,1){40}{$\cdot$}

\end{picture}
\end{center}
For example, in the matrix
$$
\left(\begin{array}{ccc} v_{11} & v_{12} & v_{13} \\ v_{21} &
v_{22} & v_{23} \end{array}\right) + \left(\begin{array}{ccc} 0 &
1 & \lambda^0 \\ 1 & 0 & 0\end{array}\right)
\left(\begin{array}{ccc} x_{11} & x_{12} & x_{13} \\ x_{21} &
x_{22} & x_{23}\\ x_{31} & x_{32} & x_{33} \end{array}\right)
$$
 the linear forms at the first row are $$v_{11}+ x_{21}+\lambda^0
x_{31},\ \ v_{12}+ x_{22}+\lambda^0 x_{32},\ \ v_{13}+ x_{23}+\lambda^0 x_{33}.$$

{\bf Lemma 8.5.2}\  Let $Q: \P\stackrel{c}{\longrightarrow} \Q$ be
a quiver,  and the diagram
$$(e)\qquad  \begin{CD} 0 @>>>\rK @>\iota_{\P}>>\rK^2
@>\pi_{\P}>>\rK @>>>0\\ && @V\alpha VV @VV\beta V @VV\alpha V \\
0 @>>>\rK^m @>>\iota_{\Q}>\rK^{2m} @>>\pi_{\Q}>\rK^m @>>>0
\end{CD}
$$ be an exact sequence over the path algebra $\rK Q$.

(1) If $\alpha=(0\ \cdots\  0 \ 1)$, then
$\beta=\left(\begin{array}{ccccc} 0&\cdots& 0&1&0\\ 0&\cdots&0&0&1
\end{array}\right)_{2\times 2m} $, and $\iota_{\P}=(0\ 1)$, $\pi_{\P}={1\choose 0}$
 under some suitable basis.

(2) If $\alpha=(0\ \cdots\  0\ 0)$, then $\beta=0$ or
$\beta=\left(\begin{array}{ccccc} 0&\cdots&0&0&1\\ 0&\cdots&0&0&0
\end{array}\right)_{2\times 2m}$, and $\iota_{\P}=(0\ 1)$, $\pi_{\P}=
{1\choose 0}$ under some suitable basis. \hfill$\square$

{\bf Lemma 8.5.3}\  Let $Q=\hspace{1cm} {\unitlength 1mm
 \begin{picture}(25, 4) \put(0.00,0.00){\circle*{1.00}}
\put(20.00,0.00){\circle*{1.00}} \put(1,0){\vector(1,0){18}}
 \put(0.00,1.50){\makebox{$\P$}}
\put(20.00,1.50){\makebox{$\Q$}} \put(8.50,1.50){\makebox{$c$}}
\put(-5.00,0.00){\oval(5,5)[t]} \put(-5.00,0.00){\oval(5,5)[bl]}
\put(-6.00,-3.00){\vector(3,1){3}} \put(-10,0){$\lambda$}
\end{picture}}$ be a quiver, and $\rK Q$ be the corresponding path algebra.
Let $\mathscr{C}$ be a full subcategory of $\mod kQ$ consisting of
modules $M$, such that the eigenvalue of $M(\lambda)$ is
$\lambda^0$.
\begin{center} \unitlength 1mm
\begin{picture}(80, 45)
\put(0,0){$\rK^{2m}$} \put(10, 1){\vector(1,0){60}}
\put(71,0){$\rK^{m}$} \put(3,15){$\rK^2$} \put(10,
16){\vector(1,0){60}} \put(73,15){$\rK$}
\put(4,14){\vector(0,-1){10}} \put(74,14){\vector(0,-1){10}}

\put(36,25){$\rK^{2m}$} \put(39, 40){$\rK^2$}
\put(35.00,43.00){\oval(5,5)[t]} \put(35.0,43.00){\oval(5,5)[bl]}
\put(33.50,40.00){\vector(3,1){3}} \put(24,44){$L(\lambda)$}

\put(-1.00,18.00){\oval(5,5)[t]} \put(-1.0,18.00){\oval(5,5)[bl]}
\put(-2.50,15.00){\vector(3,1){3}} \put(-13,19){$E(\lambda)$}

\put(79.00,18.00){\oval(5,5)[t]} \put(79.0,18.00){\oval(5,5)[br]}
\put(80.5,15.00){\vector(-3,1){3}} \put(82,19){$S(\lambda)$}

\put(40,39){\vector(0,-1){10}}

\put(41,39){\vector(3,-2){30}} \put(41,24){\vector(3,-2){30}}
\put(39,39){\vector(-3,-2){30}}

\put(39,24){\vector(-3,-2){31}}

\put(35, -4){$\pi_\Q={{I_m}\choose 0 }$} \put(34,
11){$\pi_\P={1\choose 0}$} \put(75, 8){$0$} \put(-4, 8){$E(c)$}
\put(55,32){${1\choose 0}$} \put(56, 6){${I_m\choose 0}$}
\put(20,30){$\varphi'_{\P}$} \put(19,8){$\varphi'_{\Q}$}
\put(41,33){$0$}

\end{picture}
\end{center}
 The above diagram in $\mathscr{C}$
consists of the modules $S$, $E$, $L$, where $$S_{\P}=\rK,\  S_{\Q}=\rK^m,\ S(\lambda)=(\lambda^0), \ S(c)=0;$$
$$E_{\P}=\rK^2, E_{\Q}=\rK^{2m},  \ E(\lambda)=\lambda^0 I_2\ \ {\rm or} \ \ J_2(\lambda^0),  \mbox{ and }\ E(c)
=\left(\begin{array}{cccc} 0&\cdots&0&1\\
0&\cdots&0&0\end{array}\right);$$ $$L_{\P}=\rK^2,\  L_{\Q}=\rK^{2m}, \ L(\lambda)=J_2(\lambda^0), \ L(c)=0;$$
and morphisms
$$
\pi=\left({1\choose 0}, {I_m \choose 0}\right): E\rightarrow S, \
\varphi=\left({1\choose 0}, {I_m \choose 0}\right): L\rightarrow S
$$ in $\rK Q$-mod.  Then

(1) the morphism $\varphi$ is not a split epimorphism.

(2) there does not exist any morphism $\varphi': L\rightarrow M$,
such that $\varphi'\pi=\varphi$.

 {\bf Proof.} (1) is obvious by  considering ${1\choose 0}:
 J_2(\lambda)\rightarrow J_1(\lambda)$.

(2) Suppose we have the contrary. Then
$\varphi'_{\P}\pi_{\P}=\varphi_{\P}$ gives
$\varphi'_{\P}=\left(\begin{array}{cc} 1&b\\ 0&d
\end{array}\right)$.
On the other hand
 $$0=L(c)\varphi'_Q=\varphi'_{\P} E(c)=\left(\begin{array}{cc}
1&b\\ 0&d \end{array}\right)\left(\begin{array}{cccc}
0&\cdots&0&1\\ 0&\cdots&0&0
\end{array}\right)=\left(\begin{array}{cccc} 0&\cdots&0&1\\
0&\cdots&0&0 \end{array}\right),$$
 which means a contradiction.
\hfill$\square$

\subsection{Some lemmas in category theory}

\kg {\bf Lemma 8.6.1}\ Let $Q_0=$
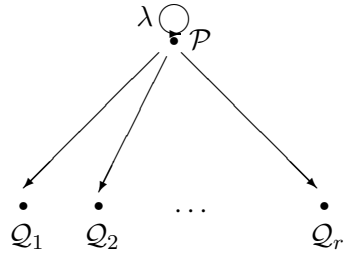
\begin{center}
\begin{picture}(80,40)
\put(32,42){\circle{4}} \put(32.3,40){\vector(1,0){0.5}}
\put(32,39){\circle*{1}}

\put(30,37.5){\vector(-1,-1){18}} \put(31,37){\vector(-1,-2){9}}
\put(33,37.5){\vector(1,-1){18}}
 \put(12,17){\circle*{1}}  \put(52,17){\circle*{1}}
\put(22,17){\circle*{1}} \put(10,12){$\Q_1$} \put(20,12){$\Q_2$}
\put(50,12){$\Q_r$}
   \put(34,38){$\P$} \put(32,15.5){\makebox{$\cdots$}}
\put(27,41){\makebox{$\lambda$}}

\end{picture}
\end{center}\vskip -1cm
be a quiver consisting of vertices $\P, \Q_1, \cdots, \Q_r$; a
loop $\lambda$ at $\P$, and $n_j$ arrows $a_{j1}, \cdots,
a_{jn_j}$ from $\P$ to $\Q_j$. Let $\Gamma_0=\rK Q_0$ be the path
algebra. We define a full subcategory $\mathscr{C}_0$ of
$\Gamma_0$-mod consisting of modules $M$, such that
$M_{\P}=\rK^m$, $M_{\Q_j}=\rK^{mn_j}$; $M(\lambda)=W$, which is a
Weyr matrix of eigenvalue $\lambda^0$ of size $m$ under a fixed
basis $B$ of $M_{\P}$; and
$$
M(a_{jl})=(0\ \cdots \ 0\ \underbrace{I_m\ 0\ \cdots \
0}_l)_{1\times n_j}
$$
under a basis $B_{j1}, \cdots, B_{jn_j}$ of $M_{\Q_j}$. Then for
any $M,L\in\mathscr{C}_0$, a morphism  $\varphi=(\varphi_{\P}; \
\varphi_{\Q_1}, \cdots,\\ \varphi_{\Q_r}): M\rightarrow L$ is
given by
$$
\varphi_{\P}=Z, \quad \varphi_{\Q_j}=\left( \begin{array}{cccc} Z&
& & \\ & Z & &  \\ && \ddots & \\ &&& Z
\end{array}\right)_{n_j\times n_j}
$$
under the basis $B$ and $ B_{j1}, \cdots, B_{jn_j}$ for $j=1,
\cdots, r$, where $M(\lambda)Z=ZL(\lambda)$. \hfill$\square$

\medskip

{\bf Lemma 8.6.2} \ Let $Q$ be a quiver obtained from $Q_0$ by adding one arrow $c: \P\rightarrow \Q$, where
either $${\rm (i)}\ \Q\notin \{\Q_1, \cdots , \Q_r \},\ \ {\rm or\  (ii)}\ \Q=\Q_{j_0}\ {\rm for\  some}\ 1\leq
j_0\leq r,$$ and $\Gamma=kQ$. We define a full subcategory $\mathscr{C}$ of $\Gamma$-mod as follows.

 Let $M_{\P}=\rK^m$ have a $\rK$-basis $B$, and
$M_{\Q_j}=\rK^{m\cdot n_j}$ have a $\rK$-basis $(B_{j1}, \cdots,
B_{jn_j})$ given by Lemma 8.6.1, where $1\leq j\leq r$ in case
(i), or
 $j\neq j_0$ in case (ii). And let  $M_{\Q}=\rK^{n}$ in case
 (i), $M_{\Q_{j_0}}=\rK^{n}\oplus\rK^{m\cdot n_{j_0}} $ in case (ii),
 where $\rK^{m\cdot n_{j_0}}$ has a $\rK$-basis $(B_{j_01}, \cdots,
 B_{j_0n_{j_0}})$. Moreover $M(\lambda)=W$ of eigenvalue $\lambda^0$,
  and $M(a_{jl})$ are given by Lemma
 8.6.1, except
 $$
M(a_{j_0l})=(0\ |\ 0\  \cdots \ 0\ \underbrace{I_m\ 0\ \cdots \
0}_l)_{1\times (n_{j_0}+1)}
$$
in case (ii), where $0$ before the vertical line is an $m\times n$ matrix. And $$M(c)=P\ \text{ in case (i),} \
\ M(c)=(P\mid 0\cdots 0)_{1\times (n_{j_0}+1)}\ \text{ in case (ii).}$$ Then

(1) $\forall M,L\in \mathscr{C}$, a morphism $$\varphi=(\varphi_{\P}; \varphi_1, \cdots, \varphi_{r},
\varphi_{\Q})\ \ {\rm or}\ \ (\varphi_{\P}, \varphi_1, \cdots, \varphi_{r})$$ from $M$ to $L$ under the fixed
basis is given by
$$
\varphi_{\P}=Z \mbox{ satisfying } M(\lambda)Z=ZL(\lambda);\
\varphi_{\Q_l}=\left(\begin{array}{cccc} Z&&&\\ &Z&&\\ &&\ddots&\\
&&&Z
\end{array}\right)$$
 where $1\leqslant l\leqslant r$ in case (i), $l\neq j_0$ in case (ii);  and
  $$\varphi_{\Q}=U
\mbox{ in case(i)};\   \varphi_{j_0}=\left(\begin{array}{c|ccc} U&*&\cdots&*\\ \hline &Z&&\\
&&\ddots&\\&&&Z
\end{array}\right) \mbox{ in case (ii)},$$ where $L(c)=P$ and $P U=ZP'$ in case (i),
 $L(c)=(P'\mid 0\cdots 0)$  with $P U=ZP'$ and $P*=0$ in case (ii).

(2) $\mathscr{C}$ is an exact category having the exact structure inherited from $\Gamma$-mod.

{\bf{Proof.}} (1) is obvious.

(2)  We only prove (E2) of 6.1. If $(e) \quad 0\longrightarrow M
\stackrel{\iota}{\longrightarrow}E\stackrel{\pi}{\longrightarrow}L\longrightarrow 0 $ is a conflation of
$\mathscr{C}$, and $\varphi: L'\rightarrow L$ is a morphism of $\mathscr{C}$,  then a pullback
$$
\begin{CD}
0 @>>> M @>\iota'>> E' @>\pi'>> L' @>>> 0\\
&& @VidVV @VV\varphi'V @VV\varphi V\\
0 @>>> M @>>\iota> E @>>\pi> L @>>> 0
\end{CD}
$$
of $\Gamma$-mod gives the desired deflation $\pi'$ of
$\mathscr{C}$. In fact the eigenvalue of $E'(\lambda)$ must be
$\lambda^0$, since $E'$ is a submodule of $E\oplus L'$. \hfill
$\square$

\medskip
{\bf Lemma 8.6.3} Let  $Q=\hspace{1cm} {\unitlength 1mm
 \begin{picture}(25, 4) \put(0.00,0.00){\circle*{1.00}}
\put(20.00,0.00){\circle*{1.00}} \put(1,0){\vector(1,0){18}}
 \put(0.00,1.50){\makebox{$\P$}}
\put(20.00,1.50){\makebox{$\Q$}} \put(8.50,1.50){\makebox{$c$}}
\put(-5.00,0.00){\oval(5,5)[t]} \put(-5.00,0.00){\oval(5,5)[bl]}
\put(-6.00,-3.00){\vector(3,1){3}} \put(-10,0){$\lambda$}
\end{picture}}$ be a quiver,  $\Gamma=kQ$ be the
 path algebra, and
 $$(e)\quad 0\rightarrow
S\rightarrow E\rightarrow S\rightarrow 0$$
  be the exact sequence
given by Lemma 7.6.2. We define a full subcategory $\mathscr{D}$
of $\Gamma$-mod consisting of modules $M$, such that
$M_{\P}=\rK^{m_{\P}}$; $M_Q=\rK^{m_{\Q}}$; and the eigenvalue of
$M(\lambda)$ is a constant $\lambda^0$. Then

 (1) $\mathscr{D}$ is an exact category having the
 exact structure inherited from $\Gamma$-mod.

(2) (e) is not an almost split sequence in $\mathscr{D}$.

{\bf Proof.} (1) is a special case  of Lemma 8.6.2.

(2) The almost split sequences in $\mathscr{D}$ has also the property of item (2) of Lemma 7.6.1. In fact if
$$0\longrightarrow M\stackrel{\iota}{\longrightarrow} N\stackrel{\pi}{\longrightarrow} L\longrightarrow 0$$
 is an
almost split conflation of $\mathscr{D}$, $K$ is a non-zero submodule of $M$, then $K$ and $ M/K$ are contained
in $\mathscr{D}$. The following diagram of $\Gamma$-mod
$$
\begin{CD}
0@>>> M@>\iota>> N @>\pi>> L @>>>0\\
& & @V\psi VV @VV\psi V @VVidV\\
0@>>> M/K@>>> N/\iota(K) @>>> L @>>>0
\end{CD}
$$
is still in $\mathscr{D}$. Thus the second sequence  splits in $\mathscr{D}$. Therefore Lemma 7.6.2 tells us
that (e) is not an almost split sequence of $\mathscr{D}$. The proof is completed. \hfill$\square$

\medskip
{\bf Lemma 8.6.4} Let $\mathscr{C}$ and $\mathscr{D}$ be two exact categories. Suppose  there exists a functor
$F: \mathscr{C}\rightarrow \mathscr{D}$, which is a representation equivalence (i.e. dense, full, and reflects
isomorphisms). If $(e)\quad M\stackrel{\iota}{\longrightarrow} N\stackrel{\pi}{\longrightarrow} L$ is an almost
split conflation of $\mathscr{C}$ such that $F(M)$ (resp. $F(L)$) is non-injective (resp. non-projective),
 then $F(e)$ is an almost split conflation in $\mathscr{D}$.

{\bf Proof.} Since $F$ reflects isomorphisms, $F(\pi)$ is not a
retraction. Suppose that $\varphi': K'\rightarrow F(L)$ is a
morphism in $\mathscr{D}$, which is not a split epimorphism. Then
there exists some morphism $\varphi: K\rightarrow L $ in
$\mathscr{C}$ such that $F(K)=K'$, $F(\varphi)=\varphi'$, since
$F$ is dense and full. And $\varphi$ is not a split epimorphism in
$\mathscr{C}$, otherwise $\varphi'$ would be so in $\mathscr{D}$.
Then there exists some $\psi: K\rightarrow L$ with
$\psi\pi=\varphi$, so $F(\psi)F(\pi)=\varphi'$,  thus $F(\pi)$ is
a right almost split morphism. Similarly, we can prove that
$F(\iota)$ ia a left almost split morphism.  And $F(e)$ is an
almost split conflation of $\mathscr{D}$ by the definition given
in 6.1. \hfill $\square$

\newpage

\begin{center}
\section{Added columns}
\end{center}

\kg This section is still devoted to presenting some technical preparations in order to prove our main theorem
in case MW5 of Theorem 5.6.1.

\bigskip

\subsection{The structure of added columns}

\kg Let $(\K,\M,H)$ be a bimodule problem with a reduction sequence $(*)$ of freely parameterized triples
starting from $(N_0, R, \n)$, such that the end term of the corresponding sequence $(***)$ is in the case of MW5
of Theorem 5.6.1. Then $\lambda$ appears for the first time in $\mathfrak{A}^s$ by corollary 5.6.1, and we have
Formula $(42)$ of 5.4. If $h(\lambda)\ne 0$ in the Formula (43), then a sequence of regularizations and a loop
reduction for $a_i$ gives a parameter $\lambda_{\gamma}=\nu$.

{\bf Definition 9.1.1} Denote the  matrix index of $\lambda$ in $H^{s}_{\n^{s}}$
 by $(\alpha(0),\beta(0))$  ; the
index of $\nu$ in $H^{s+i}_{\n^{s+i}}$ by $(\alpha^0, \beta^0)$. Suppose  that the parameter $\nu$ is contained
in a block $G$ partitioned by $(T,\sim)$,  the index of the block-column of $G$ is $l^0$, and $l^0\in\L^0\in
T/\!\sim$. \hfill $\square$

Let $(\K, \M/\M', d)$ be a partial bimodule problem, $\mathfrak{B}$ be the corresponding partial bocs. Suppose
that $(\K, \M/\M', d)$ determines a matrix problem $(\K_1\times\K_2, \overline{\M}, \mathcal{V})$, for example
given by Proposition 8.1.1 or 8.1.2.

{\bf Remark.} Let us take
 $$\displaystyle \mathcal{N}^r=\sum_{l=0}^{r-1}\overline{N^l_{p_lq_l}} \otimes \chi^l\
 {\rm for}\  1\le r\le s,$$
  where $(\chi^l)^*$ is the first arrow of $(\K_1^l\times\K_2^l,
  \overline{\M}^l, \mathcal{V}^l)$,
 $\overline{N^0_{p_0q_0}}=\overline{N_{pq}}$, $\chi^0=\chi$.
Now we fix a $k$-basis $B'$ for the matrix problem $(\K_1\times \K_2, \overline{\M}, \mathcal{V})$ given in 8.1,
and we regard the coefficient functions $(B')^*$ as the indeterminates. In particular,
$\overline{(\xi_{\I\J}^v)^*}$, $\overline{(\xi_{\I\J}^{v_{ij}})^*}$, $(ij)=(12), (23) $ or $(31)$,
$\overline{(\xi_{\I\J}^{v_i})^*}$, $i=1,2,3,$ and $\overline{(E_{\I})^*}$ are regarded as $n_{\I}\times n_{\J}$
and $n_{\I}\times n_{\I}$ matrices with entries being indeterminates (see 8.1). Fix an index $r$, $R_1^r, R_2^r,
V^r$ can be obtained from $R_1, R_2, V$ restricted by the matrix equation:
$$R_1\mathcal{N}^r\equiv_{\preceq(p_{r-1}, q_{r-1})}
V+\mathcal{N}^rR_2.$$ the equation system II$^{\prime}$ is given by the above matrix equation and the equation
$x^3_{p_rq_r}=0$ is given by (the $(p_r,q_r)$-block of $\mathcal{N}^rR_2-R_1\mathcal{N}^r)=0$.

\medskip

{\bf Definition 9.1.2} Consider the matrix equation $R_1N^r=V+N^rR_2$. Denote by $eq_{(\alpha, \beta)}$ the
$(\alpha, \beta)$-entry, and by $Eq_{(\alpha, \beta)}$ the set of entries indexed before $(\alpha, \beta)$ under
the matrix order of Definition 2.3.1. \hfill $\square$

\medskip

{\bf Lemma 9.1.1} Suppose that $\n^r_1\times\n_2^r$ has all the components to be $1$,  and $N_{p_rq_r}$ is the
first non-zero block of $N^r$ with the  matrix index $(\alpha_r, \beta_r)$. Then

(1) $eq_{(\alpha_r, \beta_r)}$ is the equation $x^3_{p_rq_r}=0$ in
 $(\K_1^r \times \K_2^r,
\overline{\M}^r, H^r)$; $Eq_{(\alpha_r, \beta_r)}$ consists of the equation system II$^{\prime}$ of Definition
8.1.1 and the condition on the diagonal elements determined by $(T^r_1\times T^r_2, \sim^r_1\times\sim^r_2)$.

(2) $eq_{(\alpha_r, \beta_r)}$ is a linear combination of some equations in $E_{q(\alpha^r, \beta_r)}$, if and
only if the equation $x^3_{p_rq_r}=0$ is a linear combination of some equations in II$^{\prime}$ of 8.1.1 in
$(\K_1^r\times\K_2^r, \overline{\M}^r, \mathcal{V}^r)$.

 Let $(\triangle)$ be a
reduction sequence of 8.1 starting from $( N, R_1\times R_2, \n_1\times \n_2)$, with $G$ being also a block
partitioned by $(T_1\times T_2, \sim_1\times\sim_2)$, such that
\begin{itemize}
\item[{\bf S01}]
  the end term $\mathfrak{B}^s$ of the corresponding
sequence of partial bocses is  in the case of MW5 of Theorem 5.6.1;
 \item[{\bf S02}] $l^0\in L^0\in T_2/\!\sim_2$, such that $\L^0\cap
 T_1=\emptyset$.
\end{itemize}
Note that since $\mathfrak{A}^s$, as well as $\mathfrak{B}^s$, are both local, we have that $\forall\, \I\in
T/\!\sim$, $\I\mid_{T_1}\ne \emptyset$ or $\I\mid_{T_2}\ne \emptyset$.

Under the assumptions S01 and S02, we define a new size vector
 $\tilde{n}$ of $(T,\sim)$, with
$$n_{\I}=n_{\I}\ \forall\ \I \in
 T/\!\sim \setminus \{\mathcal{L}^0\}
\ {\rm and}\ \tilde{n}_{\mathcal{L}^0}=n_{\mathcal{L}^0}+1.$$

Let $\wt{\mathcal{N}}^s$ be obtained from $\mathcal{N}^s$ by adding a column to the left of each $l$-th block
column consisting of $0$'s  for all $l\in \L^0$. Suppose $\L^0=\{l_1<l_2<\cdots<l_{\L^0}\}$, we  denote by
$$ o_1<o_2<\cdots<o_{l_{\L^0}}$$ the indices of the added columns.
Let $$\displaystyle \wt T=T\bigcup \{o_1,o_2,\cdots,o_{l_{\L^0}} \}.$$
 In fact, $\wt{\mathcal{N}}^s$
is just a direct sum of $\mathcal{N}^s$ and a zero representation of size vector $\m$, where  $m_{\L^0}=1$,
$m_{\I}=0,\ \forall \ \I\in T/\!\!\sim\setminus \{\L^0\}$ according to the partition of $(T, \sim)$. We will
prove in the next subsection, that there exists a reduction sequence $(\wt{\bigtriangleup})$ based on
$(\bigtriangleup)$:
$$(\wt{\bigtriangleup})\quad (\wt{N}, \wt{R}_1\times \wt{R}_2,
\wt{\n}_1\times \wt{\n}_2),
  \cdots, (\wt{N}^r, \wt{R}_1^r\times \wt{R}_2^r,
\wt{\n}_1^r\times \wt{\n}_2^r),  \cdots,  (\wt{N}^s, \wt{R}_1^s\times \wt{R}_2^s, \wt{\n}_1^s\times \wt{\n}_2^s)
$$

Next  we may compare the entries in  the following two  matrix equations:
$$R_1N^r=V+N^rR   \mbox{  and  } \wt{R}_1\wt{\mathcal{N}}^r=
\wt{V}+\wt{\mathcal{N}}^r\wt{R}_2.$$ Given any pair of indices $(\alpha, \beta)$ of $N^r$, denote by
$\wt{eq}_{(\alpha, \beta)}$ the equation given by the $(\alpha, \beta)$-entry and by $\wt{eq}_{(\alpha, o_l)}$
the equation given by the $(\alpha, o_l)$-entry in the second matrix equation.

{\bf Lemma 9.1.2} (1) $eq_{(\alpha, \beta)}$ and $\widetilde{eq}_{(\alpha, \beta)}$ are the same;  their right
hand sides given by ${V}+{{\mathcal N}}^r{R}_2$ and $\wt{V}+\wt{{\mathcal N}}^r\wt{R}_2$ contain only the
indeterminates at the $\beta$-th column of $R_2, V$ and $\wt{R}_2$, $\wt{V}$ respectively.

(2) The left side of the equation $\widetilde{eq}_{(\alpha, o_l)}$
given by $\widetilde R_1\wt {\mathcal N}^r$ equals zero; and the
right side of $\widetilde{eq}_{(\alpha, o_l)}$ given by
$\wt{V}+\wt{{\mathcal N}}^r\wt{R}_2$ contains only the
indeterminates at the $o_l-$th column of $\wt R_2$ and $\wt
V$.\hfill $\Box$

\subsection{The reduction sequence $(\wt{\bigtriangleup})$}

\kg We may write the matrix problem $(\K_1\times \K_2, \M, \mathcal{V})$ instead of $(\K_1\times \K_2,
\overline{\M}, \mathcal{V})$ for simplicity.

{\bf Proposition 9.2.1} There exists a reduction sequence $(\wt{\bigtriangleup})$ given in 9.1, such that the
added columns of $\wt{\cal N}^s$ consist of zero's and $\emptyset$'s.

 \bcen
\unitlength 1mm
\begin{picture}(120,30)
\put(-10,0){\framebox(15,30)} \put(-7,0){\line(0,1){30}} \put(-7,18){\line(1,0){12}} \put(-7,25){\line(1,0){12}}
\put(-4,20){$G$}

\put(20,0){\framebox(30,20)[l]{\dashbox{1}(3,20)}} \put(30,10){\framebox(10,8)}

\put(60,0){\framebox(30,20)[l]{\dashbox{1}(3,20)}} \put(60,10){\framebox(10,8)}

\put(100,0){\framebox(30,20)} \put(103, 0){\line(0,1){20}} \put(100,10){\framebox(10,8)}

\put(18,25){\makebox{$(i,l)$-block:}}

\put(-9, -5){Figure 1} \put(32, -5){Figure 2} \put(70, -5){Figure 3} \put(112, -5){Figure 4}
\end{picture}

\end{center}
\vskip 3mm

{\bf Proof.} When $r=0$, we have $\wt{T}=T\cup\{o_1,o_2,\cdots,o_{l_{\L^0}}\}$ and  the numbers of equivalent
classes of $T/\!\sim$ and of $\wt{T}/\wt{\sim}$ are the same, and the equation system II$^{\prime}$ at
$(\P_0,\Q_0)$ for $(\K_1\times\K_2,\M, \mathcal{V})$ and at $(\wt \P_0,\wt \Q_0)$ for
$(\wt{\K}_1\times\wt{\K}_2, \wt{\M}, \wt{\mathcal{V}})$ are also the same.  Suppose the $r$-th reduction has
been done and we obtain a triple $(\wt{N}^r, \wt{R}_1^r\times\wt{R}_2^r, \wt{\n}_1^r\times \wt{\n}_2^r)$, such
that $\wt{\mathcal{N}}^r$ is obtained from $\mathcal{N}^r$ by adding a column at the left of the $l$-th block
column consisting of  $0$'s and $\emptyset$'s for all $l\in\L^0$.

(1) If the number of the equivalent classes of $\wt{T}^r/\! \wt{\sim}^r$ equals that of $T^r/\!\sim^r$, the
equation system \textrm{II} at $(\P_r, \Q_r)$ in Definition 8.1.1 and the equation $x_{p_rq_r}^3=0$  for
$(\K_1^r\times\K_2^r, \M^r, \mathcal{V}^r)$ and those at $(\wt \P_r,\wt \Q_r)$ for
$(\wt{\K}_1^r\times\wt{\K}_2^r, \wt{\M}^r, \wt{\mathcal{V}}^r)$ are the same. Thus the reduction type of both of
them are the same. Suppose $N^r_{p_rq_r}$, the first block of $N^r$, belongs to the $(i,j)$-block according to
the partition of $(T,\sim)$. If $j\notin \L^0$, or $ j\in \L^0$ but $N_{p_rq_r}^r$ does not intersect the first
column of the $j$-th block column of $N^r$, then $\wt{n}_{p_r}^r=n_{p_r}^r$, $\wt{n}_{q^r}^r=n_{q^r}^r$, and we
set $\overline{\wt{N}^r_{p_rq_r}}=\overline{N^r_{p_rq_r}}$ (see figure 2).  If $j=l\in{\L^0}$ and $N_{p_rq_r}^r$
intersects the first column of  the $l$-th  block column of $N^r$, then
 $\wt{n}_{p_r}^r=n_{p_r}^r$, but $\wt{n}_{q_r}^r=n_{q_r}^r+1$ (see
figure 3). Then $p_r\nsim q_r$ in $(T^r,\sim^r)$, otherwise we would have $n^r_{p_r}=n^r_{q_r}$ and
$\wt{n}^r_{p_r}=\wt{n}^r_{q_r}$, which is a contradiction. Now we give the $(r+1)$-th reduction for
$({\wt{N}}^r, \wt{R}_1^r\times \wt{R}_2^r, \wt{\n}_1^r\times \wt{\n}_2^r )$ based on
 that for $(N^r, R_1^r\times R_2^r,\n_1^r\times\n_2^r)$.

 { Regularization}. we set
$\overline {\wt{N}^r_{p_rq_r}}=\emptyset$;

{ Edge reduction}. we set $\overline{{\wt{N}^r}_{p_rq_r}}=\left(\begin{array}{cc}{\large 0}_{n_{p_r}\times 1} &
\overline{{N}^r_{p_rq_r}}\end{array}\right)=\left(
\begin{array}{ccc} 0&0&I_d\\ 0&0&0 \end{array}\right) $;

{ Loop reduction.} it can not occur,  because of $p_r\nsim q_r$.

Therefore $\wt{\mathcal{N}}^{r+1}$ is obtained from $\mathcal{N}^{r+1}$ by adding a column at the left of each
$l$-th block column consisting of $0$'s or $\emptyset$'s for all $l\in\L^0$ in this  case.

(2) If  $d=n_{q_{r_0}}^{r_0}$ in the preceding edge reduction for some step $r_0$, and
$$
\overline{\wt{N}^r_{p_rq_r}}=\left(\unitlength=1mm
\begin{picture}(17,5)
\mput(12,-3)(0,1.5){6}{\line(0,1){1}} \put(1,-2.5){$0$}\put(1,2.5){$0_{d\times 1}$}
\put(15,-2.5){$0$} \put(15,2.5){$I_d$}
\end{picture}
\ \ \right),
$$
then a new equivalent class $\Q=\{o_1,o_2,\cdots,o_{l_{\L_0}}\}$ takes place at this stage. Thus
$$\wt{T}^{r_0+1}/\!\wt{\sim}^{r_0+1}=(T^{r_0+1}/\!\sim^{r_0+1})\cup
\{\Q\}.$$
 We stress that the situation must occur, since
$\mathfrak{B}^s$ is local.

(3) Suppose $r>r_0$. If $j\notin \L^0$, or $j\in\L^0$ but $N^r_{p_rq_r}$ does not intersect the first column of
the $j$-th block column of $N^r$, then $\wt{n}^r_{\P_r}=n^r_{\P_r}$, $\wt{n}^r_{\Q_r}=n^r_{\Q_r}$. We set
$\overline{\wt{N}^{\, r}_{p_rq_r}}=\overline{{N}^r_{p_rq_r}}$, since the equation systems \textrm{II}$^{\prime}$
are the same at $(\P_r,\Q_r)$ for both $(\K_1^r\times\K_2^r, \M^r, \mathcal{V}^r)$ and
$(\wt{\K}_1^r\times\wt{\K}_2^r, \wt{\M}^r, \wt{\mathcal{V}}^r)$ by Lemma 9.1.2. Otherwise
$\wt{n}_{p_r}^r=n_{p_r}^r$, $\wt{n}_{q_r}^r=n_{q_r}^r$, $\wt{n}_{o_l}=\wt{n}_{\Q}^r=1$ (see figure 4). The
$(r+1)$-th reduction for $\wt{N}_{p_rq_r}^r$ is divided into two reductions based on the reduction for
$N_{p_rq_r}^r$. Since  equation system \textrm{II}$'$ at $(\P_r,\Q_r)$ are the same  in both
$(\K_1^r\times\K_2^r, \M^r, \mathcal{V}^r)$ and $(\wt{\K}_1^r\times\wt{\K}_2^r, \wt{\M}^r, \wt{\mathcal{V}}^r)$,
we set
$$
 {\overline{{\wt{N}^r}_{p_rq_r}} }=(0 \ | \
 \overline{{N}^r_{p_rq_r}}\ )\ \ {\rm or}\ \ \overline{{\wt{N}^r}_{p_rq_r}}=(\emptyset \ | \
  \overline{{N}^r_{p_rq_r}}\ )
  $$
according to the equation $x^3_{p_ro_l}=0$ is , or is not a linear combination of the equations  in
\textrm{II}$'$ of 8.1.1 at $(\P_r,\Q)$ for
 $(\wt{\K}_1^r\times\wt{\K}_2^r, \wt{\M}^r, \wt{\mathcal{V}}^r)$.

Thus $\wt{\mathcal{N}}^{r+1}$ is obtained from $\mathcal{N}^{r+1}$ by adding a column at the left of  the $l$-th
block column consisting of $0$'s and $\emptyset$'s for all $l\in\L^0$ in this case.

Finally, $(\wt{\bigtriangleup})$ is a reduction sequence by induction on $r$.  The proof is completed.
\hfill$\square$

\medskip

{\bf Corollary 9.2.1}\ $(\wt{T}^s, \wt{\sim}^s)$ possesses two equivalent classes: one is the unique class $\P$
of $( T^{s},\sim^{s})$, another one is $\Q=\{o_1,o_2,\cdots, o_{l_\L^0}\}$.

Suppose that  for any $\ i,j\in T$,  $R^{s}_{ij}$ stands for the $(i,j)$-block of $R_2^{s}$ partitioned by
$T_2/\sim_2$, and $\wt{R}^{s}_{ij}$  for the $(i,j)$-block of $\wt{R}_2^{s}$ partitioned of  by
$\wt{T}_2/\!\wt{\sim}_2$ .

{\bf Corollary 9.2.2} (1) $\wt{R}_1^r=R_1^r$ for $r=0,1, \cdots, s $.

(2) $\forall\ i, j\in T_2\setminus\L^0$, $\wt{R}^{s}_{ij}=R^{s}_{ij}$; and $\forall \,l_1,l_2\in\L^0$,
$\wt{R}^{s}_{l_1j}$ has a subblock $R^{s}_{l_1j}$ at the lower $n_{\L^0}\times n_{\J}$ part, $\wt{R}^{s}_{il_2}$
has a subblock $R^{s}_{il_2}$ at the right $n_{\I}\times n_{\L^0}$ part, $\wt{R}^{s}_{l_1l_2}$ has a subblock
$R^{s}_{l_1l_2}$ at the lower-right $n_{\L^0}\times n_{\L^0}$ part (see figure 5).

(3) The diagonal block $\wt{R}^{s}_{ll}$ shown in figure 6 has  free entries at the top-row and zero's at the
 left-column below the top-row.

$$
\unitlength 1mm
\begin{picture}(100,30)

\put(10, 5){\framebox(21,21)} \put(10,19){\line(1,0){21}} \put(17,12){\line(1,0){14}}
\put(24,26){\line(0,-1){21}} \put(17,26){\line(0,-1){14}} \put(13,21){$1$} \put(19.5,21){$2$}
\put(26.5,21){$3$}\put(19.5,14){$4$} \put(26,14){$5$} \put(26.5,7){$6$} \put(8,21){$i$} \put(19,27){$l_2$}
\put(26,27){$j$} \put(31.5,14){$l_1$}

\put(35,16){\line(1,0){3}} \put(37.8,15){$\rightsquigarrow$}

\put(45,3){\framebox(24,24)} \put(45,20){\line(1,0){24}} \put(52,10){\line(1,0){17}}
\put(62,27){\line(0,-1){24}} \put(52,27){\line(0,-1){17}} \put(55,27){\line(0,-1){17}}
\put(52,17){\line(1,0){17}} \multiput(52.5,25.5)(0,-1){17}{$\cdot$} \multiput(53.5,25.5)(0,-1){17}{$\cdot$}
\multiput(55,17)(1,0){13}{$\cdot\cdot$} \multiput(55,18)(1,0){13}{$\cdot\cdot$}

 \put(48,22){$1$} \put(57.5,22){$2$}
\put(64.5,22){$3$}\put(57.5,12){$4$} \put(64.5,12){$5$} \put(64.5,5){$6$}

\put(80,3){\framebox(24,24)} \put(80,23){\line(1,0){24}} \put(84,3){\line(0,1){24}} \put(81,24){$*$}
\put(86,24){$**\cdots *$} \put(81,12){$0$} \put(90,12){$R_{ll}^{s}$}
\end{picture}
$$

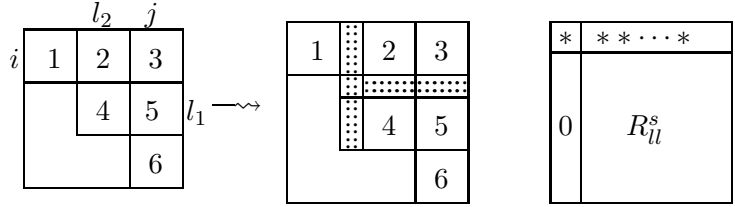
\hskip 40mm  figure 5: from $R_2^s$ to $\wt{R}_2^s$ \hskip 28mm figure 6 \hfill \vspace{2mm}

{\bf Proof.} We first prove the assertions (1) and (2). By  lemma 9.1.2, $\emptyset$'s, appearing in the added
column, only influence the added column of $\wt{R}_2$,  do not influence $R_1\times R_2 $ at all.

(3) The elements $*$ in $\wt{R}^s_{ll}$ of figure 6 are not involved in any equations, therefore they are
independent of all the others. The elements $0$ at the left column of $\wt{R}^s_{ll}$ are given by case (2) of
the proof of Proposition 9.2.1, since the $*$ at the upper-left corner of figure 6 must determine a single
equivalent class $\Q$ at this stage. \hfill$\square$

\ {\bf Example 1.}\ Consider Example of 7.6. Let $\K_1=\wt\Lambda, \K_2=\wt\Lambda,$\ $\M=rad \wt\Lambda$ and
${\cal V}=0$. We obtain a matrix problem $ ({\cal K}_1\times{\cal K}_2, \M,{\cal V})$ given by Proposition
8.1.1 with an index set $(T,\sim)$, such that $T/\!\sim=\{\I_1,\I_2\}$, where $\I_1=\{1,\cdots,5\},\
\I_2=\{6,\cdots,10\}$. Let $T_1=\I_1$, $T_2=\I_2$, then ${\cal K}_1$ is partitioned by $T_1$, ${\cal K}_2$ by
$T_2$ respectively. Define a size vector $\n$, such that $n_1=\cdots=n_5=2, n_6=\cdots=n_{10}=2.$ Then after a
sequence of reductions, we obtain a bocs $\frak B^s=\frak A^s$ (not properly partial!)  given by
$$
A=\left(\begin{array}{cc} 1&0\\0& 1 \end{array} \right),\ \ B=\left(\begin{array}{cc} 0&1\\ 0&0 \end{array}
\right),\ \ C=\left(\begin{array}{cc} c_3&c_4\\ \lambda& c_2
\end{array}\right),\ \ D=\left(\begin{array}{cc} d_3&d_4\\
d_1&d_2
\end{array} \right),
$$
which  satisfies MW5 of Theorem 5.6.1. Continue the reductions, then $\displaystyle
 C=\left(\begin{array}{cc}
\emptyset &\emptyset\\ \lambda& \emptyset
\end{array}\right)$  yields a local bocs $\frak A^{s'}$. Let
$\displaystyle D=\left(\begin{array}{cc} \emptyset&\emptyset\\
\eta&\emptyset
\end{array} \right)$. Thus $G=D$, and  $10$, the column index of $G$,
 belongs to $\I_2$, such that $\I_2\cap T_1=\emptyset$.
  Therefore the matrix problem $ ({\cal K}_1\times{\cal
K}_2, \M,{\cal V})$ satisfies S01 and S02 of 9.1.
 Define a new size vector $\wt\n$, such that $\wt n_1=\cdots=\wt n_5=2$,
 $\wt n_6=\cdots=\wt n_{10}=3$,
 and we add a column to the left of  the $l$-th block column for $l=6,\cdots,10$ as follows:
$$\wt A=\left(\begin{array}{ccc} 0&1&0\\0& 0&1 \end{array} \right),\ \ \wt B=\left(\begin{array}{ccc}
0&0&1\\0& 0&0 \end{array} \right),\ \ \wt C=\left(\begin{array}{ccc} c_0&c_3&c_4\\ 0& \lambda & c_2
\end{array}\right),\ \ \wt D=\left(\begin{array}{ccc} d'_0&d_3&d_4\\
d_0&d_1&d_2
\end{array} \right).
$$
Then we obtain a bocs $\wt{ \frak B}^s.$ If we continue the reductions up to $\wt C'= \left(\begin{array}{ccc}
\emptyset &\emptyset&\emptyset\\0&\lambda&\emptyset\end{array} \right) $, and obtain a layered bocs $\wt{\frak
B}^{s'}$ having two vertices :
\begin{center}
\unitlength 1mm
 \begin{picture}(30, 6) \put(10.00,-1.00){\circle*{1.00}}
\put(25.00,-1.00){\circle*{1.00}} \put(5.00,0.00){\oval(5,5)[t]} \put(5.00,0.00){\oval(5,5)[bl]}
\put(4.00,-3.00){\vector(3,1){3}} \put(11,-1){\line(1,0){12}} \put(23,-1){\vector(1,0){1}}
\put(10.00,3.00){\makebox{$\P$}} \put(25.00,3.00){\makebox{$\Q$}} \put(0.00,0.0){\makebox{$\lambda$}}
\put(15.50,1.00){\makebox{$d_0$}}
\end{picture}.
\end{center}

{\bf Example 2 .}  Consider  Example 4 in 2.1 and 4.3. A sequence of reductions leads to a minimal bocs below in
figure 7, where $\lambda^0$ is a fixed eigenvalue. Let us set the shadowed part $(1, \nu)$ free, then we obtain
a matrix problem $(\K_1\times \K_2, \M,{\cal V})$ in figure 8, where $T=\{1,2,3\}$ and $2\sim 3$, $T_1=\{1\}$,
$T_2=\{2,3\}$;
$$K_1=\{(s) |\  \forall\, s\in k\},\ \ K_2=\left\{\Big(\begin{array}{cc}u& 0\\0
&u\end{array}\Big)\Big|\  \forall\, u\in\, k\right\};$$ $\M=\{(m_1,m_2)|\forall \  m_1, m_2\in k\};\ \ {\cal
V}=\{(0,0)\}.$ Denote the equivalent class $\{1\}$ by $\P$, and $\{2,3\}$ by $\I$. Define a size vector $\n$,
such that $n_{\P}=1$, $n_{\I}=1$, then $\n_1=(1)$, $\n_2=(1,1)$.  Let $$N=(f_1,f_2),\ R_1=(x), \
R_2=\left(\begin{array}{cc} y& 0\\ 0& y\end{array}\right),\ V=(0,\ 0),$$ then $(N, R_1\times R_2,
\n_1\times\n_2)$ is a triple of matrix problem $(\K_1\times K_2, \M, \cal V)$ , which possesses an equation
 $$(x)(f_1,
f_2)= (f_1, f_2)\left(\begin{array}{cc} y& 0\\ 0& y\end{array}\right).$$
 Let us define a new size vector $\wt{\n}$,
such that $\wt{n}_{\P}=n_{\P}=1$, $\wt{n}_{\I}=n_{\I}+1=2$.  Then we add a column to the left at each
block-column of the matrix ${\cal N}=(0,0)$. The  structure in 9.1  gives a triple of matrices $(\wt{N},
\wt{R}_1\times \wt{R}_2, \wt{\n}_1\times\wt{\n}_2)$, where $\wt{N}=(f_{11}\ f_{12}\ \vdots \ f_{21}\ f_{22})$,
$\wt{R}_1=R_1=(x)$,
$$
\wt{R}_2= \left(
\begin{array}{cc|cc} y_{11}& y_{12} & & \\ y_{21}& y_{22} & &\\
\hline & & y_{11}& y_{12}\\ & & y_{21}& y_{22} \end{array} \right).
$$
After an edge reduction $(\overline{f_{11}\ f_{12}})=(0\ 1)$ in case of $\overline{f_1}=(1)$, we have
$\wt{\mathcal{N}}'=(0\ 1\ \vdots\ 0 \ 0)$,
 $\wt{N}'=(0\ 0\ \vdots\ c \ *)$, where $c=f_{21}, *=f_{22}$, $\wt{R}'_1=(x)$,
$$
\wt{R}'_2= \left( \begin{array}{cc|cc} y_{11}& y_{12} & & \\ 0& x & &\\ \hline & & y_{11}& y_{12}\\ & &  0& x
\end{array} \right).
$$
There are two equivalent classes $\P,\Q$ in $(\wt{N}', \wt{R}'_1\times \wt{R}'_2,\wt{\n}'_1\times\wt{\n}'_2)$
(see figure 9).
\begin{center}
\unitlength=1mm $\left(\begin{array}{c}
\begin{picture}(13,30)
\put(0,0){\dashbox{1}(12,30)} \mput(0,6)(2,0){6}{\line(1,0){1}} \mput(0,12)(2,0){6}{\line(1,0){1}}
\mput(0,18)(2,0){6}{\line(1,0){1}} \mput(0,24)(2,0){6}{\line(1,0){1}} \mput(6,6)(0,2){12}{\line(0,1){1}}
\mput(0,25.5)(0,1.5){3}{\line(1,0){12}}

\put(2,2){$0$} \put(8,2){$1$} \put(2,8){$1$} \put(8,8){$\emptyset$} \put(2,14){$1$} \put(8,14){$1$}
\put(2,20){$1$} \put(8,20){$\lambda^0$} \put(2,26){$1$} \put(8,26){$\nu$} \put(1,-12){Figure 7}
\end{picture}
\end{array}\right)$ \hskip 20mm
$\left(\begin{array}{c}
\begin{picture}(13,30)
\put(0,0){\dashbox{1}(12,30)} \mput(0,6)(2,0){6}{\line(1,0){1}} \mput(0,12)(2,0){6}{\line(1,0){1}}
\mput(0,18)(2,0){6}{\line(1,0){1}} \mput(0,24)(2,0){6}{\line(1,0){1}} \mput(6,6)(0,2){12}{\line(0,1){1}}
\put(2,2){$0$} \put(8,2){$1$} \put(2,8){$1$} \put(8,8){$\emptyset$} \put(2,14){$1$} \put(8,14){$1$}
\put(2,20){$1$} \put(8,20){$\lambda^0$} \put(2,26){$f_1$} \put(8,26){$f_2$} \put(-10,26){$\P\{$}
\put(0,30){$\overbrace{\hskip 12mm}$}
\put(-11,10){$\I \left\{\begin{array}{l} \\ \\ \\ \\ \\
\end{array}\right.$} \put(5,33){$\I$} \put(-1,-12){Figure 8}
\end{picture}
\end{array}\right)$ \hskip 20mm
$\left(\begin{array}{c}
\begin{picture}(24,45)
\put(0,0){\dashbox{1}(24,45)} \mput(0,10)(2,0){12}{\line(1,0){1}} \mput(0,20)(2,0){12}{\line(1,0){1}}
\mput(0,30)(2,0){12}{\line(1,0){1}} \mput(0,40)(2,0){12}{\line(1,0){1}} \mput(12,10)(0,2){18}{\line(0,1){1}}
\put(2,2){$0$} \put(7,2){$0$} \put(14,2){$0$} \put(19,2){$1$} \put(2,6){$0$} \put(7,6){$0$} \put(14,6){$1$}
\put(19,6){$0$} \put(2,12){$0$} \put(7,12){$1$} \put(14,12){$\emptyset$} \put(19,12){$\emptyset$}
\put(2,16){$1$} \put(7,16){$0$} \put(14,16){$\emptyset$} \put(19,16){$\emptyset$} \put(2,22){$0$}
\put(7,22){$1$} \put(14,22){$0$} \put(19,22){$1$} \put(2,26){$1$} \put(7,26){$0$} \put(14,26){$1$}
\put(19,26){$0$} \put(2,32){$0$} \put(7,32){$1$} \put(15,32){$0$} \put(19,32){$\lambda^0$} \put(2,36){$1$}
\put(7,36){$0$} \put(14,36){$\lambda^0$} \put(20,36){$0$} \put(2,42){$0$} \put(7,42){$1$} \put(14,42){$c$}
\put(19,42){$*$}

\put(8,-5){Figure 9}

\end{picture}
\end{array}\right)$
\end{center}
\vskip 3mm

This example satisfies the assumption  S02 of 9.1, but does not satisfy S01. The example  can be used to help to
understand the structure of 9.1 and also Formula (54) of 8.4.

\subsection{The differentials in partial bocs $\wt{\mathfrak{B}}^s$}

\kg In this subsection we will calculate the differentials of the solid arrows of $\wt{\mathfrak{B}}^s$.

{\bf Lemma 9.3.1}\ Let $(\alpha(0),\beta(0)) $ be the matrix index of $\lambda$, and $(\alpha^0,\beta^0)$ that
of $\nu$ given by Definition 9.1.1. Then $\alpha(0)>\alpha^0$ under the order of Definition 2.3.1.

{\bf Proof.} Since $(\alpha(0), \beta(0))\prec(\alpha^0, \beta^0)$,  $\alpha(0)\geq \alpha^0$. If
$\alpha(0)=\alpha^0$, then any solid arrow $a_l$ with $\lambda \prec a_l \preceq \nu $ has row index $\alpha^0
$. Therefore the summands of $\delta(a_l)$
  involve only the terms of $v$ or
$\lambda^e v$ for some dotted arrows $v$ and positive integer $e$. Thus  we obtain a contradiction to $\wt{w}\ne
0$, $(\lambda-\mu)|\wt{h}(\lambda,\mu)$, since such a $\wt{h}(\lambda, \mu)$ must contain some summands
$v\lambda^e$. Therefore $\alpha(0)>\alpha^0$ as desired.  \hfill$\square$

 Suppose  that the partial bocs
$\wt{\mathfrak{B}}^s$ has a partial layer
$$\wt{L}^s=(\wt{\Gamma}';\; \wt{\omega};\; a_1, \cdots,  a_j,
a_{j+1},\cdots,a_{j+m}; c_1,\, \cdots, c_r, c^0,\, d_1, \cdots, d_n;\, w,\, u,\, v),$$
 where $\ind\wt{\Gamma}'=\{\P,\Q\}$,
$\wt{\Gamma}'(\P,\P)=\rK[\lambda]$ and $\lambda$ appears for the first time in $\wt{\mathfrak{B}}^s$,
$\wt{\Gamma}'(\Q,\Q)=\rK$; $a_l:\, \P\longrightarrow \P$, $l=1,2,\cdots,j$, are solid loops appearing before the
$\alpha^0$-th row of $G$ according to the matrix order of Definition  2.3.1, and $a_{j+1}, \cdots, a_{j+m} $ are
solid loops at the $\alpha^0$-th row of $G$; $c_1,\cdots, c_r,\, c^0: \P\rightarrow \Q$ are solid edges sitting
at the added columns up to the $\alpha^0$-th row of $G$, and the last one $c^0$ is sitting at the intersection
of the $o_{l^0}$-th column and the $\alpha^0$-th row; $d_1, \cdots, d_n$ are all  the solid arrows after
$a_{j+m}$;  $w: \P\rightarrow \P$, $u : \P\rightarrow \Q$ are dotted arrows (see the figures below).

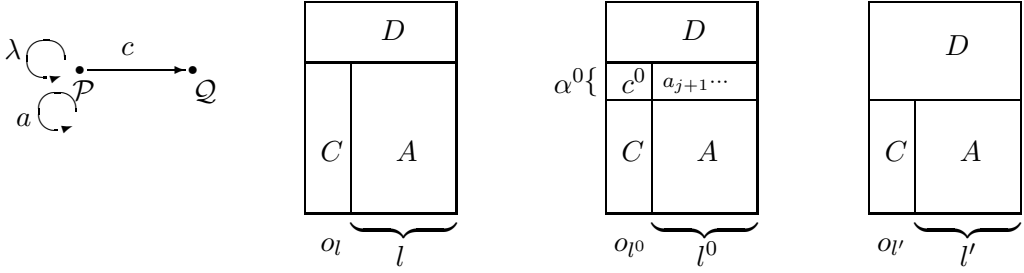
\unitlength 1mm
\begin{center}
\begin{picture}(140,30)
\unitlength 1mm  \put(10.00,19.00){\circle*{1.00}} \put(25.00,19.00){\circle*{1.00}}
\put(11,19){\line(1,0){12}} \put(23,19){\vector(1,0){1}} \put(9.00,15.00){\makebox{$\P$}}
\put(25.00,15.00){\makebox{$\Q$}} \put(0.00,20.50){\makebox{$\lambda$}} \put(15.50,21.00){\makebox{$c$}}
\put(5.50,20.50){\oval(5,5)[t]} \put(5.50,20.5){\oval(5,5)[bl]} \put(4.50,17.25){\vector(3,1){3}}
\put(7.00,13.50){\oval(5,5)[t]} \put(7.0,13.50){\oval(5,5)[bl]} \put(6.00,10.50){\vector(3,1){3}}
\put(1.50,11.50){\makebox{$a$}}

\put(40,0){\framebox(20,28)} \put(46,20){\line(0,-1){20}} \put(40,20){\line(1,0){20}} \put(50,23){\makebox{$D$}}
\put(52,7){\makebox{$A$}} \put(42,7){\makebox{$C$}} \put(46,-1){\makebox{$\underbrace{\hskip 14mm}$}}
\put(52,-7){\makebox{$l$}} \put(42,-5){\makebox{$o_{l}$}}


\put(80,0){\framebox(20,28)} \put(86,20){\line(0,-1){20}} \put(80,20){\line(1,0){20}}
\put(80,15){\line(1,0){20}} \put(90,23){\makebox{$D$}} \put(92,7){\makebox{$A$}} \put(82,16){\makebox{$c^0$}}
\put(82,7){\makebox{$C$}} \put(87.5,16.8){\makebox{$\s a_{j+1} \cdots$}}
\put(86,-1){\makebox{$\underbrace{\hskip 14mm}$}} \put(77,16.5){\makebox{$\{$}}
\put(73,16){\makebox{$\alpha^0$}} \put(92,-7){\makebox{$l^0$}}

\put(81,-5){\makebox{$o_{l^0}$}}

\put(115,0){\framebox(20,28)} \put(121,15){\line(0,-1){15}} \put(115,15){\line(1,0){20}}
\put(125,21){\makebox{$D$}} \put(127,7){\makebox{$A$}} \put(117,7){\makebox{$C$}}
\put(121,-1){\makebox{$\underbrace{\hskip 14mm}$}}
\put(127,-7){\makebox{$l'$}}
\put(116,-5){\makebox{$o_{l'}$}}

\end{picture}
\end{center}
\vspace{2mm} \noindent
 Where $l,l'\in\L^0$,  $a_{j+1}, \cdots,
a_{j+m}$ belong to the area $A$.

\hskip 8mm\setcounter{equation}{56}

Lemma 9.1.1 and 9.1.2 enables us to calculate the differentials of the solid arrows of $\wt{\mathfrak{B}}^s$.
Let
\begin{equation}
 \left\{
\begin{array}{ccl}
\dz(a_1)^0 &=& h_{11}(\lambda,\mu)w_1 \\
\cdots   &  & \qquad\cdots \\
\dz(a_j)^0 &=& h_{j1}(\lambda,\mu)w_1 +\cdots+h_{jj}(\lambda,\mu)w_j
\end{array} \right.
\end{equation}
be given by b Formula (42) of 5.4 up to $j$, then $h_{ll}(\lambda,\lambda)\ne 0$, (see part A in the figure). On
the other hand, part $C$ in the figure gives

\begin{equation}
 \left\{
\begin{array}{ccl}
\dz(c_1)^0 &=& f_{11}(\lambda)u_1 \\
\dz(c_2)^0 &=& f_{21}(\lambda)u_1 +f_{22}(\lambda)u_2 \\
\cdots   &  & \qquad\cdots \\
\dz(c_r)^0 &=& f_{r1}(\lambda)u_1 +f_{r2}(\lambda)u_2+\cdots+f_{rr}(\lambda)u_r
\end{array} \right.
\end{equation}
where $u_{r(1)}=0, \cdots, u_{r(\beta)}=0$ for some indices $1\leq r(1)<\cdots<r(\beta)\leq r $, and the dotted
arrows $u_l$, $\forall \ l\in I^0$,   are linearly independent, where $I^0=\{1,2,\cdots,r\}\setminus\{r(1),
\cdots,r(\beta)\}$.

 If
\begin{equation}
  \wt{R}^s_{ll}=\left(\begin{array}{c|cccc} x_{00}& x_{01}
& x_{02} & \cdots & x_{0m} \\ \hline &&&& \\  0& & R^s_{ll}& &
\\ &&&&
\end{array}\right), \mbox{  then  }
 \left\{
\begin{array}{ccl}
\dz(a_{j+1}) &=& c^0 x_{01}+\cdots \\
\dz(a_{j+2})&=& c^0 x_{02}+ \cdots \\
\cdots   &  & \cdots \\
\delta(a_{j+m}) &=& c^0 x_{0m}+\cdots
\end{array} \right.
\end{equation}

 Define a polynomial:
 \begin{equation}  \qquad h^0(\lambda,\mu)=\prod_{l=1}^jc_l(\lambda)
 h_{ll}(\lambda,\mu)\cdot \prod_{l\in I^0}f_l(\lambda).\end{equation} (see
 Formula (30) of 4.3 and Formula (43) of 5.4.)
 Let $h^0(\lambda,\lambda)\neq 0$,
we set $a_l=\emptyset$, for $l=1,\cdots,j$; $c_l=\emptyset$ for $l\in I^0 $, $c_{r(l)}=0$ for
$l=1,\cdots,\beta$,  then  we obtain an induced partial bocs $\wt{\mathfrak{B}}^{s+j+r}$.  We will prove in the
next subsection that $\delta(c^0)=0$ in $\wt{\mathfrak{B}}^{s+j+r}$.  Let $c^0=1$; $a_{j+l}=\emptyset$, and
$x_{o_l}+\cdots =0$ in Formula (59) for  $l=1,\cdots,m$. Then we obtain an induced partial bocs
$\wt{\mathfrak{B}}^{s+J}$ with $J=j+r+1+m$, which is parameterized, but not necessarily freely parameterized in
general, because there is no guarantee that $h^0(\lambda,\mu)$ satisfies item (3) of Definition 4.2.1.

{\bf Lemma 9.3.2}\ $\wt{\mathfrak{B}}^{s+J}=(\wt{\Gamma}^{S+J}, \wt{\Omega}^{s+J})$ is an induced bocs of
$\wt{\mathfrak{B}}=(\wt{\Gamma}, \wt{\Omega})$, such that

(1) $\wt{\mathfrak{B}}^{s+J}$ is local;

(2) $\wt{\mathfrak{B}}^{s+J}$ possess a unique parameter $\lambda$ with the domain $h^0(\lambda,\mu)\ne 0$ given
in Formula (60), which is not necessarily freely parameterized;

(3) $\wt{\Gamma}^{S+J}$ is freely generated by $d_1,d_2,\cdots,d_n$ over $k[\lambda, \mu, h^0(\lambda,
\mu)^{-1}]$;

(4) $\wt{\Omega}^{s+\J}$ is obtained from $\wt{\Omega}$ by the restriction $w _l=0$ for $1\leqslant l\leqslant
j$, (Formula(57)), $u_l=0$ for any $l\in I^0$ (Formula(58)), and $x_{ol}+\cdots =0$ for $1\leqslant l \leqslant
m$ (Formula(59)). \hfill $\square$

Let $T_1, T_2$ be the index set of $\mathfrak{B}$, and let $\mathcal{S}=\{1,2,\cdots,n_1 \}$, if the number of
the rows of matrix $N$  is $n_1$. We define a set of pairs of indices
$$\{(\alpha,l)\mid \alpha\in\mathcal{S}, l\in T_2\}.$$ And we give an order
in the set, such that $(\alpha_1, l_1)\prec(\alpha_2,l_2)$, provided $\alpha_1>\alpha_2$, or $\alpha_1=\alpha_2$
but $l_1<l_2$, (see the following two examples).

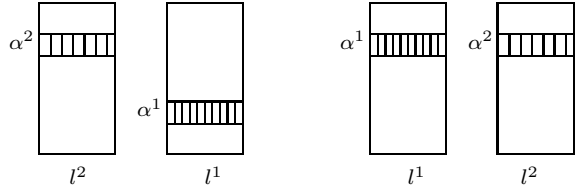
\unitlength 1mm
\begin{center}
\begin{picture}(80,25)
\unitlength 1mm \put(0,0){\framebox(10,20)} \put(0,13){\line(1,0){10}} \put(0,16){\line(1,0){10}}
\put(-4,14){$\s \alpha^2$} \put(1.5,13){\line(0,1){3}} \put(3,13){\line(0,1){3}}\put(4.5,13){\line(0,1){3}}
\put(6,13){\line(0,1){3}}\put(7.5,13){\line(0,1){3}} \put(9,13){\line(0,1){3}} \put(4,-4){$\s l^2$}

\put(17,0){\framebox(10,20)} \put(17,4){\line(1,0){10}} \put(17,7){\line(1,0){10}} \put(13,5){$\s \alpha^1$}
\put(18,4){\line(0,1){3}} \put(19,4){\line(0,1){3}}\put(20,4){\line(0,1){3}}
\put(21,4){\line(0,1){3}}\put(22,4){\line(0,1){3}} \put(23,4){\line(0,1){3}}\put(24,4){\line(0,1){3}}
\put(25,4){\line(0,1){3}}\put(26,4){\line(0,1){3}} \put(22,-4){$\s l^1$}

\put(44,0){\framebox(10,20)} \put(44,13){\line(1,0){10}} \put(44,16){\line(1,0){10}} \put(40,14){$\s \alpha^1$}
\put(45,13){\line(0,1){3}} \put(46,13){\line(0,1){3}}\put(47,13){\line(0,1){3}}
\put(48,13){\line(0,1){3}}\put(49,13){\line(0,1){3}} \put(50,13){\line(0,1){3}}\put(51,13){\line(0,1){3}}
\put(52,13){\line(0,1){3}}\put(53,13){\line(0,1){3}} \put(49,-4){$\s l^1$}

\put(61,0){\framebox(10,20)} \put(61,13){\line(1,0){10}} \put(61,16){\line(1,0){10}} \put(57,14){$\s \alpha^2$}
\put(62.5,13){\line(0,1){3}} \put(64,13){\line(0,1){3}}\put(65.5,13){\line(0,1){3}}
\put(67,13){\line(0,1){3}}\put(68.5,13){\line(0,1){3}} \put(70,13){\line(0,1){3}} \put(64,-4){$\s l^2$}

\end{picture}
\end{center}

We specify a pair of indices $(\alpha^0,l^0)$ for $\nu$ given by Definition 9.1.1. And we also specify a pair
$(\alpha(0), l(0))$ for $\lambda$, where $\alpha(0)$ is the row-index of $\lambda$ in the matrix;  $l(0)$ the
index of the  block-column partitioned  by $(T,\sim)$, such that $\lambda$ sits in that block.

\medskip
 Next we use the
procedure of Lemma 5.3.1 once more for  the local partial bocs $\wt{\mathfrak{B}}^{s+J}$ to get a sequence of
parameters $\lambda=\lambda_{10}, \lambda_{11}, \cdots,\lambda_{1\gamma^1} $.

{\bf Proposition 9.3.1} (1) If $\gamma^1=0$, we have
\begin{equation}
\left\{\begin{array}{ccl}
\delta(d_1)^0&=&g_{11}(\lambda, \mu)v_1\\
\delta(d_2)^0&=&g_{21}(\lambda, \mu)v_1+g_{22}(\lambda, \mu)v_2\\
\cdots & &\cdots\qquad\cdots  \\
\delta(d_n)^0&=&g_{n1}(\lambda, \mu)v_1+g_{n2}(\lambda, \mu)v_2+\cdots+g_{nn}(\lambda, \mu)v_n
\end{array} \right.
\end{equation}
where $g_{ll}(\lambda, \lambda)\neq 0$ (see Formula (40) and (41) in 5.4, and we use $\lambda, \mu$ instead of
$\nu,\kappa$).

(2) If $\gamma^1=1$, denote $\lambda_{11}$ by $\nu_1$. Then $(\alpha^1,l^1)\succ(\alpha^0,l^0)$.

(3) If $\gamma^1\geq 2$, denote $\lambda_{1,\gamma^1-1}$ by $\lambda_1$, $\lambda_{1\gamma^1}$ by $\nu_1$
respectively. Then $(\alpha(1),l(1))\succ(\alpha(0),l(0))$.

{\bf Proof.} (2) If $\gamma^1=1$, then $\nu_1$ is located outside of $G$ or  above the $\alpha^0$-row of $G$  by
Formula (59). Thus $(\alpha^1,l^1)\succ(\alpha^0,l^0)$.

(3)  If $\gamma^1\ge 2$, then $\lambda_1$ is located  outside of $G$ or  above the $\alpha^0$-row of $G$  still
by Formula (59), thus $(\alpha(1), l(1))\succ (\alpha^0,l^0)$. On the other hand,
$(\alpha^0,l^0)\succ(\alpha(0),l(0))$ by Lemma 9.3.1. Consequently
$(\alpha(1),l(1))\succ(\alpha(0),l(0))$.\hfill$\square$

\subsection{The calculation of $\delta(c^0)$}

\kg In this subsection we will prove that
 $\delta(c^0)=0$ in $\wt{\mathfrak{B}}^{s+j+r}$ under
the assumption S01 and S02 of structure 9.1.

Given any $l\in \mathcal{L}^0$, we define a matrix index $\beta_l\in l$, such that the distance from the
$\beta_l$-th column to the left edge of  the $l$-th block-column is the same as the distance from the
$\beta^0$-th column to the left edge of $l^0$-th block-column. And $\beta_{l^0}=\beta^0$.
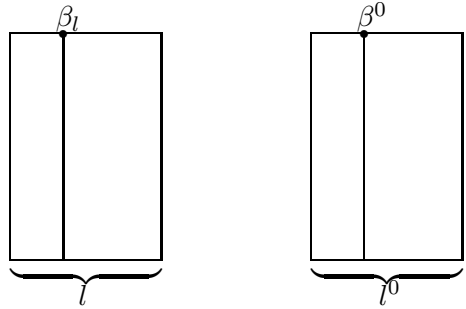
\begin{center}\unitlength=1mm
\begin{picture}(60, 32)
\put(0,0){\framebox(20,30)} \put(7,0){\line(0,1){30}} \put(6, 31){$\beta_l$} \put(7, 30){\circle*{1}} \put(0,
-1){$\underbrace{\hskip 20mm}$} \put(9,-6){$l$}

\put(40,0){\framebox(20,30)} \put(47,0){\line(0,1){30}} \put(46, 31){$\beta^0$} \put(47, 30){\circle*{1}}
\put(40, -1){$\underbrace{\hskip 20mm}$} \put(49,-6){$l^0$}
\end{picture}
\end{center}\vskip 5mm

Recall the matrix equation $R_1\mathcal{N}^s=V+ \mathcal{N}^s R_2$
 from 9.1.  If $eq$ is an equation, $eq^R$ stands for the
right-hand side of $eq$ given by $V+\mathcal{N}^sR_2$, and if $Eq$ is a set of equations, $Eq^R$ stands for the
set of the right-hand side of the equations in  $Eq$. Recall from  Lemma 9.1.2, that $eq^R_{(\alpha, \beta)}$
contains only the indeterminates at the $\beta$-th column of $R_2$ and $V$ for any pair of indices $(\alpha,
\beta)$.

\medskip
{\bf Lemma 9.4.1} Let $(\alpha^0,\beta^0)$ be the index of $\nu$ given in Definition 9.2.1. Then
$eq^R_{(\alpha^0,\beta^0)}$
  is a linear combination of the forms of $Eq^R $, where $Eq= Eq_{(\alpha_0,\beta_0)}$ is
  given by Definition 9.1.2.

{\bf Proof.} (1)  We set $\lambda=\mu$, then $eq_{(\alpha^0,\beta^0)}^{c}$ is a linear combination of the
equations of $Eq^c$ in $\mathfrak{B}^{s+i-1}$ by Lemma 9.1.1 and Formula (42) of 5.4, where $eq^c$ stands for
the equation $eq$ when we do not distinguish the multiplication of the parameters from left or right.

(2) Let $X_1$ be the set of  indeterminates at the   $l$-th block column of $R_2$ and $V$ for any   $l\in \L^0$,
and $X_2$ be that at the $j$-th block column of $R_1$, $R_2$ and $V$ for any $j\in T\setminus\L^0$. Thus
$X_1\stackrel{\cdot}{\cup} X_2$ is the complete set of the indeterminates. Deleting $X_2$ from
$eq_{(\alpha^0,\beta^0)}^{c}$ and $Eq^c$, we obtain an equation $eq_{(\alpha^0,\beta^0)}^{'c}$ and a set of
equations $Eq^{'c}$. It is straightforward that $eq_{(\alpha^0,\beta^0)}^{'c}$ is the same linear combination of
the equations in $Eq^{' c}$.

(3) The left-hand sides of $eq_{(\alpha^0,\beta^0)}^{'c}$ and the
equations in  $Eq^{'c}$ are all zero by assumption S02  of   9.1.
And the right-hand side of $eq_{(\alpha_0,\beta_0)}^{'c}$  is just
$eq_{(\alpha^0,\beta^0)}^R$, and that of any equation in $Eq^{'c}$
is either a form in   $Eq^R$ containing the indeterminates at the
$\beta$-th column for some $\beta\in l \in \L^0$ or zero.
Therefore $eq_{(\alpha_0,\beta_0)}^R$ is a linear combination of
the forms of $Eq^R$. The proof is completed.\hfill$\square$

\medskip
{\bf Lemma 9.4.2} In the matrix equation $\wt{R}_1\wt{\mathcal{N}}^s=\wt{V}+ \wt{\mathcal{N}}^s\wt{R}_2$,
$\wt{eq}_{(\alpha^0,o_{l^0})}^R$ is a linear combination of the forms in $\wt{Eq}^R$, where
$\wt{Eq}=\wt{Eq}_{(\alpha^0,o_{l^0})} $ is given by Definition 9.1.2.

{\bf Proof.} (1) Since for any pair of indices $(\alpha,\beta)$ of $\mathcal{N}^s$, $eq_{(\alpha,\beta)}$ is the
same as $\wt{eq}_{(\alpha,\beta)}$ by Lemma 9.1.2,  it is also the case for the right-hand sides of the
equations. Because $eq^R_{(\alpha^0,\beta^0)}$ is a linear combination of the forms in $Eq^R$ by Lemma 9.4.1,
so is $\wt{eq}^R_{(\alpha^0,\beta^0)}$ in $\wt{Eq}^R$.

(2)  For any $ \alpha\ge \alpha^0$, and $ l\in \L^0$, the linear form  of the $(\alpha, o_l)$-th entry of
$\wt{V}+\wt{\mathcal{N}}^{s}\wt{R}_2$ in the indeterminates at the $o_l$-th column of $\wt{V}$ and $\wt{R}_2$ is
the same as the linear form of the $(\alpha, \beta_l)$-th entry in the indeterminates at the $\beta_l$-th column
by Lemma 8.5.1, where $\beta_l\in l$ given by the above diagram.

Combining (1) and (2),  $\wt{eq}_{(\alpha^0,o_{l^0})}^R$
  is a  linear combination of the forms in $\wt{Eq}^R$ in the
  indeterminates of the $o_l$-th column of $\wt{R}_2$ and $\wt{V}$
  for $l\in\L^0$. \hfill $\square$

\medskip
 {\bf Proposition 9.4.1} $\wt{eq}_{(\alpha^0,o_{l^0})}$ is
 a linear combination of the equations in  $\wt{Eq}$.

 {\bf Proof.} Since the left-hand sides of all the equations $\wt{eq}_{(\alpha,o_l)}$
 for  $\alpha\ge \alpha^0$ and  $l\in \L^0$ equal zero,
$\wt{eq}_{(\alpha^0,o_{l^0})}$ is a linear combination of equations of $\wt{Eq}$ by Lemma 9.4.2. \hfill
$\square$

\medskip
{\bf Corollary 9.4.1}  $\delta(c^0)=0$ in $\wt{\mathfrak{B}}^{s+j+r}$.

{\bf Proof.} By Proposition 9.4.1 and Lemma 9.1.1. \hfill$\square$

\subsection{The induction on $(\alpha^{\eta}, l^{\eta})$}

 \kg Suppose we are given $\mathfrak{B}^s$, $\wt{\mathfrak{B}}^s$,
$\wt{\mathfrak{B}}^{s+J}$ in 9.3, then $\wt{\mathfrak{B}}^{s+J}$ is local by Lemma 9.3.2. Assume that when we
perform the procedure 5.3 starting from $\wt{\mathfrak{B}}^{s+J}$, $\gamma^1=1$ in Proposition 9.3.1, and
$\wt{\mathfrak{B}}^{s+J}$
 possesses   the following formulae of 5.4

\begin{equation} \left\{\begin{array}{lcll}
\delta(a_1^{1})^0&=&h^{1}_{11}(\lambda, \mu)w^{1}_1&\\
\cdots & & \cdots\qquad\cdots &\\
\delta(a^{1}_{i^{1}-1})^0 &=& h^{1}_{i^1-1,1}(\lambda,\mu)w^{1}_1&+\cdots+
h^{1}_{i^{1}-1,i^{1}-1}(\lambda,\mu)w^{1}_{i^{1}-1}\\
\delta(a^{1}_{i^{1}})^0&=&h^{1}_{i^1,1}(\lambda,\mu)w^{1}_1& +\cdots+
h^{1}_{i^{1},i^{1}-1}(\lambda,\mu)w^{1}_{i^{1}-1} \quad + \wt{h}^{1}(\lambda,\mu)\wt{w}^{1}
\end{array} \right.
\end{equation}
  where $a_l^1=d_l$ for $l=1,\cdots,i^1$, $\wt{w}^1$ is
linearly independent of $w^1_1, \cdots, w^1_{i^1-1}$, and $(\lambda-\mu) \mid \wt{h}^1(\lambda,\mu)$,
$(\lambda-\mu)^2 \nmid\wt{h}^1(\lambda,\mu)$.  Let
\begin{equation}
 h^1(\lambda, \mu)=
h^0(\lambda,\mu)\prod\limits_{l=1}^{j^1}c_l^1(\lambda)h^1_{ll}(\lambda, \mu)\cdot\prod\limits_{l \in
{I^1}}f_{l}^1(\lambda)
\end{equation}
For any fixed $\lambda^0\in\rK$ with $h^1(\lambda^0, \lambda^0)\neq 0$,  set $a_l^1=\emptyset$,
$l=1,\cdots,(i^1-1)$, $a^1_{i^1}=\nu_1$, then
$$
\left\{\begin{array}{lcl}
\delta(d_1^1)^0&=&g^1_{11}(\nu_1, \kappa_1)v^1_1\\
\cdots & & \cdots\qquad\cdots \\
\delta(d^1_{n^1-1})^0 &=& g^1_{n^1-1,1}(\nu_1, \kappa_1)v^1_1+\cdots+
g^1_{n^1-1,n^1-1}(\nu_1, \kappa_1)v^1_{n^1-1}\\
\delta(d^1_{n^1})^0&=&g^1_{n^1,1}(\nu_1, \kappa_1)v^1_1\quad+\cdots+ g^1_{n^1,n^1-1}(\nu_1,
\kappa_1)v^1_{n^1-1}\quad + g^1_{n^1,n^1}(\nu_1, \kappa_1)v^1_{n^1}
\end{array} \right.
$$
where $d^1_l=d_{i^1+l}$, $n^1=n-i^1$, (see Formula (40) of 5.4), and $\sigma_l^1(\nu_1)$ are given by  (41) of
5.4, $g^1_{ll}(\nu_1,\kappa_1)\in \rK[\nu_1,\kappa_1,\sigma_l^1(\nu_1)^{-1}, \sigma_l^1(\kappa_1)^{-1}]$ are all
invertible.

 Suppose $\nu_1$ is contained in a block$G^1$ partitioned by
$(\wt{T}, \wt{\sim})$, and $l^1$, the index of the block-column of $G^1$, belongs to $\L^1$. We propose the
following definition and assumption.

\medskip
{\bf Definition 9.5.1} A local  bocs $\mathscr{A}$ (not necessarily layered) is said to have generalized MW5, if
$\mathscr{A}$ has a unique parameter $\lambda$ with a domain $g(\lambda, \mu)\neq 0$, such that Formula (42) of
Lemma 5.4.1 holds, and for any fixed $\lambda^0\in k$ with $g(\lambda^0, \lambda^0)\neq 0$, all the $g_{ll}(\nu,
\kappa)\in k[\nu, \kappa, \sigma_l(\nu)^{-1}, \sigma_l(\kappa)^{-1}]$ are invertible in Formula (41). \hfill
$\square$

Our assumptions are:

{\bf S11}\ $\gamma^1=1$, the end term $\wt{\mathfrak{B}}^{s+J}$ of the sequence of partial bocses is in the case
of generalized MW5.

{\bf S12} $l^1\in \L^1\in \wt{T}_2/\wt{\sim}_2$, such that $\L^1\cap \wt{T}_1=\phi.$

 From now on, we  write
$(~^1{N}^r, ~^1{R}_1^r\times ~^1{R}_2^r, ~^1{\n}_1^r\times ~^1{\n}_2)^r$ instead of $(\wt{N}^r,
\wt{R}_1^r\times\wt{R}_2^r, \wt{\n}_1^r\times\wt{\n}_2^r)$ and sequence $(~^1\bigtriangleup)$ instead of
$(\wt{\bigtriangleup})$. We will construct a new sequence $(~^2\triangle)$ based on $(~^1\triangle)$ similarly
to structure 9.1  as follows.

First we define a new size vector $~^2{\n}$ of $({~^1}{T}, {~^1}{\sim})$, such that
$~^2{n}_{\L^1}={~^1}{n}_{\L^1}+1$, and $~^2{n}_{\I}={~^1}{n}_{\I}$, $\forall \, I\in\!~^1T/\!~^1\!\sim \setminus
\{\L^1\}$. Then we  add a column consisting of $0$'s and $\emptyset$'s at the left of the $l$-th block-column
for each $l\in \L^1$  into $~^1\mathcal{N}^s$. Thus $~^2T=~^1T \cup \Q^1$ where $\Q^1=\{o^1_1, o_2^1, \cdots,
o_{l_{\L_1}}^1\}$ is the set of the indices  of the added columns of $\L^1=\{l_1^1, l_2^1, \cdots,
l_{\L^1}^1\}$. And $~^2T^s$  has either equivalent classes $\P, \Q, \Q_1$ with $Q_1=Q^1$ when $\L^1\ne\L^0$; or
$\P, \Q_1$ with $\Q_1=\Q \cup \Q^1$ when $\L^1= \L^0$. The formula (57) and (59) of 9.3 do not change in
$~^2{\mathfrak{B}}^s$. When $\L^1\ne \L^0$, Formulae (58) remains,  (see Figure 1 below); when $\L^1=\L^0$,
$c_l, u_l, l=1, \cdots, r$ and $c^0$ change to $1\times 2$ matrices with the unchanged differentials of $c_l$
given in Formulae (58), (see Figure 2).
\begin{center}\unitlength=1mm
\begin{picture}(100, 50)
\put(0,10){\framebox(20,40)} \put(0,25){\line(1,0){20}} \put(0,40){\line(1,0){20}} \put(5,5){\line(0,1){20}}
\put(20,5){\line(0,1){4.5}} \put(0,20){\line(1,0){5}} \put(10,6){\vector(-1,0){5}} \put(15,6){\vector(1,0){5}}
 \put(-5,21){$\alpha^0$} \put(1, 21){$c^0$} \put(1, 14){$C$}
 \put(10, 16){$A$} \put(8, 30){$A^1$} \put(8, 43){$D^1$}
  \put(11,5){$l^0$}

\put(30,10){\framebox(20,40)} \put(35,20){\line(1,0){15}} \put(30,40){\line(1,0){20}} \put(35,5){\line(0,1){35}}
\put(50,5){\line(0,1){4.5}} \put(30,35){\line(1,0){5}} \put(40,6){\vector(-1,0){5}} \put(45,6){\vector(1,0){5}}
 \put(26,37){$\alpha^1$} \put(31, 36){$c^1$} \put(30.5, 23){$C^1$}
 \put(40, 15){$A$} \put(40, 30){$A^1$} \put(38, 43){$D^1$}
  \put(41,5){$l^1$}

\put(70,10){\framebox(25,40)} \put(70,25){\line(1,0){25}} \put(70,40){\line(1,0){25}} \put(80,6){\line(0,1){19}}
\put(75,25){\line(0,1){15}}\put(70,20){\line(1,0){10}} \put(70,35){\line(1,0){5}} \put(95,2){\line(0,1){8}}
\put(75,2){\line(0,1){4}}

\put(66,37){$\alpha^1$} \put(71, 36){$c^1$} \put(70.5, 29){$C^1$}
 \put(66,21){$\alpha^0$} \put(75, 21){$c^0$} \put(74, 14){$C$}
 \put(87, 16){$A$} \put(83, 30){$A^1$} \put(78, 43){$D^1$}
\put(87,6){$l^0$} \put(84,2){$l^1$} \put(83,4){\vector(-1,0){8}} \put(87,4){\vector(1,0){8}}
\put(86,8){\vector(-1,0){6}} \put(90,8){\vector(1,0){5}}

\put(18,-2){\mbox{figure 1}} \put(78,-2){\mbox{figure 2}}

\end{picture}
\end{center}
 Denote by $c^1$ the  $(\alpha^1, o_{l^1}^1)$-entry; and by
$C^1$ the set of entries at $o_l^1$-th column before $c^1$ (and after $c^0$ in the case of $\L^1=\L^0$). Let
$A^1\subset D$ up to the $(\alpha^1,l^1)$-th row block. The differentials of the solid arrows before the
$(\alpha^1, l^1)$-th row in $A^1$, are given by Formula (62) and we may denote it by $(57)^1$; Those at the
$(\alpha^1, l^1)$-th row of $A^1$ are parallel to  Formula (59) of 9.3, we may denote it by $(59)^1$. And we
also have Formula $(58)^1$ for the entries in the area $C^1$. Let $h^1(\lambda,\lambda)\neq 0$, we set $a\in
A\cup A^1$ before $c^1$ to be $\emptyset$; $c\in C\cup C^1$ to be $\emptyset$ or $0$; $c^0=(1)$ when $\L^1\ne
\L^0$, or $c^0=(0\ 1)$ when $\L^1=\L^0$.
 We conclude that $\delta(c^1)=0$ in the
induced partial bocs $~^2\mathfrak{B}^{s+J+j^1+r^1}$ under the assumption S11, S12 similarly to Corollary 9.4.1.
Let $c^1=(1)$,
 set the  solid arrows at the $(\alpha^1, l^1)$-th row of $A^1$ to be $\phi$,
  then we obtain an induced bocs
$~^2{\mathfrak{B}}^{s+J+J^1}$, where $J^1=j^1+r^1+1+m^1$. The local bocs $~^2{\mathfrak{B}}^{s+J+J^1}$  may not
be freely parameterized in general, because $h^1(\lambda,\mu)$  does not necessarily satisfy the item (3) of
Definition 4.2.1. Then we perform procedure 5.3 starting from the local bocs ${}^2{\mathfrak{B}}^{s+J+J^1}$ once
again, and obtain a sequence of parameters $\lambda=\lambda_{20},\lambda_{21},\cdots,\lambda_{2\gamma^2}$. If
$\gamma^2=1$ and  $~^2{\mathfrak{B}}^{s+J+J^1}$ satisfies the assumption S21, S22 parallel to S11, S12,  we may
continue to use the structure of 9.1 once more.

Now suppose we have  a chain of indices:
$$ (\alpha^0,l^0)\prec (\alpha^1,l^1)\prec\ \cdots\ \prec
(\alpha^{p-1},l^{p-1})\prec (\alpha^p,l^p),$$ such that $\gamma^1=1, \cdots, \gamma^{p-1}=1$, and  we have a
chain of index sets
$$T\subseteq ~^1T\subseteq ~^2T \subseteq \cdots \subseteq
~^{p}T$$ such that $~^qT=T\cup \Q\cup \Q^1\cup \cdots \cup \Q^{q-1}$ for $0\leq q \leq p$, where
$\Q^{q-1}=\{o_1^{q-1}, o_2^{q-1}, \cdots o_{\L^{q-1}}^{q-1}\}$ is the set of the indices of the added columns if
$\L^{q-1}=\{l_1^{q-1}, l_2^{q-1}, \cdots, l_{\L^{q-1}}^{q-1}\}$. And $~^qT^s$ possesses the equivalent classes
$\P$ and $\Q_{q1},\cdots,\Q_{qr_q}$, which are unions of some $\Q, \Q^1, \cdots, Q^{q-1}$. Moreover we have a
chain of reduction sequences:
$$(\bigtriangleup), \,
({~^1}{\bigtriangleup}), \, (~^2{\bigtriangleup}), \, \cdots, (~^p{\bigtriangleup});$$ a chain of the
corresponding sequence of partial bocses
$$
(***), ({~^1}{***}), (~^2{***}), \cdots, (~^p{***});
$$
and a chain of the end terms of the preceding sequences:
$$\mathfrak{B}^s,
{~^1}{\mathfrak{B}}^{s}, ~^2{\mathfrak{B}}^{s}, \cdots, ~^p{\mathfrak{B}}^{s};$$
 such that all the local partial bocses
$$\mathfrak{B}^s,
{~^1}{\mathfrak{B}}^{s+J}, ~^2{\mathfrak{B}}^{s+J+J^1}, \cdots,
~^{(p-1)}{\mathfrak{B}}^{s+\sum_{\eta=0}^{(p-2)}J^{\eta}}$$ satisfy the assumption  S01, S02; S11, S12; S21,
S22; $\cdots$; S(p-1)1, S(p-1)2 inductively ,  where $J^{\eta}=j^{\eta}+r^{\eta}+1+m^{\eta}$.  We perform the
procedure 5.3 for the local  bocs $~^p{\mathfrak{B}}^{s+\sum_{\eta=0}^{p-1} J^{\eta}}$ induced from
$~^p\mathfrak{B}^s$, and obtain a sequence of parameters $\lambda=\lambda_{p0}, \lambda_{p1},\cdots, \lambda_{p
\gamma^{p}}$. If $\gamma^p=1$, let
\begin{equation} h^p(\lambda, \mu)=h^{p-1}(\lambda,
\mu)\prod\limits_{l=1}^{j^p}c_l^p(\lambda) h^p_{ll}(\lambda, \mu)\prod\limits_{l \in I^{p}} f_{l}^{p}(\lambda);
\end{equation}  If $~^p{\mathfrak{B}}^{s+\sum_{\eta=0}^{p-1}
J^{\eta}}$ still satisfies the assumption SP1, SP2, and we set $\nu_p=\lambda_{p1}$, then $(\alpha^p,l^p)\succ
(\alpha^{p-1},l^{p-1})$.
 Adding a column to the
left of  the $l$-th block column of $~^p{\mathcal{N}}^s$ for all $l\sim l^p$, we obtain a bocs
$~^{(p+1)}{\mathfrak{B}}^{s}$. Such a procedure must stop at some stage since the index set $\{(\alpha^q,l^q)\}$
is bounded by $(n_1,t_2)$, (see the explanation below Lemma 9.3.2).

{\bf Proposition 9.5.1}$\quad$  We perform the procedure of 5.3 starting from bocs
$~^p{\mathfrak{B}}^{s+\sum_{\eta=0}^{p-1}J^{\eta}}$. If $\gamma^p\ne 1$, or $\gamma^p=1$ but
$~^p{\mathfrak{B}}^{s+\sum_{\eta=0}^{p-1}J^{\eta}}$ does not satisfy the assumption SP1 and SP2, then we must
have one of the following possibilities.

{\bf Q1.} $\gamma^{p}=0$, then Formula (60) of 9.3 becomes
$$
\left\{\begin{array}{ccll}
\delta(d^{p-1}_1)^0&=&g^{p-1}_{11}(\lambda, \mu)v_1^{p-1}&\\
\delta(d^{p-1}_2)^0&=&g^{p-1}_{21}(\lambda, \mu)v_1^{p-1}
&+\ g^{p-1}_{22}(\lambda, \mu)v_2^{p-1}\\
\cdots & &\cdots\qquad\cdots  \\
\delta(d^{p-1}_{n^{p-1}})^0&=&g^{p-1}_{n^{p-1}1}(\lambda, \mu)v^{p-1}_1&+\ g^{p-1}_{n^{p-1}2}(\lambda,
\mu)v_2^{p-1}+\cdots+g^{p-1}_{n^{p-1}n^{p-1}}(\lambda, \mu)v_{n^{p-1}}^{p-1}
\end{array} \right.
$$ and
\begin{equation}
\sigma_l^{p-1}(\lambda)=h^{p-1}(\lambda, \lambda ) \prod\limits_{l\in I^p} f_l^p(\lambda) \prod\limits_{q=1}^{l}
c_q^{p-1}(\lambda),  \mbox{ write } \sigma^{p-1}_{n_{p-1}}(\lambda)=\sigma^{p-1}(\lambda),
\end{equation}
 (see  (41) of 5.3), where
$v_1^{p-1},v_2^{p-1},\cdots,v_{n^{p-1}}^{p-1}$ are linearly independent, and $g_{ll}^{p-1}(\lambda,\lambda)\neq
0$.

(1) There exists some $\lambda^0\in k$ with $\sigma_l^{p-1}(\lambda^0)\neq 0$ such that $g_{ll}^{p-1}(\lambda,
\mu)\in k[\lambda, \mu, $ $\sigma_l^{p-1}(\lambda)^{-1} \sigma_l^{p-1}(\mu)^{-1}]$ for $1\leq l \leq e-1$ are
all invertible, but $g_{ee}^{p-1}(\lambda, \mu)\in k[\lambda, \mu, \sigma_e^{p-1}(\lambda)^{-1}$ $
\sigma_e^{p-1}(\mu)^{-1}]$ is not invertible.

(2) $\forall \, \lambda^0\in k $ with $\sigma_l^{p-1}(\lambda^0)\neq 0$, $g_{ll}^{p-1}{(\lambda, \mu)}\in
k[\lambda, \mu, \sigma_l^{p-1}(\lambda)^{-1} \sigma_l^{p-1}(\mu)^{-1}]$ are all invertible for $1\leq l \leq
n^{p-1}$.

{\bf Q2.} $\gamma^{p}=1$, then  Formula (62)of 9.5  becomes
$$
\left\{\begin{array}{lcl}
\delta(a_1^{p})^0&=&h^{p}_{11}(\lambda, \mu)w^{p}_1\\
\cdots & & \cdots\qquad\cdots \\
\delta(a^{p}_{i^{p}-1})^0 &=& h^{p}_{i^{p}-1,1}(\lambda,\mu)w^{p}_1+\cdots+
h^{p}_{i^{p}-1,i^{p}-1}(\lambda,\mu)w^{p}_{i^{p}-1}\\
\delta(a^{p}_{i^{p}})^0&=&h^{p}_{i^{p},1}(\lambda,\mu)w^{p}_1\quad+\cdots+
h^{p}_{i^{p},i^{p}-1}(\lambda,\mu)w^{p}_{i^{p}-1} \quad + \wt{h}^{p}(\lambda,\mu)\wt{w}^{p}
\end{array} \right.
$$
where $\wt{w}^{p}$ is a linear combination of $w_1^p, \cdots, w^p_{i^p-1}$ or $\wt{w}^{p}$ is linearly
independent of $w_1^p, \cdots, w^p_{i^p-1}$, but  $(\lambda-\mu)^2\mid \wt{h}^{p}(\lambda,\mu)$.

 Moreover for any fixed  $\lambda^0\in\rK$ with
$h^{p}(\lambda^0, \lambda^0)\neq 0$, if we  set $a_l^{p}=\emptyset$, $l=1,\cdots, i^{p}-1$, and
$a_{i^{p}}^{p}=\nu_{p}$, and write $d_l^{p}=a_{i^p+l}^p$, then we have the following formula:
$$
\left\{\begin{array}{ccl}
\delta(d_1^{p})^0&=&g^{p}_{11}(\nu_p, \kappa_p)v^{p}_1\\
\cdots & & \cdots\qquad\cdots \\
\delta(d^{p}_{n^{p}})^0&=&g^{p}_{n^{p},1}(\nu_p, \kappa_p)v^{p}_1+ \cdots +
 g^{p}_{n^{p},n^{p}}(\nu_p,
\kappa_p)v^{p}_{n^{p}}
\end{array} \right.
$$  Let
\begin{equation*} \sigma_l^{p}(\nu_p)=h^p(\nu_p, \nu_p) \prod\limits_{q=1}^l c_q^p(\nu_p), \mbox{ write }
\sigma^p(\nu_p)=\sigma^p_{n_p}(\nu_p) \end{equation*}
be given by Formula (64) of 9.5 and (41) of 5.4, then
$\forall \, \lambda^0\in\rK$ with $h^{p}(\lambda^0, \lambda^0)\neq 0$, $g^{p}_{ll}(\nu_p,\kappa_p)\in
k[\nu_p,\kappa_p, \sigma_l^p(\nu_p)^{-1}\sigma_l^p(\kappa_p)^{-1}]$ are all invertible.

{\bf Q3.} $\gamma^{p}\geq 2$, then we have a sequence of parameters $\lambda=\lambda_{p0}, \lambda_{p1},\cdots,
\lambda_{p\gamma^{p}} $. If we set $\lambda_1=\lambda_{p,\gamma^{p}-1}$ having a pair of indices
$(\alpha(1),l(1))$, then $(\alpha(1),l(1))>(\alpha(0),l(0))$.

\subsection{Non-homogeneous partial bocs $~^p{\mathfrak{B}}^s$}

\kg We will prove that the partial bocs $~^p{\mathfrak{B}}^s$ given in Q1 and Q2 of  Proposition 9.5.1 is not
homogeneous in this subsection.

 In case of Q1, let $\lambda^0 \in \rK$ with $\sigma_e^{p-1}(\lambda^0) \neq
0$ in item (1), or $\sigma^{p-1}(\lambda^0) \neq 0$ in item (2), (see Formula (65)).  We define an object $S \in
R(~^p{\mathfrak{B}}^s)$, such that $S_\P=\rK$; $S_{Q^q}=\rK$ for $0\leq  q < p$; $S(\lambda)=(\lambda^0)$,
$S(a)=\emptyset$, $\forall\ a\in \displaystyle \bigcup_{\eta=0}^{p-1}A^{\eta}$; $S(c)=\emptyset$ or $0$,
$\forall\ c\in \displaystyle \bigcup_{\eta=0}^{p-1}C^\eta$; \ $S(c^\eta)=(0 \cdots 1)$ for $0<\eta \le p-1$;
$S(d)=\emptyset$, $\forall\ d\in D^{p-1}$.

{\bf Lemma 9.6.1}\ \ Suppose we are given a partial bocs $~^p{\mathfrak{B}}^s$ satisfying Q1 of Proposition
9.5.1 and the object $S\in R(~^p\mathfrak{B}^s)$ defined above. If $(e):$ $S \stackrel{\iota}{\longrightarrow} E
\stackrel{\pi}{\longrightarrow}S$ is an almost split conflation in $R(~^p{\mathfrak{B}}^s)$, then $E_\P=\rK^2$;
$E_{\Q^q}=\rK^2$;
$$
E(\lambda)=J_2(\lambda^0); \ E(a)=\emptyset;\ E(c)=\emptyset\ \ {\rm  or}\ \ 0,\ \ E(c^\eta)=(0\ \cdots\ 0\
I_2),
$$
 and
$\iota, \pi$ are given by Corollary 6.4.1.

{\bf Proof. } The eigenvalue of $E$ is $\lambda^0$, since the sequence $0\longrightarrow
S\stackrel{\iota_0}{\longrightarrow} E \stackrel{\pi_0}{\longrightarrow} S\longrightarrow 0$ is exact over
$\rK[\lambda, \sigma_e^{p-1}(\lambda)^{-1}]$ or $k[\lambda, \sigma^{p-1}(\lambda)^{-1} ]$

 Suppose we have  some  $c_i=c^q$ or $c_i\in C^q$ with $S(c_i)\ne
 \emptyset$, such  that $E(a)=\emptyset$,
$\forall\ a\in \bigcup_{\eta=0}^{p-1} A^{\eta}$ before $c_i$;
 $E(c)=\emptyset$ or $0$, $\forall\ c\in \bigcup_{\eta=0}^{p-1}
 C^{\eta}$ before $c_i$; $E(c^{\eta})=( 0\ \cdots \ 0 \
I_2 )$ for any $c^{\eta}$ before $c_i$. We claim that $E(c_i)=(\begin{array}{cccc} 0& \cdots & 0&
I_2\end{array})$ when $S(c_i)=(\begin{array}{cccc} 0& \cdots & 0& 1\end{array})$ by item (1) of Lemma 8.5.2; and
$$E(c_i)=0\ \ {\rm or}\ \ \left(\begin{array}{ccccc} 0& \cdots & 0&0& 1\\  0& \cdots & 0& 0&0 \end{array} \right)$$ when
$S(c_i)=0$ by item (2) of Lemma 8.5.2. In the second case, we define a full subcategory $\mathscr{C}$ of
$R(~^p\mathfrak{B}^s)$ consisting of the objects $M$ such that $M_{\P}=k^m$, $=M_{\Q^q}=k^m$, $0\leq q< p$,  the
eigenvalue of $M(\lambda)$ is $\lambda^0$, $M(a)=\emptyset$, $M(c)=\emptyset$ or $0$, $M(c^{\eta})=(0, \cdots,
0, I_m)$ for all the $a,c,c^{\eta}$ before $c_i$. We take
 an object $L\in\mathscr{C}$, such that $L_{\P}=k^2$,
 $L_{\Q^q}=k^2$, $0\leq q< p$;
$L(\lambda)=J_2(\lambda^0)$; $L(a)=\emptyset$, $\forall\ a\in
\bigcup_{\eta=0}^{p-1} A^{\eta}$ before $c_i$; $L(c)=\emptyset$ or
$0$, $\forall c\in \bigcup\limits_{\eta=0}^{p-1}C^{\eta}$ before
$c_i$; $L(c^{\eta})=(0, \cdots, 0, I_2)$ for any $c^{\eta}$ before
$c_i$; and $L(c_i)=0$. We also define a morphism $\varphi:
L\rightarrow S$, such that $\varphi_{\P}={1\choose 0}$,
$\varphi_{\Q^q}={1\choose 0}$ for $0\leq q< p$. If $(e)$ is an
almost split conflation in $R(~^p\mathfrak{B}^s)$, it is so in
$\mathscr{C}$. But Lemma 8.5.3 leads to a contradiction  in
case of $$E(c_i)=\left(\begin{array}{ccccc}  0& \cdots & 0&0& 1\\
0& \cdots & 0& 0&0 \end{array} \right).$$ Thus $E(c_i)=0$. Our conclusion on $E(a)$, $E(c)$ and $E(c^{\eta})$
follows by induction.

 Finally, $E(\lambda)$ must equal to
$J_2(\lambda^0)$, otherwise $E(\lambda)=\lambda^0 I_2$ would lead to $(e)$ to be split. The proof is finished.
\hfill$\square$

{\bf Proposition 9.6.1}\  The partial bocs $~^p{\mathfrak{B}}^s$ given by condition Q1 of Proposition 9.5.1 is
not homogeneous.

{\bf Proof.} We take an infinite list of objects $S_{\lambda^0}$ defined in the beginning of 9.6. If
$R(~^p\mathfrak{B}^s)$ is homogeneous, there exists a cofinite subset $D_0\subset k\setminus\{\mbox{the roots of
} \sigma^{p-1}(\lambda)\}$ (see Formula (65) of 9.5), such that $S_{\lambda^0}$ is homogeneous for any
$\lambda^0\in D_0$. Fix such a $\lambda^0$ we denote $S_{\lambda^0}$ by $S$ for simplicity. If $(e)$ $S
\stackrel{\iota}{\rightarrow}
  E \stackrel{\pi}{\rightarrow} S$  is an almost split conflation starting and
ending at $S$ in $R(~^p{\mathfrak{B}}^s)$, then  Lemma 9.6.1 gives the structure of $E$.

(1) In case of item (1) of Q1, we define an object $L\in R(~^p{\mathfrak{B}}^s)$, such that $L_{\P}=k^2$,
$L_{\Q^q}=k^2$;
$$L(\lambda)=\left(\begin{array}{cc}\lambda^0& 0\\
0&\mu^0\end{array}\right)\ \ {\rm with}\ \ g_{ee}^{p-1}(\lambda^0, \mu^0)=0,\ \
\sigma^{p-1}(\lambda^0)\sigma_e^{p-1}(\mu^0)\ne 0;$$ $L(a)=0$, $\forall a\in \cup_{\eta=0}^{p-1}A^{\eta}$;
$L(c)=0$, $\forall c\in \cup_{\eta=0}^{p-1}C^{\eta}$; $L(c^{\eta})=(0, \cdots, 0, I)$ for $0\leq \eta \leq p-1$;
$L(d)=0$, $\forall d\in D^{p-1}$. We also define a morphism $\varphi: L\rightarrow S$, such that
$\varphi_{\P}={1\choose 0}$, $\varphi_{\Q^q}={1\choose 0}$. Then a contradiction appears similarly to the proof
of Proposition 7.4.1.

(2) In case of item (2) of Q1, we recall that
 $~^pT^s$ has the equivalent classes $\P; \Q_1, \cdots,
\Q_{r}$, where $\P$ is the unique equivalent class of $\mathfrak{B}^s$, $\Q_1, \cdots, \Q_{r}$ are some unions
of $\Q, \Q^1, \Q^2, \cdots, $ $ \Q^{p-1}$ respectively. Suppose $\Q_j=\cup_{l=1}^{n_j}\Q^{j_l}$.  We define  a
full subcategory $\overline{\mathscr{C}}_0$ of $R(~^{(p-1)}\mathfrak{B}^s)$ consisting of objects $M$ (see Lemma
8.6.1) such that $M_{\P}=k^m, M_{\Q^q}=k^m$, $0\leq q< p$; $M(\lambda)$ has  eigenvalue $\lambda^0$ with
$\sigma^{p-1}(\lambda^0)\ne 0$; $M(a)=\emptyset$, $\forall a\in \cup_{\eta=0}^{p-1}A^{\eta}$; $M(c)=\emptyset$
or $0$, $\forall c\in \cup_{\eta=0}^{p-1}C^{\eta}$; $M(c^{\eta})=(0, \cdots, 0, I^m)$ for $0\leq \eta \leq p-2$.
Moreover we define a full subcategory $\overline{\mathscr{C}}$ of $R(~^p{\mathfrak{B}}^s)$ (see Lemma 8.6.2)
according to the following two cases:

(i) $\Q^{p-1}\in ~^pT^s/~^p\sim^s$ is an equivalent class;

(ii) $\Q^{p-1}\subset \Q_{j_0}$ in $ ~^pT^s/~^p\sim^s$ for some $1\leq j_0\leq r$.

\noindent Then $\forall M\in \overline{\mathscr{C}}$, $M_{\P}=k^m$; $M_{\Q_j}=k^{mn_j}$ for $1\leq j\leq r$ in
case (i), or $1\leqslant j \leqslant r$,  $j\ne j_0$ in case (ii); $M_{\Q^{p-1}}=k^n$ in case (i) or
$M_{\Q_{j_0}}=k^{n+mn_j}$ in case (ii). And $M(\lambda)$ has eigenvalue $\lambda^0$ with
$\sigma^{p-1}(\lambda^0)\ne 0$; $M(a)=\emptyset$, $\forall a\in \cup_{\eta=0}^{p-1}A^{\eta}$; $M(c)=\emptyset$
or $0$, $\forall c\in \cup_{\eta=0}^{p-1}C^{\eta}$; $M(c^{\eta})=(0, \cdots, 0, I^m)$ for $0\leq \eta \leq p-2$;
$M(c^{p-1})=P$.

Let $\mathscr{C}$ be the category defined in Lemma 8.6.2, we define a functor
$$
F_1: \overline{\mathscr{C}}\rightarrow \mathscr{C} \quad \mbox{with  } F_1(M)=M'
$$
such that $M'_{\P}=k^m$; $M'_{\Q_j}=k^{mn_j}$ for $1\leq j\leq r$ in case (i) or $1\leqslant j \leqslant r$,
$j\ne j_0$ in case (ii); $M'_{\Q^{p-1}}=k^n$ in case (i) or $M'_{\Q_{j_0}}=k^{n+mn_j}$ in case (ii). And
$M'(\lambda)=M(\lambda)$; $M'(a_{jl})=(M(c^{j_l})\mid 0)_{1\times n_{j_l}}$ for $1\leq j\leq r$ in case (i) or
$1\leq j\leq r$, $j\ne j_0$ in case (ii); $M'(c)=P$ in case (i), or $M'(c)=(P\mid 0)_{1\times (n_{j_0}+1)}$ in
case (ii). It is obvious that $F_1$ is a representation equivalence. In fact, the morphisms in
$\overline{\mathscr{C}}$  go to the diagonal parts according to the partition of $(~^pT^s/~^p\sim^s)$ under the
action of $F_1$.

Let $\mathscr{D}$ be the category defined in Lemma 8.6.3. Then we have a representation equivalence $$F_2:
\mathscr{C}\rightarrow \mathscr{D},$$ such that $F_2(M')=M''$ with $M''_{\P}=k^m$ and $M''_{\Q}=k^n$;
$M''(\lambda)=M(\lambda)$; $M''(c)=P$. In fact, the morphisms in $\mathscr{C}$  go to the diagonal blocks under
the action of $F_2$, i.e. ``$*$" of  $\varphi_{j_0}$ in case (ii) go to zero in Lemma 8.6.2.

Thus $$F=F_2F_1: \overline{\mathscr{C}}\rightarrow \mathscr{D}$$ is a representation equivalence. If $(e)$ given
in the beginning of the proof was an almost split conflation in $R(~^p{\mathfrak{B}}^s)$, then so was in
$\overline{\mathscr{C}}$. Therefore $F(e)$  would be an almost split conflation in $\mathscr{D}$ by Lemma 8.6.4.
Which contradicts to Lemma 8.6.3.

The proof is finished, i.e.  $~^p\mathfrak{B}^s$ is not homogeneous in both cases of Q1 as desired.
\hfill$\square$

\medskip

{\bf Proposition 9.6.2} The partial bocs $~^p{\mathfrak{A}}^s$ given by condition Q2 of  Proposition  9.5.1 is
not homogeneous.

{\bf Proof.}  We fix some $\lambda^0\in k$ with $h^{p}(\lambda^0, \lambda^0 )\neq 0$ given in  Formula (64) of
9.5. Let $\lambda=\lambda^0$, $a_l^{p}=\emptyset$, $l=1,\cdots,i^{p}-1$, $a_{i^p}^p=\nu^p$ then the partial bocs
$~^p{\mathfrak{B}}^{s'}_{(\lambda^0)}$ with $s'=s+(\sum\limits_{\eta=0}^{p-1} J^{\eta})+(i^p+n^p)$ induced from
$~^p{\mathfrak{B}}_{(\lambda^0)}^{s+\sum_{\eta=0}^{p-1} J^{\eta}}$ is minimal.  We construct an infinite list of
the objects $\{S'_{\nu^0_{p}}\mid \nu^0_p\in \rK \mbox{ with } \sigma^{p}(\nu^0_{p})\neq 0 \}\subset
R(~^p{\mathfrak{B}}^{s'}_{(\lambda^0)}),$ such that
 $(S'_{\nu^0_p})_{\P}=k$; $(S'_{\nu^0_p})(\nu_p)=(\nu^0_p)$;
 $(S'_{\nu^0_p})(d_l^p)=\emptyset$ for $1\leq l\leq n_p$.
  If
$$\vartheta: R(~^p{\mathfrak{B}}^{s'}_{(\lambda^0)})\rightarrow R(~^p{\mathfrak{B}}^{s})$$ is the reduction
functor, then $S_{\nu^0_p}=\vartheta (S'_{\nu^0_p})$ provide  an infinite list of the indecomposables over
$~^p\mathfrak{B}^s$. If $~^p\mathfrak{B}^s$ is homogeneous, there exists a cofinite subset $D_0\subseteq
k\setminus \{\mbox{the roots of } \sigma^p(\nu_p)\}$, such that $S_{\nu^0_p}$ is homogeneous for any $\nu^0_p\in
D_0$. Denote $S_{\nu^o_p}'$ by $S'$ for simplicity. Let    $(e'):$ $S'\rightarrow E' \rightarrow S'$ be the
almost split conflations starting and ending at $S'$ given by Proposition 6.6.1. Then $(e):  S\rightarrow E
\rightarrow S$, the image of $(e')$ under $\vartheta$, is an almost split conflation of $R(~^p\mathfrak{B}^s)$
by Theorem 7.2.1. Construct an object
 $L\in R(~^p{\mathfrak{B}^s})$ such that $L_{\P}=k^2$, $L_{\Q^q}=k^2$;
  $L(\lambda)=J_2(\lambda^0)$;
$L(a)=\emptyset$ for any $ a\in \cup_{\eta=0}^{p-1} A^{\eta}$;
 $L(c)=\emptyset$ or $0$, $\forall\ c\in \cup_{\eta=0}^{p-1} C^{\eta}$; $L(c^{\eta})=
 (0\ \cdots\ 0\ I_2)$ for $\eta=1,\cdots, p-1$;
$L(a_l^{p})=\emptyset$ for $l=1,\cdots,(i^{p}-1)$, $L(a_{i^{p}}^{p})=\nu_p^0 I_2$; and $L(d_l^{p})=\emptyset$
for $l=1,\cdots,n^{p}$, then $L$  is indecomposable.

\vspace{0.5cm} \hspace{4.2cm}
\begin{picture}(80,32) \unitlength=1mm
\put(10, 0){$\rK^2$} \put(17, 1){\vector(1,0){30}} \put(49, 0){$\rK$} \put(34, 20){\vector(1,-1){15}} \put(29.5,
21){$\rK^2$} \qbezier[18](28,20)(21,13)(14,6) \put(14,6){\vector(-1,-1){1}}

\put(25,25){\oval(5,5)[t]} \put(25,25){\oval(5,5)[bl]} \put(23,22){\vector(3,1){3}}

 \put(37,25){\oval(5,5)[t]} \put(37,25){\oval(5,5)[br]} \put(39.00,22.00){\vector(-3,1){3}}

\put(7,6){\oval(5,5)[t]} \put(7,6){\oval(5,5)[bl]} \put(5,3){\vector(3,1){3}}

\put(7,0){\oval(5,5)[t]} \put(7,0){\oval(5,5)[bl]} \put(6,-3){\vector(3,1){3}}

\put(55,6){\oval(5,5)[t]} \put(55,6){\oval(5,5)[br]} \put(57,3){\vector(-3,1){3}}

\put(55,0){\oval(5,5)[t]} \put(55,0){\oval(5,5)[br]} \put(56.00,-3){\vector(-3,1){3}}

\put(-2, 8){$\s E(\lambda)$} \put(-3, 1){$\s E(\nu_p)$} \put(58, 6){$\s S(\lambda)$}\put(58, 0){$\s S(\nu_p)$}
\put(17, 28){$\s L(\lambda)$} \put(36, 28){$\s L(\nu_p)$} \put(39, 15){$\s {1\choose 0}$} \put(30, -2){$\s
{1\choose 0}$} \put(20, 15){$\s \varphi'$}
\end{picture}

\vspace{0.4cm} \noindent A contradiction follows similarly to the proof of Proposition 7.5.1 for
$\mathfrak{B}_4$. Therefore $~^p{\mathfrak{B}}^s$ is not homogeneous as desired. \hfill$\square$

\newpage

\bcen\section{Bipartite bimodule Problems}\ecen

\subsection{Bipartite property}

\kg We stress that according to Proposition 7.3.1, 7.3.2, 7.4.1 and 7.5.1,  the only possible exceptional
situation for a minimally wild bimodule problem being homogeneous is given by  MW5 of Theorem 5.6.1. On the
other hand Example (4) of 6.5 shows that $Mat(\K,\M)$ can be strongly homogeneous in this case. Thus we are
forced to focus on some special class of bimodule problems, in which $P_1(\Lambda)$ are included.

{\bf Definition 10.1.1}\ A bimodule problem $(\K,\M, H)$ is said
to be {\it bipartite}, if
\begin{itemize}
\item[\textrm{I}.] $T=T_1\bigcup T_2 $, where $T_1=\{1, 2,\cdots,t_1 \}$, $ T_2=\{t_1+1, t_1+2, \cdots,
t_1+t_2\}$, $ T/\sim \,=T_1/\sim_1\,{\bigcup}\, T_2/\sim_2$.

\item[\textrm{II}.] $\K=\Bigg\{\left(\begin{array}{cc} S_1&0\\
0&S_2\end{array}\right)\bigg\}$, where $ S_r$ are $t_r\times t_r$
upper triangular matrices, such that $ s^r_{ii} =s^r_{jj},$ if $
i\sim_r  j$,  and if $i<j$, $s_{ij}^r$ satisfies
 the following $r$-th equation for $r=1,2$,
$$\begin{array}{rcl} \sum\limits_{\I^1\ni i<j\in
\J^1}c_{ij}^{1l_1}x_{ij}&=&0,\\[+4ex]
0&=&\sum\limits_{\I^2\ni i<j\in \J^2}c_{ij}^{2l_2}y_{ij},\end{array}$$ where $ 1\leq l_r\leq q^r_{\I\J}$  for
some $q^r_{\I\J}\in \mathbb{N}$, and for each pair $$ (\I^r,\J^r)\in (T_r/\sim_r)\times (T_r/\sim_r).$$

\item[\textrm{III}.]  $\M=\left\{\left(\begin{array}{cc} 0&A\\ 0&0
\end{array}\right) \right\}$ where $ A $ are $ t_1\times t_2$ matrices
 satisfying the equations
 $$
\sum_{(i,j)\in \I^1\times\J^2}d_{ij}^lz_{ij}=0,$$ $1\leq l\leq
q_{\I\J} \mbox{ for some }q_{\I\J}\in \mathbb{N}, \mbox{ and for
each pair } (\I^1,\J^2)\in  (T_1/\sim_1)\times (T_2/\sim_2).
$

\item[\textrm{IV}.] $H=0$.

\item[\textrm{V}.] The row-indices of the free entries of $N_0$
(see Formula (19) of 3.6) are pairwise different.
\end{itemize}

 For the sake of convenience we sometimes denote
 $( \K, \M, H=0)$ by $(\K_1\times\K_2, \M,H=0)$, where  $ \K_1\times \K_2=
 \{(S_1)\times (S_2)\}$ and $\M=\{(A)\}$. And we denote a triple
$$\left(
\left(\begin{array} {cc} 0&N\\ 0&0
\end{array}\right),\ \ \left(\begin{array}{cc} R_1&0\\ 0&R_2
\end{array} \right), \n\right)
$$
 given in  Formula (6) of 2.6 by $(N, R_1\times R_2, \n_1\times\n_2)$.

The differential biquiver corresponding to a bipartite bimodule problem  is as follows. Let
$$T_1/\!\!\sim_1=\{\I_1, \cdots, \I_{s_1}\},\ \   T_2/\!\!\sim_2=\{\J_1, \cdots , \J_{s_2}\},$$ then we draw $s_1$
vertices on the top and $s_2$ vertices on the bottom. The solid arrows are all from top to bottom, and the
dotted arrows given by $\K_1$ are sitting at the top, and those by $\K_2$ are sitting at the bottom, see
Examples in 4.3.

Given a finite-dimensional algebra $\Lambda$, the bimodule problem $P_1(\Lambda)$, or equivalently $(\K,\M,
H=0)$  with $\K=\wt{\Lambda}\times \wt{\Lambda}$, $\M=\rad\wt{\Lambda}$ is obviously bipartite. Namely, let
$T_1=\{1',2',\cdots, t'\}$ according to the row-indices of $\wt{\Lambda}$ from top to bottom, and let
$T_2=\{1,2,\cdots, t\}$ according to the column-indices of $\wt {\Lambda}$ from left to right. Then $\forall \
\J\in T_2/\sim_2 $, with $$\J=\{j_1<j_2<\cdots <j_{\J}\},$$ the $j_{\J}$-th column of $N_0$ consists of free
entries and zeros (see examples of 2.1), conversely any free entries of $N_0$ must sit at $j_{\J}$-column for
some $\J\in T_2/\sim_2$. In particular,  all the  row-indices of the free entries are pairwise different.

We claim that for the simplicity of the proof, we stick to have
the hypothesis \textrm{V} in the definition, but it may not be
essential for our purpose.

\subsection{Locations}

\kg Assume that our original bimodule problem $(\K,\M,H=0)$ is wild and bipartite, which has a reduction
sequence $(*)$, such that the end term of $(***)$ is in case MW5 of Theorem 5.6.1. Then Lemma 9.3.1 shows the
following possibilities of the locations of $\lambda$ and $\nu$ with respect to the  block $\overline{G}$
containing $\nu$ partitioned by $(T,\sim)$.
 \bcen
\unitlength=1mm
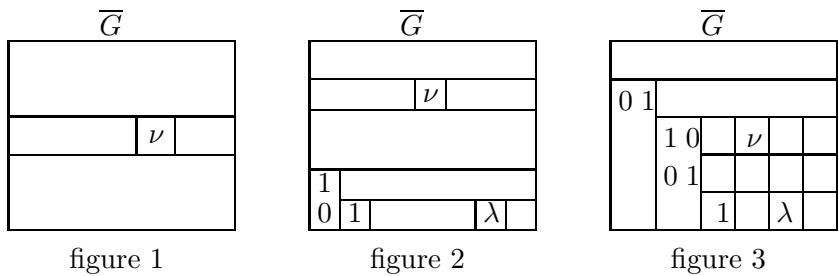
\begin{picture}(100,30)
\put(0,0){\framebox(30,25)} \put(0,10){\framebox(30,5)}
\put(17,10){\framebox(5,5)} \put(18.5,11.5){$\nu$}

\put(40,0){\framebox(30,25)} \put(40,16){\framebox(30,4)}
\put(54,16){\framebox(4,4)} \put(55,17){$\nu$}
\put(40,8){\line(1,0){30}} \put(44,0){\line(0,1){8}}
\put(44,4){\line(1,0){26}} \put(45,1){$1$} \put(41,1){$0$}
\put(48,0){\line(0,1){4}} \mput(62,0)(4,0){2}{\line(0,1){4}}
 \put(41,5){$1$}
 \put(63,1){$\lambda$}

\put(80,0){\framebox(30,25)} \put(80,20){\line(1,0){30}}
\put(86,15){\line(1,0){24}} \mput(92,5)(0,5){2}{\line(1,0){18}}
\put(86,0){\line(0,1){20}} \put(86,0){\line(0,1){15}}
\mput(92,0)(4.5,0){4}{\line(0,1){15}}
 \put(94,1){$1$}
 \put(102,1){$\lambda$}
 \put(98,11){$\nu$} \put(87,11){$1\;0$}
\put(87,6){$0\;1$} \put(81,16){$0\;1$}

\put(12,26){$\overline {G}$}
\put(52,26){$\overline{G}$}\put(92,26){$\overline{G}$}

\put(8,-5){figure 1} \put(48,-5){figure 2}\put(88,-5){figure 3}
\end{picture}\vskip 5mm
\ecen

 {\bf Definition 10.2.1} \cite{XZ}. Let $(\K,\M,H=0)$ be a bimodule problem
 having
 a reduction sequence $(*)$ of  parameterized triples, (not necessarily freely). If
 $$\overline{N^r_{p_rq_r}}=\left(
 \begin{array}{cc} 0&I\\ 0&0 \end{array}\right)\ \ {\rm or}\ \ W$$ given by edge or
 loop reduction respectively, where
 $W=\oplus_{j=1}^{\alpha} W_{\lambda^j}$ is a Weyr matrix defined in 2.5 such that
$$W_{l-1,l}^j=\left(\begin{array}{c} I_{m_l^j}\\ 0 \end{array}
\right) _{m_{l-1}^j\times m_l^j}.$$
 Then all the $1$'s appearing in
$I$ and $I_{m_l^j}$, $ l=2, \cdots, d_j$, $j=1, \cdots, \alpha$,
 are called links of $(\K^s, \M^s, H^s)$.

{\bf Location 1.} There is no  link in $\overline{G}$, or all the
links in $\overline{G}$ are lower than $\nu$ (see figure 1 and 2).

{\bf Location 2.} There exists at least one link  in
$\overline{G}$, which is parallel to or higher than $\nu$ (see
figure 3).

{\bf Examples.} Consider again the Examples in 9.2. Then Example 1
gives  $$\overline{G}=G=\left(\begin{array}{cc} \emptyset&\emptyset\\
\nu&\emptyset\end{array}\right)$$ which is in the case of Location 1.

The figure 7 of Example 2  may be used to illustrate the Location 2, but it  has only one parameter $\nu$.

\subsection{Links}

\kg Recall that for any size vector $\n$ of $(T,\sim)$ its
dimension $d=\dim\, \n =\sum_{\I\in T/\sim} n_{\I}$ (see 2.2).
Denote by $\tau(H_{\n})$ the number of links in $H_{\n}$. Then the
following proposition can be regarded as a generalization of
Theorem A in \cite{XZ}.

{\bf Proposition 10.3.1} Let
$$
 (N^r,R^r,\underline{n}^r),\cdots,(N^s,R^s,\underline{n}^s)
$$
be a reduction sequence of  parameterized triples with $r\geq 1$.
Then for any $s\geq r$,
$$
{\rm dim}\ \underline{n}^r+\tau (H_{\underline{n}^r}^r)={\rm dim}\
\underline{n}^s+\tau (H_{\underline{n}^s}^s)
$$

{\bf Proof. } We use induction on $s$. $s=r$ is trivial.  Suppose
the proposition is true for $s$, now we consider the $(s+1)$-th
reduction.

Regularization. Clearly, dim $\underline{n}^s=$ dim
$\underline{n}^{s+1}$ and $\tau(H^s)=\tau (H^{s+1})$ by 3.1.

Edge reduction. Suppose $\overline{N_{p_sq_s}^s}=\left(\begin{array}{cc}0&I_d\\0&0\end{array}\right)$, then dim
$\underline{n}^{s+1}$= dim ${\underline{n}^s}-d$ by 3.1, but $\tau(H^{s+1})=\tau (H^s)+d$, hence the formula
holds for $(s+1)$.

Loop reduction. Assume that $\overline{N_{p_sq_s}^s}=W$ is similar to $\bigoplus_{l=1}^\alpha
\bigoplus_{i=1}^{d_l}J_i(\lambda_l)^{e_l^i}$, then $\tau(W)=\sum_{l=1}^\alpha \sum_{i=1}^{d_l} e_l^i(i-1)$. On
the other hand, dim$\underline{n}^{s+1}=$ dim ${\underline{n}^s}-\tau(W)$ by 3.1 (see also [XZ, Lemma 5]). On
the other hand $\tau(H^{s+1})=\tau(H^s)+\tau(W)$. Therefore the formula still holds for $(s+1)$.\hfill$\square$

We write $\mathcal{N}_r^s=\displaystyle
\sum_{l=r}^{s-1}\overline{N_{p_lq_l}^l}\otimes \rho^l$, and
$\tau(\mathcal{N}_r^s)$ stands for the number of links in
$\mathcal{N}_r^s$. It is clear that
$$\tau(\mathcal{N}_r^s)=\tau(H_{\underline{n}^s}^s)-\tau(H_{\underline{n}^r}^r).$$

{\bf Corollary 10.3.1 } Let $
(N^r,R^r,\underline{n}^r),\cdots,(N^s,R^s,\underline{n}^s) $
 be a reduction sequence of  parameterized triples given by one of the
three reductions of 3.1. Then $(N^s, R^s, \n^s)$ is local if and
only if ${\rm dim}\
 \underline{n}^r=\tau(\mathcal{N}_r^s)+1$.

{\bf Proof.} $R^s$ is a local ring, if and only if $\dim
\,\n^s=1$, if and only if $\dim\,\n^r + \tau(H^r_{\n^r})=1
+\tau(H^s_{\n^s})$, if and only if $\dim\,
\n^r=\tau(\mathcal{N}^s_r)+1$. \hfill$\square$

\subsection{The non-homogeneous property in case of MW5}

\hspace*{0.65cm}{\bf Theorem 10.4.1 }\ Let $(\K,\M,H=0)$ be a
bipartite bimodule problem. If $(\K,\M,H)$ has a reduction
sequence $(*)$ of minimal size, such that the end term
$\mathfrak{A}^s$ of the corresponding sequence $(***)$ belongs to
MW5 of Theorem 5.6.1, then $Mat(\K,\M)$ is not homogeneous.

{\bf Proof.} In case of Location 1. We use the structure  of
adding column given in 9.1  starting  from the beginning term
$(N,R_1\times R_2,\underline{n}_1\times \underline{n}_2)$ (see
Definition 10.1.1). Since
 all the indices of the block-columns of $N$
belong to $T_2$, and $T_1\cap T_2=\emptyset$,  the assumptions S02 of  Structure 9.1 and S12, $\cdots$, S(p-1)2
 of 9.5 always take place. We stress that the triple
 $(N, R_1\times R_2, \n_1\times\n_2)$ corresponds to the  whole bocs
 $\mathfrak{B}=\mathfrak{A}$
 (not properly partial!) and $\overline{G}=G$ in this case. We use induction both on the
 index  pairs
 $(\alpha^\eta,,l^\eta)$ of  $\nu_\eta$ and the index
pairs $(\alpha(\zeta),l(\zeta))$ of $\lambda_\zeta$.

Going back to 9.5, a sequence of layered bocses $\mathfrak{B}^s=\mathfrak{A}^s$,
${~^1}\mathfrak{B}^s={~^1}{\mathfrak{A}}^s $, ${~^2}\mathfrak{B}^s={~^2}{\mathfrak{A}}^s$, $\cdots$,
 constructed under the assumptions Sq1 and Sq2 for $q=0, 1 \cdots$ must stop
at some stage, say $p$ since the pairs $(\alpha^\eta,l^\eta)$ are
bounded by $(n_1,  t_2)$, where $n_1$ stands for the number of the
rows of matrix $N$, and $t_2$ the number of the elements of $T_2$.
If $~^p{\mathfrak{A}}^{s}$ satisfies Q1 or Q2 of Proposition
9.5.1, then $~^p{\mathfrak{A}}^s$ is not homogeneous by
Proposition 9.6.1 or 9.6.2  respectively. If
$~^p{\mathfrak{A}}^{s}$ is in case Q3 of Proposition 9.5.1, and we
set $\lambda_1=\lambda_{p,\gamma^p-1}$ with a pair of indices
$(\alpha(1),l(1))\succ(\alpha(0),l(0))$, then we start a new round
of the preceding procedure from $\lambda_1$. Inductively, such
round must stop at some stage, since the pairs
$(\alpha(\zeta),l(\zeta))$ are also bounded by $(n_1,t_2)$. Our
conclusion for Location 1 follows from double induction.

In case of Location 2. $\overline{G}$ is of the following shape, and we will focus on the shadowed rectangle.
$\overline G \supsetneq G$ in this case. \vspace{1cm}

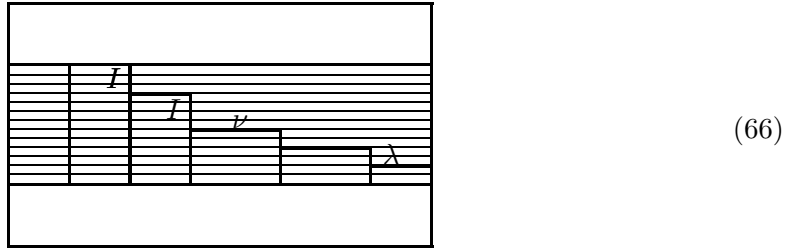
\begin{equation}
\begin{array}{c}\unitlength=0.8mm
\begin{picture}(70,40)
\thicklines\put(0,0){\framebox(70,40)} \put(0,10){\line(1,0){70}}
\put(0,30){\line(1,0){70}} \put(10,10){\line(0,1){20}}
\put(20,10){\line(0,1){20}} \put(30,10){\line(0,1){15}}
\put(45,10){\line(0,1){9}} \put(60,10){\line(0,1){6}}
\put(20,25){\line(1,0){10}} \put(30,19){\line(1,0){15}}
\put(45,16){\line(1,0){15}} \put(60,13){\line(1,0){10}}
\put(16,26){$I$} \put(26,21){$I$}\put(16,26){$I$}
\put(37,19.5){$\nu$} \put(62,13.5){$\lambda$}

\thinlines \mput(0,11.8)(0,1.5){12}{\line(1,0){70}}
\end{picture}
\end{array}
\end{equation}

 Suppose that a triple $(N^r_0, R^r, \n^r)$ (see the notation in the proof of Corollary 8.1.1) is the term of the
sequence $(*)$ having $N^r_{p_rq_r}$ as the first block of  the shadowed part in Figure (66). Then Proposition
8.1.2 and Condition $\textrm{V}$ of Definition 10.1.1 determines a matrix problem $(\K_1^r\times
\K_2^r,\overline{\M}^r, \mathcal{V}^r)$ where $$\K_1^r=\{(s)|\forall s\in k \};\ \ \K_2^r=\{S_{jj}^r\},$$ here
$S_{jj}^r$ is the $j$-th diagonal block of $\K^r$ partitioned by $(T, \sim)$ with $j$ being the column index of
$\overline{G}$; $\mathcal{V}=\text{ the shadowed part of}\ \{SH_{\n^r}-H_{\n^r}S | \, \forall S\in \K^r\}.$
Furthermore the matrix problem $(\K_1^r\times \K_2^r, \overline{\M}^r, \mathcal{V}^r)$ corresponds to a partial
bocs $\mathfrak{B}^r$ of $\mathfrak{A}^r$, which must be one sided.

We first consider the induced local bocs $\mathfrak{A}^r_{\P}$ at $\P$ obtained by deletion  with a layer
$$L_{\P}=(\Gamma'_{\P}; \omega_{\P}; b_1,\cdots, b_m, b_{m+1}, \cdots, b_{n}; u)$$ such that $ b_1,\cdots, b_m$
belong to $\mathfrak{B}^r$  but  $b_{m+1}, \cdots, b_{n}$ do not. Then a partial bocs $\mathfrak{B}_{\P}^r$ is
obtained from $\mathfrak{A}^r_{\P}$ restricted to the shadowed part.

(1) If any reduction sequence starting from the induced bocs $\mathfrak{A}^{r}_{\P}$ meets configuration (1) of
Theorem 5.6.1, then the size of the local triple at $\P$ is small than that of $(N^r, R^r, \n^r)$ since the
latter one is not local. We have a contradiction to the minimal size assumption. Therefore $\mathfrak{B}_{\P}^r$
can not satisfy $P3$ or $P4$ and must be  in case P1 or P2 of 8.2.

(2) If  $\mathfrak{B}^{r}_{\P}$ is in case P2, then
$\mathfrak{B}^r$ is given by Figure (46) of 8.2. Suppose that
$\delta(\overline{b})$ satisfies Formula (52) of 8.3,  then the
conflation $(e)$ constructed in Proposition 8.3.1 can  also be
regarded  as  an almost split conflation of $R(\mathfrak{A}^r)$.
Therefore the same proposition tells that $\mathfrak{A}^r$ is not
homogeneous. If $\delta(\overline{b})$ satisfies Formula (51) and
$\delta(a)^0$ with respect to $\overline{b}$ are given by
Proposition 8.3.2, then $\mathfrak{A}^r$ is not homogenous. In
fact the induced partial bocs of $\mathfrak{B}^r$ constructed in
the same proposition determines   a global induced bocs of
$\mathfrak{A}^r$, which is also in the case of Proposition 7.3.2.

(3) If $\mathfrak{B}_{\P}^r$ is in case P2, and $\delta(\overline{b})$ satisfies Formula (51) of 8.3, moreover
$\delta(a)^0$ with respect to $\overline{b}$ satisfy Formula (53) of 8.3. Then going back to Structure  9.1, let
$N^r$ be   the shadowed part of $N^r_0$. Let  $p\in\P\in T^r/\sim^r$ be the block-row index of $N^r$, then
$T^r_1/\sim_1^r=\{\P\}$ and $T^r_2/\sim_2^r=\{\P, \I_1, \cdots, \I_i\}$. $N^r$ is partitioned by
$(n_{\P})\times\n_2$, $R_1=(x_{ij})_{n_{\P}\times n_{\P}}$, $R_2$ is partitioned by $\n_2\times\n_2$. Therefore
the assumptions S02 of the structure 9.1 and  S12, $\cdots$, S(p-1)2  of 9.5 take place by Proposition 8.4.1,
since  $\overline{B}$, the block determined by $\overline{b}$, will not contain any parameter after any
reductions.

(4) If $\mathfrak{B}^r_{\P}$ is in case P1 of 8.2, then the
assumption S02 of 9.1 and S12, $\cdots$, S(p-1)2 of 9.5 always
take place.

Next we perform the structure  9.1 of adding columns in cases (3) and (4) starting from the triple $(N^r,
R_1^r\times R_2^r, \n_1^r\times\n_2^r)$.  Going back to 9.1, 9.2 and 9.3,  it is clear that $\lambda$, $\nu$
both belong to $\mathfrak{B}^s$ and $\Gamma'(\P,\P)=k[\lambda]$.  Then we add columns to the partial bocs
$\mathfrak{B}^s$ and obtain a partial bocs ${~^1}{\mathfrak{B}}^s$. ${~^1}{\mathfrak{B}}^s$ determines a global
bocs ${~^1}{\mathfrak{A}}^r$, which contains ${~^1}{\mathfrak{B}}^s$ as the first part. Since
${~^1}{\mathfrak{B}}^{s+J}$ given in Lemma 9.3.1 is local, ${~^1}{\mathfrak{A}}^{s+J}$, having
${~^1}{\mathfrak{B}}^{s+J}$ as the first part, must be local. In fact let $\mathcal{N}_r^s=\displaystyle
\sum_{l=r}^{s-1}\overline{N_{p_lq_l}^l}\otimes \rho^l$ in the sequence $(\bigtriangleup)$ of 9.1, and
${~^1}{N}_r^{s+J}=\displaystyle \sum_{l=r}^{s+J-1}\overline{{~^1}\mathcal{N}_{p_lq_l}^l} \otimes \rho^l$ in the
sequence $({~^1}{\bigtriangleup})$. Then ${\rm dim}_k \n^r=\tau(\mathcal{N}_r^s)+1$  according to Corollary
10.3.1 since $\mathfrak{A}^s$ is local. On the other hand,
$\tau({~^1}{\mathcal{N}}_r^{s+J})=\tau(\mathcal{N}_r^s)+1$, and ${\rm dim}_k ~^1\n^r={\rm dim}_k \n^r+1$, thus
${\rm dim}_k {~^1}\n^r=\tau({~^1}{\mathcal{N}}_r^{s+J})+1$. Then ${~^1}{\mathfrak{A}}^{s+J}$ being local follows
from Corollary 10.3.1 once again. Therefore we are able to perform the procedure 5.3 for
${~^1}{\mathfrak{A}}^{s+J}$ and obtain parameters
$\lambda=\lambda_{10},\lambda_{11},\cdots,\lambda_{1\gamma^1}$. If $\gamma^1=0$, we have the case Q1 of
Proposition 9.5.1 for the global bocs $~^1\mathfrak{A}^s$, which is not homogeneous by proposition 9.6.1.  If
$\gamma^1>0$ and  $\lambda_{1\gamma^1}$ locates outside of the shadowed part, then we reach the case of Location
1. Now we assume that
\begin{itemize}
\item[(i)] $\lambda_{1\gamma^1}$ locates inside of the shadowed part
and
 \item[(ii)] ${~^1}{\mathfrak{A}}^s$ (not only ${~^1}{\mathfrak{B}}^s$!)
satisfies  S11  of 9.5.
\end{itemize}
we stress that $~^1\mathfrak{A}^s$ also satisfies S12 in both
cases (3) and (4) as $~^1\mathfrak{B}^s$ does.
 Continue  to use  the structure  of adding columns under the  assumptions
(i) and (ii), we obtain by induction a sequence of partial bocses
$$
\mathfrak{B}^s,\ {~^1}{\mathfrak{B}}^s,\ ~^2{\mathfrak{B}}^s,\
\cdots,\ ~^p{\mathfrak{B}}^s,$$ (see 9.5), as well as a sequence of
global bocses
$$ \mathfrak{A}^s,\
{~^1}{\mathfrak{A}}^s,\ ~^2{\mathfrak{A}}^s,\ \cdots,\
~^p{\mathfrak{A}}^s.
$$
Since the pairs of indices $(\alpha^\eta,l^\eta)$ are bounded by $(n_\P^r,t_2^r)$, where $t_2^r$ stands for the
number of elements of $T_2^r$, the procedure must stop at some stage, say $p$. If $\lambda_{p \gamma^p}$ locates
outside of the shadowed part, we reach the case of Location 1. Now suppose $\lambda_{p \gamma^p}$ locates inside
the shadowed part. If $~^p{\mathfrak{A}}^{s}$ is in the condition  Q1 or Q2 of Proposition 9.5.1, then
$~^p{\mathfrak{A}}^s$ is not homogeneous by Proposition 9.6.1 or 9.6.2. If $~^p{\mathfrak{A}}^{s}$ is in the
condition of Q3 of Proposition 9.5.1, then we set $\lambda_1=\lambda_{p,\gamma^p-1}$. If $\lambda_1$ is outside
of the shadowed part, we obtain the case of Location 1. Otherwise,  the pair of indices of $\lambda_1$,
$(\alpha(1),l(1))\succ(\alpha(0),l(0))$ in the shadowed part. Thus we start a new round of the above procedure
for $\lambda_1$. Inductively such round must stop at some stage, since the pairs $(\alpha(\zeta),l(\zeta))$ are
also bounded by $(n_{\P}^r,t_2^r)$. Finally our conclusion for Location 2 follows from double induction.
\hfill$\square$

{\bf Remark.} All the discussion in this subsection also applies to MW4 of Theorem 5.6.1, which allows us to
treat MW4 and MW5 in a unified way. But because of the particularity of MW5, we prefer to treat them separately.
On the other hand we have already presented an easy proof for MW4 in Proposition 7.5.1.

\subsection{The main theorem}

\kg {\bf Key Theorem 10.5.1}\ Let $(\K,\M,H=0)$ be a bipartite
bimodule problem of wild representation type. Then Mat$(\K,\M)$ is
not homogeneous.

{\bf Proof.} Since $(\K,\M,H=0)$ is of wild representation type, there must exist a reduction sequence $(*)$  of
minimal size, such that the end term $\mathfrak{A}^s$ of the corresponding sequence $(***)$ is in one of the
cases MW1----MW5 by Theorem 5.6.1. $\mathfrak{A}^s$ is not homogeneous in case MW1, MW2, MW3, MW4 by Proposition
7.3.1, 7.3.2, 7.4.1, 7.5.1 respectively. Thus $\rm{Mat}(\K,\M)$ is not homogeneous by Corollary 7.2.1. And
$\rm{Mat}(\K, \M)$ is not homogeneous if $\mathfrak{A}^s$ is in case MW5 by Theorem 10.4.1. \hfill $\square$

{\bf Corollary 10.5.1}\ If $P_1(\Lambda) $ is of wild type, then
$P_1(\Lambda)$ is not homogeneous.

{\bf Theorem 10.5.2}\  Let $\Lambda$ be a finite-dimensional
$\rK$-algebra of representation  wild type, then $\Lambda$-mod is
not homogeneous.

{\bf Proof.} Since there exists an one-to-one correspondence between the almost split sequences of $\Lambda$-mod
and the almost split conflations of $P_1(\Lambda)$ except finitely many by Theorem 6.2.1,  the non-homogeneous
property of $P_1(\Lambda)$ implies the same property of $\Lambda$-mod. \hfill$\square$

Our Main theorem 1.3.2 mentioned in the introduction follows from Theorem 10.5.2 immediately.

\bigskip

\newpage
{\bf Acknowledgement} \ The first author expresses her thanks to
R. Bautista. She learnt this problem  from him and received his
preprint on the problem $14$ years ago. She would like to thank
C.M. Ringel,  D. Simson and S. Liu for their concern for this
work. She also thanks Professor Huang Weiming, Hu Yongjian and
Zhao Xu'an in Beijing Normal University for their helpful
conversations. Some ideas of the manuscript  has been discussed
with Drozd, Crawley-Boevey and Assashiba during ICRA  \textrm{XI}
held in Mexico in August 2004. Especially Crawley-Boevey proposed
many valuable suggestions at the Asia-link Conference in Beijing
in May 2005. In particular, both authors are grateful to Chen Bo,
Zhao Deke and Pan Jun, the students of Beijing Normal University,
for detailed discussions on this paper at our algebra seminars.

\medskip

The first author is supported by the Important Project 19331030 of NNSFC, the Research Found for the Doctoral
Program of Ministry of Education of China, the Cultivation Fund of the key Scientific and Technical innovation
project of Ministry of Education of China, China-UK joint project of the Royal Society (No. 15262), and AsiaLink
project of the European Community (ASI/B7-301/98/679-11); the second author is supported by NNSFC 10201004 and
10426014.

\newpage
\addcontentsline{toc}{section}{\protect\numberline{} References}

\noindent Zhang Yingbo

\noindent Department of Mathematics, Beijing Normal University,
 Beijing 100875, P.R.China

\noindent Email: Zhangyb@bnu.edu.cn

\bigskip
\noindent Xu Yunge

\noindent  Faculty of Mathematics $\&$ Computer Science, Hubei
University,
 Wuhan 430062, P.R.China

\noindent Email: xuy@hubu.edu.cn

\end{document}